\documentclass[12pt,a4paper]{article}
\usepackage[utf8]{inputenc}
\usepackage[T1]{fontenc}
\usepackage{lmodern}
\usepackage[shortlabels]{enumitem}
\usepackage{amsmath,amssymb,amsthm,mathtools}
\usepackage{mathrsfs}
\usepackage{array}
\usepackage[margin=1in]{geometry}

\usepackage{cite}
\usepackage{hyperref}
\hypersetup{
	pdfdisplaydoctitle=true,
	pdftitle={Packaging Systems: Intrinsic Structure, Recovery Ranks, and Globalization Theorems},
	pdfauthor={Rongchao Ma},
	colorlinks=true,
	linkcolor=blue,
	citecolor=blue,
	urlcolor=blue,
}
\usepackage{xurl}

\theoremstyle{plain}
\newtheorem{theorem}{Theorem}[section]
\newtheorem{lemma}[theorem]{Lemma}
\newtheorem{proposition}[theorem]{Proposition}
\newtheorem{corollary}[theorem]{Corollary}

\newtheorem{definition}[theorem]{Definition}
\newtheorem{example}[theorem]{Example}
\newtheorem{convention}[theorem]{Convention}

\theoremstyle{remark}
\newtheorem{remark}[theorem]{Remark}

\author{
	Rongchao Ma \\
	\textit{Department of Physics, University of Alberta} \\
	\textit{Edmonton, Alberta, Canada} \\
}

\date{}

\begin{document}

\setlength{\emergencystretch}{2em}

\title{Packaging Systems: Intrinsic Structure, Recovery Ranks, and Globalization Theorems}

\maketitle

\begin{abstract}
We define packaged structures through carriers, internal units, incidence, candidate spaces, and selected families. A realized packaging system additionally retains a context and formation datum related to its output by a possibly partial, multivalued realization relation. We establish comparison-relative recoverability and non-collapse criteria and an intrinsic refinement-completion theory.

For atlas specializations, exact extension becomes a specialization-relative packaging completion problem. Under additional recovery and elimination axioms, we establish general globalization sufficiency criteria and sharp classifications for three classes of compatible finite atlases on a fixed cover $\mathcal H$. A geometric obstruction rank $\nu(\mathcal H)$ yields thresholds $2k+1<\nu(\mathcal H)$ for empty-preserving closures of local Carath\'eodory number at most $k\ge1$; $d<\nu(\mathcal H)$ for complete functional-dependency theories of local minimum Armstrong size at most $d\ge3$; and $3<\nu(\mathcal H)$ for transitive relations, including equivalence relations.

For nonempty finite label sets and empty-preserving closures, minimum commutant realization dimension equals minimum Armstrong size. We also determine intrinsic recovery ranks and construct canonical clique completions under bounded conformality without chordality. Explicit obstructions establish sharpness; a transport hypothesis is essential for the general recovery--obstruction converse. Complexity bounds apply locally, not necessarily to global extensions.
\end{abstract}



\tableofcontents

\section{Introduction}
\label{SEC:IntroductionA}

We develop an axiomatic theory of packaged structures and packaging systems.
A packaged structure records package carriers, internal units, incidence, local candidate spaces, and selected candidate families.
A realized packaging system also retains the context and formation datum associated with a particular packaged output.
Its realization relation need not be total or single-valued.
Within this framework, we study intrinsic completion, recoverability, and comparison-relative non-collapse, and establish exact-globalization criteria for several classes of local structures.

The terminology is motivated by our earlier work on packaged entangled states and symmetry packaging \cite{Ma2017,MaSPI2025,MaSPII2025}.
Related constructions arise in field theory, lattice-gauge simulation, and gauge-invariant quantum information
\cite{MaFieldTheory2025,MaLatticeGauge2025,MaGaugeQI2025}.
Those works use specific physical or algebraic formation laws.
Here we define the general framework and study which structural and local-to-global results it supports.

\subsection{Packaged structures and the globalization problem}

A realized packaging system is a triple
\[
\mathbb P=(\Xi,P,\mathfrak P),
\qquad
\Xi\vdash P\rightsquigarrow\mathfrak P,
\]
where \(\Xi\) is a packaging context, \(P\in\operatorname{Form}_\Xi\) is a formation datum, and the realization judgment means \((P,\mathfrak P)\in\mathsf J_\Xi\).
In a packaging context, a formation datum may realize no packaged structure or several distinct packaged structures; different data may also realize the same structure.
A realized system requires an actual realization, whereas its packaged endpoint can be studied independently of the formation datum.
Restriction and gluing laws are introduced only in the relevant specializations.

The intrinsic structure of a packaged output raises questions independent of its formation mechanism.
For refinement-packaged structures, selected local refinements may agree exactly on every pairwise overlap of package interiors without admitting a global refinement with the prescribed local traces.
This leads to the intrinsic completion problem.

The later globalization theorems concern additional mathematical specializations.
In particular, closure and functional-dependency atlases are defined on covers of descriptor-label sets, which need not coincide with package interiors.
Their restriction and compatibility laws belong to the corresponding specializations and are not primitive axioms of packaging.
For the finite globalization results, let $V$ be a finite label set and let $\mathcal H=(I_j)_{j\in J}$ be a finite cover of $V$.

The universal exact-globalization problem fixes a cover $\mathcal H$ and a class of local structures, and asks whether every compatible atlas in that class has a global structure whose restriction to each original chart is exactly the prescribed local structure.
Only the local atlas varies; the cover and the local class are fixed.

The obstruction rank \(\nu(\mathcal H)\) measures how the cover geometry can prevent exact globalization.
The primal graph joins two labels whenever they belong to a common chart.
We call a clique \emph{guarded} if all its labels occur together in one chart.
If the primal graph has an induced cycle of length at least four, then \(\nu(\mathcal H)=3\).
Otherwise, \(\nu(\mathcal H)\) is the size of the smallest unguarded clique, with value \(\infty\) if every clique is guarded.
The globalization theorems relate this rank to the complexity of the local structures.

\subsection{Principal results}

We first establish what the packaging axioms allow and what they do not force.
Deterministic formation is a proper special case, as is refinement packaging on any nonempty structural domain.
We also give exact tests for recoverability and for failure of a proposed reduction.
In the access--commutant example, the commutant determines its represented source but need not be an algebra quotient of that source.

For refinement-packaged structures, we construct an intrinsic completion operator and a completion space \(\mathcal G(\mathfrak P)\).
A static isomorphism induces a bijection between completion spaces.
We give two finite examples with the same package--unit incidence and matching types of local selections.
In both examples, the selected refinements agree on every pairwise overlap.
Yet one example has exactly one global completion, while the other has none.

We classify universal exact globalization for three classes of local mathematical structures.
For a class $\mathscr K$, write $\operatorname{UG}_{\mathscr K}(\mathcal H)$ when every compatible $\mathscr K$-atlas on $\mathcal H$ admits an exact global extension.
The three classes are specified as follows:
\begin{itemize}
	\item $\operatorname{UG}_{\operatorname{car}\le k}$ concerns finite empty-preserving closures of Carath\'eodory number at most $k$;
	
	\item $\operatorname{UG}_{\operatorname{arm}\le d}$ concerns complete functional-dependency theories of minimum Armstrong size at most $d$, allowing empty premises;
	
	\item $\operatorname{UG}_{\mathrm{trans}}$ concerns transitive binary relations.
\end{itemize}
The same transitive criterion holds for preorders, partial orders, strict partial orders, and equivalence relations.
The complexity bounds apply to local structures; a global extension need not satisfy them.

\begin{theorem}[Globalization obstruction spectrum]
	\label{THM:PackagingIntroGlobalizationSpectrum}
	For every finite cover $\mathcal H$, every integer $k\ge1$, and every integer $d\ge3$,
	\[
	\begin{aligned}
		\operatorname{UG}_{\operatorname{car}\le k}(\mathcal H)
		&\Longleftrightarrow 2k+1<\nu(\mathcal H),\\
		\operatorname{UG}_{\operatorname{arm}\le d}(\mathcal H)
		&\Longleftrightarrow d<\nu(\mathcal H),\\
		\operatorname{UG}_{\mathrm{trans}}(\mathcal H)
		&\Longleftrightarrow 3<\nu(\mathcal H).
	\end{aligned}
	\]
	Every failed inequality admits an explicit finite obstruction supported on primitive cover geometry.
\end{theorem}

The globalization classifications arise within restriction-equipped atlas specializations of the packaging framework. The atlas realization bridge (Theorem~\ref{THM:PackagingAtlasRealizationBridge}) represents compatible atlases as packaging systems and identifies their exact extensions with a specialization-relative completion space. Its nonemptiness requires further mathematical structure: Section~\ref{SEC:SpectrumA} formulates recovery and elimination axioms that give general globalization sufficiency, and isolates the obstruction-transport hypotheses needed for sharp converses. The three displayed thresholds are established by verifying suitable recovery mechanisms and constructing explicit obstructions in the respective classes; they are not consequences of the primitive realization axioms alone.

Proposition~\ref{PROP:PackagingIntroTwoSeparations} distinguishes extension of an individual atlas from universal exact globalization and shows that the universal answer depends on the class of local structures.
Native dimension two is exceptional: universal exact globalization holds on every finite cover (Proposition~\ref{THM:PackagingNativeDimensionTwoClassification}).
In higher dimensions, a global commutant realization can require more dimensions than its local realizations; Proposition~\ref{PROP:PackagingNativeStrictGlobalDimensionJump} exhibits an increase from three to four dimensions.

The first two thresholds differ because the two closure classes recover composed implications in different ways.
For a nonempty finite attribute set and an empty-preserving closure \(c\), we prove
\[
d_{\mathrm{comm}}(c)=\operatorname{arm}(c),
\]
where the left-hand side is the minimum dimension of a commuting self-adjoint realization and the right-hand side is the minimum Armstrong size.

For the native class $\mathbf{CommCl}_d$, Theorem~\ref{THM:PackagingNativeGeneralResolution} establishes original-chart recovery for composable local implications whose combined support is a primal clique, provided every primal clique of size at most $d$ is guarded.
The recovered implication has a nonempty premise drawn from the original premises, of size at most $d-1$, and is valid in an original chart.
Thus, unlike ordinary Horn composition \cite{Liberatore2025DefiniteHornForgetting}, the theorem controls both premise size and the location of the recovered consequence.

For \(d\ge3\), the intrinsic recovery ranks satisfy
\[
\operatorname{ipr}(\mathbf{CommCl}_d)=d,
\qquad
\operatorname{ipr}(\mathscr C_{d-1})=2d-1,
\]
where \(\mathscr C_{d-1}\) is the class of finite empty-preserving closures of Carath\'eodory number at most \(d-1\).
Thus bounding the minimum Armstrong size yields a lower recovery rank than bounding the Carath\'eodory number alone.

For locally propagation-recoverable closure classes, controlled propagation gives clique-locality under guarding and chordality: the global closure restricts to each primal clique as the closure generated by the original rules supported there (Theorem~\ref{THM:PackagingControlledPropagationGlobalization}).
Without chordality, bounded conformality still yields compatible least clique fillers that may be adjoined simultaneously without changing the globally generated closure (Theorem~\ref{THM:PackagingCanonicalCliqueCompletion}).
The transitive case instead uses composition along paths, with cycle-based obstructions (Theorem~\ref{THM:PackagingTransitiveUniversalAmalgamation}).

More generally, the intrinsic propagation--recovery rank, defined through rooted circuits, gives a sufficient criterion for universal exact globalization.
For a restriction-stable class $\mathscr C$ of finite empty-preserving closures with finite rank $p:=\operatorname{ipr}(\mathscr C)\ge3$, the inequality $p<\nu(\mathcal H)$ implies $\operatorname{UG}_{\mathscr C}(\mathcal H)$ (Theorem~\ref{THM:PackagingIntrinsicRecoverySufficiency}).
A converse requires additional obstruction-transport hypotheses, verified by explicit constructions for the bounded-Carath\'eodory and bounded-Armstrong classes.
Theorem~\ref{THM:PackagingTransportNecessityCounterexample} provides a counterexample showing that restriction and relabeling stability alone do not suffice for the converse.

\subsection{Relation to previous work}

Our results draw on closure systems, database dependencies, hypergraphs, and amalgamation.
The comparisons below distinguish the mathematical objects being extended and the quantifiers governing extension, rather than terminology alone.
The primitive packaged coordinates admit a faithful relational presentation; our non-collapse claims therefore concern specified comparisons, not absolute non-representability.

Sheaf-type gluing and local equivalence relations provide established approaches to compatible local data
\cite{MacLaneMoerdijk1992Sheaves,
	Rosenthal1982LocalEquivalence,
	Rosenthal1984SheavesLocalEquivalence,
	KockMoerdijk1996LocalEquivalence}.
By contrast, the primitive packaging framework specifies package-side structure and a realization relation without assuming restriction maps; these enter only in appropriate specializations.

Finite closure systems, implicational bases, rooted circuits, and abstract convexity have a substantial theory
\cite{Birkhoff1967Lattice,
	GanterWille1999FCA,
	CaspardMonjardet2003ClosureSystems,
	NakamuraKashiwabara2013RootedCircuits}.
Gluing closure and convexity structures has also been studied
\cite{Sierksma1981ConvexityUnions,DegreefFourneau1984GluingsClosures}.
For a specified finite closure atlas, the least generated global closure and an exact extension criterion are considered in
\cite{Kim2026ClosureAtlases}.
Our bounded-Carath\'eodory theorem instead fixes the cover and classifies when every compatible atlas satisfying the local complexity bound has an exact global extension.

Functional dependencies, Armstrong relations, and partition semantics are classical
\cite{BeeriDowdFaginStatman1984Armstrong,
	MannilaRaiha1986DesignByExample,
	Lee1983AlgebraicRelational,
	CosmadakisKanellakisSpyratos1986PartitionSemantics}.
Local and global consistency under functional dependencies has been studied for database schemes
\cite{Sagiv1983GloballyConsistent,
	LaverMendelzonGraham1983CyclicFD,
	SagivShmueli1989FDAcyclicity,
	Wang1992BoundedDatabaseSchemes}.
Functional dependence in locally consistent \(K\)-relations is studied in \cite{BarlagEtAl2026LocallyConsistentKRelations}.
Our bounded-Armstrong theorem classifies universal exact extension for compatible complete dependency theories subject to a local minimum-Armstrong-size bound.

Related local-to-global methods occur in constraint consistency, relational elimination, and consistency witnesses for annotated relations
\cite{VanBeekDechter1997Tightness,
	DechterVanBeek1997RelationalConsistency,
	AtseriasKolaitis2025Monoids,
	AtseriasKolaitis2025Witnesses,
	AtseriasKolaitis2026Gamma}.
Our recovery invariant measures how much of a clique must lie in one original chart before a composed consequence can be recovered there.

For transitive structures, related results arise in model-theoretic amalgamation \cite{Lipparini2021SAPU,Otto2024AmalgamationPatterns}.
Our theorem classifies finite covers on which every compatible transitive atlas has an exact global extension.
The cover criteria use chordality and bounded conformality, notions familiar from hypergraph acyclicity
\cite{BeeriFaginMaierYannakakis1983Acyclic,
	Fagin1983Degrees,
	BraultBaron2016HypergraphAcyclicity,
	Otto2012HighlyAcyclic}.
Classical $\alpha$-acyclic database hypergraphs are characterized by chordality and full conformality; in the three principal exact-globalization classifications here, chordality remains necessary but only class-dependent bounded conformality is required.
The semantic objects and extension questions differ from relational-instance consistency.

\subsection{Organization}

The paper proceeds in three stages:
\begin{enumerate}
    \item \emph{Axiomatic foundations.}
    Section~\ref{SEC:FoundationsPackagingSystems} develops the primitive packaging axioms, formation and realization relations, comparison-relative recovery and non-collapse criteria, and examples of distinct native formation laws.

    \item \emph{Intrinsic structure and atlas specialization.}
    Section~\ref{SEC:PackagingIntrinsicStructureA} establishes intrinsic refinement completion, coherence, and invariance results. Section~\ref{SEC:PackagingControlledGlobalization} realizes compatible atlases as packaging systems and proves controlled propagation, canonical clique completion, and sharp globalization for bounded-Carath\'eodory closures.

    \item \emph{Further specializations and synthesis.}
    Section~\ref{SEC:ArmstrongNativeA} relates minimum commutant realization dimension to minimum Armstrong size. Section~\ref{SEC:NativeRecoveryA} proves native clique resolution and the bounded-Armstrong globalization classification. Section~\ref{SEC:PackagingCycleGeometry} develops transitive and equivalence-relation amalgamation. Section~\ref{SEC:SpectrumA} formulates recovery-equipped packaging classes, develops intrinsic recovery ranks and obstruction transport, and synthesizes the sharp class-specific results.
\end{enumerate}

The primitive packaging axioms, intrinsic completion properties, and additional restriction and recovery laws of the specializations play distinct logical roles throughout.
In particular, the atlas bridge records realized packaging systems, while exact globalization requires the further hypotheses verified in the respective classes.

\section{Foundations of packaging systems}
\label{SEC:FoundationsPackagingSystems}

This section defines the objects used throughout the paper: packaged structures, packaging contexts, formation data, and realized packaging systems.
We keep two questions separate: whether an output is allowed by a context and whether a formation datum actually realizes it.

The main foundational result, Theorem~\ref{THM:PackagingSelectionFreedom}, shows that the possible selections on a fixed scaffold form a Boolean algebra, while the formation relation can still be arbitrary.
Further results identify proper specializations and make precise which comparisons can show non-collapse.
Minimal-cover selection and access--commutant formation give concrete examples of different formation laws.
Together, these results supply the language needed for the later completion and globalization theorems.

\subsection{Minimal structural vocabulary}
\label{SUBSEC:MinimalStructuralVocabulary}

Before defining packaging, we specify the kinds and sizes of objects it can contain.
A structural domain is a set of possible carriers; each carrier can have local candidate data.
These definitions organize the data but impose no maps, quotients, or comparisons of information.
This lets us distinguish properties of the general axioms from assumptions added in later examples.

\begin{convention}[Set-theoretic size]
	\label{CONV:PackagingSetTheoreticSize}
	Fix Grothendieck universes \(\mathcal U\in\mathcal V\).
	A \(\mathcal U\)-small set belongs to \(\mathcal U\).
	A \(\mathcal U\)-class is a \(\mathcal V\)-small collection of \(\mathcal U\)-small objects, not necessarily itself in \(\mathcal U\).
	Unless stated otherwise, structural domains, carrier and package families, internal-unit sets, descriptor index sets and codomains, and sets carrying refinements are \(\mathcal U\)-small.
	Contexts, reconstruction-regime	collections, and class-valued descriptor admissibility assignments are \(\mathcal V\)-small.
	The local candidate spaces remain \(\mathcal U\)-small even when the regime collections are larger.
	For \(I\in\mathcal U\), \(\mathbf{Ref}(I)\) is \(\mathcal U\)-small.
	No construction takes a tagged union over a proper \(\mathcal U\)-class.
\end{convention}

\begin{definition}[Structural domain]
	\label{DEF:StructuralDomainPackaging}
	A \textbf{structural domain} \(\mathscr C\) is a specified \(\mathcal U\)-small set of mathematical objects.
	No common morphisms, topology, algebraic operations, incidence, or dynamics are assumed unless supplied by a specialization.
\end{definition}

\begin{definition}[Structural carrier]
	\label{DEF:StructuralCarrier}
	A \textbf{structural carrier} is an object \(K\in\mathscr C\) that may support package-side structure.
	A context specifies which carriers are admissible, through \(\operatorname{Car}_{\Xi}\subseteq\mathscr C_{\Xi}\).
	A carrier need not itself be a subset of the internal-unit set or a quotient class.
\end{definition}

\begin{definition}[Local candidate datum]
	\label{DEF:PackagingLocalCandidateDatum}
	Let \(K\) be a structural carrier and let \(I\in\mathcal U\) be a specified local unit set.
	A \textbf{local candidate datum} on the pair \((K,I)\) consists of a \(\mathcal U\)-small set
	$
	\mathsf L(K;I)
	$
	together with a distinguished subset
	\[
	\mathsf R(K;I)\subseteq\mathsf L(K;I).
	\]
	The set \(\mathsf L(K;I)\) is the \textbf{local candidate space}, and \(\mathsf R(K;I)\) is the \textbf{selected candidate family}.	
	The set \(\mathsf L(K;I)\) need not have an order, algebra, topology, probability law, or categorical structure.
	We also assume no representation by maps, partitions, or equivalence relations unless a specialization supplies one.
	The candidate space and the selected family are both permitted to be empty.
\end{definition}

The general candidate coordinate of a packaged structure will use only Definition~\ref{DEF:PackagingLocalCandidateDatum}.
Later we specialize the local candidate space to a refinement poset on the package interior.
Thus normalized refinements, quotient representatives, and their kernels are not part of the definition of an arbitrary local candidate datum.

Descriptor families are introduced for later reconstruction equipment.
Their relation to refinement candidates belongs to the refinement specialization and is not required by the general candidate coordinate.

\subsection{Packaged structures and packaging contexts}
\label{SUBSEC:PackagingCoreContextOutput}

The basic theory separates the data retained by a packaged output from the conditions governing its realization.
A \emph{packaged structure} specifies carriers, units, incidence, candidates, and selections.
A \emph{packaging context} specifies admissible outputs and the formation data that realize them.

\subsubsection{Packaged structure \texorpdfstring{$\mathfrak P$}{P}}

Axioms~\textup{(S1)}--\textup{(S5)} specify carriers, internal units, incidence, local candidate spaces, and selected families without imposing coherence between packages. Lemma~\ref{LEM:PackagingDataSize} establishes the required set-theoretic size of the model class.

\begin{definition}[Packaged-structure axioms and models]
	\label{DEF:PackagingStructureAxioms}
	\label{DEF:GeneralPackaging}	
	Fix a structural domain \(\mathscr C\).
	A candidate package-side datum over \(\mathscr C\) is a five-tuple
	\[
	\mathfrak P
	=
	\bigl(
	\operatorname{Pkg}_{\mathfrak P},
	E_{\mathfrak P},
	\triangleleft_{\mathfrak P},
	\mathsf L_{\mathfrak P},
	\mathsf R_{\mathfrak P}
	\bigr).
	\]
	The \textbf{packaged-structure axiom system}
	$
	\mathsf{Ax}_{\mathrm{PkgStr}}(\mathscr C)
	$
	consists of the following five conditions:
	\begin{enumerate}[label=(S\arabic*),leftmargin=3em]
		\item \textbf{Carrier axiom.}
		\[
		\operatorname{Pkg}_{\mathfrak P}
		\subseteq\mathscr C.
		\]
		The members of \(\operatorname{Pkg}_{\mathfrak P}\)	are the package carriers.
		
		\item \textbf{Internal-unit axiom.}
		\[
		E_{\mathfrak P}\in\mathcal U.
		\]
		Thus the package-side internal-unit set is \(\mathcal U\)-small.
		
		\item \textbf{Incidence axiom.}
		\[
		\triangleleft_{\mathfrak P}
		\subseteq
		E_{\mathfrak P}
		\times\operatorname{Pkg}_{\mathfrak P}.
		\]
		For \(K\in\operatorname{Pkg}_{\mathfrak P}\), its interior is
		\[
		I_K
		=
		\operatorname{Int}_{\mathfrak P}(K)
		:=
		\{e\in E_{\mathfrak P}:
		e\triangleleft_{\mathfrak P}K\}.
		\]
		
		\item \textbf{Local-candidate axiom.}
		\[
		\mathsf L_{\mathfrak P}
		=
		\bigl(
		\mathsf L_{\mathfrak P}(K)
		\bigr)_{K\in\operatorname{Pkg}_{\mathfrak P}}
		\]
		is a family of \(\mathcal U\)-small sets indexed by \(\operatorname{Pkg}_{\mathfrak P}\), where \(\mathsf L_{\mathfrak P}(K)\) is the declared local candidate space
		attached to \((K,I_K)\).
		
		\item \textbf{Selection axiom.}
		\[
		\mathsf R_{\mathfrak P}
		=
		\bigl(
		\mathsf R_{\mathfrak P}(K)
		\bigr)_{K\in\operatorname{Pkg}_{\mathfrak P}}
		\]
		is a family satisfying
		\[
		\mathsf R_{\mathfrak P}(K)
		\subseteq
		\mathsf L_{\mathfrak P}(K)
		\qquad
		(K\in\operatorname{Pkg}_{\mathfrak P}).
		\]
	\end{enumerate}
	\medskip
	A \textbf{packaged structure} is a model of \textup{(S1)}--\textup{(S5)}; write
	\[
	\mathbf{PkgDat}_{\mathscr C}=
	\operatorname{Mod}\bigl(\mathsf{Ax}_{\mathrm{PkgStr}}(\mathscr C)\bigr),
	\qquad \mathscr C_{\mathfrak P}:=\mathscr C.
	\]
	The ambient domain $\mathscr C$ is part of the type, not a sixth coordinate.
	The same five-tuple with different ambient domains is regarded as differently typed (equivalently, domain-tagged).
\end{definition}

\begin{definition}[Package and internal unit]
	\label{DEF:PackageInternalUnit}
	Let \(\mathfrak P\) be a packaged structure.
	A \textbf{package} is a carrier \(K\in\operatorname{Pkg}_{\mathfrak P}\), together with its interior \(I_K\), local candidate space
	\(\mathsf L_{\mathfrak P}(K)\),	and selected candidate family
	\(\mathsf R_{\mathfrak P}(K)\).	
	An element \(e\in E_{\mathfrak P}\)	is a \textbf{package-side internal unit}.
\end{definition}

\begin{remark}[Scope of the primitive axioms]
	The axioms permit empty package families, empty or overlapping interiors, uncovered internal units, and empty candidate or selected families.
	They impose no compatibility condition between packages and do not require a global completion.
	Such conditions may be introduced by specialization.
\end{remark}

\begin{lemma}[Size of the packaged-data space]
	\label{LEM:PackagingDataSize}
	Let \(\mathscr C\) be a \(\mathcal U\)-small structural domain.
	Then
	$
	\mathbf{PkgDat}_{\mathscr C}
	$
	is a \(\mathcal V\)-small set.
\end{lemma}

\begin{proof}
	Because \(\mathscr C\in\mathcal U\), the collection
	$
	\mathcal P(\mathscr C)
	$
	of package-carrier sets is \(\mathcal U\)-small.
	For each \(\operatorname{Pkg}\subseteq\mathscr C\) and each \(E\in\mathcal U\), the possible incidence relations form the \(\mathcal U\)-small set
	\[
	\mathcal P(E\times\operatorname{Pkg}).
	\]
	Since \(\mathcal U\in\mathcal V\), the collection of all \(\mathcal U\)-small sets is \(\mathcal V\)-small.
	Thus the possible indexed families
	\[
	\bigl(\mathsf L(K)\bigr)_{K\in\operatorname{Pkg}},
	\qquad
	\bigl(\mathsf R(K)\bigr)_{K\in\operatorname{Pkg}},
	\qquad
	\mathsf R(K)\subseteq\mathsf L(K),
	\]
	form a \(\mathcal V\)-small collection.
	Combining the five coordinates gives a \(\mathcal V\)-small set of packaged structures over \(\mathscr C\).
\end{proof}

\subsubsection{Packaging context \texorpdfstring{$\Xi$}{Xi}}

Having fixed what counts as a packaged output, we specify the environment in which that output can be formed.
A packaging context declares a structural domain, admitted carriers, formation data, and a primitive realization relation.
Keeping this relation separate from structural compatibility is what allows partial and multivalued formation in the ensuing system axioms.

\begin{definition}[Packaging context $\Xi$]
	\label{DEF:PackagingContext}	
	A \textbf{packaging context} is a quadruple
	\[
	\Xi
	=
	\bigl(
	\mathscr C_\Xi,
	\operatorname{Car}_\Xi,
	\operatorname{Form}_\Xi,
	\mathsf J_\Xi
	\bigr)
	\]
	that satisfies the following conditions:	
	\begin{enumerate}[label=(P\arabic*),leftmargin=3em]	
		\item $\mathscr C_\Xi$ is a structural domain,
		and
		\[
		\operatorname{Car}_\Xi
		\subseteq
		\mathscr C_\Xi
		\]
		is the set of admissible carriers.
		
		\item
		\[
		\operatorname{Form}_\Xi\in\mathcal V
		\]
		is a set, called the \textbf{formation-data space}.
		
		\item With
		\[
		\mathbf{PkgDat}_\Xi
		=
		\left\{
		\mathfrak Q
		\in
		\mathbf{PkgDat}_{\mathscr C_\Xi}:
		\operatorname{Pkg}_{\mathfrak Q}
		\subseteq
		\operatorname{Car}_\Xi
		\right\},
		\]
		the \textbf{realization relation} is
		\[
		\mathsf J_\Xi
		\subseteq
		\operatorname{Form}_\Xi
		\times
		\mathbf{PkgDat}_\Xi.
		\]		
	\end{enumerate}
\end{definition}

\subsection{Foundational compatibility}
\label{SUBSEC:PackagingCompatibilityRealization}

An output can be allowed by a context without being realized by any formation datum.
We call the first property \emph{compatibility} and define it here.
Compatibility is necessary for realization but does not imply it.
This use of the word is different from agreement of local data on chart overlaps, introduced later.

\begin{definition}[Compatibility of context and packaged structure]
	\label{DEF:PackagingContextCompatibility}	
	A packaging context $\Xi$ and a packaged structure $\mathfrak P$ are \textbf{compatible} if
	\[
	\mathscr C_{\mathfrak P}
	=
	\mathscr C_\Xi,
	\qquad
	\operatorname{Pkg}_{\mathfrak P}
	\subseteq
	\operatorname{Car}_\Xi.
	\]
	Equivalently,
	\[
	\mathfrak P\in\mathbf{PkgDat}_\Xi.
	\]
\end{definition}

\begin{remark}[Compatibility and realization]	
	For every packaging context $\Xi$,
	\[
	(P,\mathfrak P)\in\mathsf J_\Xi
	\quad\Longrightarrow\quad
	P\in\operatorname{Form}_\Xi
	\quad\text{and}\quad
	\mathfrak P\in\mathbf{PkgDat}_\Xi.
	\]	
	The converse does not hold in general.
	In particular, a context may have a nonempty formation-data space and compatible packaged structures while
	\[
	\mathsf J_\Xi=\varnothing.
	\]	
	Thus compatibility is necessary for realization, but does not imply realizability.
	The realization relation supplies the additional condition distinguishing admissible pairs from realized formations.
\end{remark}

\subsection{Formation datum \texorpdfstring{$P$}{P} and realization}
\label{SUBSEC:Packaging}

A formation datum realizes an admissible endpoint precisely when the pair belongs to the context's realization relation.
The relation need not be functional, allowing deterministic and multivalued formation as distinct cases.

\begin{definition}[Formation datum]
	\label{DEF:PackagingFormationDatum}
	Let $\Xi$ be a packaging context.
	A \textbf{formation datum} over $\Xi$ is an element
	$
	P\in\operatorname{Form}_\Xi.
	$
\end{definition}

\begin{definition}[Realization judgement]
	\label{DEF:PackagingRealizationJudgement}
	Let \(P\in\operatorname{Form}_\Xi\) and let
	$
	\mathfrak P\in\mathbf{PkgDat}_\Xi.
	$
	\begin{itemize}
		\item The \textbf{realization judgement}
		\[
		{
			\Xi\vdash P\rightsquigarrow\mathfrak P
		}
		\]
		means
		$
			(P,\mathfrak P)\in\mathsf J_\Xi.
		$	
		The abbreviated notation
		$
		P:\Xi\rightsquigarrow\mathfrak P
		$
		has the same meaning.
		
		\item A pair \((P,\mathfrak P)\) satisfying this relation is a \textbf{realized formation}.
		A formation datum having no realized endpoint does not by itself constitute a packaging system.
	\end{itemize}
\end{definition}

\subsection{Packaging system \texorpdfstring{$\mathbb P$}{P}}

The preceding definitions are assembled into the realized packaging system \(\mathbb P=(\Xi,P,\mathfrak P)\).
Its axioms require both contextual compatibility and an actual realization judgement, so a merely admissible pair is not yet a realized system.
Retaining all three coordinates makes source--endpoint comparisons and different native formation laws mathematically expressible within one formalism.

\begin{definition}[Packaging-system axioms and models]
	\label{DEF:PackagingSystemAxioms}
	\label{DEF:PackagingSystem}
	The \textbf{packaging-system axiom system}, denoted
	$
	\mathsf{Ax}_{\mathrm{PkgSys}},
	$
	consists of the following conditions on a triple
	\[
	\mathbb P=(\Xi,P,\mathfrak P).
	\]
	\begin{enumerate}[label=(Y\arabic*),leftmargin=3em]
		\item \textbf{Context axiom.}
		\(\Xi\) is a packaging context in the sense of Definition~\ref{DEF:PackagingContext}.
		
		\item \textbf{Formation axiom.}
		$
		P\in\operatorname{Form}_\Xi.
		$
		
		\item \textbf{Endpoint axiom.}
		\(\mathfrak P\) is a model of
		$
		\mathsf{Ax}_{\mathrm{PkgStr}}(\mathscr C_\Xi)
		$
		and is compatible with the context, equivalently
		$
		\mathfrak P\in\mathbf{PkgDat}_\Xi.
		$
		
		\item \textbf{Realization axiom.}
		$
		(P,\mathfrak P)\in\mathsf J_\Xi,
		$
		equivalently
		$
		\Xi\vdash P\rightsquigarrow\mathfrak P.
		$
	\end{enumerate}
	A model of \(\mathsf{Ax}_{\mathrm{PkgSys}}\) is called a realized packaging system.
\end{definition}

\begin{definition}[Outcome families and realizability]
	\label{DEF:PackagingGeneralRealizability}
	For a formation datum and a compatible target, respectively, put
	\[
	\operatorname{Out}_\Xi(P)=\{\mathfrak Q:(P,\mathfrak Q)\in\mathsf J_\Xi\},
	\qquad
	\operatorname{FormFib}_\Xi(\mathfrak Q)
	=\{P:(P,\mathfrak Q)\in\mathsf J_\Xi\}.
	\]
	The \textbf{realized-output family} is
	\[
	\operatorname{Real}_\Xi
	=\{\mathfrak Q\in\mathbf{PkgDat}_\Xi:
	\operatorname{FormFib}_\Xi(\mathfrak Q)\ne\varnothing\}.
	\]
	The context is \textbf{deterministic} if every outcome family has at most one member, and \textbf{total on formation data} if every such family is nonempty.
	These are additional properties, not universal axioms.
\end{definition}

\subsection{Primitive theory and specialization principles}
\label{SUBSEC:PackagingPrimitiveAxioms}

The primitive theory combines the packaged-structure and packaging-system axioms.
An \emph{axiomatic specialization} imposes additional conditions on formation laws or local data, separating general consequences from those specific to a chosen class.

\begin{definition}[Primitive axiom system of packaging]
	\label{DEF:PackagingPrimitiveAxioms}
	The primitive theory \(\mathsf{Ax}_{\mathrm{Pkg}}\) has two parts.
	Axioms~\textup{(S1)}--\textup{(S5)} define packaged outputs.
	Axioms~\textup{(Y1)}--\textup{(Y4)} define realized systems \((\Xi,P,\mathfrak P)\) built from those outputs.
	A triple is a realized model exactly when it satisfies both parts.
\end{definition}

\begin{definition}[Axiomatic specialization of packaging]
	\label{DEF:PackagingAxiomaticSpecialization}
	Let \(\Sigma\) be any additional family of mathematically specified conditions on packaging contexts, formation data, realization relations, packaged structures, or their interaction.
	We say that:
	\begin{itemize}
		\item The theory
		\[
		\mathsf{Ax}_{\mathrm{Pkg}}+\Sigma
		\]
		is an \textbf{axiomatic specialization of packaging}.
		
		\item It is a \textbf{proper specialization} when there exists a realized packaging system satisfying \(\mathsf{Ax}_{\mathrm{Pkg}}\) that does not satisfy \(\Sigma\).
	\end{itemize}
\end{definition}

\begin{remark}[Primitive non-axioms]
	The primitive theory does not require a total or single-valued realization relation.
	It does not require formation data to compose, a designated source map, or formation by a map, quotient, kernel, functor, or channel.
	It also assumes no information ordering, common structure on candidate spaces, or global reconstruction.	
	Any such condition is additional structure and therefore belongs to a specialization rather than to the definition of general packaging.
\end{remark}

\begin{example}[Map-based, quotient, and functorial specializations]
	\label{EX:PackagingClassicalSpecializations}
	The primitive packaging axioms admit the following specializations.	
	\begin{enumerate}
		\item \emph{Map-based realization.}
		Let \(X,Y\) be \(\mathcal U\)-small sets and let \(f:X\to Y\) be a map.
		Choose the structural domain and admitted carriers to be \(Y\).
		For each \(y\in Y\), let \(\mathfrak P_y\) be the one-carrier packaged structure with carrier \(y\), empty internal-unit set, and empty candidate and selected families.
		These structures satisfy the primitive axioms and are pairwise distinct.
		Taking
		\[
		\operatorname{Form}_\Xi=X,
		\qquad
		\mathsf J_\Xi=\{(x,\mathfrak P_{f(x)}):x\in X\}
		\]
		gives a total deterministic realization relation.
		For a partial map \(f:X\rightharpoonup Y\), use its graph in the same construction, allowing formation data outside its domain to have no realized endpoint.
		
		\item \emph{Quotient realization.}
		Let \(X\) be \(\mathcal U\)-small and \(\sim\) an equivalence relation on \(X\).
		Take the formation datum \(P=(X,\sim)\) and the quotient map
		\[
		q_\sim:X\twoheadrightarrow X/{\sim}.
		\]
		Choose a structural domain admitting the quotient map \(q_\sim\) as a carrier, and a packaged endpoint \(\mathfrak P_{q_\sim}\) containing that carrier.
		For example, a singleton formation-data space and the realization relation
		\[
		\mathsf J_\Xi=\{((X,\sim),\mathfrak P_{q_\sim})\}
		\]
		define such a context.
		Requiring that \(q_\sim\) satisfy the quotient universal property specifies a quotient specialization.
		For an algebraic quotient, the equivalence relation must also be an appropriate congruence.
		
		\item \emph{Functorial realization.}
		Let \(\mathcal C,\mathcal D\) be \(\mathcal U\)-small categories and
		\(F_0:\operatorname{Ob}(\mathcal C)\to\operatorname{Ob}(\mathcal D)\) an object assignment.
		The map-based construction represents \(F_0\) by a deterministic realization relation.
		A functorial specialization additionally specifies an action \(F_1\) on morphisms, compatible with sources and targets, and imposes
		\[
		F_1(\mathrm{id}_X)=\mathrm{id}_{F_0(X)},
		\qquad
		F_1(g\circ f)=F_1(g)\circ F_1(f).
		\]
		These conditions make \((F_0,F_1)\) a functor.
		They require additional categorical structure or formation data, rather than following from deterministic realization alone.
	\end{enumerate}
\end{example}

\subsection{Selection geometry and freedom of realization}
\label{SUBSEC:PackagingSelectionFreedom}

Fix the carriers, internal units, incidence, and candidate spaces.
The remaining choice is which candidates to select in each package.
Theorem~\ref{THM:PackagingSelectionFreedom} shows that these choices form a Boolean algebra.
But the axioms place no corresponding order rule on which selections a formation datum can realize.
The theorem therefore describes both a structural consequence and a limit of the primitive axioms.

\begin{definition}[Packaging scaffold and selection order]
	\label{DEF:PackagingSelectionScaffold}	
	Fix a structural domain $\mathscr C$.
	A \emph{packaging scaffold} is a quadruple
	\[
	\mathcal D=(\mathcal K,E,\triangleleft,\mathsf L)
	\]
	satisfying axioms \textup{(S1)}--\textup{(S4)} of
	Definition~\ref{DEF:PackagingStructureAxioms},
	with the selected-family coordinate omitted.	
	Write
	\[
	\operatorname{Sel}(\mathcal D)
	=
	\prod_{K\in\mathcal K}\mathcal P(\mathsf L(K)).
	\]	
	For $R\in\operatorname{Sel}(\mathcal D)$,
	denote the corresponding packaged structure by
	\[
	\mathfrak P(\mathcal D,R)
	=
	(\mathcal K,E,\triangleleft,\mathsf L,R).
	\]	
	The \emph{selection order} is defined by
	\[
	R\le_{\mathrm{sel}}T
	\quad\Longleftrightarrow\quad
	R(K)\subseteq T(K)
	\quad(K\in\mathcal K).
	\]	
	This order compares selections, not refinement candidates inside $\mathsf L(K)$.
\end{definition}

\begin{theorem}[Selection algebra and unrestricted realization]
	\label{THM:PackagingSelectionFreedom}	
	Let $\mathcal D$ be any packaging scaffold over
	$\mathscr C$, and put
	\[
	\Omega_{\mathcal D}
	=
	\bigsqcup_{K\in\mathcal K}
	\bigl(\{K\}\times\mathsf L(K)\bigr).
	\]	
	Then the following hold:	
	\begin{enumerate}[label=(\roman*),leftmargin=2.5em]	
		\item The selection order is canonically a complete Boolean lattice:
		\[
		\operatorname{Sel}(\mathcal D)
		\cong
		\mathcal P(\Omega_{\mathcal D}),
		\]
		under the correspondence
		\[
		R\longmapsto
		\{(K,\ell):\ell\in R(K)\}.
		\]		
		Its meets, joins, and complements are computed
		coordinatewise by intersection, union, and complementation in $\mathsf L(K)$.		
		When $\Omega_{\mathcal D}\ne\varnothing$, its atoms are the single-candidate selections.
		
		\item Let $F\in\mathcal V$ be any set of formation data, and let
		$
		\mathcal T
		\subseteq
		F\times\operatorname{Sel}(\mathcal D)
		$
		be any relation.
		
		There exists a packaging context $\Xi_{\mathcal T}$ with formation-data space $F$ and realization relation exactly
		\[
		\mathsf J_{\Xi_{\mathcal T}}
		=
		\left\{
		(p,\mathfrak P(\mathcal D,R)):
		(p,R)\in\mathcal T
		\right\}.
		\]		
		Consequently,
		\[
		\operatorname{Out}_{\Xi_{\mathcal T}}(p)
		=
		\left\{
		\mathfrak P(\mathcal D,R):
		(p,R)\in\mathcal T
		\right\}
		\qquad(p\in F).
		\]
	
		\item If $\Omega_{\mathcal D}$ contains at least two elements, there is a context and a formation datum whose realized selections are closed under neither binary meets nor binary joins in $\operatorname{Sel}(\mathcal D)$.	
	\end{enumerate}	
	Thus, although $\operatorname{Sel}(\mathcal D)$ is a complete Boolean lattice, the primitive packaging axioms impose neither lattice-closure requirements on realized selections nor any order-monotonicity requirement on the realization relation.
\end{theorem}

\begin{proof}
	\leavevmode
	\begin{enumerate}[label=(\roman*),leftmargin=2.5em]
		\item A selected family $R$ specifies exactly the subset
		\[
		\{(K,\ell):\ell\in R(K)\}
		\]
		of the tagged disjoint union.	
		Conversely, any subset of that union determines	one and only one subset of each fibre $\mathsf L(K)$.	
		These operations are inverse and preserve inclusion.
		All stated lattice operations and the characterization of atoms follow from the power-set description.
	
		\item Define
		\[
		\mathscr C_{\Xi_{\mathcal T}}=\mathscr C,
		\qquad
		\operatorname{Car}_{\Xi_{\mathcal T}}=\mathcal K,
		\qquad
		\operatorname{Form}_{\Xi_{\mathcal T}}=F,
		\]
		and give the context the displayed realization relation.	
		Every $\mathfrak P(\mathcal D,R)$ satisfies
		\textup{(S1)}--\textup{(S5)} and is compatible with these admitted carriers.	
		The size conditions follow from the universe convention.	
		Hence $\Xi_{\mathcal T}$ is a packaging context, and the formula for its outcome families is immediate.	
		If $\mathcal T$ is nonempty, every related pair	gives a realized packaging system.
		If it is empty, the context legitimately has no	realized systems.
	
		\item Choose two distinct elements $u,v$ of $\Omega_{\mathcal D}$, and let $R_u,R_v$ correspond under \textup{(i)} to $\{u\}$ and $\{v\}$.
	
		Apply \textup{(ii)} with
		\[
		F=\{p\},
		\qquad
		\mathcal T=
		\{(p,R_u),(p,R_v)\}.
		\]	
		The two realized selections have meet corresponding	to $\varnothing$ and join corresponding to $\{u,v\}$, neither of which is realized by $p$.	
		Their realized set is also neither an up-set nor a down-set.	
	\end{enumerate}
	This proves the assertion.
\end{proof}

\begin{remark}[Scope of the primitive packaging axioms]
	\label{REM:PackagingFoundationalStrength}
	\leavevmode
	\begin{itemize}[leftmargin=2.1em]		
		\item Theorem~\ref{THM:PackagingSelectionFreedom}
		characterizes the freedom permitted by the primitive
		packaging axioms on a fixed scaffold.
		It is not a globalization theorem.
		
		\item The complete Boolean lattice
		$\operatorname{Sel}(\mathcal D)$ arises from independent	choices of selected candidates.
		Nevertheless, realization relations remain unrestricted, and realized-selection families need not be closed under meets or joins.
		
		\item The primitive axioms impose neither overlap
		compatibility nor exact-gluing requirements.
		They also prescribe no intrinsic direction of
		information loss or gain
		(Proposition~\ref{THM:PackagingRecoverabilityNeutrality}).
		Compatibility and globalization results require
		additional hypotheses, as specified in the
		subsequent specializations.		
	\end{itemize}	
	The power-set and relational constructions underlying
	Theorem~\ref{THM:PackagingSelectionFreedom} are classical.
	Their role here is to characterize precisely the
	freedom allowed by the primitive packaging axioms.
\end{remark}

\subsection{Generality of the realization relation}
\label{SUBSEC:PackagingRealizationGenerality}

A formation datum need not have an outcome, and it may have more than one.
Proposition~\ref{THM:PackagingGeneralDeterminism} shows that deterministic and total deterministic formation are successively more restrictive cases.
Explicit examples prove both inclusions are strict.
This is a difference between relations and single-valued maps, not a claim that formation cannot be encoded by a relation.

\begin{proposition}[Determinism and totality as proper specializations]
	\label{THM:PackagingGeneralDeterminism}	
	Let $\Xi$ be a packaging context, and put
	\[
	\operatorname{Dom}(\mathsf J_\Xi)
	=
	\{P\in\operatorname{Form}_\Xi:
	\operatorname{Out}_\Xi(P)\ne\varnothing\}.
	\]	
	Then:	
	\begin{enumerate}[label=(\roman*),leftmargin=2.5em]		
		\item $\Xi$ is deterministic if and only if
		$\mathsf J_\Xi$ is the graph of a unique map
		\[
		\operatorname{ev}_\Xi:
		\operatorname{Dom}(\mathsf J_\Xi)
		\longrightarrow
		\mathbf{PkgDat}_\Xi.
		\]
		
		\item $\Xi$ is total and deterministic if and only if
		$\mathsf J_\Xi$ is the graph of a map
		\[
		\operatorname{ev}_\Xi:
		\operatorname{Form}_\Xi
		\longrightarrow
		\mathbf{PkgDat}_\Xi.
		\]
		
		\item Let
		$\mathbf{PkgCtx}_{\mathrm{gen}}$,
		$\mathbf{PkgCtx}_{\mathrm{det}}$, and
		$\mathbf{PkgCtx}_{\mathrm{tdet}}$
		denote the classes of all packaging contexts,
		deterministic contexts, and total deterministic
		contexts, respectively.
		Then
		\[
		\mathbf{PkgCtx}_{\mathrm{tdet}}
		\subsetneq
		\mathbf{PkgCtx}_{\mathrm{det}}
		\subsetneq
		\mathbf{PkgCtx}_{\mathrm{gen}}.
		\]
		The corresponding classes of realized packaging
		systems, classified by their contexts, satisfy
		the same strict inclusions.		
	\end{enumerate}
\end{proposition}

\begin{proof}
	We prove the three assertions separately.	
	\begin{enumerate}[label=(\roman*),leftmargin=2.5em]		
		\item \emph{Determinism.}
		By definition, $\Xi$ is deterministic precisely
		when every nonempty outcome family $\operatorname{Out}_\Xi(P)$ is a singleton.		
		Thus, each $P\in\operatorname{Dom}(\mathsf J_\Xi)$
		determines a unique endpoint $\operatorname{ev}_\Xi(P)$, and
		\[
		\mathsf J_\Xi
		=
		\left\{
		(P,\operatorname{ev}_\Xi(P)):
		P\in\operatorname{Dom}(\mathsf J_\Xi)
		\right\}.
		\]
		Conversely, the graph of a map is right-unique.
		
		\item \emph{Totality.}
		A deterministic context is total on formation data precisely when
		\[
		\operatorname{Dom}(\mathsf J_\Xi)
		=
		\operatorname{Form}_\Xi.
		\]
		Thus total deterministic realization is equivalent to evaluation by a map defined on the entire formation-data space.
		
		\item \emph{Strictness.}
		Choose a structural domain $\mathscr C=\{K\}$ and admit $K$ as a carrier.		
		Define two packaged structures $\mathfrak Q_0$ and $\mathfrak Q_1$ with the same package family $\{K\}$,
		internal-unit set $\{e\}$, incidence relation $\{(e,K)\}$, and local candidate space $\{0\}$, but with selected candidate families
		\[
		\mathsf R_{\mathfrak Q_0}(K)=\varnothing,
		\qquad
		\mathsf R_{\mathfrak Q_1}(K)=\{0\}.
		\]
		Both satisfy the packaged-structure axioms, and $\mathfrak Q_0\ne\mathfrak Q_1$.
		
		Fix distinct formation data $p_0,p_1$.
		Consider three contexts with structural domain $\mathscr C$ and admitted-carrier set $\{K\}$, specified by
		\[
		\begin{array}{c|c|c}
			& \operatorname{Form} & \mathsf J\\
			\hline
			\Xi_{\mathrm{tdet}}
			& \{p_0\}
			& \{(p_0,\mathfrak Q_0)\}\\
			\Xi_{\mathrm{det}}
			& \{p_0,p_1\}
			& \{(p_0,\mathfrak Q_0)\}\\
			\Xi_{\mathrm{gen}}
			& \{p_0\}
			& \{(p_0,\mathfrak Q_0),
			(p_0,\mathfrak Q_1)\}
		\end{array}
		\]		
		The first context is total and deterministic.
		The second is deterministic but not total, whereas the third is nondeterministic.		
		Each context has a realized formation.
		Hence both inclusions are strict for packaging contexts and, by classifying	realized packaging systems according to their contexts, also for realized systems.		
	\end{enumerate}
\end{proof}

\subsection{Refinement specialization}
\label{SUBSEC:PackagingRefinementSpecialization}

So far, local candidate spaces may be arbitrary sets.
We now take them to be refinements of package interiors, ordered by how finely they distinguish units.
Equivalence relations give canonical representatives of these refinements, and descriptor maps provide another description.
Proposition~\ref{PROP:PackagingRefinementProperSpecialization} shows that this condition excludes some general packaged structures.
The completion theory of the next section applies to this narrower class.
Equivalence relations and partition lattices used here are classical \cite{Ore1942Equivalence,Birkhoff1967Lattice,DaveyPriestley2002Lattices}.

\begin{definition}[Normalized labeled refinement]
	\label{DEF:NormalizedLabeledRefinement}
	A \textbf{normalized labeled refinement} of a set \(I\) is a surjection \(\pi:I\twoheadrightarrow J\).
	Its distinction relation is	\(\operatorname{Dist}(\pi)=\{(x,y)\in I^2:\pi(x)\ne\pi(y)\}\).
	It is \textbf{trivial} when \(|J|\le1\), and \textbf{nontrivial} otherwise.
\end{definition}

\begin{definition}[Relabeling equivalence]
	\label{DEF:RelabelingEquivalenceRefinements}
	Refinements \(\pi_i:I\twoheadrightarrow J_i\) are \textbf{relabeling-equivalent}, written \(\pi_1\cong\pi_2\), when \(\pi_2=b\circ\pi_1\) for a bijection \(b:J_1\to J_2\).
	Write \([\pi]\) for its canonical representative, specified below.
\end{definition}

\begin{definition}[Refinement order]
	\label{DEF:RefinementOrderPackaging}
	The refinement poset \(\mathbf{Ref}(I)\) consists of the canonical representatives of these classes, with \([\pi_0]\preceq[\pi_1]\) when \(\pi_0=h\circ\pi_1\) for some map \(h\).
	Thus \(\pi_1\) is finer and \(\pi_0\) is a coarsening.
\end{definition}

\medskip\noindent\textbf{Canonical representatives and size.}
For set-theoretic purposes, put
\[
\operatorname{Eq}(I)=\{R\subseteq I^2:R\text{ is an equivalence relation}\},
\qquad q_R:I\twoheadrightarrow I/R,
\]
and implement
\[
\mathbf{Ref}(I)=\{q_R:R\in\operatorname{Eq}(I)\},
\qquad [\pi]=q_{\ker\pi}.
\]
Write
\[
\bot_I:=q_{I\times I}.
\]
This is the least element of \(\mathbf{Ref}(I)\).
It is the trivial refinement (with empty codomain when \(I=\varnothing\)).
Thus the notation \([\pi]\) designates a canonical representative, not the literal collection of all surjections obtained by relabeling its codomain.
This implements the preceding equivalence:
when two normalized maps have equal kernels, \(b(\pi_0(x))=\pi_1(x)\) defines a well-defined bijection between their codomains, uniquely determined by surjectivity.
Conversely, bijective relabeling preserves the kernel.
These arguments also cover \(I=\varnothing\), whose unique normalized refinement has empty codomain.

Since \(\operatorname{Eq}(I)\subseteq\mathcal P(I^2)\) is \(\mathcal U\)-small, so is the displayed family of canonical quotient maps.
Moreover,
\[
q_R\preceq q_T\quad\Longleftrightarrow\quad T\subseteq R:
\]
if \(T\subseteq R\), the map \([x]_T\mapsto[x]_R\) gives the required factorization.
The converse follows by comparing kernels of a factorization.
This convention gives the stated refinement order and justifies the transport and pullback formulas used below.
This convention codes only local refinement data.
It does not identify a packaging formation datum, context, or system with quotient formation.

\medskip
\noindent\textbf{Canonical refinement represented by an arbitrary map.}
If
$
f:I\longrightarrow Y
$
is any map, let
$
\bar f:I\twoheadrightarrow f(I)
$
denote its corestriction to its attained image, and define
\[
[f]
:=
[\bar f]
=
q_{\ker f}
\in\mathbf{Ref}(I).
\]
Thus the refinement represented by a map depends only on its kernel and is unchanged by unused codomain labels.
In particular, the notation \([v_a]\) and \([D_{\mathcal A}]\) used for descriptor maps below is well defined even when those maps are not surjective.

\begin{definition}[Descriptor family]
	\label{DEF:DescriptorFamilyPackaging}
	A \textbf{descriptor family} on a carrier/set pair \((K,I)\) is an indexed family \(\mathcal A=(f_i:I\to Y_i)_{i\in A}\), with \(A,Y_i\in\mathcal U\).
	Its joint descriptor is
	\[
	D_{\mathcal A}:I\longrightarrow\prod_{i\in A}Y_i,
	\qquad D_{\mathcal A}(x)=(f_i(x))_{i\in A}.
	\]
	Restricting to \(B\subseteq A\) gives a \textbf{subfamily}.
	The empty family has the unique map \(D_{\varnothing}:I\to\{*\}\).
\end{definition}

\begin{definition}[Refinement-packaged structure]
	\label{DEF:PackagingRefinementPackagedStructure}	
	A packaged structure $\mathfrak P$ over
	$\mathscr C$ is a
	\textbf{refinement-packaged structure} if
	\[
	\mathsf L_{\mathfrak P}(K)
	=
	\mathbf{Ref}(I_K)
	\qquad
	(K\in\operatorname{Pkg}_{\mathfrak P}).
	\]	
	Denote the class of such structures by
	\[
	\mathbf{PkgDat}^{\mathrm{ref}}_{\mathscr C}
	\subseteq
	\mathbf{PkgDat}_{\mathscr C}.
	\]
\end{definition}

\noindent\textbf{Refinement conventions.}
The local candidate spaces inherit the refinement order $\preceq$ and the least elements
$
\bot_{I_K}=q_{I_K\times I_K}.
$
The selected families
\[
\mathsf R_{\mathfrak P}(K)
\subseteq
\mathbf{Ref}(I_K)
\]
are arbitrary subsets.
They need not be down-sets or contain $\bot_{I_K}$.
Since the local candidate spaces are determined by the incidence relation, we may suppress $\mathsf L_{\mathfrak P}$ and use the abbreviated notation
\[
\mathfrak P
=
\bigl(
\operatorname{Pkg}_{\mathfrak P},
E_{\mathfrak P},
\triangleleft_{\mathfrak P},
\mathsf R_{\mathfrak P}
\bigr)
\]
when no ambiguity arises.

\begin{proposition}[Refinement packaging is a proper specialization]
	\label{PROP:PackagingRefinementProperSpecialization}
	For every structural domain \(\mathscr C\),
	\[
	\mathbf{PkgDat}^{\mathrm{ref}}_{\mathscr C}
	\subseteq
	\mathbf{PkgDat}_{\mathscr C}.
	\]
	If \(\mathscr C\neq\varnothing\), the inclusion is strict:
	\[
	{
		\mathbf{PkgDat}^{\mathrm{ref}}_{\mathscr C}
		\subsetneq
		\mathbf{PkgDat}_{\mathscr C}.
	}
	\]
\end{proposition}

\begin{proof}
	The inclusion is immediate from	Definition~\ref{DEF:PackagingRefinementPackagedStructure}.	
	Suppose \(\mathscr C\neq\varnothing\), and choose \(K\in\mathscr C\).
	Let
	\[
	\operatorname{Pkg}_{\mathfrak P}=\{K\},
	\qquad
	E_{\mathfrak P}=\{0,1\},
	\]
	with both internal units incident with \(K\), so that
	\[
	I_K=\{0,1\}.
	\]
	Choose a three-element local candidate space
	\[
	\mathsf L_{\mathfrak P}(K)=\{a,b,c\},
	\qquad
	\mathsf R_{\mathfrak P}(K)=\{a,b,c\}.
	\]
	This is a packaged structure by Definition~\ref{DEF:GeneralPackaging}.
	
	On the other hand, \(\mathbf{Ref}(\{0,1\})\) has exactly two members, namely the indiscrete and discrete refinements.
	Hence
	\[
	\mathsf L_{\mathfrak P}(K)
	\neq
	\mathbf{Ref}(I_K),
	\]
	so \(\mathfrak P\) is not refinement-packaged.
\end{proof}

\begin{definition}[Refinement-compatible endpoint]
	\label{DEF:PackagingRefinementCompatibleEndpoint}
	A compatible pair \(\Xi,\mathfrak P\) is \textbf{refinement-compatible} when
	\[
	\mathfrak P
	\in
	\mathbf{PkgDat}^{\mathrm{ref}}_{\mathscr C_\Xi}.
	\]
	A packaging system
	$
	(\Xi,P,\mathfrak P)
	$
	is a \textbf{refinement-packaging system} when its realized endpoint \(\mathfrak P\) is refinement-packaged.
\end{definition}

\subsection{Source--endpoint recoverability and comparison neutrality}
\label{SUBSEC:PackagingRecoverabilityNeutrality}

A realized formation links a formation datum to an output, but it does not say which one contains more information.
To ask whether one determines the other, we first specify how source data and endpoints are compared.
Proposition~\ref{THM:PackagingRecoverabilityNeutrality} shows that the primitive axioms imply neither direction of recoverability, even when both comparisons use equality.

\begin{definition}[Comparison relations on realized systems]
	\label{DEF:PackagingRealizedSystemComparisonSpace}	
	Let $\Xi$ be a packaging context.
	Identifying a realized packaging system
	$(\Xi,P,\mathfrak P)$ with its realized pair
	$(P,\mathfrak P)$, put
	\[
	\operatorname{Sys}(\Xi):=\mathsf J_\Xi.
	\]	
	A \textbf{source comparison relation} and an
	\textbf{endpoint comparison relation} are specified binary relations
	\[
	\mathcal S_\Xi
	\subseteq
	\operatorname{Form}_\Xi^2,
	\qquad
	\mathcal E_\Xi
	\subseteq
	\mathbf{PkgDat}_\Xi^2,
	\]
	respectively.	
	They induce relations $\widehat{\mathcal S}_\Xi$ and
	$\widehat{\mathcal E}_\Xi$ on $\operatorname{Sys}(\Xi)$ by
	\[
	\begin{aligned}
		(P_0,\mathfrak P_0)
		\mathrel{\widehat{\mathcal S}_\Xi}
		(P_1,\mathfrak P_1)
		&\iff
		P_0\mathrel{\mathcal S_\Xi}P_1,
		\\[4pt]
		(P_0,\mathfrak P_0)
		\mathrel{\widehat{\mathcal E}_\Xi}
		(P_1,\mathfrak P_1)
		&\iff
		\mathfrak P_0
		\mathrel{\mathcal E_\Xi}
		\mathfrak P_1,
	\end{aligned}
	\]
	for all	$(P_i,\mathfrak P_i)\in\operatorname{Sys}(\Xi)$,
	$i=0,1$.	
\end{definition}

The comparison relations need not be equivalence relations. Here ``source comparison'' refers to a declared comparison of formation data and does not require a distinguished source coordinate.

\begin{definition}[Comparison-relative recoverability]
	\label{DEF:PackagingComparisonRelativeRecoverability}
	In the setting of Definition~\ref{DEF:PackagingRealizedSystemComparisonSpace}:	
	\begin{enumerate}[label=(\roman*),leftmargin=2.5em]		
		\item The declared source distinction is \textbf{recoverable from the endpoint comparison} when
		\[
		{
			\widehat{\mathcal E}_{\Xi}
			\subseteq
			\widehat{\mathcal S}_{\Xi}.
		}
		\]
		Equivalently, if two realized systems are endpoint-indistinguishable, then they are also source-indistinguishable.
		
		\item The declared endpoint distinction is \textbf{determined by the source comparison} when
		\[
		{
			\widehat{\mathcal S}_{\Xi}
			\subseteq
			\widehat{\mathcal E}_{\Xi}.
		}
		\]
		
		\item Source and endpoint comparisons are \textbf{comparison-equivalent} on realized systems when
		\[
		{
			\widehat{\mathcal S}_{\Xi}
			=
			\widehat{\mathcal E}_{\Xi}.
		}
		\]
		
		\item There is a \textbf{source-distinction collapse at the endpoint}
		when
		\[
		\widehat{\mathcal S}_{\Xi}
		\subsetneq
		\widehat{\mathcal E}_{\Xi}.
		\]
		
		\item There is \textbf{endpoint-side extra discrimination} when
		\[
		\widehat{\mathcal E}_{\Xi}
		\subsetneq
		\widehat{\mathcal S}_{\Xi}.
		\]
		
		\item The two comparisons are \textbf{incomparable} when neither relation contains the other.
	\end{enumerate}	
	The terms ``collapse'' and ``extra discrimination'' are relative to the declared comparison relations.
	They do not impose an intrinsic numerical, entropic, dimensional, or probabilistic notion of information.
\end{definition}

\begin{proposition}[Source--endpoint comparison neutrality]
	\label{THM:PackagingRecoverabilityNeutrality}
	The primitive packaging axioms impose no universal inclusion between
	\[
	\widehat{\mathcal S}_{\Xi}
	\qquad\text{and}\qquad
	\widehat{\mathcal E}_{\Xi}.
	\]
	Even when the underlying source and endpoint comparisons are equality relations, packaging contexts exist realizing each of the four regimes
	\[
	\widehat{\mathcal S}_{\Xi}
	=
	\widehat{\mathcal E}_{\Xi},
	\]
	\[
	\widehat{\mathcal S}_{\Xi}
	\subsetneq
	\widehat{\mathcal E}_{\Xi},
	\]
	\[
	\widehat{\mathcal E}_{\Xi}
	\subsetneq
	\widehat{\mathcal S}_{\Xi},
	\]
	and
	\[
	\widehat{\mathcal S}_{\Xi},
	\widehat{\mathcal E}_{\Xi}
	\quad\text{incomparable}.
	\]	
	Consequently, general packaging implies neither source-distinction collapse nor endpoint-side extra discrimination.
	In particular, packaging has no intrinsic information-loss or information-gain direction.
\end{proposition}

\begin{proof}
	Take distinct compatible endpoints $\mathfrak Q_0,\mathfrak Q_1$ and use equality as both source and endpoint comparison.
	On realized systems, the relations
	\[
	\{(p_0,\mathfrak Q_0),(p_1,\mathfrak Q_1)\},\quad
	\{(p_0,\mathfrak Q_0),(p_1,\mathfrak Q_0)\},\quad
	\{(p,\mathfrak Q_0),(p,\mathfrak Q_1)\}
	\]
	realize respectively equality, source-distinction collapse, and endpoint-side extra discrimination.
	Adding $(r,\mathfrak Q_0)$ to the third construction, with $r\ne p$, gives one pair related only by the source comparison and another related only by the endpoint comparison, hence incomparability.
\end{proof}

\subsection{Comparison-based structural non-collapse}
\label{SUBSEC:PackagingComparisonNoncollapse}

A reduced description may treat two packaged objects as the same even though a chosen invariant distinguishes them.
We call such a pair a \emph{separation witness}.
Lemma~\ref{THM:PackagingNativeNoncollapse} gives a precise test for these witnesses.
The conclusion concerns the chosen reduction; it does not rule out faithful descriptions of packaging in other mathematical languages.

\begin{definition}[Packaging comparison and separation witness]
	\label{DEF:PackagingComparisonRelation}	
	Let $\mathcal X$ be a specified collection of packaging objects.
	\begin{itemize}
		\item A \textbf{packaging comparison} on $\mathcal X$
		is a specified binary relation
		$
		\mathcal C\subseteq\mathcal X\times\mathcal X.
		$	
		
		\item Let
		$
		J:\mathcal X\longrightarrow\mathcal Z
		$
		be a retained datum, and let $\simeq_{\mathcal Z}$ be an equivalence relation on $\mathcal Z$.	
		Define the \textbf{$J$-agreement relation}
		\[
		\operatorname{Agr}_J
		:=
		\left\{
		(X_0,X_1)\in\mathcal X^2:
		J(X_0)\simeq_{\mathcal Z}J(X_1)
		\right\}.
		\]	
		The comparison $\mathcal C$ \textbf{retains $J$}
		if
		$
		\mathcal C\subseteq\operatorname{Agr}_J.
		$	
		A \textbf{$J$-separation witness against $\mathcal C$} is a pair
		$
		(X_0,X_1)
		\in
		\mathcal C\setminus\operatorname{Agr}_J.
		$
	\end{itemize}	
\end{definition}

\medskip
\noindent\textbf{Relative faithfulness.}
If $\mathcal X$ carries a specified equivalence relation
$\simeq_{\mathcal X}$, a packaging comparison $\mathcal C$ is \textbf{faithful relative to $\simeq_{\mathcal X}$} if
$
\mathcal C\subseteq\simeq_{\mathcal X}.
$

\begin{lemma}[Retained-invariant non-collapse criterion]
	\label{THM:PackagingNativeNoncollapse}
	Let \(\mathcal C\) and \(J\) be as above.
	\begin{enumerate}[label=(\roman*),leftmargin=2.5em]
		\item A \(J\)-separation witness exists if and only if
		$
		\mathcal C\not\subseteq\operatorname{Agr}_J.
		$
		
		\item If
		$
		\mathcal C\subseteq\mathcal D,
		$
		then every \(J\)-separation witness against \(\mathcal C\) is also a separation witness against \(\mathcal D\).
		
		\item Suppose \(J\) is invariant under a declared equivalence \(\simeq_{\mathcal X}\), that is,
		\[
		X_0\simeq_{\mathcal X}X_1
		\quad\Longrightarrow\quad
		J(X_0)\simeq_{\mathcal Z}J(X_1).
		\]
		Then every \(J\)-separation witness	\(X_0,X_1\) satisfies
		$
		X_0\not\simeq_{\mathcal X}X_1.
		$
		In particular, a comparison relating such a pair is not faithful relative to \(\simeq_{\mathcal X}\).
	\end{enumerate}
\end{lemma}

\begin{proof}
	\leavevmode
	\begin{itemize}
		\item Part \textup{(i)} follows from the definition of \(\operatorname{Agr}_J\):
		a separation witness exists exactly when
		\(\mathcal C\not\subseteq\operatorname{Agr}_J\).
		
		\item For \textup{(ii)}, suppose \(\mathcal C\subseteq\mathcal D\).
		Every pair related by \(\mathcal C\) is	also related by \(\mathcal D\).
		Hence every separation witness against \(\mathcal C\) is also one against \(\mathcal D\).
		
		\item For \textup{(iii)}, suppose	\(X_0\simeq_{\mathcal X}X_1\).
		Invariance of \(J\) gives \(J(X_0)\simeq_{\mathcal Z}J(X_1)\), which contradicts the separation-witness condition.
		Thus \(X_0\not\simeq_{\mathcal X}X_1\), and any comparison relating this pair is not faithful.
	\end{itemize}	
\end{proof}

The criterion is deliberately relative.
It does not assert that packaging objects are intrinsically distinct merely because a chosen comparison fails to retain a chosen datum.
Nor must the comparison arise from a map, quotient, kernel, functor, channel, or any other particular representation mechanism.

\begin{remark}[Comparison-relative structural non-collapse]
	\label{REM:PackagingRepresentabilityScope}	
	Proposition~\ref{PROP:PackagingRelationalPresentation}
	gives a faithful many-sorted representation of the
	bare packaged-structure coordinates.
	
	Corollary~\ref{COR:PackagingIntrinsicCompletionNoncollapseA} identifies a completion invariant that distinguishes structures under a coarser incidence-and-local-type comparison.
	
	Proposition~\ref{PROP:PackagingAccessCommutantNonQuotient} provides a separate non-factorization witness
	for native access--commutant formation.
	
	Together these results give explicit and distinct
	criteria for comparison-relative non-collapse
	and for non-quotient formation.	
\end{remark}

\subsection{Reduction, reconstruction, and non-factorization}

For a reduction given by a map, reconstruction of retained data is equivalent to constancy on the map's fibres (Proposition~\ref{PROP:PackagingReductionFactorizationCriterion}).
The bare packaged coordinates also admit a faithful many-sorted relational presentation (Proposition~\ref{PROP:PackagingRelationalPresentation}).
Non-representability claims must therefore specify the class of reductions being excluded.

\begin{definition}[Reduction and reconstruction]
	\label{DEF:PackagingShadowReconstruction}	
	Let $\mathcal X$ be a collection of packaging objects.	
	A \textbf{reduction} is a specified map
	\[
	U:\mathcal X\longrightarrow\mathcal Y,
	\qquad
	\mathcal Y=U(\mathcal X).
	\]	
	Let
	\[
	J:\mathcal X\longrightarrow\mathcal Z
	\]
	be a retained datum.
	We distinguish two notions of reconstruction:	
	\begin{enumerate}[label=(\roman*),leftmargin=2.5em]		
		\item \textbf{Exact reconstruction.}
		The datum $J$ is \textbf{exactly reconstructible through $U$} if there exists a map
		\[
		\overline J:
		\mathcal Y\longrightarrow\mathcal Z
		\]
		such that
		\[
		J=\overline J\circ U.
		\]
		
		\item \textbf{Reconstruction up to equivalence.}
		Suppose $\mathcal Y$ and $\mathcal Z$ carry equivalence relations $\simeq_{\mathcal Y}$ and $\simeq_{\mathcal Z}$, with canonical quotient maps $q_{\mathcal Y}$ and $q_{\mathcal Z}$.
		The datum $J$ is \textbf{reconstructible through $U$ up to equivalence} if there exists a map
		\[
		\widehat J:
		\mathcal Y/{\simeq_{\mathcal Y}}
		\longrightarrow
		\mathcal Z/{\simeq_{\mathcal Z}}
		\]
		such that
		\[
		q_{\mathcal Z}\circ J
		=
		\widehat J\circ q_{\mathcal Y}\circ U.
		\]		
	\end{enumerate}
\end{definition}

A reduction is a specified comparison map, not a primitive realization mechanism.
Neither reduction nor reconstruction is required by the general packaging axioms.

\begin{proposition}[Factorization and non-reconstruction criterion]
	\label{PROP:PackagingReductionFactorizationCriterion}
	Let
	\[
	U:\mathcal X\longrightarrow\mathcal Y=U(\mathcal X),
	\qquad
	J:\mathcal X\longrightarrow\mathcal Z.
	\]	
	Then \(J\) is exactly reconstructible through \(U\) if and only if
	\[
	U(x_0)=U(x_1)
	\quad\Longrightarrow\quad
	J(x_0)=J(x_1)
	\]
	for all \(x_0,x_1\in\mathcal X\).	
	More generally, with declared equivalences \(\simeq_{\mathcal Y}\) and \(\simeq_{\mathcal Z}\), \(J\) is reconstructible through \(U\) up to equivalence if and only if
	\[
	U(x_0)\simeq_{\mathcal Y}U(x_1)
	\quad\Longrightarrow\quad
	J(x_0)\simeq_{\mathcal Z}J(x_1).
	\]	
	Consequently, a pair satisfying
	\[
	U(x_0)\simeq_{\mathcal Y}U(x_1),
	\qquad
	J(x_0)\not\simeq_{\mathcal Z}J(x_1)
	\]
	is an obstruction to reconstruction through \(U\).	
	If, in addition, \(\mathcal X\) carries an equivalence
	\(\simeq_{\mathcal X}\) and \(J\) is invariant under it, then the same witness implies
	\[
	x_0\not\simeq_{\mathcal X}x_1.
	\]
\end{proposition}

\begin{proof}
	Suppose first that
	\[
	J=\overline J\circ U.
	\]
	Then $U(x_0)=U(x_1)$ implies
	$J(x_0)=J(x_1)$.
	
	Conversely, if $J$ is constant on every
	fibre of $U$, define
	\[
	\overline J(U(x)):=J(x).
	\]
	Since $\mathcal Y=U(\mathcal X)$, this defines a map on all of $\mathcal Y$.
	Fibre-constancy makes it well defined, and $J=\overline J\circ U$.
	
	For reconstruction up to equivalence, consider
	\[
	q_{\mathcal Y}\circ U:
	\mathcal X
	\longrightarrow
	\mathcal Y/{\simeq_{\mathcal Y}}.
	\]
	This map is surjective.
	
	The stated implication is precisely the condition that $q_{\mathcal Z}\circ J$	is constant on its fibres.
	Hence
	\[
	\widehat J([U(x)]_{\simeq_{\mathcal Y}})
	:=
	[J(x)]_{\simeq_{\mathcal Z}}
	\]
	defines a well-defined map satisfying
	\[
	q_{\mathcal Z}\circ J
	=
	\widehat J\circ q_{\mathcal Y}\circ U.
	\]	
	Conversely, such a factorization immediately implies the stated condition.	
	The obstruction and nonequivalence assertions follow by contraposition and invariance of $J$, respectively.
\end{proof}

\begin{proposition}[Tagged relational presentation of the bare packaged data]
	\label{PROP:PackagingRelationalPresentation}
	Let
	\[
	\mathfrak P
	=
	\bigl(
	\operatorname{Pkg}_{\mathfrak P},
	E_{\mathfrak P},
	\triangleleft_{\mathfrak P},
	\mathsf L_{\mathfrak P},
	\mathsf R_{\mathfrak P}
	\bigr)
	\]
	be a packaged structure.
	Put
	\[
	\mathsf C_{\mathfrak P}
	=
	\bigsqcup_{K\in\operatorname{Pkg}_{\mathfrak P}}
	\bigl(\{K\}\times\mathsf L_{\mathfrak P}(K)\bigr).
	\]
	On the three sorts
	\[
	\operatorname{Pkg}_{\mathfrak P},
	\qquad
	E_{\mathfrak P},
	\qquad
	\mathsf C_{\mathfrak P},
	\]
	define
	\[
	\operatorname{Inc}(e,K)
	\quad\Longleftrightarrow\quad
	e\triangleleft_{\mathfrak P}K,
	\]
	\[
	\operatorname{Own}((K,\ell),L)
	\quad\Longleftrightarrow\quad
	K=L,
	\]
	and
	\[
	\operatorname{Sel}((K,\ell),L)
	\quad\Longleftrightarrow\quad
	K=L
	\text{ and }
	\ell\in\mathsf R_{\mathfrak P}(K).
	\]	
	The tagged presentation determines the five coordinates of \(\mathfrak P\).
	Viewed as an abstract many-sorted relational structure, it determines these coordinates up to isomorphism of the corresponding sorts and candidate fibres.	
\end{proposition}

\begin{proof}
	Use three sorts: package carriers, internal units, and the tagged disjoint union $\bigsqcup_K(\{K\}\times\mathsf L_{\mathfrak P}(K))$.
	The incidence, ownership, and selection relations recover respectively $\triangleleft_{\mathfrak P}$, the fibres $\mathsf L_{\mathfrak P}(K)$, and the selected subsets $\mathsf R_{\mathfrak P}(K)$.
	Hence the displayed	many-sorted relational datum reconstructs the five bare coordinates.
	Any	additional internal structure carried by the candidate spaces is not part of this bare relational signature.
\end{proof}

\subsection{Derived package-incidence geometry}

Package interiors determine two ordinary combinatorial objects: a hypergraph whose edges are the interiors and a graph recording which packages overlap.
These objects retain incidence but omit local candidate spaces, selections, and formation data.
Later we use related cover geometry for charts of descriptor labels, which need not be package interiors.

\begin{definition}[Derived package-incidence geometry]
	\label{DEF:PackagingIncidenceHypergraph}	
	Let $\mathfrak P$ be a packaged structure.
	Define
	\[
	E^{\mathrm{act}}_{\mathfrak P}
	=
	\bigcup_{K\in\operatorname{Pkg}_{\mathfrak P}} I_K,
	\qquad
	\operatorname{Pkg}^{+}_{\mathfrak P}
	=
	\{K\in\operatorname{Pkg}_{\mathfrak P}:I_K\ne\varnothing\}.
	\]	
	The \textbf{derived package-incidence geometry} of $\mathfrak P$ consists of the following two structures:	
	\begin{enumerate}[label=(\roman*),leftmargin=2.5em]		
		\item The \textbf{package-incidence hypergraph}
		\[
		\mathcal H_{\mathfrak P}
		=
		\left(
		E^{\mathrm{act}}_{\mathfrak P},
		\{I_K:K\in\operatorname{Pkg}^{+}_{\mathfrak P}\}
		\right),
		\]
		in which repeated interiors are represented by a single hyperedge.
		
		\item The \textbf{package-overlap graph}
		$G^{\mathrm{ov}}_{\mathfrak P}$, whose vertex set is $\operatorname{Pkg}^{+}_{\mathfrak P}$,	with distinct vertices $K,L$ adjacent	if and only if
		\[
		I_K\cap I_L\ne\varnothing.
		\]		
	\end{enumerate}
\end{definition}

Both structures are determined by the package-incidence coordinates of $\mathfrak P$.
They do not retain the local candidate spaces or selected candidate families, nor the formation datum or context of a realized packaging system.

\subsection{Native witness: multivalued cover formation}
\label{SUBSEC:PackagingMultivaluedCoverFormation}

We now give a concrete example in which one datum has several outcomes.
The datum lists available sets, and the outcomes are their inclusion-minimal covers.
Proposition~\ref{PROP:PackagingMultivaluedCoverLargeFibres} produces \(3^m\) outcomes for every \(m\ge1\).
This proves that multivalued formation occurs naturally, not only because the axioms permit arbitrary relations.

Fix a finite nonempty set $E$ and put
$
\mathscr C_E=\mathcal P(E).
$
For each family $\mathcal F\subseteq\mathcal P(E)$, define the packaged structure
\[
\mathfrak P_{\mathcal F}
=
\bigl(
\mathcal F,
E,
\triangleleft_{\mathcal F},
\mathsf L_{\mathcal F},
\mathsf R_{\mathcal F}
\bigr),
\]
where
\[
e\triangleleft_{\mathcal F}K
\quad\Longleftrightarrow\quad
e\in K,
\]
and
\[
\mathsf L_{\mathcal F}(K)
=
\mathsf R_{\mathcal F}(K)
=
\{*_K\}
\qquad(K\in\mathcal F).
\]
Every such $\mathfrak P_{\mathcal F}$ satisfies the packaged-structure axioms, whether or not $\mathcal F$ covers $E$.

\begin{example}[Minimal-cover packaging]
	\label{EX:PackagingMultivaluedCover}	
	Define the packaging context $\Xi_E^{\mathrm{cov}}$ by
	\[
	\mathscr C_{\Xi_E^{\mathrm{cov}}}
	=
	\operatorname{Car}_{\Xi_E^{\mathrm{cov}}}
	=
	\mathcal P(E),
	\]
	\[
	\operatorname{Form}_{\Xi_E^{\mathrm{cov}}}
	=
	\mathcal P(\mathcal P(E)).
	\]	
	Its realization relation is
	\[
	\mathsf J_{\Xi_E^{\mathrm{cov}}}
	=
	\left\{
	(\mathcal A,\mathfrak P_{\mathcal F}):
	\begin{array}{l}
		\mathcal F\subseteq\mathcal A
		\subseteq\mathcal P(E),\\[2pt]
		\displaystyle
		\bigcup_{K\in\mathcal F}K=E,\\[4pt]
		\displaystyle
		\bigcup_{K\in\mathcal F\setminus\{L\}}K
		\ne E
		\quad(L\in\mathcal F)
	\end{array}
	\right\}.
	\]	
	Thus a formation datum $\mathcal A$ specifies the available carriers, and
	\[
	\operatorname{Out}_{\Xi_E^{\mathrm{cov}}}
	(\mathcal A)
	=
	\left\{
	\mathfrak P_{\mathcal F}:
	\mathcal F\subseteq\mathcal A
	\text{ is an inclusion-minimal cover of }E
	\right\}.
	\]	
	Since $E$ is finite, this outcome family is nonempty if and only if
	\[
	\bigcup_{K\in\mathcal A}K=E.
	\]
\end{example}

\begin{proposition}[Multivalued realization with unbounded outcome families]
	\label{PROP:PackagingMultivaluedCoverLargeFibres}
	\leavevmode
	\begin{enumerate}
		\item The context $\Xi_E^{\mathrm{cov}}$ is nondeterministic whenever $|E|\ge2$.
			
		\item For every integer $m\ge1$, there exist a finite set $E_m$ and a formation datum
		\[
		\mathcal A_m
		\in
		\operatorname{Form}_{\Xi_{E_m}^{\mathrm{cov}}}
		\]
		such that
		\[
		\left|
		\operatorname{Out}_{\Xi_{E_m}^{\mathrm{cov}}}
		(\mathcal A_m)
		\right|
		=
		3^m.
		\]
	\end{enumerate}
\end{proposition}

\begin{proof}
	\leavevmode
	\begin{enumerate}
		\item Suppose $|E|\ge2$ and choose $a\in E$.	
		For
		$
		\mathcal A
		=
		\{E,\{a\},E\setminus\{a\}\},
		$
		the two families
		\[
		\mathcal F_0=\{E\},
		\qquad
		\mathcal F_1=
		\bigl\{\{a\},E\setminus\{a\}\bigr\}
		\]
		are distinct inclusion-minimal covers of $E$.	
		Hence
		$
		\mathfrak P_{\mathcal F_0}
		\ne
		\mathfrak P_{\mathcal F_1},
		$
		and both belong to $\operatorname{Out}_{\Xi_E^{\mathrm{cov}}}
		(\mathcal A)$.
		Thus the context is nondeterministic.
		
		\item For the second assertion, let
		\[
		E_m
		=
		\bigsqcup_{i=1}^{m}
		\{a_i,b_i,c_i\},
		\]
		and let $\mathcal A_m$ consist of the three two-element subsets of each triple $\{a_i,b_i,c_i\}$.
		
		An inclusion-minimal cover of $E_m$ contains exactly two of the three available	subsets from each triple.		
		Conversely, every independent choice of two subsets from each triple gives an inclusion-minimal cover of $E_m$.
		
		Each triple therefore contributes exactly three choices, independently of the others.
		Since distinct carrier families determine distinct packaged structures,
		\[
		\left|
		\operatorname{Out}_{\Xi_{E_m}^{\mathrm{cov}}}
		(\mathcal A_m)
		\right|
		=
		3^m.
		\]
	\end{enumerate}
\end{proof}

The following access--commutant example contrasts deterministic realization with the multivalued minimal-cover mechanism.

These are properties of particular specializations, not requirements imposed on general packaging systems.

\subsection{Native witness: access--commutant formation}
\label{SUBSEC:PackagingAccessCommutantGuidingModel}

Our second example is deterministic: a concrete algebra is sent to its commutant.
In finite dimension, taking the commutant again recovers the represented source.
Yet Proposition~\ref{PROP:PackagingAccessCommutantNonQuotient} shows that this formation need not be an algebra quotient.
The example separates recovery from quotient formation and leads to the dependency theory studied later.

\begin{example}[Access--commutant packaging]
	\label{EX:PackagingAccessCommutant}	
	Fix a nonzero finite-dimensional complex Hilbert space $\mathcal H$, and let $\mathscr C_{\mathcal H}$ be the set of concrete unital $*$-subalgebras of $B(\mathcal H)$.
	
	For $A\in\mathscr C_{\mathcal H}$, let $\operatorname{Und}(A)$ denote its underlying set.
	Define the one-carrier refinement-packaged structure
	\[
	\mathfrak P_A
	=
	\bigl(
	\{A\},
	\operatorname{Und}(A),
	\operatorname{Und}(A)\times\{A\},
	\mathsf R_A
	\bigr),
	\]
	using the abbreviated notation of
	Definition~\ref{DEF:PackagingRefinementPackagedStructure},
	where
	\[
	\mathsf L_A(A)
	=
	\mathsf R_A(A)
	=
	\mathbf{Ref}(\operatorname{Und}(A)).
	\]	
	Define the packaging context $\Xi_{\mathcal H}^{\mathrm{com}}$ by
	\[
	\mathscr C_{\Xi_{\mathcal H}^{\mathrm{com}}}
	=
	\operatorname{Car}_{\Xi_{\mathcal H}^{\mathrm{com}}}
	=
	\mathscr C_{\mathcal H},
	\]
	\[
	\operatorname{Form}_{\Xi_{\mathcal H}^{\mathrm{com}}}
	=
	\{(\mathrm{com},M):M\in\mathscr C_{\mathcal H}\},
	\]
	and
	\[
	\mathsf J_{\Xi_{\mathcal H}^{\mathrm{com}}}
	=
	\left\{
	((\mathrm{com},M),\mathfrak P_{M'}):
	M\in\mathscr C_{\mathcal H}
	\right\},
	\]
	where
	\[
	M'
	=
	\{T\in B(\mathcal H):
	Tm=mT\text{ for every }m\in M\}.
	\]	
	Thus
	\[
	\Xi_{\mathcal H}^{\mathrm{com}}
	\vdash
	(\mathrm{com},M)
	\rightsquigarrow
	\mathfrak P_{M'}
	\qquad(M\in\mathscr C_{\mathcal H}).
	\]	
	This context is total and deterministic.
	Its realization law is commutant formation, a mechanism specific to this specialization.
\end{example}

\begin{proposition}[Recovery from the represented commutant]
	\label{PROP:PackagingAccessCommutantRecovery}	
	Let $\mathcal H$ be a nonzero finite-dimensional complex Hilbert space, and let $M\subseteq B(\mathcal H)$ be a concrete	unital $*$-subalgebra.	
	If $N=M'$, then
	$
	N'=M.
	$	
	Consequently, the represented endpoint $N\subseteq B(\mathcal H)$, together with the ambient Hilbert space, determines the represented	source algebra.
	Commutant formation is an involution on	$\mathscr C_{\mathcal H}$.
\end{proposition}

\begin{proof}
	The finite-dimensional bicommutant theorem
	\cite[Theorem~3.2.1]{Jones2010VonNeumannNotes} gives
	\[
	M''=M.
	\]
	Since $N=M'$, we obtain
	\[
	N'=M''=M.
	\]
	The commutant of a concrete unital $*$-subalgebra is again a concrete unital	$*$-subalgebra, proving the involution assertion.
\end{proof}

\begin{corollary}[Preservation of represented source distinctions]
	\label{COR:PackagingAccessCommutantNoSourceCollapse}	
	In the uniform context $\Xi_{\mathcal H}^{\mathrm{com}}$, for all $M_0,M_1\in\mathscr C_{\mathcal H}$,
	\[
	M_0=M_1
	\quad\Longleftrightarrow\quad
	M_0'=M_1'.
	\]	
	Hence equality of represented sources and equality of represented endpoints induce the same comparison relation on realized systems.
\end{corollary}

\begin{proof}
	The forward implication is immediate.
	Conversely, $M_0'=M_1'$ implies
	\[
	M_0=M_0''=(M_0')'=(M_1')'=M_1''=M_1.
	\]
\end{proof}

This recovery uses the concrete algebra retained as the packaged carrier inside $B(\mathcal H)$.
It does not require the source coordinate of the formation datum.
Nor does it assert reconstruction from the bare incidence and refinement coordinates alone.

\begin{proposition}[Commutant formation of generated sources]
	\label{PROP:PackagingAccessCommutantGeneration}	
	Let $(M_i)_{i\in A}$ be a small family of concrete unital $*$-subalgebras of $B(\mathcal H)$, and put
	\[
	M=\operatorname{alg}^{*}
	\left(\bigcup_{i\in A}M_i\right),
	\]
	where the generated algebra is understood to be unital.	
	Then
	\[
	M'=\bigcap_{i\in A}M_i'.
	\]	
	This identity holds on an arbitrary Hilbert	space.
	If $\mathcal H$ is finite-dimensional, then
	\[
	\left(\bigcap_{i\in A}M_i'\right)'=M.
	\]	
	For $A=\varnothing$, the conventions are
	\[
	M=\mathbb C I_{\mathcal H},
	\qquad
	\bigcap_{i\in\varnothing}M_i'=B(\mathcal H).
	\]
\end{proposition}

\begin{proof}
	An operator commutes with the generated unital $*$-algebra precisely when it commutes with every $M_i$.
	Hence
	\[
	M'=\bigcap_{i\in A}M_i'.
	\]	
	In finite dimension, taking commutants and applying the bicommutant identity gives
	\[
	\left(\bigcap_{i\in A}M_i'\right)'
	=M''=M.
	\]
	The empty-family case follows from the stated conventions.
\end{proof}

Thus forming a commutant after generating an algebra from several sources gives the intersection of their commutants.
This is a native compatibility law of the access--commutant specialization.

A represented commutant formation $M\rightsquigarrow N$, where $N=M'$, is called an \textbf{algebra-quotient formation} if there exists a surjective complex-linear algebra homomorphism
\[
q:M\twoheadrightarrow N.
\]

\begin{proposition}[Access--commutant formation need not be algebra-quotient formation]
	\label{PROP:PackagingAccessCommutantNonQuotient}	
	There exist a nonzero finite-dimensional complex Hilbert space $\mathcal H$ and $M\in\mathscr C_{\mathcal H}$ such that
	\[
	\Xi_{\mathcal H}^{\mathrm{com}}
	\vdash
	(\mathrm{com},M)
	\rightsquigarrow
	\mathfrak P_{M'},
	\]
	but no surjective complex-linear algebra homomorphism
	\[
	M\twoheadrightarrow M'
	\]
	exists.
\end{proposition}

\begin{proof}
	Take
	\[
	\mathcal H=\mathbb C^2\otimes\mathbb C^3,
	\qquad
	M=M_2(\mathbb C)\otimes I_3.
	\]
	Then
	\[
	M'=I_2\otimes M_3(\mathbb C),
	\]
	and
	\[
	\dim_{\mathbb C}M=4,
	\qquad
	\dim_{\mathbb C}M'=9.
	\]
	Hence no surjective complex-linear map $M\to M'$ exists, and therefore no surjective complex-linear algebra homomorphism exists.
\end{proof}

A larger algebra dimension for $M'$ does not mean that the output carries more information than the source.
Exact represented-source recovery still holds by
Proposition~\ref{PROP:PackagingAccessCommutantRecovery}.

\begin{remark}[Endpoint representation versus native formation]
	\label{REM:PackagingAccessCommutantEndpointRepresentation}	
	In the non-quotient example above, write $N=M'$.
	The packaged endpoint $\mathfrak P_N$ has the single nonempty interior $\operatorname{Und}(N)$.	
	Its incidence therefore admits the one-block partition representation,	equivalently the fibre partition of a constant map or the canonical quotient by the indiscrete equivalence relation on $\operatorname{Und}(N)$.	
	These are representations of endpoint incidence.
	They do not identify the native formation $M\rightsquigarrow N$
	with an algebra quotient of $M$.
\end{remark}

Thus the two native formation laws exhibit different behaviors within the same primitive axioms.

We next study the intrinsic completion theory of refinement-packaged structures, where global coherence depends on selected local
refinements and their interactions across package overlaps.

\section{Intrinsic structure of refinement-packaged outputs}
\label{SEC:PackagingIntrinsicStructureA}

Selected local refinements do not always come from one global equivalence relation.
This section asks exactly when they do.
The main theorem, Theorem~\ref{THM:PackagingIntrinsicRefinementCompletionA}, constructs a completion operator and characterizes the local selections that extend to a global refinement.
It also shows that the completion space is preserved by static isomorphism.
Corollary~\ref{COR:PackagingIntrinsicCompletionNoncollapseA} gives two finite structures with matching incidence and locally matching refinements, but with different completion behavior.
These conclusions depend only on the packaged output, not on how it was formed.

\subsection{Intrinsic refinement completion}

We ask which selected refinements on overlapping interiors come from a single equivalence relation on all active units.
Theorem~\ref{THM:PackagingIntrinsicRefinementCompletionA} answers the question with an intrinsic operator: its fixed points are exactly the realizable local refinement tuples.
The theorem also gives a least completion for each such tuple and proves invariance under static isomorphism.

\begin{definition}[Static isomorphism of refinement-packaged structures]
	\label{DEF:PackagingStaticIsomorphismA}	
	Let $\mathfrak P$ and $\mathfrak P'$ be	refinement-packaged structures over the same structural domain.	
	A \textbf{static isomorphism}
	\[
	(F,\phi):\mathfrak P\longrightarrow\mathfrak P'
	\]
	consists of bijections
	\[
	F:\operatorname{Pkg}_{\mathfrak P}
	\longrightarrow
	\operatorname{Pkg}_{\mathfrak P'},
	\qquad
	\phi:E_{\mathfrak P}
	\longrightarrow E_{\mathfrak P'},
	\]
	satisfying the following conditions:	
	\begin{enumerate}[label=(\roman*),leftmargin=2.5em]		
		\item \textbf{Incidence preservation.}
		For every $e\in E_{\mathfrak P}$ and $K\in\operatorname{Pkg}_{\mathfrak P}$,
		\[
		e\triangleleft_{\mathfrak P}K
		\quad\Longleftrightarrow\quad
		\phi(e)\triangleleft_{\mathfrak P'}F(K).
		\]
		
		\item \textbf{Selected-refinement preservation.}
		For each $K\in\operatorname{Pkg}_{\mathfrak P}$,
		let
		\[
		\phi_K:
		I_K\xrightarrow{\sim}I'_{F(K)}
		\]
		be the restriction of $\phi$, which is a bijection by \textup{(i)}.		
		Define the induced pullback
		\[
		\phi_K^*:
		\mathbf{Ref}(I'_{F(K)})
		\longrightarrow
		\mathbf{Ref}(I_K)
		\]
		by
		$
		\phi_K^*(q_R)
		=
		q_{(\phi_K\times\phi_K)^{-1}(R)}.
		$		
		Then
		\[
		\phi_K^*
		\bigl(
		\mathsf R_{\mathfrak P'}(F(K))
		\bigr)
		=
		\mathsf R_{\mathfrak P}(K),
		\]
		where $\phi_K^*$ acts elementwise on subsets of the refinement poset.		
	\end{enumerate}
\end{definition}

Static isomorphism preserves the package--unit incidence and selected local refinements.
It does not require preservation of additional algebraic, topological, or other native structure carried by the package carriers.
Thus static isomorphism is a derived comparison of refinement-packaged outputs, not a primitive requirement of general packaging or realization.

\begin{definition}[Intrinsic refinement-completion space]
	\label{DEF:PackagingIntrinsicCompletionA}	
	Let $\mathfrak P$ be a refinement-packaged structure, and put
	\[
	E^{\mathrm{act}}_{\mathfrak P}
	=
	\bigcup_{K\in\operatorname{Pkg}_{\mathfrak P}} I_K.
	\]	
	The \textbf{local-refinement space} is
	\[
	\mathsf{Loc}(\mathfrak P)
	=
	\prod_{K\in\operatorname{Pkg}_{\mathfrak P}}
	\mathbf{Ref}(I_K),
	\]
	equipped with the coordinatewise refinement order.	
	The \textbf{local trace map}
	\[
	\operatorname{Tr}_{\mathfrak P}:
	\operatorname{Eq}(E^{\mathrm{act}}_{\mathfrak P})
	\longrightarrow
	\mathsf{Loc}(\mathfrak P)
	\]
	is defined by
	\[
	\operatorname{Tr}_{\mathfrak P}(H)
	=
	\left(
	q_{H\cap I_K^2}
	\right)_{K\in\operatorname{Pkg}_{\mathfrak P}}.
	\]	
	The \textbf{intrinsic refinement-completion space} of $\mathfrak P$ is
	\[
	\mathcal G(\mathfrak P)
	=
	\left\{
	H\in\operatorname{Eq}
	(E^{\mathrm{act}}_{\mathfrak P}):
	\operatorname{Tr}_{\mathfrak P}(H)
	\in
	\prod_{K\in\operatorname{Pkg}_{\mathfrak P}}
	\mathsf R_{\mathfrak P}(K)
	\right\}.
	\]
\end{definition}

\begin{definition}[Generated equivalence and coherence]
	\label{DEF:PackagingGeneratedEquivalenceA}	
	Let $\mathfrak P$ be a refinement-packaged structure, and let
	\[
	r=(r_K)_{K\in\operatorname{Pkg}_{\mathfrak P}}
	\in\mathsf{Loc}(\mathfrak P).
	\]	
	The \textbf{generated equivalence relation}
	of $r$ is
	\[
	G_{\mathfrak P}(r)
	=
	\left(
	\bigcup_{K\in\operatorname{Pkg}_{\mathfrak P}}
	\ker r_K
	\right)^{\mathrm{eq}},
	\]
	where equivalence closure is taken on $E^{\mathrm{act}}_{\mathfrak P}$.
	
	Define the \textbf{local completion operator}
	\[
	\Gamma_{\mathfrak P}:
	\mathsf{Loc}(\mathfrak P)
	\longrightarrow
	\mathsf{Loc}(\mathfrak P)
	\]
	by
	$
	\Gamma_{\mathfrak P}(r)
	=
	\operatorname{Tr}_{\mathfrak P}
	\bigl(G_{\mathfrak P}(r)\bigr).
	$	
	A tuple $r\in\mathsf{Loc}(\mathfrak P)$ is called:	
	\begin{enumerate}[label=(\roman*),leftmargin=2.5em]		
		\item \textbf{selected} if
		\[
		r\in
		\prod_{K\in\operatorname{Pkg}_{\mathfrak P}}
		\mathsf R_{\mathfrak P}(K);
		\]
		
		\item \textbf{coherent} if
		\[
		\Gamma_{\mathfrak P}(r)=r.
		\]		
	\end{enumerate}
\end{definition}

\noindent\textbf{Notation.}
When $\mathfrak P$ is fixed, we suppress the subscript and write
\[
G(r),\qquad
\Gamma(r),\qquad
\operatorname{Tr}(H)
\]
for the corresponding constructions.

\begin{theorem}[Intrinsic completion operator and exact coherent traces]
	\label{THM:PackagingIntrinsicRefinementCompletionA}
	Let \(\mathfrak P\) be a refinement-packaged structure.
	\begin{enumerate}[label=(\roman*),leftmargin=2.5em]
		\item The map
		\[
		\Gamma_{\mathfrak P}:
		\mathsf{Loc}(\mathfrak P)
		\longrightarrow
		\mathsf{Loc}(\mathfrak P)
		\]
		is order-preserving and idempotent, and satisfies
		$
		\Gamma_{\mathfrak P}(r)
		\preceq
		r
		$
		coordinatewise for every \(r\in\mathsf{Loc}(\mathfrak P)\).
		Thus \(\Gamma_{\mathfrak P}\) is an interior operator on the coordinatewise refinement poset.
		
		\item The image of the local trace map is exactly the set of coherent local refinement tuples:
		\[
		\operatorname{Tr}_{\mathfrak P}
		\left(
		\operatorname{Eq}(E^{\mathrm{act}}_{\mathfrak P})
		\right)
		=
		\operatorname{Fix}(\Gamma_{\mathfrak P}).
		\]		
		Consequently, the traces of global completions are exactly the coherent selected tuples:
		\[
		\operatorname{Tr}_{\mathfrak P}
		\bigl(\mathcal G(\mathfrak P)\bigr)
		=
		\operatorname{Fix}(\Gamma_{\mathfrak P})
		\cap
		\prod_K\mathsf R_{\mathfrak P}(K).
		\]		
		In particular,
		\[
		\mathcal G(\mathfrak P)\ne\varnothing
		\quad\Longleftrightarrow\quad
		\operatorname{Fix}(\Gamma_{\mathfrak P})
		\cap
		\prod_K\mathsf R_{\mathfrak P}(K)
		\ne\varnothing.
		\]
		
		\item Let
		\[
		r\in
		\prod_K\mathsf R_{\mathfrak P}(K)
		\]
		be coherent.
		Then
		$
		G_{\mathfrak P}(r)
		\in
		\mathcal G(\mathfrak P)
		$
		and the complete fibre of the trace map over \(r\) is
		\[
		\operatorname{Tr}_{\mathfrak P}^{-1}(r)
		=
		\left\{
		H\in
		\operatorname{Eq}
		\bigl(E^{\mathrm{act}}_{\mathfrak P}\bigr):
		\begin{array}{l}
			G_{\mathfrak P}(r)\subseteq H,\\[2pt]
			H\cap I_K^2=\ker r_K
			\text{ for every }K
		\end{array}
		\right\}.
		\]
		In particular, \(G_{\mathfrak P}(r)\) is the unique least global completion with local trace \(r\), ordered by inclusion of equivalence relations.
		
		\item Let
		$
		(F,\phi):
		\mathfrak P\longrightarrow\mathfrak P'
		$
		be a static isomorphism.
		It induces an order isomorphism
		\[
		\Phi^*:
		\mathsf{Loc}(\mathfrak P')
		\longrightarrow
		\mathsf{Loc}(\mathfrak P)
		\]
		defined by
		\[
		\bigl(\Phi^*r'\bigr)_K
		=
		\phi_K^*
		\bigl(
		r'_{F(K)}
		\bigr),
		\]
		and
		\[
		\Phi^*\circ\Gamma_{\mathfrak P'}
		=
		\Gamma_{\mathfrak P}\circ\Phi^*.
		\]
		It also induces an order isomorphism, with respect to inclusion of equivalence relations,
		\[
		\mathcal G(\mathfrak P')
		\longrightarrow
		\mathcal G(\mathfrak P),
		\qquad
		H'
		\longmapsto
		(\phi\times\phi)^{-1}(H')
		\big|_{(E^{\mathrm{act}}_{\mathfrak P})^2}.
		\]
		Hence
		$
		g(\mathfrak P)
		:=
		|\mathcal G(\mathfrak P)|
		$
		is a static-isomorphism invariant.
		If \(E^{\mathrm{act}}_{\mathfrak P}\) is finite, then \(g(\mathfrak P)\) is a finite numerical invariant.
	\end{enumerate}
\end{theorem}

\begin{proof}
	We prove the four assertions separately.	
	\begin{enumerate}[label=(\roman*),leftmargin=2.5em]		
		\item \emph{Interior-operator properties.}
		Fix
		\[
		r=(r_K)_K\in\mathsf{Loc}(\mathfrak P).
		\]
		For every package $K$,
		\[
		\ker r_K
		\subseteq
		G_{\mathfrak P}(r)\cap I_K^2.
		\]
		Since
		\[
		q_R\preceq q_T
		\quad\Longleftrightarrow\quad
		T\subseteq R,
		\]
		we obtain
		\[
		\Gamma_{\mathfrak P}(r)\preceq r.
		\]
		
		Next, suppose $r\preceq s$.
		Then
		\[
		\ker s_K\subseteq\ker r_K
		\qquad\text{for every }K,
		\]
		and hence
		\[
		G_{\mathfrak P}(s)
		\subseteq
		G_{\mathfrak P}(r).
		\]
		Taking local traces gives
		\[
		\Gamma_{\mathfrak P}(r)
		\preceq
		\Gamma_{\mathfrak P}(s).
		\]
		Thus $\Gamma_{\mathfrak P}$ is order-preserving.
		
		For idempotence, observe that
		\[
		\ker r_K
		\subseteq
		\ker\Gamma_{\mathfrak P}(r)_K
		\subseteq
		G_{\mathfrak P}(r)
		\]
		for every $K$.
		Taking unions and generated equivalence closures yields
		\[
		G_{\mathfrak P}(r)
		\subseteq
		G_{\mathfrak P}
		\bigl(\Gamma_{\mathfrak P}(r)\bigr)
		\subseteq
		G_{\mathfrak P}(r).
		\]
		Therefore
		\[
		G_{\mathfrak P}
		\bigl(\Gamma_{\mathfrak P}(r)\bigr)
		=
		G_{\mathfrak P}(r).
		\]
		Applying the local trace map gives
		\[
		\Gamma_{\mathfrak P}^2(r)
		=
		\Gamma_{\mathfrak P}(r).
		\]
		Hence $\Gamma_{\mathfrak P}$ is an interior operator on $\mathsf{Loc}(\mathfrak P)$.
		
		\item \emph{Coherent traces and selected completions.}
		Let
		\[
		H\in\operatorname{Eq}(E^{\mathrm{act}}_{\mathfrak P}),
		\qquad
		r=\operatorname{Tr}_{\mathfrak P}(H).
		\]
		For every package $K$,
		\[
		\ker r_K=H\cap I_K^2.
		\]
		Since $H$ contains every local kernel,
		\[
		G_{\mathfrak P}(r)\subseteq H.
		\]
		Hence
		\[
		\ker r_K
		\subseteq
		G_{\mathfrak P}(r)\cap I_K^2
		\subseteq
		H\cap I_K^2
		=
		\ker r_K,
		\]
		so
		\[
		\Gamma_{\mathfrak P}(r)=r.
		\]
		This proves
		\[
		\operatorname{Im}(\operatorname{Tr}_{\mathfrak P})
		\subseteq
		\operatorname{Fix}(\Gamma_{\mathfrak P}).
		\]
		
		Conversely, if
		$r\in\operatorname{Fix}(\Gamma_{\mathfrak P})$,
		then
		\[
		\operatorname{Tr}_{\mathfrak P}
		\bigl(G_{\mathfrak P}(r)\bigr)
		=
		\Gamma_{\mathfrak P}(r)
		=
		r.
		\]
		Thus every fixed point is a global trace, and
		\[
		\operatorname{Tr}_{\mathfrak P}
		\bigl(
		\operatorname{Eq}(E^{\mathrm{act}}_{\mathfrak P})
		\bigr)
		=
		\operatorname{Fix}(\Gamma_{\mathfrak P}).
		\]
		
		By definition, a global equivalence relation
		belongs to $\mathcal G(\mathfrak P)$ exactly
		when its trace is selected. Therefore
		\[
		\operatorname{Tr}_{\mathfrak P}
		\bigl(\mathcal G(\mathfrak P)\bigr)
		=
		\operatorname{Fix}(\Gamma_{\mathfrak P})
		\cap
		\prod_K\mathsf R_{\mathfrak P}(K).
		\]
		The nonemptiness criterion follows.
		
		\item \emph{Completion fibres and leastness.}
		Fix a coherent selected tuple
		\[
		r\in\prod_K\mathsf R_{\mathfrak P}(K).
		\]
		By \textup{(ii)},
		\[
		G_{\mathfrak P}(r)\in\mathcal G(\mathfrak P),
		\qquad
		\operatorname{Tr}_{\mathfrak P}
		\bigl(G_{\mathfrak P}(r)\bigr)=r.
		\]
		
		Let
		\[
		H\in
		\operatorname{Eq}(E^{\mathrm{act}}_{\mathfrak P}).
		\]
		The condition
		\[
		\operatorname{Tr}_{\mathfrak P}(H)=r
		\]
		is equivalent to
		\[
		H\cap I_K^2=\ker r_K
		\qquad\text{for every }K.
		\]
		
		Since $H$ is an equivalence relation, these equalities imply
		\[
		G_{\mathfrak P}(r)\subseteq H.
		\]
		Conversely, the stated local equalities, together with this inclusion, give	$\operatorname{Tr}_{\mathfrak P}(H)=r$.
		
		Hence
		\[
		\operatorname{Tr}_{\mathfrak P}^{-1}(r)
		=
		\left\{
		H\in\operatorname{Eq}
		(E^{\mathrm{act}}_{\mathfrak P}):
		\begin{array}{l}
			G_{\mathfrak P}(r)\subseteq H,\\[2pt]
			H\cap I_K^2=\ker r_K
			\quad\text{for every }K
		\end{array}
		\right\}.
		\]		
		Because $G_{\mathfrak P}(r)$ belongs to this fibre and is contained in every other member, it is the unique least global completion with trace $r$.
		
		\item \emph{Invariance under static isomorphism.}
		Let
		\[
		(F,\phi):
		\mathfrak P\longrightarrow\mathfrak P'
		\]
		be a static isomorphism.
		
		Incidence preservation gives
		\[
		\phi(E^{\mathrm{act}}_{\mathfrak P})
		=
		E^{\mathrm{act}}_{\mathfrak P'}.
		\]
		Thus $\phi$ restricts to a bijection
		\[
		\phi_{\mathrm{act}}:
		E^{\mathrm{act}}_{\mathfrak P}
		\longrightarrow
		E^{\mathrm{act}}_{\mathfrak P'}.
		\]
		
		For
		$H'\in\operatorname{Eq}
		(E^{\mathrm{act}}_{\mathfrak P'})$,
		define
		\[
		\Psi(H')
		=
		(\phi_{\mathrm{act}}\times
		\phi_{\mathrm{act}})^{-1}(H').
		\]
		
		Pullback by a bijection commutes with restriction, unions of relations, and generated equivalence closure.
		Consequently, for every $r'\in\mathsf{Loc}(\mathfrak P')$,
		\[
		\Psi\bigl(G_{\mathfrak P'}(r')\bigr)
		=
		G_{\mathfrak P}(\Phi^*r'),
		\]
		and for every global equivalence relation $H'$,
		\[
		\operatorname{Tr}_{\mathfrak P}
		\bigl(\Psi(H')\bigr)
		=
		\Phi^*
		\bigl(\operatorname{Tr}_{\mathfrak P'}(H')\bigr).
		\]
		
		Combining these identities yields
		\[
		\begin{aligned}
			\Phi^*
			\bigl(\Gamma_{\mathfrak P'}(r')\bigr)
			&=
			\Phi^*
			\bigl(
			\operatorname{Tr}_{\mathfrak P'}
			(G_{\mathfrak P'}(r'))
			\bigr)\\
			&=
			\operatorname{Tr}_{\mathfrak P}
			\bigl(
			G_{\mathfrak P}(\Phi^*r')
			\bigr)\\
			&=
			\Gamma_{\mathfrak P}(\Phi^*r').
		\end{aligned}
		\]
		Therefore
		\[
		\Phi^*\circ\Gamma_{\mathfrak P'}
		=
		\Gamma_{\mathfrak P}\circ\Phi^*.
		\]
		
		Each $\phi_K^*$ is an order isomorphism of refinement posets, and $F$ is a bijection of package families.
		Hence $\Phi^*$ is an order isomorphism of local-refinement spaces.
		
		Moreover, selected-refinement preservation gives
		\[
		\Phi^*
		\left(
		\prod_{K'}
		\mathsf R_{\mathfrak P'}(K')
		\right)
		=
		\prod_K\mathsf R_{\mathfrak P}(K).
		\]
		Together with trace compatibility, this shows that
		\[
		H'\in\mathcal G(\mathfrak P')
		\quad\Longleftrightarrow\quad
		\Psi(H')\in\mathcal G(\mathfrak P).
		\]
		
		Pullback by $\phi_{\mathrm{act}}^{-1}$ provides the inverse map.
		Thus $\Psi$ restricts to a bijection
		\[
		\mathcal G(\mathfrak P')
		\xrightarrow{\sim}
		\mathcal G(\mathfrak P).
		\]
		
		Consequently,
		\[
		g(\mathfrak P)=g(\mathfrak P').
		\]
		If $E^{\mathrm{act}}_{\mathfrak P}$ is finite, then
		$\operatorname{Eq} (E^{\mathrm{act}}_{\mathfrak P})$
		is finite, so $g(\mathfrak P)$ is a	finite numerical invariant.
		Since inverse images under bijections preserve and reflect inclusion, the induced bijection is an order isomorphism of completion posets.	
	\end{enumerate}
\end{proof}

\subsection{Intrinsic non-collapse of packaged structures}
\label{SUBSEC:PackagingIntrinsicNoncollapse}

Local agreement on overlaps does not always guarantee a global refinement.
Corollary~\ref{COR:PackagingIntrinsicCompletionNoncollapseA} gives two finite structures with the same incidence and matching local refinement types.
Both have agreement on every pairwise overlap, but one has a global completion and the other has none.
This is a concrete non-collapse result relative to the stated local comparison.

\begin{corollary}[Intrinsic non-collapse despite non-singleton overlap compatibility]
	\label{COR:PackagingIntrinsicCompletionNoncollapseA}
	There exist finite refinement-packaged structures
	$
	\mathfrak P_+,
	~
	\mathfrak P_-
	$
	with the following properties.
	\begin{enumerate}[label=(\roman*),leftmargin=2.5em]
		\item They have exactly the same package--unit incidence skeleton.
		
		\item For every corresponding package \(K\), the selected-candidate families of \(\mathfrak P_+\) and \(\mathfrak P_-\) are isomorphic by transport along a bijection of the corresponding interiors.
		
		\item Every nonempty pairwise chart overlap has two elements, and within each structure the selected local refinements have exactly the same restriction on that overlap.
		
		\item Their intrinsic completion spaces have different cardinalities:
		\[
		|\mathcal G(\mathfrak P_+)|=1,
		\qquad
		|\mathcal G(\mathfrak P_-)|=0.
		\]
	\end{enumerate}
	Therefore, even these three pieces of local information do not decide whether a global completion exists: incidence, independent local candidate types, and exact agreement on overlaps.
They also do not determine the packaged structure up to static isomorphism.
\end{corollary}

\begin{proof}
	Fix the structural domain
	$
	\mathscr C=\{A,B,C\}
	$
	and six pairwise disjoint two-element sets
	\[
	\begin{aligned}
		A_0&=\{a_0,a_1\},&
		B_0&=\{b_0,b_1\},&
		C_0&=\{c_0,c_1\},\\
		P_0&=\{p_0,p_1\},&
		Q_0&=\{q_0,q_1\},&
		R_0&=\{r_0,r_1\}.
	\end{aligned}
	\]	
	Put
	\[
	E=A_0\sqcup B_0\sqcup C_0
	\sqcup P_0\sqcup Q_0\sqcup R_0,
	\]
	and define the package interiors by
	\[
	\begin{aligned}
		I_A&=A_0\sqcup B_0\sqcup P_0,\\
		I_B&=B_0\sqcup C_0\sqcup Q_0,\\
		I_C&=C_0\sqcup A_0\sqcup R_0.
	\end{aligned}
	\]
	
	For pairwise disjoint blocks $X,Y,Z$, write $XY\mid Z$ for the refinement whose equivalence classes are $X\cup Y$ and $Z$.	
	Define two refinement-packaged structures $\mathfrak P_+$ and $\mathfrak P_-$ with the same package family $\mathscr C$, internal-unit
	set $E$, and incidence determined by the displayed interiors.	
	Their selected candidate families are
	\[
	\begin{aligned}
		\mathsf R_{\mathfrak P_+}(A)
		&=\{A_0P_0\mid B_0\},&
		\mathsf R_{\mathfrak P_-}(A)
		&=\{A_0B_0\mid P_0\},\\
		\mathsf R_{\mathfrak P_+}(B)
		&=\{B_0Q_0\mid C_0\},&
		\mathsf R_{\mathfrak P_-}(B)
		&=\{B_0C_0\mid Q_0\},\\
		\mathsf R_{\mathfrak P_+}(C)
		&=\{C_0R_0\mid A_0\},&
		\mathsf R_{\mathfrak P_-}(C)
		&=\{C_0R_0\mid A_0\}.
	\end{aligned}
	\]	
	We verify the four assertions.	
	\begin{enumerate}[label=(\roman*),leftmargin=2.5em]		
		\item \emph{Common incidence skeleton.}		
		Both structures have the same package family $\{A,B,C\}$, internal-unit set $E$, and package interiors $I_A,I_B,I_C$.		
		Their incidence relations therefore coincide.
		
		This proves \textup{(i)}.

		\item \emph{Chartwise selected-candidate isomorphism.}		
		Each selected refinement has exactly two equivalence classes, of cardinalities $4$ and $2$.
		
		More explicitly, on $I_A$ a bijection interchanging $B_0$ and $P_0$ elementwise, while fixing $A_0$, transports
		\[
		A_0P_0\mid B_0
		\quad\text{to}\quad
		A_0B_0\mid P_0.
		\]
		
		Similarly, on $I_B$ a bijection interchanging $C_0$ and $Q_0$ elementwise, while fixing $B_0$, transports
		\[
		B_0Q_0\mid C_0
		\quad\text{to}\quad
		B_0C_0\mid Q_0.
		\]
		
		On $I_C$, the corresponding selected refinements are identical.		
		Thus the selected-candidate families are isomorphic independently on every corresponding package.
		
		This proves \textup{(ii)}.

		\item \emph{Exact pairwise overlap compatibility.}		
		The pairwise intersections are
		\[
		I_A\cap I_B=B_0,
		\qquad
		I_B\cap I_C=C_0,
		\qquad
		I_C\cap I_A=A_0.
		\]		
		Each overlap contains exactly two elements.		
		In both $\mathfrak P_+$ and $\mathfrak P_-$, the adjacent selected refinements restrict to the indiscrete equivalence relation	on the corresponding overlap.		
		Consequently, within each structure, the selected refinements agree exactly on every pairwise overlap.
		Moreover, the corresponding overlap traces of the two structures coincide.
		
		This proves \textup{(iii)}.

		\item \emph{Different intrinsic completion spaces.}
		\begin{itemize}
			\item First consider $\mathfrak P_+$.		
			Let $H_+$ be the equivalence relation on $E$ whose equivalence classes are
			\[
			A_0\cup P_0,
			\qquad
			B_0\cup Q_0,
			\qquad
			C_0\cup R_0.
			\]		
			Its restriction to each package interior is exactly the prescribed selected refinement.
			Hence
			$
			H_+\in\mathcal G(\mathfrak P_+).
			$		
			We claim that $H_+$ is the unique completion.		
			Indeed, every global completion must identify all elements within each of the three displayed classes, because these identifications are forced by the selected local refinements.		
			Thus every global completion contains $H_+$.		
			Conversely, merging any two distinct $H_+$-classes would force at least one of
			\[
			A_0\sim B_0,\qquad
			B_0\sim C_0,\qquad
			C_0\sim A_0.
			\]		
			These identifications are forbidden by the prescribed refinements on $I_A$, $I_B$, and $I_C$, respectively.		
			Therefore no additional identification is possible, and
			$
			\mathcal G(\mathfrak P_+)=\{H_+\}.
			$
			
			\item Now consider $\mathfrak P_-$.		
			The selected refinement on $I_A$ forces
			$
			A_0\sim B_0,
			$
			while the selected refinement on $I_B$ forces
			$
			B_0\sim C_0.
			$		
			Transitivity of any global equivalence relation would therefore imply
			$
			A_0\sim C_0.
			$		
			However, the selected refinement on $I_C$ has equivalence classes
			$
			C_0\cup R_0
			$
			and
			$
			A_0,
			$
			and therefore requires $A_0$ and $C_0$ to remain distinct.		
			This contradiction shows that
			$
			\mathcal G(\mathfrak P_-)=\varnothing.
			$		
			Consequently,
			$
			|\mathcal G(\mathfrak P_+)|=1,
			\quad
			|\mathcal G(\mathfrak P_-)|=0.
			$				
		\end{itemize}	
		This proves \textup{(iv)}.
	\end{enumerate}	
	Finally, by Theorem~\ref{THM:PackagingIntrinsicRefinementCompletionA},
	the cardinality of the intrinsic completion	space is invariant under static isomorphism.	
	Since the two cardinalities differ,	$\mathfrak P_+$ and $\mathfrak P_-$ are	not statically isomorphic.
\end{proof}

The examples differ globally even though every nonempty overlap has two elements and the local restrictions agree exactly.
Moreover, since $\mathfrak P_+$ and $\mathfrak P_-$ have the same incidence coordinates, their local-refinement spaces and completion operators
coincide:
\[
\mathsf{Loc}(\mathfrak P_+)
=
\mathsf{Loc}(\mathfrak P_-),
\qquad
\Gamma_{\mathfrak P_+}
=
\Gamma_{\mathfrak P_-}
=:\Gamma.
\]
The difference lies in the selected candidate families.
Specifically,
\[
\operatorname{Fix}(\Gamma)
\cap
\prod_K\mathsf R_{\mathfrak P_+}(K)
\ne\varnothing,
\]
whereas
\[
\operatorname{Fix}(\Gamma)
\cap
\prod_K\mathsf R_{\mathfrak P_-}(K)
=
\varnothing.
\]
Thus neither the common incidence geometry nor the independent chartwise selected-candidate types determine global completion feasibility.
What matters is which refinements were selected, even when the completion operator itself is the same.

\begin{remark}[Relation to classical local equivalence-relation problems]
	\label{REM:PackagingLocalEquivalenceLiterature}
	Local-to-global questions for equivalence relations have classical antecedents.
	In particular, local equivalence relations and the problem of recovering global equivalence relations from local data were studied by
	Rosenthal~\cite{Rosenthal1982LocalEquivalence,Rosenthal1984SheavesLocalEquivalence} and Kock--Moerdijk~\cite{KockMoerdijk1996LocalEquivalence}.	
	Nonextendibility of compatible local equivalence relations is known.
	Theorem~\ref{THM:PackagingIntrinsicRefinementCompletionA} represents generated equivalences and local traces by an interior operator on the refinement space.
	The completion space depends on the selected candidate families; its cardinality and order structure are invariant under static isomorphism.
	This formulation yields the comparison-relative non-collapse witness above.
\end{remark}

The construction gives a separation witness for the comparison-based criterion of Section~\ref{SUBSEC:PackagingComparisonNoncollapse}.
The two packaged structures agree on their incidence skeleton, independent local selected-candidate types, and pairwise overlap traces, but have different values of the intrinsic invariant
\[
g(\mathfrak P)=|\mathcal G(\mathfrak P)|.
\]
Thus the completion invariant is intrinsic to the refinement-packaged output, while its failure to be determined by coarser data is a comparison-relative non-collapse result.

\section{Controlled propagation and exact globalization}
\label{SEC:PackagingControlledGlobalization}

This section characterizes universal exact globalization of closure operators on finite covers of descriptor labels, which need not coincide with package interiors.
For local Carath\'eodory number at most \(k\), Theorem~\ref{THM:PackagingGeneralCaratheodoryGluing} gives the sharp chordality and \((2k+1)\)-conformality criterion.
The proof uses controlled elimination and the recovery bounds of Theorems~\ref{THM:PackagingControlledPropagationGlobalization} and~\ref{THM:PackagingBoundedDeterminationPropagationRank}.
Theorem~\ref{THM:PackagingCanonicalCliqueCompletion} also constructs canonical clique completions without assuming chordality.
The complexity bound applies to the local closures, not necessarily to their global extension.

\subsection{Atlas compatibility and exact extension}
\label{SUBSEC:PackagingAtlasCompatibility}

Admissibility in a packaging context does not imply agreement between local charts.
Atlas compatibility requires equality of restrictions on chart overlaps.
Theorem~\ref{THM:PackagingAtlasRealizationBridge} represents each compatible atlas as a packaging system, without asserting that a global exact extension exists.

\begin{definition}[Atlas compatibility and exact extension]
	Let
	$
	\mathcal H
	$
	be a cover of a label set \(V\), and suppose that each chart \(I\in\mathcal H\) carries a local structure \(s_I\) admitting restriction to subsets.
	\begin{enumerate}
		\item The atlas
		\[
		\mathcal A=(s_I)_{I\in\mathcal H}
		\]
		is \emph{compatible} when the local structures agree on every overlap:
		\[
		{
			s_I|_{I\cap J}=s_J|_{I\cap J}
			\qquad
			(I,J\in\mathcal H).
		}
		\]
		
		\item A global structure \(S\) on \(V\) is an \emph{exact extension} of the atlas when
		\[
		{
			S|_I=s_I
			\qquad
			(I\in\mathcal H).
		}
		\]
	\end{enumerate}
\end{definition}

Atlas compatibility is necessary for exact extension but is not, in general, sufficient.
Indeed, if \(S\) is an exact extension, then
$
s_I|_{I\cap J}
=
S|_{I\cap J}
=
s_J|_{I\cap J}.
$
The converse is the central question of this part of the paper.
Compatible local structures may still fail to extend.
A consequence derived through other charts can return to an original chart even though that chart does not allow it.

For closure atlases this definition becomes
$
c_I|_{I\cap J}=c_J|_{I\cap J},
$
and compatibility has the following basic consequence:
a local implication whose complete support lies in an overlap has the same truth value in either chart.
This transport principle will be used repeatedly in propagation, filling, recovery, and obstruction arguments.

\medskip
\noindent\textbf{Realizing atlases as packaging systems.}

The packaging of an atlas can be made explicit at both levels of the primitive axiom system.
In particular, realizing the atlas as a packaging system and globally extending its local structures are distinct assertions.

Fix a finite indexed cover
\[
\mathcal H=(I_j)_{j\in J}
\]
of a finite set $V$.
Let $\mathscr S(I)$, for $I\subseteq V$, be
$\mathcal U$-small sets of structures equipped with identity- and composition-preserving restrictions.
Let
\[
\mathscr K(I)\subseteq\mathscr S(I)
\]
be the allowed local class, which need not contain every possible global extension.
Write
\[
\operatorname{Atl}_{\mathscr K}(\mathcal H)
=
\left\{
(s_j)_{j\in J}:
\begin{array}{l}
	s_j\in\mathscr K(I_j),\\
	s_i|_{I_i\cap I_j}=s_j|_{I_i\cap I_j}
	\quad(i,j\in J)
\end{array}
\right\}.
\]
Write $\operatorname{UG}_{\mathscr K}(\mathcal H)$
for the statement that every atlas in
$\operatorname{Atl}_{\mathscr K}(\mathcal H)$
has an exact extension in $\mathscr S(V)$.

\begin{theorem}[Atlas realization and exact-completion bridge]
	\label{THM:PackagingAtlasRealizationBridge}	
	With the preceding data, define for each
	\[
	\mathcal A=(s_j)_{j\in J}
	\in\operatorname{Atl}_{\mathscr K}(\mathcal H)
	\]
	the packaged structure
	\[
	\mathfrak P_{\mathcal A}
	=
	(J,V,\triangleleft,\mathsf L,\mathsf R_{\mathcal A}),
	\]
	where
	\[
	v\triangleleft j
	\quad\Longleftrightarrow\quad
	v\in I_j,
	\]
	and
	\[
	\mathsf L(j)=\mathscr S(I_j),
	\qquad
	\mathsf R_{\mathcal A}(j)=\{s_j\}.
	\]	
	Then the following hold.	
	\begin{enumerate}[label=(\roman*),leftmargin=2.5em]		
		\item There is a total deterministic packaging context $\Xi_{\mathcal H,\mathscr S,\mathscr K}$
		specified by
		\[
		\begin{aligned}
			\mathscr C_{\Xi_{\mathcal H,\mathscr S,\mathscr K}}
			&=J,\\
			\operatorname{Car}_{\Xi_{\mathcal H,\mathscr S,\mathscr K}}
			&=J,\\
			\operatorname{Form}_{\Xi_{\mathcal H,\mathscr S,\mathscr K}}
			&=\operatorname{Atl}_{\mathscr K}(\mathcal H),\\
			\mathsf J_{\Xi_{\mathcal H,\mathscr S,\mathscr K}}
			&=
			\left\{
			(\mathcal A,\mathfrak P_{\mathcal A}):
			\mathcal A\in
			\operatorname{Atl}_{\mathscr K}(\mathcal H)
			\right\}.
		\end{aligned}
		\]		
		Each compatible atlas therefore determines a
		realized packaging system, and the assignment
		\[
		\mathcal A\longmapsto\mathfrak P_{\mathcal A}
		\]
		is injective.
		
		\item Define the specialization-relative completion
		space by
		\[
		\mathcal G_{\mathscr S}(\mathfrak P_{\mathcal A})
		=
		\left\{
		S\in\mathscr S(V):
		S|_{I_j}\in
		\mathsf R_{\mathcal A}(j)
		\text{ for every }j\in J
		\right\}.
		\]		
		Then $\mathcal A$ has an exact extension in
		$\mathscr S(V)$ if and only if
		\[
		\mathcal G_{\mathscr S}(\mathfrak P_{\mathcal A})
		\ne\varnothing.
		\]		
		Equivalently, the universal property for the specified local class is
		\[
		\operatorname{UG}_{\mathscr K}(\mathcal H)
		\quad\Longleftrightarrow\quad
		\forall\mathcal A\in
		\operatorname{Atl}_{\mathscr K}(\mathcal H),
		\quad
		\mathcal G_{\mathscr S}(\mathfrak P_{\mathcal A})
		\ne\varnothing.
		\]		
	\end{enumerate}
\end{theorem}

\begin{proof}	
	The declared coordinates of	$\mathfrak P_{\mathcal A}$
	satisfy \textup{(S1)}--\textup{(S5)}.
	
	The displayed context admits all its package carriers, and its realization relation is the graph of an evaluation defined on the whole formation-data space.	
	Hence the context is total and deterministic.
	
	If
	\[
	\mathfrak P_{\mathcal A}
	=
	\mathfrak P_{\mathcal B},
	\]
	then
	\[
	\{s_j\}=\{t_j\}
	\]
	for every $j$, so $s_j=t_j$ and
	$\mathcal A=\mathcal B$.	
	This proves \textup{(i)}.
	
	By construction, $S\in\mathscr S(V)$ belongs to the displayed completion space precisely when
	\[
	S|_{I_j}=s_j
	\qquad(j\in J).
	\]	
	This is exactly the definition of an exact atlas extension.	
	Universal quantification over the allowed compatible atlases gives the final equivalence, proving \textup{(ii)}.
\end{proof}

\begin{proposition}[Selection meets and mixed global completions]
	\label{PROP:PackagingSelectionCompletionLaws}	
	Keep the preceding cover and restriction system, but allow arbitrary selected families
	\[
	R(j)\subseteq\mathscr S(I_j)
	\]
	over its fixed scaffold.	
	Put
	\[
	\mathcal G_{\mathscr S}(R)
	=
	\left\{
	S\in\mathscr S(V):
	S|_{I_j}\in R(j)
	\text{ for all }j\in J
	\right\}.
	\]	
	For every family $(R^\alpha)_{\alpha\in A}$	of selections,
	\[
	\mathcal G_{\mathscr S}
	\left(
	\bigwedge_{\alpha\in A}R^\alpha
	\right)
	=
	\bigcap_{\alpha\in A}
	\mathcal G_{\mathscr S}(R^\alpha),
	\]
	where meets are taken in the selection lattice of Theorem~\ref{THM:PackagingSelectionFreedom}.	
	Also,
	\[
	\bigcup_{\alpha\in A}
	\mathcal G_{\mathscr S}(R^\alpha)
	\subseteq
	\mathcal G_{\mathscr S}
	\left(
	\bigvee_{\alpha\in A}R^\alpha
	\right),
	\]
	and this inclusion can be strict, even for two compatible singleton selections.
\end{proposition}

\begin{proof}	
	A global structure belongs to the completion space for a coordinatewise intersection if and only if every local trace lies in every corresponding selected family.
	This proves the meet identity.
	
	The join inclusion follows from the	coordinatewise unions.	
	For strictness take
	\[
	V=\{a,b\},
	\qquad
	I_1=\{a\},
	\qquad
	I_2=\{b\},
	\]
	and
	\[
	\mathscr S(I)=\{0,1\}^I
	\]
	with ordinary coordinate restriction.
	
	Let $R^0$ select the constant-zero function	on each chart, and let $R^1$ select the	constant-one function on each chart.
	
	Both singleton atlases are compatible, and
	\[
	\mathcal G_{\mathscr S}(R^0)
	=
	\{(0,0)\},
	\qquad
	\mathcal G_{\mathscr S}(R^1)
	=
	\{(1,1)\}.
	\]
	
	But $R^0\vee R^1$ permits either value independently on each chart, so
	\[
	\mathcal G_{\mathscr S}(R^0\vee R^1)
	=
	\{0,1\}^{\{a,b\}},
	\]
	which also contains the mixed assignments $(0,1)$ and $(1,0)$.
\end{proof}

The distinction between an atlas and its completion is especially important for the nonextendible compatible examples established in Theorem~\ref{THM:PackagingSharpDependencyGluing}
and in the later native obstruction constructions.
In each case the packaged atlas has a realized formation in the context above, while its specialization-relative completion space is empty.

\begin{lemma}[Transfer of supported consequences]
	\label{LEM:PackagingCompatibleConsequenceTransfer}
	Let
	$
	(c_I)_{I\in\mathcal H}
	$
	be a compatible closure atlas on a finite cover.
	Let \(I,J\in\mathcal H\), let
	\[
	A\subseteq I\cap J,
	\qquad
	b\in I\cap J.
	\]
	Then
	\[
	{
		b\in c_I(A)
		\quad\Longleftrightarrow\quad
		b\in c_J(A).
	}
	\]
	More generally, the same conclusion holds when either side is replaced by the restriction of a chart closure to any set containing \(A\cup\{b\}\).
\end{lemma}

\begin{proof}
	Compatibility gives
	\[
	c_I(A)\cap(I\cap J)
	=
	c_J(A)\cap(I\cap J).
	\]
	Since \(b\in I\cap J\), membership of \(b\) is the same on the two sides.
	If a chart closure is restricted to a set containing \(A\cup\{b\}\), then that restriction has value \(c(A)\cap L\) on \(A\), so membership of \(b\) is unchanged.
\end{proof}

\subsection{Exact gluing of descriptor dependencies on simplex boundaries}
\label{SUBSEC:PackagingDependencyGluing}

A simplex-boundary cover consists of all sets obtained by deleting one label from a fixed finite set.
It is a useful test case: implications from different charts can combine to force a consequence forbidden by another chart.
Theorem~\ref{THM:PackagingDependencyCoatomCriterion} gives an exact test for one atlas.
Theorem~\ref{THM:PackagingSharpDependencyGluing} gives a sharp threshold for \emph{all} compatible atlases with nonempty-premise rules of bounded width.
These results isolate the obstruction used in the general-cover theorem.

Here each chart is a set of descriptor labels, not necessarily a set of internal package units.
Their local data are dependency closures rather than prescribed descriptor functions.
For a descriptor vocabulary
\[
\mathscr V=(v_a:U\to Y_a)_{a\in V},
\]
write $D_A$ for the joint descriptor indexed by $A\subseteq V$ and define
\[
\operatorname{cl}_{\mathscr V}(A)
=
\{b\in V:[v_b]\preceq[D_A]\}.
\]
Thus $\operatorname{cl}_{\mathscr V}$ records which descriptor labels are determined by the observations indexed by $A$.
The extension question concerns the closure operators themselves.
A realization is allowed to choose a new common set of units and new codomains.
No fixed-unit, alphabet-bound, or native-formation constraint is imposed.

A \textbf{closure operator} on a finite set $V$ is an extensive, monotone, idempotent map
\[
c:2^V\longrightarrow 2^V.
\]
An implication $A\to b$ is valid when $b\in c(A)$ and nontrivial when $b\notin A$.
An \textbf{implicational basis} is a family of valid implications generating $c$ by forward chaining.
Such a basis always exists, for example by taking all valid nontrivial implications.

A model of an implicational basis is a subset closed under its rules.
Closure operators are ordered pointwise by inclusion.
Empty preservation, $c(\varnothing)=\varnothing$, is assumed only when explicitly stated.
For classical closure-system and implicational-system background, see
\cite{Wille1982ConceptLattices,
	GanterWille1999FCA,
	CaspardMonjardet2003ClosureSystems}.

The finite descriptor realization below uses the classical semantics of functional dependencies and Armstrong relations \cite{BeeriDowdFaginStatman1984Armstrong}. We include the construction and its proof to specify the finite-size convention used later.

\begin{lemma}[Finite descriptor realization: the Armstrong construction]
	\label{LEM:PackagingArmstrongRealization}	
	Every closure operator $c:2^V\to2^V$ on a finite set $V$ admits a finite descriptor realization
	\[
	\mathscr V=(v_a:U\to Y)_{a\in V}
	\]
	such that
	\[
	\operatorname{cl}_{\mathscr V}=c,
	\qquad
	|U|=|Y|=|\operatorname{Fix}(c)|.
	\]
\end{lemma}

\begin{proof}
	Put $\mathcal F=\operatorname{Fix}(c)$,	and define
	\[
	U=\{u_F:F\in\mathcal F\},
	\qquad
	Y=\{0\}\cup
	\{\lambda_F:F\in\mathcal F\setminus\{V\}\},
	\]
	where all symbols are distinct.
		
	For each $a\in V$, set
	\[
	v_a(u_F)=
	\begin{cases}
		0,&a\in F,\\
		\lambda_F,&a\notin F.
	\end{cases}
	\]
	
	For $F,G\in\mathcal F$, the coordinates	on which the two rows agree are exactly
	\[
	\{a\in V:v_a(u_F)=v_a(u_G)\}
	=
	\begin{cases}
		V,&F=G,\\
		F\cap G,&F\ne G.
	\end{cases}
	\]
	This agreement set is closed under $c$,	since fixed sets of a closure operator are closed under intersections.
	
	Let $A\subseteq V$ and $b\in V$.
	By definition,
	$
	b\in\operatorname{cl}_{\mathscr V}(A)
	$
	if and only if every pair of rows agreeing on $A$ also agrees on $b$.
	
	If $b\in c(A)$, every such pair agrees on $b$, because its agreement set is closed and contains $A$.
	Hence
	$
	c(A)\subseteq
	\operatorname{cl}_{\mathscr V}(A).
	$
	
	Conversely, if $b\notin c(A)$, put $F=c(A)$.
	Then $F\ne V$, and the rows $u_V,u_F$ agree on $A$ but disagree on $b$:
	$
	v_b(u_V)=0,
	\quad
	v_b(u_F)=\lambda_F.
	$
	Thus
	$
	\operatorname{cl}_{\mathscr V}(A)
	\subseteq c(A).
	$
	
	Therefore
	$
	\operatorname{cl}_{\mathscr V}(A)=c(A)
	\quad(A\subseteq V).
	$
	
	Finally, since $V\in\mathcal F$, we have
	\[
	|U|=|\mathcal F|,
	\qquad
	|Y|=1+|\mathcal F\setminus\{V\}|
	=|\mathcal F|.
	\]
\end{proof}

\begin{remark}[Existence of a vocabulary versus extension on fixed units]
	\label{REM:PackagingArmstrongUnitScope}
	The realized vocabulary is compatible with the primitive packaging model:
	choose a carrier $K$, a target with interior $U$ and all refinements as candidates, and admit the presentation classes of vocabulary subfamilies.
	A singleton formation datum whose realization relation contains this target supplies a full system.
	This proves existence in a chosen context.
	For example, the identity closure on labels $\{a,b\}$ cannot be realized on a two-point unit set:
	empty-preservation forces both descriptors to be nonconstant, hence injective, so each would determine the other.
	The lemma realizes that closure on a different finite set.
\end{remark}

\begin{lemma}[Closure restriction and least generated extension]
	\label{LEM:PackagingClosureRestrictionGeneration}	
	Let $V$ be a finite set.	
	\begin{enumerate}[label=(\roman*),leftmargin=2.5em]		
		\item \textbf{Closure restriction.}
		If $c:2^V\to2^V$ is a closure operator and $J\subseteq V$, then
		\[
		c|_J(A):=c(A)\cap J
		\qquad(A\subseteq J)
		\]
		is a closure operator on $J$, with
		\[
		\operatorname{Fix}(c|_J)
		=
		\{F\cap J:F\in\operatorname{Fix}(c)\}.
		\]
		
		\item \textbf{Least generated extension.}
		Let $\mathcal H$ be a finite cover of $V$, and let $c_I$ be a closure operator on	each $I\in\mathcal H$.
		
		Choose an implicational basis $T_I$ of each $c_I$, and let $C_*$ be the	closure on $V$ generated by
		\[
		T=\bigcup_{I\in\mathcal H}T_I.
		\]		
		Then
		\[
		\operatorname{Fix}(C_*)
		=
		\left\{
		F\subseteq V:
		F\cap I\in\operatorname{Fix}(c_I)
		\text{ for every }I\in\mathcal H
		\right\}.
		\]		
		Consequently, $C_*$ is independent of the chosen local bases and is
		the least closure $C$ on $V$ satisfying
		\[
		c_I(A)\subseteq C(A)
		\qquad
		(I\in\mathcal H,\ A\subseteq I).
		\]
		
		\item \textbf{Exact extension criterion.}
		An exact global extension exists if and only if
		\[
		C_*|_I=c_I
		\qquad(I\in\mathcal H).
		\]		
		When this holds, $C_*$ is the least exact global extension.		
	\end{enumerate}
\end{lemma}

\begin{proof}
	We prove the three assertions separately.	
	\begin{enumerate}[label=(\roman*),leftmargin=2.5em]		
		\item \emph{Closure restriction.}
		Extensivity and monotonicity of	$c|_J$ follow from those of $c$.		
		For $A\subseteq J$, monotonicity and idempotence of $c$ give
		\[
		\begin{aligned}
			(c|_J)^2(A)
			&=
			c(c(A)\cap J)\cap J\\
			&\subseteq
			c(c(A))\cap J\\
			&=
			c(A)\cap J.
		\end{aligned}
		\]
		The reverse inclusion follows from extensivity.
		Thus $c|_J$ is idempotent.		
		If $F\in\operatorname{Fix}(c)$, then
		\[
		F\cap J
		\subseteq
		c(F\cap J)\cap J
		\subseteq
		c(F)\cap J
		=
		F\cap J.
		\]
		Hence $F\cap J$ is fixed by $c|_J$.		
		Conversely, if $B\in\operatorname{Fix}(c|_J)$, then $F:=c(B)$ is fixed by $c$ and
		\[
		B=c(B)\cap J=F\cap J.
		\]
		This proves the fixed-set identity.
		
		\item \emph{Least generated extension.}
		A subset $F\subseteq V$ is a model of $T$ if and only if its restriction $F\cap I$ is a model of $T_I$ for every $I\in\mathcal H$.		
		Since $T_I$ generates $c_I$, this is equivalent to
		\[
		F\cap I\in\operatorname{Fix}(c_I)
		\qquad(I\in\mathcal H).
		\]
		The models of $T$ are precisely	the fixed sets of $C_*$, proving
		the displayed characterization.
		
		This characterization depends only on the local closures, not on their
		chosen bases.
		Hence $C_*$ is independent of those bases.
		
		Every implication valid in $c_I$ is valid in $C_*$, so
		\[
		c_I(A)\subseteq C_*(A)
		\qquad
		(I\in\mathcal H,\ A\subseteq I).
		\]
		
		Conversely, let $C$ be any closure on $V$ satisfying these inclusions.
		Every rule in $T_I$ is then valid in $C$, for every $I$.
		Thus every rule of $T$ is valid	in $C$.
		
		By minimality of the closure generated by $T$,
		\[
		C_*(A)\subseteq C(A)
		\qquad(A\subseteq V).
		\]
		Hence $C_*$ is the least closure preserving every local consequence.
		
		\item \emph{Exact extension criterion.}
		By \textup{(ii)}, for every $I\in\mathcal H$ and $A\subseteq I$,
		\[
		c_I(A)
		\subseteq
		C_*(A)\cap I
		=
		(C_*|_I)(A).
		\]
		
		Suppose an exact extension $C$ exists.
		Then
		\[
		C|_I=c_I
		\qquad(I\in\mathcal H).
		\]
		Since $C$ preserves every local consequence, \textup{(ii)} gives
		$C_*\le C$ pointwise.
		Therefore
		\[
		c_I
		\le C_*|_I
		\le C|_I
		=
		c_I.
		\]
		Hence $C_*|_I=c_I$ for every chart.
		
		Conversely, if these equalities	hold, then $C_*$ itself is an
		exact global extension.
		
		Its leastness among exact extensions follows from \textup{(ii)}.		
	\end{enumerate}
\end{proof}

Fix a label set \(V\) of size \(n\ge3\) and the simplex-boundary cover
\[
V_x=V\setminus\{x\}\qquad(x\in V).
\]
A family of closure operators \(c_x:2^{V_x}\to2^{V_x}\) is \emph{compatible} when
\begin{equation}
	\label{EQ:PackagingDependencyCompatibility}
	c_x(A)\cap V_y=c_y(A)\cap V_x
	\quad(A\subseteq V_x\cap V_y).
\end{equation}
An \emph{exact extension} is a closure operator \(C\) on \(V\) satisfying \(C(A)\cap V_x=c_x(A)\) for all \(A\subseteq V_x\).
Choose an implicational basis \(\mathcal B_x\) for each \(c_x\), with rules \(A\to b\), where \(A\subseteq V_x\) and \(b\in V_x\setminus A\).
Empty premises are allowed here.
Write $C_*$ for the closure generated by the union of these bases.
Lemma~\ref{LEM:PackagingClosureRestrictionGeneration} proves its leastness and the exact-extension criterion.
Compatibility transfers a valid local rule to every chart containing its premise and conclusion, by applying \eqref{EQ:PackagingDependencyCompatibility} on their common support.
The next proof follows a derivation until it first produces a new conclusion inside the original chart.

\begin{lemma}[Two-rule obstruction on a simplex boundary]
	\label{LEM:PackagingTwoRuleObstruction}	
	Let $V$ be finite with $|V|\ge3$, and put
	$
	V_x=V\setminus\{x\}
	~ (x\in V).
	$
	Let $(c_x)_{x\in V}$ be a compatible closure atlas on this cover.
	Choose an arbitrary implicational basis	$\mathcal B_x$ for each $c_x$, and put
	\[
	T=\bigcup_{x\in V}\mathcal B_x.
	\]	
	The following are equivalent:	
	\begin{enumerate}[label=(\roman*),leftmargin=2.5em]		
		\item The atlas admits no exact	global extension.
		
		\item There exist distinct $x,y\in V$ and two implications in $T$,
		\[
		A\to x,
		\qquad
		B\cup\{x\}\to y,
		\]
		where
		$
		A,B\subseteq V\setminus\{x,y\},
		$
		such that
		\begin{equation}
			\label{EQ:PackagingTwoRuleSupport}
			A\cup B=V\setminus\{x,y\},
			\qquad
			c_x(A\cup B)=A\cup B.
		\end{equation}		
	\end{enumerate}
\end{lemma}

\begin{proof}
	Let $C_*$ be the closure generated by $T$.
	
	By Lemma~\ref{LEM:PackagingClosureRestrictionGeneration}, an exact extension exists precisely when
	$
	C_*|_{V_x}=c_x
	\quad(x\in V).
	$
	
	Compatibility has the following consequence:
	if a rule $P\to b$ belongs to an original local basis and
	$
	P\cup\{b\}\subseteq V_z,
	$
	then it is valid in $c_z$.
	Indeed, compatibility transfers the rule from its source chart to $V_z$.
	
	\begin{itemize}
		\item \emph{(i)$\Rightarrow$(ii).}
		
		Suppose exact extension fails.
		Then, for some $x\in V$, there is $S\subseteq V_x$ such that
		\[
		C_*(S)\cap V_x
		\not\subseteq c_x(S).
		\]
		Put
		$
		F=c_x(S).
		$
		
		Since
		$
		S\subseteq F\subseteq C_*(S),
		$
		monotonicity and idempotence give
		$
		C_*(F)=C_*(S).
		$
		Thus forward chaining from $F$ eventually derives an element outside $F$ but inside $V_x$.
		
		Stop a derivation at its first such conclusion $y\in V_x\setminus F$.
		
		Before this step, the only possible	new label is $x$.
		Moreover, $x$ must have been derived:
		otherwise the rule deriving $y$ would be supported inside $V_x$ and hence
		valid in $c_x$, contradicting $c_x(F)=F$.
		
		The first derivation of $x$ therefore uses an original rule
		\[
		A\to x,
		\qquad A\subseteq F.
		\]
		The rule producing $y$ must use $x$, since otherwise it would again be
		valid in $c_x$.
		It consequently has the form
		\[
		B\cup\{x\}\to y,
		\qquad B\subseteq F.
		\]
		
		In particular,
		\[
		A,B\subseteq F
		\subseteq V\setminus\{x,y\}.
		\]
		
		We claim that
		\[
		A\cup B\cup\{x,y\}=V.
		\]
		
		Otherwise, choose
		\[
		z\in V\setminus
		(A\cup B\cup\{x,y\}).
		\]
		Both displayed rules are supported inside $V_z$.
		By compatibility, both are valid in $c_z$, so their composition gives
		\[
		y\in c_z(A\cup B).
		\]
		
		Since
		\[
		A\cup B\cup\{y\}
		\subseteq V_x\cap V_z,
		\]
		compatibility yields
		\[
		y\in c_x(A\cup B)
		\subseteq c_x(F)=F,
		\]
		contradicting $y\notin F$.
		
		Therefore
		$
		A\cup B=V\setminus\{x,y\}.
		$
		Together with
		\[
		A\cup B\subseteq F
		\subseteq V\setminus\{x,y\},
		\]
		this gives
		$
		F=A\cup B.
		$
		
		Since $F$ is closed under $c_x$, we have
		$
		c_x(A\cup B)=A\cup B.
		$
		Thus \textup{(ii)} holds.
		
		\item \emph{(ii)$\Rightarrow$(i).}
		
		Conversely, suppose the two rules in \textup{(ii)} exist, and put
		\[
		D=A\cup B=V\setminus\{x,y\}.
		\]
		
		Starting from $D$, the first rule derives $x$, and the second then
		derives $y$.
		Hence
		$
		y\in C_*(D).
		$
		
		However,
		$
		D\subseteq V_x,
		\quad
		c_x(D)=D,
		\quad
		y\notin D.
		$
		Thus
		$
		C_*(D)\cap V_x
		\ne c_x(D).
		$
		
		By the least-generated-extension criterion, no exact global extension
		exists.
	\end{itemize}
\end{proof}

Eliminating a variable from a definite-Horn theory is a classical Horn-forgetting operation; see \cite{Liberatore2025DefiniteHornForgetting} for related elimination methods. The following lemma proves the precise model-projection identity used below.

\begin{lemma}[One-pivot definite-Horn elimination and closure projection]
	\label{LEM:PackagingNativeOnePivotHorn}
	\label{LEM:PackagingHornProjection}	
	Let $W$ be finite, and let $T$ be a finite set of nontrivial implications
	\[
	A\to b,
	\qquad
	\varnothing\ne A\subseteq W\setminus\{b\}.
	\]
	Write $C_T$ for the closure generated by $T$.	
	Fix $v\in W$, put $J=W\setminus\{v\}$, and let
	\[
	T[J]
	=
	\{A\to b\in T:A\cup\{b\}\subseteq J\}.
	\]	
	Define the \textbf{one-pivot resolvent family} $\operatorname{Res}_v(T)$ to consist of all nontrivial implications
	$
	A\cup B\to b
	$
	obtained from pairs of rules
	\[
	A\to v,
	\qquad
	B\cup\{v\}\to b
	\]
	in $T$, where
	$
	A,B\subseteq J,
	\quad
	b\in J\setminus(A\cup B).
	$
	
	Then
	\[
	\operatorname{Mod}
	\bigl(T[J]\cup\operatorname{Res}_v(T)\bigr)
	=
	\{F\setminus\{v\}:F\models T\}.
	\]	
	Consequently,
	\[
	C_T|_J
	=
	C_{\,T[J]\cup\operatorname{Res}_v(T)}.
	\]
\end{lemma}

\begin{proof}
	Put
	$
	T_v=T[J]\cup\operatorname{Res}_v(T).
	$
	\begin{itemize}
		\item \emph{Projection.}
		
		Let $F\models T$ and set $S=F\setminus\{v\}$.
		Clearly $S\models T[J]$.
		
		Consider a resolvent
		$
		A\cup B\to b
		$
		arising from
		\[
		A\to v,
		\qquad
		B\cup\{v\}\to b.
		\]
		If $A\cup B\subseteq S$, then $A\subseteq F$ forces $v\in F$, and the second rule forces $b\in F$.
		
		Since $b \in J$, we have $b\in S$.		
		Thus $S\models\operatorname{Res}_v(T)$,	and hence
		\[
		\{F\setminus\{v\}:F\models T\}
		\subseteq
		\operatorname{Mod}(T_v).
		\]
		
		\item \emph{Lifting.}
		
		Conversely, let $S\models T_v$.
		Define
		\[
		F_S=
		\begin{cases}
			S\cup\{v\},
			&\text{if some }(A\to v)\in T
			\text{ satisfies }A\subseteq S,\\
			S,
			&\text{otherwise}.
		\end{cases}
		\]
		We show that $F_S\models T$.
		
		Every rule supported in $J$ is satisfied because $S\models T[J]$.
		
		Every rule with conclusion $v$ is also satisfied by the definition of $F_S$.
		
		It remains to consider a rule
		$
		B\cup\{v\}\to b
		$
		whose premise is contained in $F_S$.
		Then $v\in F_S$ and $B\subseteq S$.
		
		By the definition of $F_S$, there is a rule
		$
		A\to v
		$
		in $T$ with $A\subseteq S$.
		
		If $b\in A\cup B$, then $b\in S$.
		Otherwise,
		$
		A\cup B\to b
		$
		belongs to $\operatorname{Res}_v(T)$.
		
		Since $S\models T_v$ and $A\cup B\subseteq S$, again $b\in S$.		
		Thus every rule of $T$ is satisfied	by $F_S$.
		
		Since
		$
		F_S\setminus\{v\}=S,
		$
		we obtain the reverse inclusion
		\[
		\operatorname{Mod}(T_v)
		\subseteq
		\{F\setminus\{v\}:F\models T\}.
		\]
		
		Therefore
		\[
		\operatorname{Mod}(T_v)
		=
		\{F\setminus\{v\}:F\models T\}.
		\]
		
		By Lemma~\ref{LEM:PackagingClosureRestrictionGeneration}, the right-hand side is $\operatorname{Fix}(C_T|_J)$, while the left-hand side is	$\operatorname{Fix}(C_{T_v})$.
		
		Hence
		$
		C_T|_J=C_{T_v},
		$
		as claimed.
	\end{itemize}
\end{proof}

\begin{theorem}[Coatom criterion and quadratic oracle test]
	\label{THM:PackagingDependencyCoatomCriterion}	
	Let $V$ be finite with $n=|V|\ge3$,	and let $(c_x)_{x\in V}$ be a compatible closure atlas on the simplex-boundary cover
	\[
	V_x=V\setminus\{x\}
	\qquad(x\in V).
	\]	
	For distinct $x,y\in V$, put
	\[
	D_{xy}=V\setminus\{x,y\},
	\qquad
	q_{xy}
	=
	\mathbf1_{\{y\in c_x(D_{xy})\}},
	\]
	and set $q_{xx}=0$.	
	An exact global extension exists if and only if
	\begin{equation}
		\label{EQ:PackagingDependencyCoatomCriterion}
		q_{yx}=1
		\quad\text{and}\quad
		\max_{z\ne y}q_{zy}=1
		\quad\Longrightarrow\quad
		q_{xy}=1
		\qquad(x\ne y).
	\end{equation}	
	Equivalently, exact extension fails if and only if there exist pairwise distinct $x,y,z\in V$ such that
	\[
	q_{xy}=0,
	\qquad
	q_{yx}=1,
	\qquad
	q_{zy}=1.
	\]	
	Given compatibility, exact extendibility is decidable using $n(n-1)$ local closure-membership queries and $O(n^2)$	additional Boolean operations.	
	When extension fails, the algorithm returns	a certificate consisting of the three displayed coatom evaluations.
\end{theorem}

\begin{proof}
	Let $C_*$ be the least closure generated by the union of the local implicational bases.	
	We first prove the equivalence with the three-entry obstruction and then establish the oracle bound.
	
	\begin{itemize}
		\item \emph{Nonextendibility implies an obstruction.}
		
		Suppose no exact extension exists.
		By Lemma~\ref{LEM:PackagingTwoRuleObstruction},	there are distinct $x,y\in V$ and original local rules
		\[
		A\to x,
		\qquad
		B\cup\{x\}\to y,
		\]
		such that
		\[
		A\cup B=D_{xy},
		\qquad
		c_x(D_{xy})=D_{xy}.
		\]
		
		Since $y\notin D_{xy}$, the latter equality
		gives
		\[
		q_{xy}=0.
		\]
		
		The first rule has
		$A\subseteq D_{xy}\subseteq V_y$.
		Compatibility transfers it to chart $V_y$,
		so
		\[
		x\in c_y(D_{xy}),
		\qquad
		q_{yx}=1.
		\]
		
		Let $V_z$ be an original chart containing
		the second rule.
		Since its support contains $x$ and $y$,
		necessarily
		\[
		z\notin\{x,y\}.
		\]
		Moreover,
		\[
		B\cup\{x\}\subseteq D_{zy}.
		\]
		The second rule therefore gives
		\[
		y\in c_z(D_{zy}),
		\qquad
		q_{zy}=1.
		\]
		
		Thus the required three-entry obstruction exists.
		
		\item \emph{An obstruction prevents exact extension.}
		
		Conversely, suppose pairwise distinct $x,y,z\in V$ satisfy
		\[
		q_{xy}=0,
		\qquad
		q_{yx}=1,
		\qquad
		q_{zy}=1.
		\]
		Put $D=D_{xy}$.
		
		Since $D$ is a coatom of $V_x$, the condition $q_{xy}=0$ gives
		$
		c_x(D)=D.
		$
		
		On the other hand, $q_{yx}=1$ implies
		\[
		x\in c_y(D)\subseteq C_*(D).
		\]
		Consequently,
		\[
		D_{zy}
		\subseteq
		D\cup\{x\}
		\subseteq
		C_*(D).
		\]
		
		Since $q_{zy}=1$, we also have
		$
		y\in c_z(D_{zy}).
		$
		Monotonicity and idempotence of $C_*$ now give
		\[
		y\in
		c_z(D_{zy})
		\subseteq
		C_*(D_{zy})
		\subseteq
		C_*(D).
		\]
		
		Hence
		\[
		y\in C_*(D)\cap V_x,
		\qquad
		y\notin c_x(D).
		\]
		Therefore
		$
		C_*|_{V_x}\ne c_x.
		$		
		By Lemma~\ref{LEM:PackagingClosureRestrictionGeneration}, no exact extension exists.
		
		We have proved that nonextendibility is equivalent to the existence of pairwise distinct $x,y,z$ with
		\[
		q_{xy}=0,
		\qquad
		q_{yx}=q_{zy}=1.
		\]
		Since $q_{xy}=0$ excludes $z=x$ from any witness to
		$\max_{z\ne y}q_{zy}=1$, this is equivalent to
		\eqref{EQ:PackagingDependencyCoatomCriterion}.
		
		\item \emph{Oracle complexity and certificates.}
		
		Evaluate $q_{xy}$ once for every ordered pair $x\ne y$, using exactly
		$
		n(n-1)
		$
		local closure-membership queries.
		
		For each $y \in V$, compute
		$
		m_y=\max_{z\ne y}q_{zy},
		$
		retaining an index $z_y$ attaining the maximum whenever $m_y=1$.		
		This requires $O(n^2)$ Boolean operations.
		
		A further scan of all ordered pairs checks whether
		$$
		q_{xy}=0,
		\quad
		q_{yx}=1,
		\quad
		m_y=1.
		$$		
		If such a pair exists, then	$q_{z_y y}=1$.
		
		Since $q_{xy}=0$, necessarily $z_y\ne x$, so $x,y,z_y$ are pairwise distinct.
		The three oracle evaluations
		\[
		q_{xy}=0,
		\qquad
		q_{yx}=1,
		\qquad
		q_{z_y y}=1
		\]
		provide the asserted obstruction certificate.
		
		If no such pair exists, then the coatom criterion holds, and an exact extension exists.
		
		The total additional Boolean complexity	is $O(n^2)$.
		Compatibility is assumed.
		Its verification is not included in this bound.
	\end{itemize}
\end{proof}

\begin{theorem}[Sharp premise bound for simplex-boundary gluing]
	\label{THM:PackagingSharpDependencyGluing}	
	Let $k\ge1$, let $V$ be finite with $n=|V|\ge3$, and consider the
	simplex-boundary cover
	\[
	V_x=V\setminus\{x\}
	\qquad(x\in V).
	\]	
	The following are equivalent:	
	\begin{enumerate}[label=(\roman*),leftmargin=2.5em]		
		\item Every compatible atlas $(c_x)_{x\in V}$ whose local closures admit implicational bases with nonempty premises of size at most $k$ has an exact global extension.
		
		\item
		$
		n\ge2k+2.
		$		
	\end{enumerate}	
	When these conditions hold, the union of any such local bases generates the least exact extension $C_*$ and is a global implicational basis with nonempty premises of size at most $k$.	
	Conversely, for every $3\le n\le2k+1$, there exists a compatible nonextendible atlas in which each local closure is generated by at most one nontrivial implication with nonempty premise of size at most $k$.
\end{theorem}

\begin{proof}
	We prove the two directions separately.
	\begin{itemize}
		\item \emph{(ii)$\Rightarrow$(i). Sufficiency.}
		
		Suppose $n\ge2k+2$.
		Let $(c_x)_{x\in V}$ be a compatible atlas as in \textup{(i)}, and choose a local basis $\mathcal B_x$ for each $c_x$ satisfying the premise bound.
		
		Let $C_*$ be the closure generated by
		\[
		T=\bigcup_{x\in V}\mathcal B_x.
		\]
		
		If exact extension failed, then
		Lemma~\ref{LEM:PackagingTwoRuleObstruction}	would give distinct $x,y\in V$ and rules in $T$ of the form
		\[
		A\to x,
		\qquad
		B\cup\{x\}\to y,
		\]
		where
		$
		A\cup B=V\setminus\{x,y\}.
		$
		
		Since the premises are nonempty and have size at most $k$,
		\[
		|A|\le k,
		\qquad
		|B|\le k-1.
		\]
		Consequently,
		\[
		\begin{aligned}
			n
			&=|A\cup B|+2\\
			&\le |A|+|B|+2\\
			&\le2k+1,
		\end{aligned}
		\]
		contradicting $n\ge2k+2$.
		
		Thus an exact extension exists.
		By Lemma~\ref{LEM:PackagingClosureRestrictionGeneration}, $C_*$ is the least exact extension.
		
		Moreover, $T$ itself is an implicational basis of $C_*$, and all its premises are nonempty and of size at most $k$.

		\item \emph{(i)$\Rightarrow$(ii). Sharpness by contrapositive.}
		
		Suppose
		$
		3\le n\le2k+1.
		$
		Choose distinct $x,y\in V$ and a partition
		\[
		V=Z\sqcup W\sqcup\{x,y\},
		\]
		where
		$
		|Z|=\min\{k,n-2\},
		~
		|W|=n-2-|Z|.
		$
		Then
		\[
		1\le |Z|\le k,
		\qquad
		0\le |W|\le k-1.
		\]
		
		Consider the two implications
		\[
		\rho_1:Z\to x,
		\qquad
		\rho_2:W\cup\{x\}\to y.
		\]
		Both have nonempty premises of size at most $k$.
		
		For each $u\in V$, let $\mathcal B_u$ consist precisely of those rules among $\{\rho_1,\rho_2\}$ whose premises and conclusions lie in $V_u$.
		Let $c_u$ be the closure generated by $\mathcal B_u$.
		
		The supports of $\rho_1$ and $\rho_2$ have union $V$.
		Hence no proper chart contains both, so
		\[
		|\mathcal B_u|\le1
		\qquad(u\in V).
		\]
		
		We first verify compatibility.
		
		Let $u,v\in V$ be distinct, put
		$
		J=V_u\cap V_v,
		$
		and fix $S\subseteq J$.
		
		Suppose
		$
		b \in c_u(S)\cap J
		$
		and
		$
		b \notin S.
		$
		Since $\mathcal B_u$ contains at most one implication, some rule
		$
		P\to b\in\mathcal B_u
		$
		satisfies $P\subseteq S$.		
		Therefore
		$
		P\cup\{b\}\subseteq J\subseteq V_v,
		$
		so the same rule belongs to $\mathcal B_v$.
		Hence $b\in c_v(S)$.
		
		By symmetry,
		\[
		c_u(S)\cap J=c_v(S)\cap J.
		\]
		Thus the atlas is compatible.
		
		We now show that it has no exact global extension.
		
		The first rule belongs to $\mathcal B_y$, and the second belongs to $\mathcal B_z$ for any $z\in Z$.
		Thus both rules occur in the union of the original local bases.
		
		Put
		$
		D=Z\cup W=V\setminus\{x,y\}.
		$
		Since neither rule is supported inside $V_x$, we have
		\[
		\mathcal B_x=\varnothing,
		\qquad
		c_x(D)=D.
		\]
		
		On the other hand, the union of the local bases contains
		\[
		Z\to x,
		\qquad
		W\cup\{x\}\to y.
		\]
		Starting from $D$, these rules successively derive $x$ and $y$.
		Hence
		$
		y\in C_*(D),
		$
		where $C_*$ is the closure generated by the union of the local bases.
		
		But
		$
		y\in V_x\setminus c_x(D),
		$
		so
		\[
		C_*(D)\cap V_x\ne c_x(D).
		\]
		
		By Lemma~\ref{LEM:PackagingClosureRestrictionGeneration}, no exact global extension exists.
		
		This constructs the required compatible counterexample for every $3\le n\le2k+1$ and proves that \textup{(i)} implies \textup{(ii)}.
	\end{itemize}
\end{proof}

\begin{corollary}[Strict hierarchy and descriptor realization]
	\label{COR:PackagingDependencyHierarchy}	
	For $k\ge1$, let $\mathcal U_k$ denote the class of finite covers $\mathcal H$ of finite label sets $V$ such that every compatible
	closure atlas $(c_I)_{I\in\mathcal H}$ admits an exact global extension whenever each $c_I$ has an implicational basis whose rules have nonempty premises of size at most $k$.	
	Then:	
	\begin{enumerate}[label=(\roman*),leftmargin=2.5em]		
		\item \textbf{Strict hierarchy.}
		For every $k\ge1$, we have
		$
		\mathcal U_{k+1}
		\subsetneq
		\mathcal U_k.
		$
		The simplex boundary on $2k+2$ labels witnesses strictness.
		
		\item \textbf{Common finite descriptor realization.}
		Let $\mathcal H\in\mathcal U_k$ and let $(c_I)_{I\in\mathcal H}$ be a compatible atlas satisfying the stated premise bound.		
		Let $C_*$ be the closure generated by the union of local implicational bases.
		Then $C_*$ is the least exact extension, and there exists a finite descriptor vocabulary
		$
		\mathscr V=(v_a:U\to Y)_{a\in V}
		$
		such that
		\[
		\operatorname{cl}_{\mathscr V}=C_*,
		\qquad
		\operatorname{cl}_{\mathscr V|_I}=c_I
		\quad(I\in\mathcal H),
		\]
		where
		$
		\mathscr V|_I=(v_a)_{a\in I}.
		$
		The vocabulary may be chosen with
		$
		|U|=|Y|=|\operatorname{Fix}(C_*)|.
		$
	\end{enumerate}
\end{corollary}

\begin{proof}
	We prove the two assertions separately.	
	\begin{enumerate}[label=(\roman*),leftmargin=2.5em]		
		\item \emph{Strict hierarchy.}
		
		Every atlas admitting local bases with nonempty premises of size at most $k$ also satisfies the bound $k+1$.
		Hence
		$
		\mathcal U_{k+1}\subseteq\mathcal U_k.
		$
				
		Let $\mathcal H_k$ be the simplex-boundary cover on $2k+2$ labels.		
		By Theorem~\ref{THM:PackagingSharpDependencyGluing}, every compatible atlas on $\mathcal H_k$ satisfying the premise bound $k$ admits an exact extension.
		Thus
		$
		\mathcal H_k\in\mathcal U_k.
		$
		
		The same theorem, applied with bound $k+1$, gives a compatible nonextendible atlas on $\mathcal H_k$, since
		$
		2k+2<2(k+1)+2.
		$
		Therefore
		$
		\mathcal H_k\notin\mathcal U_{k+1}.
		$
		The inclusion is strict.
		
		\item \emph{Common finite descriptor realization.}
		
		Let $\mathcal H\in\mathcal U_k$ and choose local bases satisfying the
		premise bound.
		
		By the definition of $\mathcal U_k$, the compatible atlas admits an exact global extension.
		
		Lemma~\ref{LEM:PackagingClosureRestrictionGeneration} therefore gives
		\[
		C_*|_I=c_I
		\qquad(I\in\mathcal H),
		\]
		with $C_*$ the least exact extension.
		
		Apply Lemma~\ref{LEM:PackagingArmstrongRealization}	to $C_*$.
		This yields a finite descriptor vocabulary
		$\mathscr V=(v_a:U\to Y)_{a\in V}$ satisfying
		\[
		\operatorname{cl}_{\mathscr V}=C_*,
		\qquad
		|U|=|Y|=|\operatorname{Fix}(C_*)|.
		\]
		
		For every chart $I\in\mathcal H$ and $A\subseteq I$, restricting the vocabulary to the labels in $I$ gives
		\[
		\begin{aligned}
			\operatorname{cl}_{\mathscr V|_I}(A)
			&=
			\operatorname{cl}_{\mathscr V}(A)\cap I\\
			&=
			C_*(A)\cap I\\
			&=
			c_I(A).
		\end{aligned}
		\]		
		Hence one finite descriptor vocabulary realizes every local dependency closure exactly.		
	\end{enumerate}
\end{proof}

The bound concerns the size of premises in a generating basis.
Iterated consequences may require more original observations.
Empty-preserving unary bases give preorders, recovering the simplex-boundary case of Theorem~\ref{THM:PackagingTransitiveUniversalAmalgamation}.
Binary bases already fail on the four- and five-label boundaries and always extend on boundaries with at least six labels.
The hierarchy therefore measures an exact change in gluing behavior as joint descriptor dependence increases.

\subsection{General covers and bounded determination}
\label{SUBSEC:PackagingGeneralDependencyGluing}

For a general finite cover, there are two possible geometric obstructions: induced cycles and cliques not contained in any chart.
We control both by limiting how many labels are needed to determine a consequence.
This leads to a propagation-recovery rank and an elimination theorem.
Theorem~\ref{THM:PackagingGeneralCaratheodoryGluing} gives the sharp globalization criterion, while Theorem~\ref{THM:PackagingCanonicalCliqueCompletion} shows what can still be completed without global chordality.

Let $\mathcal H$ be a finite cover of a finite label set $V$.
Charts may be empty or repeated, and compatibility makes repeated chart closures identical.
For the geometric statements a clique is nonempty.
The empty label set is allowed: all its closures and extension questions are trivial.
The primal graph joins two distinct labels belonging to a common chart.
A set is \emph{guarded} when it is contained in a chart.
A graph is \emph{chordal} when it has no induced cycle of length at least four.
Every chart is a clique if it is nonempty.

Viewing the cover $\mathcal H$ as a hypergraph on $V$, we say that it is \textbf{$N$-conformal} when every clique of the primal graph of cardinality at most $N$ is guarded.
This is the standard bounded-conformality terminology used in hypergraph theory.
Full conformality requires every primal clique to be guarded~\cite{Otto2012HighlyAcyclic}.
Full conformality together with chordality is the classical acyclicity condition for database hypergraphs
\cite{BeeriFaginMaierYannakakis1983Acyclic,Fagin1983Degrees}.
The results below deliberately use full chordality but only a sharp bounded level of conformality.
We use the standard chordal-graph and clique-tree theory \cite{Gavril1974Chordal,HoLee1989CountingCliqueTrees,HaraTakemura2006BoundaryCliques}.
Compatibility means $c_I|_{I\cap J}=c_J|_{I\cap J}$, and an exact extension satisfies $C|_I=c_I$.
Restriction and generation are as in Lemma~\ref{LEM:PackagingClosureRestrictionGeneration}.
The arguments below classify dependency-closure data.

A closure has \emph{Carath\'eodory number at most \(k\)} when
\begin{equation}
	\label{EQ:PackagingCaratheodoryBound}
	b\in c(S)\quad\Longrightarrow\quad
	\exists A\subseteq S\quad |A|\le k,\quad b\in c(A).
\end{equation}
This standard invariant bounds the sizes of minimal determining sets \cite[Section~2]{NourineVilmin2021ConsistentClosedSets}.
Throughout this subsection \(k\ge1\) and closures preserve the empty set.
For a vocabulary this excludes descriptor labels already determined by the empty observation.
It is an additional hypothesis, not an axiom of packaging.
Thus nontrivial determining sets are nonempty.
This bound is hereditary under restriction.
A bound on premises in a generating basis need not be.

\begin{lemma}[Determination and restrictions of bases]
	\label{LEM:PackagingHereditaryPremiseBound}
	A finite empty-preserving closure \(c\) has Carath\'eodory number at most \(k\) if and only if every restriction \(c|_S\) admits an implicational basis with nonempty premises of size at most \(k\).
\end{lemma}

\begin{proof}
	Under \eqref{EQ:PackagingCaratheodoryBound}, take all valid implications \(A\to b\) supported in \(S\), with \(1\le |A|\le k\).
	Every consequence has such a witness in its input, so these implications generate \(c|_S\).
	Conversely, let \(A\) be an inclusion-minimal determining set for \(b\notin A\).
	In the restriction to \(A\cup\{b\}\), a derivation of \(b\) from \(A\) using a basis of premise size at most \(k\) must first add \(b\) by a rule \(D\to b\) with \(D\subseteq A\) and \(|D|\le k\).
	Minimality gives \(D=A\), and therefore \(|A|\le k\).
	Now let \(b\in c(S)\).
	If \(b\in S\), the singleton \(\{b\}\) is a witness of size one.
	If \(b\notin S\), choose an inclusion-minimal \(A\subseteq S\) with \(b\in c(A)\).
	The preceding argument gives \(|A|\le k\).
	Hence \(c\) has Carath\'eodory number at most \(k\).
\end{proof}

For example, the binary basis \(\{a_1,a_2\}\to u,\ \{u,a_3\}\to b\) requires all three labels \(a_1,a_2,a_3\) to determine \(b\).
Its restriction deleting \(u\) has no binary basis.
The following theorem uses the hereditary bound \eqref{EQ:PackagingCaratheodoryBound}.

\begin{lemma}[Elimination with a prescribed clique last]
	\label{LEM:PackagingCliqueLastElimination}
	If $G$ is a finite chordal graph and $I$ is a clique or the empty set, there is an ordering of its vertices in which every vertex has clique neighborhood among the later vertices, and the vertices of $I$ occur last.
\end{lemma}

\begin{proof}
	If $V(G)=\varnothing$, the empty ordering satisfies the assertion.
	If $G$ is disconnected, process its connected components separately, taking the component containing $I$ last when $I\ne\varnothing$.
	It therefore suffices to consider a connected chordal graph.
	
	Choose a maximal clique $K_0$ containing $I$, and root a clique tree at $K_0$ \cite{Gavril1974Chordal}.
	When $I=\varnothing$, choose any root.
		
	Repeatedly take a leaf clique $K\ne K_0$ and eliminate the vertices of $K$ that do not belong to its parent clique.
		
	By the running-intersection property, each such vertex belongs to no other remaining maximal clique.
	Its remaining neighbours therefore lie in $K$ and form a clique.
		
	After removing these vertices, discard the redundant leaf clique and continue toward the root.
		
	Finally, eliminate the vertices of $K_0\setminus I$, followed by those of $I$.
	Since $K_0$ is a clique, these final eliminations also have clique later
	neighbourhoods.
	
	The resulting ordering is perfect and places every vertex of $I$ last.
\end{proof}

The globalization argument combines two ingredients: one-pivot Horn elimination from Lemma~\ref{LEM:PackagingHornProjection} and clique-last perfect elimination from
Lemma~\ref{LEM:PackagingCliqueLastElimination}.

Eliminating a label can create new Horn implications.
Chordality controls where the labels involved in those implications can occur.

The key question is whether each new resolvent is already implied by a short rule from an original chart.
If so, elimination does not need to add a new rule.

The following recovery rank measures the guarding level required to guarantee this property uniformly.
It belongs to the closure-theoretic specialization, not to the primitive packaging axioms.

\begin{definition}[Chart-local rule scheme and local propagation-recovery rank]
	\label{DEF:PackagingPropagationRecoveryScheme}
	Let
	\[
	\mathscr C=(\mathscr C(I))_{I}
	\]
	be a family of classes of finite empty-preserving closure operators, stable under restriction.
	A \textbf{chart-local rule scheme} for \(\mathscr C\) assigns to every \(c\in\mathscr C(I)\) a family
	\[
	\mathscr R(c)
	\]
	of nontrivial valid implications
	\[
	A\to b,
	\qquad
	\varnothing\ne A\subseteq I\setminus\{b\},
	\qquad
	b\in c(A),
	\]
	satisfying:
	\begin{enumerate}[label=(R\arabic*),leftmargin=3em]
		\item \textbf{local generation:}
		\(\mathscr R(c)\) generates \(c\);
		
		\item \textbf{trace locality:}
		if \(J\subseteq I\) and	\(A\cup\{b\}\subseteq J\), then
		\[
		A\to b\in\mathscr R(c)
		\quad\Longleftrightarrow\quad
		A\to b\in\mathscr R(c|_J).
		\]
	\end{enumerate}	
	Fix an integer \(m\ge2\).
	The pair \((\mathscr C,\mathscr R)\) is called \textbf{locally \(m\)-propagation-recoverable} when the following holds.
	For every finite cover \(\mathcal H\) in which every nonempty primal clique of cardinality at most \(m\) is guarded, with no chordality assumption, for every compatible atlas
	\[
	(c_I)_{I\in\mathcal H},
	\qquad
	c_I\in\mathscr C(I),
	\]
	and for every pair of original-chart rules
	\[
	A\to v\in\mathscr R(c_I),
	\qquad
	B\cup\{v\}\to b\in\mathscr R(c_J),
	\]
	with
	\[
	v,b\notin A\cup B,
	\qquad
	Q=A\cup B\cup\{v,b\}
	\]
	a primal clique, there exist an original chart \(K\in\mathcal H\) and a nonempty set
	\[
	D\subseteq A\cup B
	\]
	such that
	\[
	{D\to b\in\mathscr R(c_K).}
	\]
	Thus local recovery depends on the closure class, the chosen chart-local rule scheme, and the guarding geometry.
	Chordality is not assumed.	
	The \textbf{local propagation-recovery rank} is
	\[
	{
		\operatorname{prk}(\mathscr C,\mathscr R)
		=
		\min\{m\ge2:(\mathscr C,\mathscr R)
		\text{ is locally }m\text{-propagation-recoverable}\},
	}
	\]
	with value \(\infty\) if no such \(m\) exists.
	Recoverability is monotone in the guarding level:
	if it holds for \(m\), then it holds for every \(m'\ge m\), since the latter hypothesis guards at least the same cliques.
\end{definition}

This rank depends on which chart-local rule scheme \(\mathscr R\) we choose.
Later we remove that dependence by working with minimal closure consequences.

\begin{lemma}[Recovered resolvents are redundant]
	\label{LEM:PackagingRecoveredResolventsRedundant}
	Let \(\mathcal H\) be a finite cover of \(V\), and let
	\[
	T=\bigcup_{I\in\mathcal H}\mathscr R(c_I)
	\]
	be the union of the original-chart rule families of a compatible atlas.
	For \(W\subseteq V\), write
	\[
	T[W]=\{A\to b\in T:A\cup\{b\}\subseteq W\}.
	\]
	Let \(v\in W\).
	Assume that every nontrivial Horn resolvent across \(v\),
	\[
	A\to v,
	\qquad
	B\cup\{v\}\to b
	\qquad\text{in }T[W],
	\]
	has an original-chart recovery
	\[
	D\to b\in T[W\setminus\{v\}],
	\qquad
	\varnothing\ne D\subseteq A\cup B.
	\]
	Then projection across \(v\) creates no new necessary rule:
	\[
	{
		\{F\setminus\{v\}:F\models T[W]\}
		=
		\operatorname{Mod}(T[W\setminus\{v\}]).
	}
	\]
\end{lemma}

\begin{proof}
	By Lemma~\ref{LEM:PackagingHornProjection}, the projected model set is the model set of the retained theory \(T[W\setminus\{v\}]\) together	with all nontrivial raw resolvents \(A\cup B\to b\).
	For each such resolvent the assumed recovered rule \(D\to b\) already belongs to the retained original-chart theory, and \(D\subseteq A\cup B\).
	Hence every	model of the retained theory satisfies the raw resolvent.
	All raw	resolvents are therefore redundant, proving the equality.
\end{proof}

\begin{theorem}[Packaging controlled-propagation globalization and clique localization]
	\label{THM:PackagingControlledPropagationGlobalization}	
	Let $(\mathscr C,\mathscr R)$ be locally $m$-propagation-recoverable.
	Let $\mathcal H$ be a finite cover whose primal graph
	is chordal and whose nonempty cliques of cardinality at most $m$ are guarded.	
	For a compatible atlas
	\[
	(c_I)_{I\in\mathcal H},
	\qquad c_I\in\mathscr C(I),
	\]
	put
	\[
	T=\bigcup_{I\in\mathcal H}\mathscr R(c_I),
	\]
	and let $C_*$ be the closure generated by $T$ on $V$.
	For $W\subseteq V$, define
	\[
	T[W]
	=
	\{A\to b\in T:A\cup\{b\}\subseteq W\}.
	\]	
	Then the following hold:
	\begin{enumerate}[label=(\roman*),leftmargin=2.5em]
		\item \emph{Suffix projection.}
		For any primal clique $Q$, including $Q=\varnothing$, choose a perfect elimination ordering with $Q$ last.
		For every suffix $W$ obtained while eliminating	vertices outside $Q$,
		\[
		\{F\cap W:F\models T\}
		=
		\operatorname{Mod}(T[W]).
		\]
		Thus every projected implication can be recovered from original-chart rules supported on the current suffix, even when that suffix is not an original chart.
		
		\item \emph{Clique locality.}
		For every primal clique $Q$, including $Q=\varnothing$,
		\[
		C_*|_Q=C_{T[Q]},
		\]
		where $C_{T[Q]}$ is the closure generated by the original-chart rules supported on $Q$.
		
		\item \emph{Exact globalization and minimality.}
		The closure $C_*$ is empty-preserving and is the least global closure preserving every local consequence.
		Moreover,
		\[
		C_*|_I=c_I
		\qquad(I\in\mathcal H).
		\]
		In particular, every compatible atlas admits an	exact global extension.
	\end{enumerate}
\end{theorem}

\begin{proof}
	\leavevmode
	\begin{enumerate}[label=(\roman*),leftmargin=2.5em]		
		\item \emph{Suffix projection.}
		Fix an arbitrary primal clique $Q$, possibly empty.
		By Lemma~\ref{LEM:PackagingCliqueLastElimination}, choose a perfect elimination ordering with $Q$ last.
		
		We prove inductively that, for each suffix $W$ obtained while eliminating vertices outside $Q$,
		\[
		\{F\cap W:F\models T\}
		=
		\operatorname{Mod}(T[W]).
		\tag{*}
		\]
		This holds initially at $W=V$, since $T[V]=T$.
		
		Suppose $(*)$ holds at $W$, and let $v\in W\setminus Q$ be the next eliminated vertex.
		Consider a nontrivial Horn resolvent across $v$	from two rules in $T[W]$:
		\[
		A\to v,\qquad B\cup\{v\}\to b.
		\]
		Both parent rules are original-chart rules.
		Every label in $A\cup B\cup\{b\}$ is a remaining primal neighbor of $v$.
		Since $v$ is simplicial in $G_{\mathcal H}[W]$,	the combined support
		\[
		A\cup B\cup\{v,b\}
		\]
		is a primal clique.
		
		Local $m$-propagation-recoverability supplies an original-chart rule
		\[
		D\to b,
		\qquad
		\varnothing\ne D\subseteq A\cup B.
		\]
		Its complete support is contained in $W\setminus\{v\}$, so it belongs to $T[W\setminus\{v\}]$.
		
		Lemma~\ref{LEM:PackagingRecoveredResolventsRedundant}
		therefore gives
		\[
		\{F\setminus\{v\}:F\models T[W]\}
		=
		\operatorname{Mod}(T[W\setminus\{v\}]).
		\]
		Combining this equality with the induction hypothesis establishes $(*)$ for the next suffix.
		This proves~\textup{(i)}.
		
		\item \emph{Clique locality.}
		Taking $W=Q$ in~\textup{(i)} gives
		\[
		\{F\cap Q:F\in\operatorname{Fix}(C_*)\}
		=
		\operatorname{Mod}(T[Q]).
		\]
		By Lemma~\ref{LEM:PackagingClosureRestrictionGeneration},
		the left-hand side is $\operatorname{Fix}(C_*|_Q)$, while the right-hand side is $\operatorname{Fix}(C_{T[Q]})$.
		
		Since finite closure operators are determined by their fixed-point families,
		\[
		C_*|_Q=C_{T[Q]}.
		\]
		This includes $Q=\varnothing$, with its	unique closure.
		
		\item \emph{Exact globalization and minimality.}
		Every rule of $T$ has nonempty premise, so $C_*(\varnothing)=\varnothing$.
		Local generation and
		Lemma~\ref{LEM:PackagingClosureRestrictionGeneration}
		show that $C_*$ is the least global closure	preserving all local consequences.
		
		Let $I_0$ be any nonempty original chart.
		Since $I_0$ is a primal clique, part~\textup{(ii)} gives
		\[
		C_*|_{I_0}=C_{T[I_0]}.
		\]
		By trace locality and compatibility, every rule in $T[I_0]$ is valid in $c_{I_0}$,	while $T[I_0]$ contains the generating family $\mathscr R(c_{I_0})$.
		Consequently,
		\[
		C_{T[I_0]}=c_{I_0}.
		\]
		For an empty original chart, the equality is automatic.		
		Hence
		\[
		C_*|_I=c_I
		\qquad(I\in\mathcal H),
		\]
		proving exact globalization.
	\end{enumerate}
\end{proof}

\begin{theorem}[Canonical simultaneous clique completion without chordality]
	\label{THM:PackagingCanonicalCliqueCompletion}
	Let $(\mathscr C,\mathscr R)$ be locally 	$m$-propagation-recoverable, and let $\mathcal H$ be a finite cover of $V$ on which every nonempty primal clique of cardinality at most $m$ is guarded.
	No chordality assumption is imposed.
	
	For any compatible atlas $\mathcal A=(c_I)_{I\in\mathcal H}$,
	with $c_I\in\mathscr C(I)$, put
	\[
	T=\bigcup_{I\in\mathcal H}\mathscr R(c_I),
	\qquad C_*=C_T,
	\]
	and, for $W\subseteq V$, define
	\[
	T[W]
	=
	\{A\to b\in T:A\cup\{b\}\subseteq W\}.
	\]
	For every primal clique $M$, put
	\[
	\lambda_M=C_{T[M]},
	\]
	with $\lambda_\varnothing$ the unique closure on the empty set.	
	Then the following hold:
	\begin{enumerate}[label=(\roman*),leftmargin=2.5em]		
		\item \emph{Automatic clique fillability.}
		For every primal clique $M$ and original chart $I$,
		\[
		\lambda_M|_{M\cap I}
		=
		c_I|_{M\cap I}.
		\]
		Thus the induced atlas on $M$ admits an exact extension, and $\lambda_M$ is its least exact extension.
		No separate weak-fillability hypothesis is required.
		
		\item \emph{Restriction coherence.}
		For every pair of primal cliques $L\subseteq M$,
		\[
		\lambda_M|_L=\lambda_L.
		\]
		Consequently, for any two primal cliques $M,N$,
		\[
		\lambda_M|_{M\cap N}
		=
		\lambda_N|_{M\cap N}.
		\]
		
		\item \emph{Simultaneous conservative conformalization.}
		Adjoining all maximal primal cliques $M$, equipped with their closures $\lambda_M$, produces a compatible closure atlas on
		a fully conformal cover with the same primal graph.		
		The least global closure preserving all local consequences of the enlarged atlas is $C_*$.
		Thus adjoining the fillers does not change the least generated global closure.
		No filling order or choice of fillers is required.
		
		\item \emph{Localization in the chordal case.}
		If the original primal graph is chordal, then
		\[
		C_*|_M=\lambda_M
		\qquad\text{for every primal clique }M.
		\]		
	\end{enumerate}	
	The family $(\lambda_M)_M$ is uniquely determined by the original atlas and is independent of the	chosen admissible chart-local rule scheme.
\end{theorem}

\begin{proof}
	\leavevmode
	\begin{enumerate}[label=(\roman*),leftmargin=2.5em]		
		\item \emph{Automatic clique fillability.}
		Fix a primal clique $M$, possibly empty.
		Restrict the original atlas to $M$, obtaining the compatible atlas
		$
		\mathcal A|_M
		=
		(c_I|_{I\cap M})_{I\in\mathcal H}
		$
		on the induced cover
		$
		\mathcal H|_M
		=
		(I\cap M)_{I\in\mathcal H}.
		$	
		Restriction stability of $\mathscr C$ ensures that all induced closures remain in the local class.
		
		The induced primal graph is complete on $M$:
		every pair of vertices of $M$ occurs together in an induced original chart.
		Moreover, every nonempty clique in $M$ of cardinality at most $m$ is guarded by an induced original chart, by the guarding
		hypothesis on $\mathcal H$.		
		Thus the induced cover has a chordal primal graph and satisfies the required $m$-conformality condition, even when
		$G_{\mathcal H}$ is nonchordal.
		
		By trace locality of $\mathscr R$,
		\[
		\bigcup_{I\in\mathcal H}
		\mathscr R(c_I|_{I\cap M})
		=
		T[M].
		\]
		
		Applying
		Theorem~\ref{THM:PackagingControlledPropagationGlobalization}
		to the induced atlas shows that
		$
		\lambda_M=C_{T[M]}
		$
		is its least exact extension.
		Therefore
		\[
		\lambda_M|_{M\cap I}
		=
		c_I|_{M\cap I}
		\qquad(I\in\mathcal H),
		\]
		proving~\textup{(i)}.
		
		\item \emph{Restriction coherence.}
		Let $L\subseteq M$ be primal cliques.
		Since the induced primal graph on $M$ is complete, $L$ is a clique of that graph.
		
		Applying the clique-locality identity of
		Theorem~\ref{THM:PackagingControlledPropagationGlobalization}
		to the induced atlas gives
		\[
		\lambda_M|_L
		=
		C_{T[M][L]}
		=
		C_{T[L]}
		=
		\lambda_L.
		\]
		This includes $L=\varnothing$.
		
		For arbitrary primal cliques $M,N$,	their intersection is a primal clique.
		Hence
		\[
		\lambda_M|_{M\cap N}
		=
		\lambda_{M\cap N}
		=
		\lambda_N|_{M\cap N},
		\]
		establishing~\textup{(ii)}.
		
		\item \emph{Simultaneous conservative conformalization.}
		Since the primal graph is finite, it has finitely many maximal cliques, and every clique is contained in a maximal one.
		
		By~\textup{(i)}, each filler agrees with every original chart on their overlap.
		By~\textup{(ii)}, all fillers agree with one another on overlaps.
		The enlarged atlas is therefore compatible.
		
		Adjoining every maximal primal clique makes	the cover fully conformal.
		Since each adjoined set is already a clique	of the original primal graph, no new primal	edge is introduced.
		
		It remains to prove conservation.
		For every primal clique $M$ and $S\subseteq M$,	the inclusion $T[M]\subseteq T$ gives
		\[
		\lambda_M(S)
		=
		C_{T[M]}(S)
		\subseteq
		C_*(S)\cap M.
		\]
		Thus every consequence of every adjoined closure already holds in $C_*$.
		
		Let $C_*^+$ be the least global closure preserving all local consequences of the	enlarged atlas.
		Since the enlarged atlas contains all original charts, leastness gives
		$
		C_*\le C_*^+.
		$
		Conversely, the preceding inclusion shows that $C_*$ preserves every consequence of the enlarged atlas.
		Hence
		$
		C_*^+\le C_*.
		$
		Therefore
		\[
		C_*^+=C_*.
		\]		
		All fillers were defined simultaneously from the original atlas, so no filling order is involved.
		
		\item \emph{Localization in the chordal case.}
		If $G_{\mathcal H}$ is chordal, apply the clique-locality identity of
		Theorem~\ref{THM:PackagingControlledPropagationGlobalization}
		directly to the original atlas.
		For every primal clique $M$, this gives
		\[
		C_*|_M
		=
		C_{T[M]}
		=
		\lambda_M,
		\]
		proving~\textup{(iv)}.		
	\end{enumerate}	
	Finally, by~\textup{(i)}, each $\lambda_M$ is the least exact extension of the induced atlas on $M$.
	By Lemma~\ref{LEM:PackagingClosureRestrictionGeneration}, the least generated extension depends only on the induced local closures, not on their generating bases.
	Trace locality identifies the induced rule families with the corresponding restrictions of the original chart-local rule scheme.
	Consequently, each $\lambda_M$ depends only on the original atlas, and not on the choice of an admissible chart-local generating scheme.
\end{proof}

\begin{remark}[Hypotheses and scope of canonical clique completion]
	\label{REM:PackagingCliqueCompletionDependency}
	Theorem~\ref{THM:PackagingCanonicalCliqueCompletion} does not require weak fillability, an order of filling, or chordality.
	It does require local propagation recoverability and the stated clique-guarding hypothesis.
	For the bounded-Armstrong/native class, propagation recoverability follows from the later all-dimensional resolution theorem, not conversely.
	The added clique fillers need not belong to the original complexity-bounded class, and conformalization alone need not yield an exact global extension on a nonchordal graph.
\end{remark}

\begin{proposition}[Monotonicity of propagation-recovery rank]
	\label{PROP:PackagingPropagationRankMonotonicity}
	Suppose \(\mathscr C'(I)\subseteq\mathscr C(I)\) for every finite \(I\), and suppose a chart-local rule scheme \(\mathscr R'\) is the
	restriction of \(\mathscr R\) to \(\mathscr C'\).
	Then
	\[
	{
		\operatorname{prk}(\mathscr C',\mathscr R')
		\le
		\operatorname{prk}(\mathscr C,\mathscr R).
	}
	\]
\end{proposition}

\begin{proof}
	Every compatible $\mathscr C'$-atlas is a compatible $\mathscr C$-atlas, and by hypothesis the rule schemes agree on it.
	Hence every recovery valid for $\mathscr C$ at a given guarding level is valid for $\mathscr C'$ at the same level.
\end{proof}

\begin{theorem}[Exact propagation-recovery rank for bounded determination]
	\label{THM:PackagingBoundedDeterminationPropagationRank}
	Fix $k\ge1$.
	Let $\mathscr C_k(I)$ be the class of finite empty-preserving closures on $I$ of Carath\'eodory number at most $k$, and define
	\[
	\mathscr R_k(c)
	=
	\{A\to b:
	1\le |A|\le k,\ b\in c(A)\setminus A\}.
	\]
	Then the following hold:
	\begin{enumerate}[label=(\roman*),leftmargin=2.5em]		
		\item \emph{Chart-local rule scheme.}
		The class $\mathscr C_k$ is stable under restriction, and $\mathscr R_k$ is a	chart-local rule scheme for $\mathscr C_k$.
		
		\item \emph{Exact propagation-recovery rank.}
		The local propagation-recovery rank is
		\[
		\operatorname{prk}(\mathscr C_k,\mathscr R_k)
		=
		2k+1.
		\]		
	\end{enumerate}
\end{theorem}

\begin{proof}
	\leavevmode
	\begin{enumerate}[label=(\roman*),leftmargin=2.5em]		
		\item \emph{Chart-local rule scheme.}
		Restriction stability of $\mathscr C_k$ follows from Lemma~\ref{LEM:PackagingHereditaryPremiseBound}.
		The same lemma shows that $\mathscr R_k(c)$	generates $c$ for every $c\in\mathscr C_k(I)$.		
		For $J\subseteq I$ and $A\cup\{b\}\subseteq J$, the identity
		\[
		(c|_J)(A)=c(A)\cap J
		\]
		gives
		\[
		A\to b\in\mathscr R_k(c)
		\quad\Longleftrightarrow\quad
		A\to b\in\mathscr R_k(c|_J).
		\]
		Thus $\mathscr R_k$ satisfies trace locality and is a chart-local rule scheme.
		
		\item \emph{Exact propagation-recovery rank.}		
		\begin{itemize}
			\item \emph{Upper bound.}
			Suppose every nonempty primal clique of cardinality at most $2k+1$ is guarded.
			No chordality assumption is needed.
			
			Consider two original-chart rules
			\[
			A\to v,
			\qquad
			B\cup\{v\}\to b,
			\]
			whose combined support
			$
			Q=A\cup B\cup\{v,b\}
			$
			is a primal clique.
			
			Since both rules belong to $\mathscr R_k$,
			\[
			|A|\le k,
			\qquad
			|B|+1\le k.
			\]
			Consequently,
			\[
			|Q|
			\le k+(k-1)+2
			=2k+1.
			\]
			
			Hence $Q$ is contained in an original chart $J$.
			Compatibility transfers both parent implications
			to $c_J$, giving
			$
			b\in c_J(A\cup B).
			$
			
			By the Carath\'eodory bound and empty preservation, there exists a nonempty set
			\[
			D\subseteq A\cup B,
			\qquad
			|D|\le k,
			\qquad
			b\in c_J(D).
			\]
			Therefore
			$
			D\to b\in\mathscr R_k(c_J).
			$		
			This proves local $(2k+1)$-propagation recoverability, and hence
			\[
			\operatorname{prk}(\mathscr C_k,\mathscr R_k)
			\le 2k+1.
			\]
			
			\item \emph{Lower bound.}
			Let $Q$ contain $2k+1$ labels and consider its simplex-boundary cover
			\[
			\mathcal H
			=
			\{Q\setminus\{q\}:q\in Q\}.
			\]
			Its primal graph is complete, and every clique of cardinality at most $2k$ is guarded.			
			Choose a disjoint decomposition
			\[
			Q=Z\sqcup W\sqcup\{x,y\},
			\qquad
			|Z|=k,
			\qquad
			|W|=k-1.
			\]			
			On each facet, take the closure generated by those of the two rules
			\[
			Z\to x,
			\qquad
			W\cup\{x\}\to y
			\]
			whose complete support lies in that facet.
			A facet containing neither support carries the identity closure.			
			Since the union of the two rule supports is $Q$, each facet carries at most one rule.
			
			We first verify compatibility.
			Let $I,J$ be facets and $S\subseteq I\cap J$.
			If
			\[
			b\in c_I(S)\setminus S,
			\]
			then the unique generating rule	$P\to b$ on $I$ has $P\subseteq S$.
			Thus
			\[
			P\cup\{b\}\subseteq I\cap J,
			\]
			so the same rule is assigned to $J$.
			Hence $b\in c_J(S)$.
			
			By symmetry, the restrictions agree, proving compatibility.
			
			Each local closure is empty-preserving and has Carath\'eodory number at most $k$:
			every nontrivial new consequence follows directly from its unique generating rule, whose premise has size at most $k$.
			
			Both displayed rules belong to their respective chart-local rule families, and their combined support is the primal clique $Q$.
			
			Suppose local $2k$-propagation recovery held.
			It would supply an original-chart rule
			\[
			D\to y,
			\qquad
			\varnothing\ne D\subseteq Z\cup W.
			\]
			
			However, the facet $Q\setminus\{x\}$ contains $D\cup\{y\}$ and carries the identity closure.
			Compatibility would transfer $D\to y$ to that identity closure, a contradiction.
			
			Therefore local $2k$-propagation recoverability fails.
			By monotonicity in the guarding level,
			\[
			\operatorname{prk}(\mathscr C_k,\mathscr R_k)
			\ge 2k+1.
			\]
		\end{itemize}		
		Combining the upper and lower bounds yields
		\[
		\operatorname{prk}(\mathscr C_k,\mathscr R_k)
		=
		2k+1.
		\]		
	\end{enumerate}
\end{proof}

\begin{theorem}[Sharp chordal bounded-conformality classification for bounded-determination gluing]
	\label{THM:PackagingGeneralCaratheodoryGluing}
	For a finite cover \(\mathcal H\) and \(k\ge1\), the following are equivalent:
	\begin{enumerate}
		\item Every compatible atlas of empty-preserving closures of Carath\'eodory number at most \(k\) has an exact extension.
		
		\item The primal graph of $\mathcal H$ is chordal, and every nonempty primal clique of cardinality at most $2k+1$ is guarded.
		Equivalently, $\mathcal H$ is $(2k+1)$-conformal and its primal graph is chordal.
	\end{enumerate}
	When these conditions hold, the union of local implicational bases generates the least exact extension \(C_*\).
	It admits a basis with nonempty premises of size at most \(k\), and a finite descriptor vocabulary realizes its restrictions.
	When the geometric condition fails, there is a compatible nonextendible atlas in which each local closure has an implicational basis containing at most one nontrivial rule.
\end{theorem}

\begin{proof}
	The empty label set is immediate, so assume \(V\ne\varnothing\).
	\begin{itemize}
		\item \emph{Sufficiency.}
		Assume the primal graph is chordal and every nonempty clique of cardinality at most \(2k+1\) is guarded.
		By Theorem~\ref{THM:PackagingBoundedDeterminationPropagationRank}, \((\mathscr C_k,\mathscr R_k)\) is locally \((2k+1)\)-propagation-recoverable.
		The controlled-propagation principle, Theorem~\ref{THM:PackagingControlledPropagationGlobalization}, therefore shows that the closure \(C_*\) generated by
		\[
		T=\bigcup_{I\in\mathcal H}\mathscr R_k(c_I)
		\]
		is the least exact global extension.
		Every rule of \(T\) has nonempty premise of size at most \(k\), so this union itself supplies the asserted global bounded-premise basis.
		Lemma~\ref{LEM:PackagingArmstrongRealization} then supplies a finite descriptor vocabulary realizing \(C_*\), and hence realizing all prescribed local restrictions exactly.
		
		\emph{Compatibility of the obstruction atlases.}
		In each construction below, a fixed collection of implications is distributed to the charts according to support, and each chart receives at most one implication.
		
		Such an atlas is compatible.
		Indeed, for charts $I,J$ and $S\subseteq I\cap J$, suppose
		\[
		b\in c_I(S)\cap(I\cap J)\setminus S.
		\]
		The unique generating rule on $I$ must have the form $P\to b$ with $P\subseteq S$.
		Its full support lies in $I\cap J$,	so it is also assigned to chart $J$.
		Thus $b\in c_J(S)$, and symmetry gives
		\[
		c_I|_{I\cap J}=c_J|_{I\cap J}.
		\]

		\item \emph{Necessity.}
		
		\emph{Failure of chordality.}
		Let $v_0,v_1,\ldots,v_{m-1},v_0$, $m\ge4$, be an induced cycle.
		Place $v_i\to v_{i+1}$ for $0\le i<m-1$ on every chart containing its support, and leave all other labels inert.
		A chart meets the cycle in at most one edge, so carries at most one rule.
		If a new conclusion is visible on an overlap, its entire premise and conclusion lie there, so the other chart contains the same rule.
		This proves compatibility.
		Each closure is empty-preserving with determination bound one.
		Globally $v_0$ determines $v_{m-1}$, while a chart containing the closing edge has identity closure.
		No exact extension exists.
		
		\emph{An unguarded small clique.}
		Choose an inclusion-minimal unguarded clique $M$ with $|M|\le2k+1$.
		Singletons and edges are guarded, so $|M|\ge3$ and every proper subset of $M$ is guarded.
		Choose distinct $x,y \in M$ and take
		\[
		|Z|=\min(k,|M|-2),
		\qquad
		W=M\setminus(Z\cup\{x,y\}),
		\]
		with $Z\subseteq M\setminus\{x,y\}$, then $1\le|Z|\le k$.
		If	$|M|-2>k$, then the bound $|M|\le2k+1$ gives $|W|=|M|-2-k\le k-1$, while otherwise $W=\varnothing$.
		Thus
		\[
		M=Z\sqcup W\sqcup\{x,y\},\qquad
		1\le |Z|\le k,\quad 0\le |W|\le k-1.
		\]
		Put $Z\to x$ and $W\cup\{x\}\to y$ on every chart containing the respective support, with other labels inert.
		No chart contains both rules, since the union of their supports is $M$.
		Each rule occurs because its support is a proper subset of $M$.
		The preceding overlap argument proves compatibility.
		A closure with at most one such rule is empty-preserving and has determination bound at most $k$: an added conclusion uses that rule's premise, and a label already present uses its singleton.
		A chart containing $M\setminus\{x\}$ cannot	contain $x$ and hence carries neither rule.
		Globally the two rules derive $y$ from $Z\cup W$, contradicting exactness on that chart.
		This proves	necessity and the asserted one-rule obstruction bound.
	\end{itemize}
\end{proof}

\begin{remark}[Bounded conformality versus classical hypergraph acyclicity]
	\label{REM:PackagingBoundedConformalityAcyclicity}
	The geometric terminology in Theorem~\ref{THM:PackagingGeneralCaratheodoryGluing} is classical.
	For a cover hypergraph, full conformality means that every clique of the primal graph is contained in a chart, while \(N\)-conformality requires this only for cliques of cardinality at most \(N\)~\cite{Otto2012HighlyAcyclic}.
	Full conformality together with chordality is the standard acyclicity condition for database hypergraphs	\cite{BeeriFaginMaierYannakakis1983Acyclic,Fagin1983Degrees}.
	
	Theorem~\ref{THM:PackagingGeneralCaratheodoryGluing} does not redefine hypergraph acyclicity.
	Rather, for the	semantic class of empty-preserving closure atlases with Carath\'eodory number at most \(k\), it identifies the exact weakening of conformality that remains sufficient and necessary under chordality:
	\[
	(2k+1)\text{-conformality}.
	\]
	The propagation-rank theorem and the explicit simplex-boundary counterexamples show that the threshold cannot be lowered.
	When \(2k+1\) is at least the clique number of the primal graph, bounded conformality coincides with full conformality, and the geometric condition reduces to the classical chordal--conformal one.
\end{remark}

\section{Armstrong complexity and native commutant realization}
\label{SEC:ArmstrongNativeA}

The preceding section bounds how many labels are needed to determine another label.
The access--commutant model gives a different measure: the dimension of a commuting self-adjoint matrix realization.
The main theorem here, Theorem~\ref{THM:PackagingCommutantDimensionEqualsArmstrongSize}, proves that the smallest such dimension equals the minimum number of rows in an Armstrong relation, for nonempty empty-preserving closures.
We also prove spectral bounds and dimension-drop results.
The spectral bound provides the generic baseline against which the sharper native globalization threshold in Section~\ref{SEC:NativeRecoveryA} is compared.

\subsection{Native commutant closures and realization dimension}
\label{SUBSEC:PackagingNativeCommutantRealization}

A commuting self-adjoint tuple determines a closure on its label set: a label is determined by others when its operator belongs to the algebra they generate.
Proposition~\ref{PROP:PackagingNativeCommutantDependence} describes this closure through joint spectral coordinates.
The minimum realization dimension and its behavior under restriction are defined below.

Fix a finite label set $V$ and an integer $d\ge1$, and put
\[
\mathcal H_d=\mathbb C^d.
\]
Let
\[
\mathsf{CSA}_{d,V}
=
\left\{
\mathcal D=(D_a)_{a\in V}
\in M_d(\mathbb C)^V:
D_a^*=D_a,\quad
[D_a,D_b]=0
\quad(a,b\in V)
\right\}.
\]

Let $\mathscr C_d$ denote the structural domain of concrete unital $*$-subalgebras of $B(\mathcal H_d)$.
Define the packaging context $\Xi_{d,V}$ by
\[
\mathscr C_{\Xi_{d,V}}
=
\operatorname{Car}_{\Xi_{d,V}}
=
\mathscr C_d,
\qquad
\operatorname{Form}_{\Xi_{d,V}}
=
\mathsf{CSA}_{d,V}\times\mathcal P(V).
\]

For a formation datum
$(\mathcal D,S)\in\operatorname{Form}_{\Xi_{d,V}}$,
define
\[
M_S
=
\operatorname{alg}^{*}
(I_d,D_s:s\in S),
\qquad
N_S=M_S',
\]
and
\[
c_{\mathcal D}(S)
=
\{a\in V:N_S\subseteq\{D_a\}'\}.
\]
All commutants are taken in $B(\mathcal H_d)$, and $\operatorname{alg}^{*}$ denotes the generated unital complex $*$-algebra.

The realization relation is
\[
\mathsf J_{\Xi_{d,V}}
=
\left\{
\bigl(
(\mathcal D,S),
\mathfrak P_{N_S}
\bigr):
\mathcal D\in\mathsf{CSA}_{d,V},
\ S\subseteq V
\right\},
\]
where $\mathfrak P_{N_S}$ is the one-carrier refinement-packaged structure of Example~\ref{EX:PackagingAccessCommutant}.
Thus $\Xi_{d,V}$ is total and deterministic.
Its native formation is
\[
M_S\rightsquigarrow N_S=M_S',
\]
while $c_{\mathcal D}$ records the dependency closure induced by the entire labelled tuple $\mathcal D$.

We first identify this closure through simultaneous spectral decomposition and establish finite-dimensional realization.
We then define the minimum commutant realization dimension and compare it with the minimum Armstrong size.

These constructions are specified independently of any dependency atlas.
They provide a native class for the later globalization problem without changing the primitive packaging axioms.

\begin{proposition}[Dependency closure represented by native commutants]
	\label{PROP:PackagingNativeCommutantDependence}
	For a finite commuting self-adjoint tuple $\mathcal D$,	$c_{\mathcal D}$ is a closure on $V$.
	In a common orthonormal eigenbasis write $D_a e_u=v_a(u)e_u$, and let $\Pi_S$ be the partition of the basis indices by equality of all $v_s$, $s\in S$.
	Then
	\[
	\begin{split}
		M_S&=\bigoplus_{B\in\Pi_S}\mathbb C I_{\mathbb C^B},\qquad
		N_S=\bigoplus_{B\in\Pi_S}B(\mathbb C^B),\\
		a\in c_{\mathcal D}(S)
		&\quad\Longleftrightarrow\quad D_a\in M_S
		\quad\Longleftrightarrow\quad
		v_a\text{ is determined by }(v_s)_{s\in S}.
	\end{split}
	\]
	If $a\notin c_{\mathcal D}(S)$, some matrix unit $E_{uv}$ belongs to $N_S$ and does not commute with $D_a$.
	Conversely, every closure $c$ on finite $V$ admits such a realization	with $d=|\operatorname{Fix}(c)|$.
	In this realization
	\[
	N_S\text{ is noncommutative}\quad\Longleftrightarrow\quad c(S)\ne V.
	\]
	Consequently, for every such $S$ there is no surjective complex algebra	homomorphism $M_S\to N_S$, although the represented endpoint determines	the source by $N_S'=M_S$.
\end{proposition}

\begin{proof}
	Simultaneous diagonalization reduces the tuple to coordinate functions on the joint spectral sectors.
	The algebra $M_S$ is the algebra of diagonal functions constant on the partition $\Pi_S$, and its commutant is the direct sum of full matrix algebras on the blocks.
	Hence $N_S\subseteq\{D_a\}'$ iff $D_a\in M_S$, equivalently iff the $a$-coordinate is determined by the coordinates in $S$.
	If this fails, a matrix unit inside	one $\Pi_S$-block witnesses noncommutation.
	Moreover,
	\[
	c_{\mathcal D}(S)
	=
	\{a\in V:D_a\in M_S\},
	\qquad
	M_{c_{\mathcal D}(S)}=M_S.
	\]
	The first identity gives extensivity and monotonicity, while the second gives idempotence.
	Thus $c_{\mathcal D}$ is a closure operator.
	Conversely, Lemma~\ref{LEM:PackagingArmstrongRealization}
	provides a descriptor realization of \(c\) with exactly \(|\operatorname{Fix}(c)|\) distinct rows.
	Encoding the finitely many descriptor values as	real numbers and taking diagonal matrices gives a commuting self-adjoint realization on
	\[
	\mathbb C^{|\operatorname{Fix}(c)|}.
	\]
	
	In this realization, the full descriptor tuple distinguishes all rows.
	Hence
	\[
	c(S)=V
	\quad\Longleftrightarrow\quad
	\Pi_S\text{ is discrete}.
	\]
	Since
	\[
	N_S=\bigoplus_{B\in\Pi_S}B(\mathbb C^B),
	\]
	this is equivalent to \(N_S\) being commutative.
	If \(c(S)\ne V\), then \(N_S\) is noncommutative while \(M_S\) is commutative, so no surjective complex-algebra homomorphism \(M_S\to N_S\) exists.
\end{proof}

In general, $c_{\mathcal D}$ need not preserve the empty set.
Indeed,
\[
c_{\mathcal D}(\varnothing)
=
\{a\in V:D_a\in\mathbb C I_d\}.
\]
Empty preservation is imposed separately when defining the native dimension classes.

\begin{definition}[Native commutant realization dimension]
	\label{DEF:PackagingCommutantRealizationDimension}	
	Let $c:2^V\to2^V$ be a closure operator	on a finite label set $V$.	
	The \textbf{native commutant realization dimension} of $c$ is
	\[
	d_{\mathrm{comm}}(c)
	=
	\min
	\left\{
	d\ge1:
	\exists\,\mathcal D\in\mathsf{CSA}_{d,V}
	\text{ with }c_{\mathcal D}=c
	\right\}.
	\]	
	The minimum exists by
	Proposition~\ref{PROP:PackagingNativeCommutantDependence}.
	For $V=\varnothing$, it equals $1$.
\end{definition}

\noindent\textbf{Notation.}
For nonempty $V$ and empty-preserving $c$, write $\operatorname{car}(c)$ for the Carath\'eodory number of~\eqref{EQ:PackagingCaratheodoryBound}.

\begin{definition}[Native dimension classes]
	\label{DEF:PackagingNativeDimensionClasses}
	Let \(I\) be a finite label set and let \(d\ge1\).
	Define
	\[
	\mathbf{CommCl}_d(I)
	=
	\left\{
	c:2^I\to2^I:
	\begin{array}{l}
		c\text{ is an empty-preserving closure, and}\\
		d_{\mathrm{comm}}(c)\le d
	\end{array}
	\right\}.
	\]
	For \(I=\varnothing\), this contains the unique closure in every
	\(\mathbf{CommCl}_d(\varnothing)\).
	For nonempty \(I\),	\(\mathbf{CommCl}_1(I)=\varnothing\), since a one-row realization makes every	attribute determined by the empty set.
\end{definition}

\begin{lemma}[Restriction and spectral-sector compression]
	\label{LEM:PackagingNativeRestrictionCompression}
	Let
	$
	\mathcal D=(D_a)_{a\in I}
	$
	be a finite commuting self-adjoint tuple, and let \(c=c_{\mathcal D}\).
	\begin{enumerate}[label=(\roman*),leftmargin=2.5em]
		\item For every \(S\subseteq I\), the restricted closure \(c|_S\) is the native closure of the subtuple \((D_a)_{a\in S}\).
		In particular, restriction never increases native realization dimension.
		
		\item If the realization has exactly \(m\) joint spectral sectors, then \(c\) admits a realization on \(\mathbb C^m\).
		Therefore every realization with the smallest possible matrix dimension
		$
		r=d_{\mathrm{comm}}(c)
		$
		has exactly \(r\) joint spectral sectors.
	\end{enumerate}
\end{lemma}

\begin{proof}
	For \textup{(i)}, simply forget the coordinates indexed by \(I\setminus S\).
	By Proposition~\ref{PROP:PackagingNativeCommutantDependence}, functional
	dependence among the remaining coordinate functions is exactly \(c|_S\).
	
	For \textup{(ii)}, choose one common joint eigenvector from each joint spectral sector.
	Restrict every diagonal coordinate function to those \(m\) representatives and use the resulting diagonal matrices on \(\mathbb C^m\).
	Two selected rows agree on a set of labels exactly when the corresponding
	joint spectral sectors agree there.
	Hence all functional dependencies, and therefore the native closure, are unchanged.
	
	If a minimal-dimensional realization in \(M_r(\mathbb C)\) had fewer than
	\(r\) joint spectral sectors, the preceding compression would realize the
	same closure in smaller matrix dimension, contradicting minimality.
\end{proof}

\subsection{Armstrong complexity and exact native identification}
\label{SUBSEC:PackagingArmstrongComplexity}

The same dependency closure can be represented either by commuting self-adjoint matrices or by a finite database relation.
Theorem~\ref{THM:PackagingCommutantDimensionEqualsArmstrongSize} proves that the least matrix dimension equals the smallest number of rows in an Armstrong relation, under the stated hypotheses.
This makes the dimension bound depend on the closure itself, not on a particular matrix presentation.

\begin{definition}[Minimum Armstrong size]
	\label{DEF:PackagingMinimumArmstrongSize}
	Let $c:2^I\to2^I$ be a closure operator on a finite attribute set $I$.	
	A \textbf{finite relation} on $I$ is a finite set
	\[
	R\subseteq\prod_{a\in I}Y_a,
	\]
	for some nonempty value sets $Y_a$.
	It satisfies the dependency $X\to b$, written
	$R\models X\to b$, when
	\[
	u|_X=v|_X\quad\Longrightarrow\quad u_b=v_b
	\qquad(u,v\in R).
	\]
	The relation $R$ is an \textbf{Armstrong relation for $c$} if
	\[
	R\models X\to b
	\quad\Longleftrightarrow\quad
	b\in c(X)
	\qquad(X\subseteq I,\ b\in I).
	\]
	The \textbf{minimum Armstrong size} is
	\[
	\operatorname{arm}(c)
	=\min\{|R|:R\text{ is an Armstrong relation for }c\}.
	\]
	The minimum exists by finite descriptor realization
	(Lemma~\ref{LEM:PackagingArmstrongRealization}) \cite{BeeriDowdFaginStatman1984Armstrong}.	
	For a functional-dependency theory $\Sigma$ on $I$, put
	\[
	c_\Sigma(X)=\{b\in I:\Sigma\models X\to b\},
	\qquad
	\operatorname{arm}(\Sigma)=\operatorname{arm}(c_\Sigma).
	\]
\end{definition}

\begin{remark}[Zero-row convention for Armstrong size]
	The empty relation is permitted and satisfies every functional dependency.
	Consequently,
	\[
	\operatorname{arm}(c)=0
	\quad\Longleftrightarrow\quad c(\varnothing)=I.
	\]
	Indeed, $c(\varnothing)=I$ forces $c(X)=I$ for every
	$X\subseteq I$.
	If $I\ne\varnothing$ and $c$ is empty-preserving, neither a zero-row nor a one-row relation represents $c$, so
	$\operatorname{arm}(c)\ge2$.
	
	For $I=\varnothing$, the unique closure has
	$\operatorname{arm}(c)=0$ but $d_{\mathrm{comm}}(c)=1$.
	Thus the equality in
	Theorem~\ref{THM:PackagingCommutantDimensionEqualsArmstrongSize}
	is stated only for nonempty, empty-preserving closures.
\end{remark}

\begin{theorem}[Native commutant dimension equals minimum Armstrong size]
	\label{THM:PackagingCommutantDimensionEqualsArmstrongSize}
	Let $V$ be finite and nonempty, and let $c$ be an
	empty-preserving closure on $V$.
	Then
	\[
	d_{\mathrm{comm}}(c)=\operatorname{arm}(c).
	\]
	A minimum-dimensional native realization can be chosen diagonal, with one basis vector per row of a minimum Armstrong relation.
	Every native realization in dimension $d_{\mathrm{comm}}(c)$
	has exactly $d_{\mathrm{comm}}(c)$ joint spectral sectors.
\end{theorem}

\begin{proof}
	Put $r=d_{\mathrm{comm}}(c)$ and choose an $r$-dimensional
	commuting self-adjoint realization.
	In a common eigenbasis,	retain one row of joint eigenvalues from each distinct joint spectral sector.
	These $m\le r$ rows form a relation $R$.
	By Proposition~\ref{PROP:PackagingNativeCommutantDependence},
	row agreement represents exactly $c$.
	Hence $R$ is an	Armstrong relation and
	\[
	\operatorname{arm}(c)\le m\le r.
	\]
	
	Conversely, choose a minimum Armstrong relation $R$ with $m=\operatorname{arm}(c)$ rows.
	Encode the finitely many values used in each attribute column
	injectively by real numbers, and place the resulting column
	entries on the diagonal of a matrix $D_a\in M_m(\mathbb C)$.
	The matrices commute and are self-adjoint.
	Their dependency closure is $c$, again by
	Proposition~\ref{PROP:PackagingNativeCommutantDependence}.
	Therefore $d_{\mathrm{comm}}(c)\le m$, proving equality.
	
	The diagonal construction is minimum-dimensional.
	In any $r$-dimensional realization the first inequality gives
	$r=\operatorname{arm}(c)\le m\le r$, so the number of joint spectral sectors is exactly $r$.
\end{proof}

\begin{remark}[Calibration scope of the Armstrong--commutant identity]
	\label{REM:PackagingArmstrongCommutantCalibration}
	Theorem~\ref{THM:PackagingCommutantDimensionEqualsArmstrongSize} is an exact translation between two finite representations of the same dependency closure.
	Partition semantics for functional dependencies and the size theory of Armstrong relations are classical
	\cite{Lee1983AlgebraicRelational,CosmadakisKanellakisSpyratos1986PartitionSemantics,BeeriDowdFaginStatman1984Armstrong}.
	The identity calibrates native matrix dimension by a classical representation-minimal invariant; it does not alter the functional-dependency implication calculus.
	Section~\ref{SEC:NativeRecoveryA} uses this native parameter to obtain a sharper recovery bound than generic closure complexity alone provides.
\end{remark}

\subsection{Spectral determination and sharp dimension bounds}
\label{SUBSEC:PackagingNativeSpectralBounds}

We next ask how matrix dimension limits the size of a set needed to determine another label.
The spectral determination theorem gives a sharp bound in terms of the number of joint spectral sectors.
The dimension-drop lemmas then show that restricting to a set that does not determine every label reduces the minimum realization dimension.
Both results are used in the sharper native recovery proof.

\begin{theorem}[Native spectral bound on determination complexity]
	\label{THM:PackagingNativeCommutantDimensionBound}
	Let $\mathcal D=(D_a)_{a\in V}$ be a finite commuting
	self-adjoint tuple on $\mathbb C^d$, and let $\Pi_V$
	be its joint spectral partition.
	Put
	\[
	m=m(\mathcal D):=|\Pi_V|
	=\dim_{\mathbb C}M_V\le d,
	\qquad
	M_V=\operatorname{alg}^{*}(I_d,D_a:a\in V).
	\]
	Then:
	\begin{enumerate}[label=(\roman*),leftmargin=2.5em]
		\item For every $S\subseteq V$ and $b\in c_{\mathcal D}(S)$,	there is $A\subseteq S$ such that
		\[
		|A|\le m-1,
		\qquad b\in c_{\mathcal D}(A).
		\]
		In particular, if $m=1$, then $c_{\mathcal D}(\varnothing)=V$.
		
		\item If $V\ne\varnothing$ and $c_{\mathcal D}$ is
		empty-preserving, then $m\ge2$ and
		\[
		\operatorname{car}(c_{\mathcal D})
		\le m-1
		=\dim_{\mathbb C}M_V-1
		\le d-1.
		\]
	\end{enumerate}
\end{theorem}

\begin{proof}
	Let $\Omega$ be the set of $m$ joint spectral sectors.
	Each $D_a$ induces a real-valued function $\bar v_a$ on $\Omega$.
	For $A\subseteq V$, let	$\Pi_A$ be the partition of $\Omega$ by joint equality of the coordinates indexed by $A$.
	By Proposition~\ref{PROP:PackagingNativeCommutantDependence},
	\[
	b\in c_{\mathcal D}(A)
	\quad\Longleftrightarrow\quad
	\Pi_A\text{ refines }\Pi_{\{b\}}.
	\]
	
	Fix $S\subseteq V$ and $b\in c_{\mathcal D}(S)$.
	Starting with $A=\varnothing$, repeat the following step
	whenever $\Pi_A$ does not refine $\Pi_{\{b\}}$.
	Choose two sectors that agree on $A$ but differ on $b$.
	Since $S$ determines $b$, some $s\in S$ separates these sectors; adjoin $s$ to $A$.
	Each step strictly refines $\Pi_A$, and no partition
	of an $m$-element set admits more than $m-1$ successive strict refinements starting from one block.
	Hence the procedure terminates with
	\[
	A\subseteq S,\qquad |A|\le m-1,
	\qquad b\in c_{\mathcal D}(A).
	\]
	This also covers $m=1$, when the empty set already determines every coordinate.
	
	The spectral decomposition of $M_V$ has one scalar summand per joint spectral sector, so $m=\dim_{\mathbb C}M_V\le d$.
	If $V\ne\varnothing$ and $c_{\mathcal D}$ is
	empty-preserving, then $m\ne1$, and the bound in
	\textup{(i)} gives the asserted Carath\'eodory estimate.
\end{proof}

\begin{lemma}[Strict native dimension drop on a non-key]
	\label{LEM:PackagingNativeDimensionDrop}
	Let \(c\) be an empty-preserving native commutant closure on a finite ground set \(I\), and put
	\[
	r=d_{\mathrm{comm}}(c).
	\]
	Let
	\[
	S\subseteq I.
	\]
	If
	\[
	c(S)\ne I,
	\]
	then
	\[
	{
		d_{\mathrm{comm}}(c|_S)\le r-1.
	}
	\]	
	In particular, if \(d\ge2\), \(c\in\mathbf{CommCl}_d(I)\), and \(S\) is not a key, then
	\[
	c|_S\in\mathbf{CommCl}_{d-1}(S).
	\]
\end{lemma}

\begin{proof}
	Choose a realization of \(c\) in the minimal matrix dimension
	\[
	r=d_{\mathrm{comm}}(c).
	\]
	By Lemma~\ref{LEM:PackagingNativeRestrictionCompression}, this realization has exactly \(r\) joint spectral sectors.
	Let
	\[
	\Pi_S
	\]
	be the joint partition of those sectors induced by the labels in \(S\).
	
	If \(\Pi_S\) were discrete on all \(r\) sectors, it would refine the partition of every label, and hence
	\[
	c(S)=I,
	\]
	contrary to hypothesis.
	Thus
	\[
	|\Pi_S|\le r-1.
	\]
	
	Collapse each block of \(\Pi_S\) to one point.
	All labels of \(S\)	descend to this quotient, and functional dependence among them is unchanged.
	Hence \(c|_S\) is realized on at most \(r-1\) joint spectral sectors, so
	\[
	d_{\mathrm{comm}}(c|_S)\le r-1.
	\]
\end{proof}

\begin{theorem}[Native dimension of the identity closure]
	\label{THM:PackagingNativeIdentityDimension}
	For a finite label set $X$, let	$\operatorname{id}_X(S)=S$ for every $S\subseteq X$.
	Then
	\[
	d_{\mathrm{comm}}(\operatorname{id}_X)=|X|+1.
	\]
\end{theorem}

\begin{proof}
	If $X=\varnothing$, the assertion is the defining convention for $d_{\mathrm{comm}}$.
	Suppose	$r=|X|\ge1$ and identify $X$ with $\{1,\ldots,r\}$.
	
	For the upper bound, take the $r+1$ rows $0,e_1,\ldots,e_r$ in $\{0,1\}^r$.
	For every $i$ and every $S\subseteq X\setminus\{i\}$, the rows $0,e_i$ agree on $S$ and differ at $i$.
	Thus the induced dependency closure is the identity, and diagonal realization gives
	\[
	d_{\mathrm{comm}}(\operatorname{id}_X)\le r+1.
	\]
	
	Conversely, let a native realization act on $\mathbb C^d$.
	Keep its distinct joint spectral rows as a finite relation $R$.
	Then $|R|\le d$, and $R$ represents $\operatorname{id}_X$.
	For each $i\in X$, the failure of $X\setminus\{i\}\to i$ supplies distinct rows	$u_i,v_i\in R$ agreeing in every coordinate except $i$.
	Join them by an edge labelled $i$.
	
	The chosen edges are distinct, since the endpoints of an edge labelled $i$ differ in coordinate $i$ and agree in every other coordinate. They form a forest:
	otherwise, on any cycle, the complementary path to an edge labelled $i$ consists of edges with different labels.
	Along that path coordinate $i$ remains unchanged, contradicting the change of coordinate $i$ across the omitted edge.
	A forest with $r$ edges has at least $r+1$ vertices, so
	\[
	d\ge|R|\ge r+1.
	\]
	This proves the lower bound and hence equality.
\end{proof}

The matching formula for minimum Armstrong size is classical \cite{BeeriDowdFaginStatman1984Armstrong}.

\begin{lemma}[Native realization of the uniform circuit closure]
	\label{LEM:PackagingNativeUniformCircuit}
	Let \(X\) be a finite set with
	\[
	|X|=r\ge2.
	\]
	Define the closure
	\[
	u_X(S)
	=
	\begin{cases}
		S,
		&
		|S|\le r-2,
		\\[1ex]
		X,
		&
		|S|\ge r-1.
	\end{cases}
	\]
	Then
	\[
	{
		d_{\mathrm{comm}}(u_X)=r,
	}
	\]
	and hence \(u_X\in\mathbf{CommCl}_r(X)\).
	Moreover, every proper restriction of \(u_X\) is the identity closure.
\end{lemma}

\begin{proof}
	Identify
	\[
	X=\{1,\ldots,r\}
	\]
	and work on \(\mathbb C^r\) with basis
	\[
	e_1,\ldots,e_r.
	\]
	For \(a\in X\), define
	\[
	D_a e_j
	=
	\delta_{aj}e_j.
	\]
	Thus the \(r\) joint spectral rows are the standard basis vectors
	\[
	e_1,\ldots,e_r
	\]
	viewed as rows of a \(0\)-\(1\) table.
	
	Let \(S\subseteq X\).	
	If
	\[
	|S|\le r-2
	\]
	and \(b\notin S\), choose
	\[
	c\in X\setminus(S\cup\{b\}).
	\]
	The rows \(e_b\) and \(e_c\) agree on every label in \(S\), but differ
	at \(b\).
	Hence
	\[
	b\notin c_{\mathcal D}(S).
	\]
	
	If
	\[
	|S|\ge r-1,
	\]
	the coordinates in \(S\) distinguish all \(r\) rows:
	the row corresponding to a possibly omitted coordinate is the unique all-zero row on \(S\), while every other row has a unique \(1\).
	Hence \(S\) determines every label.	
	Therefore
	\[
	c_{\mathcal D}=u_X.
	\]
	
	Finally, if \(J\subsetneq X\), then every \(S\subseteq J\) has
	\[
	|S|\le r-1.
	\]
	The only case in which \(u_X(S)=X\) is \(S=J\) with \(|J|=r-1\), whose restriction back to \(J\) is again \(J\).
	Thus
	\[
	u_X|_J=\operatorname{id}_J.
	\]
	
	The construction gives \(d_{\mathrm{comm}}(u_X)\le r\).
	Conversely,	choose \(J\subsetneq X\) with \(|J|=r-1\).
	Restriction monotonicity and Theorem~\ref{THM:PackagingNativeIdentityDimension} give
	\[
	r=d_{\mathrm{comm}}(\operatorname{id}_J)
	=d_{\mathrm{comm}}(u_X|_J)
	\le d_{\mathrm{comm}}(u_X).
	\]
	Hence \(d_{\mathrm{comm}}(u_X)=r\).
\end{proof}

\begin{corollary}[Native saturation criterion]
	\label{COR:PackagingNativeSaturation}
	Let \(d\ge2\),
	\(c\in\mathbf{CommCl}_d(I)\),
	and \(S\subseteq I\).
	
	If
	\[
	c|_S\notin\mathbf{CommCl}_{d-1}(S),
	\]
	then
	\[
	{
		c(S)=I.
	}
	\]
	
	In particular, if
	\[
	|S|=d-1
	\]
	and
	\[
	c|_S=\operatorname{id}_S,
	\]
	then
	\[
	{
		c(S)=I.
	}
	\]
\end{corollary}

\begin{proof}
	The first assertion is the contrapositive of Lemma~\ref{LEM:PackagingNativeDimensionDrop}.
	
	For the second, Theorem~\ref{THM:PackagingNativeIdentityDimension} gives
	\[
	d_{\mathrm{comm}}(\operatorname{id}_S)
	=
	|S|+1
	=
	d.
	\]
	Thus
	\[
	\operatorname{id}_S\notin\mathbf{CommCl}_{d-1}(S),
	\]
	and the first assertion applies.
\end{proof}

\begin{corollary}[Carath\'eodory lower bound for native realization dimension]
	\label{COR:PackagingCommutantDimensionLowerBound}
	Let \(c\) be a closure on a nonempty finite label set \(V\), and assume
	\[
	c(\varnothing)=\varnothing.
	\]
	Then
	\[
	{
		d_{\mathrm{comm}}(c)
		\ge
		\operatorname{car}(c)+1.
	}
	\]
	More strongly, every native commutant realization
	\[
	c=c_{\mathcal D}
	\]
	satisfies
	\[
	{
		\dim_{\mathbb C}M_V
		\ge
		\operatorname{car}(c)+1.
	}
	\]
\end{corollary}

\begin{proof}
	Apply Theorem~\ref{THM:PackagingNativeCommutantDimensionBound} to any native realization:
	empty preservation gives at least two spectral sectors and
	\[
	\operatorname{car}(c)\le \dim_{\mathbb C}M_V-1\le d-1.
	\]
	Rearranging gives both displayed lower bounds, and minimizing over $d$ gives the assertion for $d_{\mathrm{comm}}(c)$.
\end{proof}

\begin{proposition}[Sharpness of the native dimension bound]
	\label{PROP:PackagingCommutantDimensionSharpness}
	For every integer \(k\ge1\), there exists a nonempty finite	empty-preserving closure \(c_k\) such that
	\[
	{
		\operatorname{car}(c_k)=k,
		\qquad
		d_{\mathrm{comm}}(c_k)=k+1.
	}
	\]
	Hence both inequalities
	\[
	\operatorname{car}(c_{\mathcal D})
	\le
	\dim_{\mathbb C}M_V-1
	\le
	d-1
	\]
	are sharp.
\end{proposition}

\begin{proof}
	Fix \(k\ge1\) and let
	\[
	\mathcal H=\mathbb C^{k+1}
	\]
	with orthonormal basis
	\[
	e_0,e_1,\ldots,e_k.
	\]
	Take the label set
	\[
	V_k=\{a_1,\ldots,a_k,b\}.
	\]
	
	Define commuting self-adjoint diagonal operators by
	\[
	D_{a_i}e_j
	=
	\begin{cases}
		e_j,&j=i,\\
		0,&j\ne i,
	\end{cases}
	\qquad
	1\le i\le k,
	\]
	and
	\[
	D_b e_j
	=
	j\,e_j,
	\qquad
	0\le j\le k.
	\]
	
	Let
	\[
	c_k=c_{\mathcal D}.
	\]
	None of the displayed operators is scalar, so
	\[
	c_k(\varnothing)=\varnothing.
	\]
	
	The joint values of
	\[
	a_1,\ldots,a_k
	\]
	separate all \(k+1\) basis points:
	the point \(0\) has the all-zero vector, while point \(i\) has the \(i\)-th standard basis vector.
	Hence
	\[
	b\in
	c_k(\{a_1,\ldots,a_k\}).
	\]
	
	Now fix \(j\in\{1,\ldots,k\}\).
	If \(a_j\) is omitted, the basis points \(0\) and \(j\) agree on every remaining \(a_i\), whereas
	\[
	v_b(0)=0,
	\qquad
	v_b(j)=j.
	\]
	Therefore
	\[
	b\notin
	c_k(\{a_1,\ldots,a_k\}\setminus\{a_j\}).
	\]
	Thus
	\[
	\{a_1,\ldots,a_k\}
	\]
	is an inclusion-minimal determining set for \(b\), and so
	\[
	\operatorname{car}(c_k)\ge k.
	\]
	
	The full joint spectral partition is discrete on
	\[
	\{0,1,\ldots,k\},
	\]
	so
	\[
	\dim_{\mathbb C}M_{V_k}=k+1.
	\]
	Theorem~\ref{THM:PackagingNativeCommutantDimensionBound} therefore gives
	\[
	\operatorname{car}(c_k)\le k.
	\]
	Hence
	\[
	\operatorname{car}(c_k)=k.
	\]
	
	The displayed tuple realizes \(c_k\) in dimension \(k+1\), so
	\[
	d_{\mathrm{comm}}(c_k)\le k+1.
	\]
	On the other hand, Corollary~\ref{COR:PackagingCommutantDimensionLowerBound} gives
	\[
	d_{\mathrm{comm}}(c_k)
	\ge
	\operatorname{car}(c_k)+1
	=
	k+1.
	\]
	Therefore
	\[
	d_{\mathrm{comm}}(c_k)=k+1.
	\]
\end{proof}

\subsection{Generic globalization baseline}
\label{SUBSEC:PackagingNativeGenericBaseline}

The spectral bound implies a Carath\'eodory bound for the local native closures.
Applying the generic closure theorem therefore gives an exact extension under chordality and \((2d-1)\)-conformality.
This is only a baseline: it does not use all the native spectral information.
The next section improves the threshold to \(d\) for \(d\ge3\), without requiring the global realization to have dimension at most \(d\).

\noindent\textbf{Native carrier identities.}
For a fixed commuting tuple $\mathcal D$ in a common ambient matrix algebra, the represented source and endpoint algebras satisfy
\[
N_{S\cup T}=N_S\cap N_T,
\qquad
c_{\mathcal D}(S)=c_{\mathcal D}(T)
\ \Longleftrightarrow\ M_S=M_T
\ \Longleftrightarrow\ N_S=N_T.
\]
Indeed, the commutant of the algebra generated by $M_S$ and $M_T$ is their commutant intersection.
Moreover,
\[
M_{c_{\mathcal D}(S)}=M_S,
\]
because $a\in c_{\mathcal D}(S)$ exactly when $D_a\in M_S$.
The bicommutant identity then gives the remaining equivalences.
Thus distinct closed label sets determine distinct native commutant carriers, and inclusion of closed label sets corresponds to reverse inclusion of their commutants.
These identities are preserved by simultaneous unitary conjugation.
They compare represented algebras within one ambient space, not unrelated abstract algebras.
They do not impose an information-loss ordering on packaging.

\begin{corollary}[Spectral-dimension-controlled native gluing]
	\label{COR:PackagingNativeDimensionControlledGluing}
	Let \(\mathcal H_0\) be a finite cover of a finite label set \(V\).	
	Suppose that for every nonempty chart \(I\in\mathcal H_0\), the local closure \(c_I\) is empty-preserving and is realized by a commuting self-adjoint tuple
	\[
	\mathcal D^{(I)}
	=
	(D_a^{(I)})_{a\in I}
	\]
	on a finite-dimensional Hilbert space \(\mathcal H_I\), such that
	\[
	\dim_{\mathbb C}
	\operatorname{alg}^{*}
	\bigl(
	I_{\mathcal H_I},
	D_a^{(I)}:a\in I
	\bigr)
	\le r
	\]
	for some fixed integer \(r\ge2\).
	
	If the atlas \((c_I)_{I\in\mathcal H_0}\) is compatible, the primal graph is chordal, and every nonempty clique of cardinality at most
	\[
	{
		2r-1
	}
	\]
	is guarded, then the atlas admits an exact global closure extension.
	
	In particular, suppose \(d\ge2\) and every nonempty local chart admits a native realization	on a Hilbert space of dimension at most \(d\).
	Then one may take \(r=d\).
	Consequently, chordality and guarding every	nonempty primal clique of cardinality at most
	\[
	2d-1
	\]
	are sufficient for exact extension.
	
	Moreover, the resulting exact global closure admits some finite native commutant realization.
	No assertion is made that this global realization has dimension at most \(d\), or that it preserves the originally chosen local matrices.
\end{corollary}

\begin{proof}
	By Theorem~\ref{THM:PackagingNativeCommutantDimensionBound}, every local closure has Carath\'eodory number at most \(r-1\).  Theorem~\ref{THM:PackagingGeneralCaratheodoryGluing} therefore gives exact extension under chordality and guarding through size \(2r-1\).
	If the local Hilbert-space dimension is at most \(d\), one may take \(r=d\).
	Proposition~\ref{PROP:PackagingNativeCommutantDependence} then provides a finite commuting self-adjoint realization of the global closure.
	No bound by the original local dimension is asserted.
\end{proof}

\section{Native recovery and fixed-dimensional globalization}
\label{SEC:NativeRecoveryA}

The generic closure argument does not use all the information available in a commuting self-adjoint realization.
We use that information to recover implications in the original charts.
The main result, Theorem~\ref{THM:PackagingNativeGeneralResolution}, proves this recovery in every dimension \(d\ge3\).
Theorem~\ref{THM:PackagingNativeGeneralClassification} then characterizes universal exact globalization by chordality and \(d\)-conformality; dimension two has a separate, unconditional classification.
We prove the result through spectral compression, conservative filling, and transport of finite obstructions.
A dimension-jump example shows that an exact global extension need not retain the local matrix dimension.

\subsection{Native commutant globalization under a local dimension bound}
\label{SUBSEC:PackagingFixedDimensionNativeGluing}

The dimension bound is imposed separately on each local commuting self-adjoint realization; the global closure may require a larger matrix dimension.
The obstruction rank measures the cover geometry.
Dimension two admits universal globalization on every finite cover, whereas the classification for \(d\ge3\) depends on chordality and bounded conformality.

\begin{definition}[Bounded conformality and globalization obstruction rank]
	\label{DEF:PackagingNativeObstructionRank}
	For an integer \(m\ge2\), call a finite cover \(\mathcal H\) \textbf{\(m\)-conformal} when every primal clique of cardinality at most \(m\) is guarded.
	This is the standard bounded-conformality terminology for hypergraphs.
	Full conformality means that every primal clique is guarded \cite{Otto2012HighlyAcyclic}.	
	Define the \textbf{globalization obstruction rank}
	\[
	\nu(\mathcal H)
	=
	\min\left(
	\{3:\ G_{\mathcal H}\text{ is nonchordal}\}
	\cup
	\{|K|:\ K\text{ is an unguarded primal clique}\}
	\right),
	\]
	with the convention that the minimum of the empty set is \(\infty\).
	Thus
	$
	\nu(\mathcal H)\ge3.
	$
	The value \(\nu(\mathcal H)\) is a geometric invariant of the cover.
	It does not depend on a closure atlas or on a native realization.
\end{definition}

\begin{lemma}[Obstruction-rank dictionary]
	\label{LEM:PackagingObstructionRankDictionary}
	Let \(\mathcal H\) be a finite cover and let \(m\ge3\).
	Then
	\[
	\nu(\mathcal H)>m
	\quad\Longleftrightarrow\quad
	\left\{
	\begin{array}{l}
		G_{\mathcal H}\text{ is chordal},\\
		\mathcal H\text{ is \(m\)-conformal}.
	\end{array}
	\right.
	\]
	Moreover,
	\[
	\nu(\mathcal H)=\infty
	\quad\Longleftrightarrow\quad
	\mathcal H\text{ is chordal and conformal}.
	\]
\end{lemma}

\begin{proof}
	By definition, \(\nu(\mathcal H)>m\) excludes the value \(3\) contributed by nonchordality and excludes every unguarded clique of cardinality at most \(m\).
	This is exactly chordality together with \(m\)-conformality.
	The second assertion follows by excluding unguarded cliques of every finite cardinality as well as nonchordality.
\end{proof}

To state the native globalization problem precisely, we first isolate the local class of closures admitting a commutant realization in dimension at most \(d\).
The bound is imposed chartwise: the desired global extension is exact as a closure but may use a larger finite realization dimension.
This separates local native complexity from the strictly stronger same-dimension realization problem.

Recall from Definition~\ref{DEF:PackagingNativeDimensionClasses} that, for every finite label set \(I\) and \(d\ge1\),
$
\mathbf{CommCl}_d(I)
$
denotes the class of empty-preserving closures whose native commutant realization dimension is at most \(d\).

Let \(\mathcal H\) be a finite cover of \(V\).
We write
$
\mathsf{UG}_{\mathrm{comm}}(d,\mathcal H)
$
when every compatible atlas
\[
(c_I)_{I\in\mathcal H},
\qquad
c_I\in\mathbf{CommCl}_d(I),
\]
admits an exact global closure
$
C:\mathcal P(V)\longrightarrow\mathcal P(V).
$

Because every finite closure admits a finite native commutant realization by Proposition~\ref{PROP:PackagingNativeCommutantDependence}, this is equivalently the assertion that every such fixed-\(d\) local atlas admits some finite-dimensional native global realization, with the global dimension allowed to increase.

We write
$
\mathsf{UG}^{=d}_{\mathrm{comm}}(d,\mathcal H)
$
when the stronger dimension-preserving statement holds:
for every such compatible atlas there exists one commuting self-adjoint tuple
\[
\mathcal D=(D_a)_{a\in V}
\subseteq M_d(\mathbb C)
\]
such that
\[
c_{\mathcal D}|_I=c_I
\qquad
(I\in\mathcal H).
\]
Thus
\[
{
	\mathsf{UG}^{=d}_{\mathrm{comm}}(d,\mathcal H)
	\Longrightarrow
	\mathsf{UG}_{\mathrm{comm}}(d,\mathcal H),
}
\]
but the converse is not assumed.

\begin{lemma}[Dimension padding preserves native dependency closure]
	\label{LEM:PackagingNativeDimensionPadding}
	Let $I$ be finite and let $c$ be a closure on $I$.
	Put	$r=d_{\mathrm{comm}}(c)$.
	For every integer $n\ge r$,	there exists a commuting self-adjoint tuple	$\widetilde{\mathcal D}\in\mathsf{CSA}_{n,I}$ with
	$c_{\widetilde{\mathcal D}}=c$.
	Consequently, for every $d\ge1$,
	\[
	\mathbf{CommCl}_d(I)\subseteq\mathbf{CommCl}_{d+1}(I).
	\]
\end{lemma}

\begin{proof}
	Choose a native realization $(D_a)_{a\in I}$ in dimension $r$ and a common unit eigenvector $e$, where $D_ae=\lambda_ae$. For $n=r$ retain this tuple.
	For $n>r$, put
	\[
	\widetilde D_a=D_a\oplus\lambda_a I_{n-r}
	\qquad(a\in I).
	\]
	Every added joint spectral row duplicates an existing row, so the induced equality partitions and hence all	dependency closures remain unchanged.
	The class inclusion	follows from the definition of $d_{\mathrm{comm}}$.
\end{proof}

\subsubsection{Two-dimensional base case}

Dimension two has a simple structure: every nonconstant diagonal coordinate separates the two spectral sectors.
The local closure is therefore determined by whether a set of labels is empty.
The proposition below shows that this closure globalizes on every cover, without any geometric condition.

\begin{proposition}[Two-dimensional base case]
	\label{THM:PackagingNativeDimensionTwoClassification}
	Let \(I\) be a nonempty finite label set.	
	There is exactly one empty-preserving closure in \(\mathbf{CommCl}_2(I)\), namely
	\[
	{
		u_I(\varnothing)=\varnothing,
		\qquad
		u_I(S)=I
		\quad
		(S\ne\varnothing).
	}
	\]	
	Consequently, for every finite cover \(\mathcal H\),
	\[
	{
		\mathsf{UG}^{=2}_{\mathrm{comm}}(2,\mathcal H)
	}
	\]
	and therefore also
	\[
	{
		\mathsf{UG}_{\mathrm{comm}}(2,\mathcal H).
	}
	\]	
	Thus in native dimension two no chordality or clique-guarding hypothesis is required.
\end{proposition}

\begin{proof}
	In a common eigenbasis of $\mathbb C^2$, empty preservation means that every $D_a$ is non-scalar, hence every coordinate separates the two spectral sectors.
	Therefore every nonempty set of labels determines every label, so the displayed closure is the only member of $\mathbf{CommCl}_2(I)$.
	For an atlas, use the same closure globally on $V$.
	It is realized in dimension two by taking every $D_a=\operatorname{diag}(0,1)$.
\end{proof}

\subsubsection{Three-dimensional cycle-free base case}

In dimension three, joint spectral partitions give a direct description of which single labels determine others and which pairs are incomparable.
Theorem~\ref{THM:PackagingNativeDimensionThreeClassification} uses this to classify exact globalization by chordality and guarded triangles.
The dimension-jump proposition then shows that a global extension may need four dimensions even when every local chart needs only three.

\begin{lemma}[Structure of three-dimensional native dependency closures]
	\label{LEM:PackagingNativeDimensionThreeStructure}
	Let \(I\) be a nonempty finite label set and let
	$
	c\in\mathbf{CommCl}_3(I).
	$
	For \(a,b\in I\),
	\begin{itemize}
		\item write
		$
		a\rightsquigarrow_c b
		\quad\Longleftrightarrow\quad
		b\in c(\{a\}),
		$
		
		\item write
		$
		a\perp_c b
		$
		when neither
		$
		a\rightsquigarrow_c b
		$
		nor
		$
		b\rightsquigarrow_c a
		$
		holds.	
	\end{itemize}	
	Then:	
	\begin{enumerate}		
		\item \(\rightsquigarrow_c\) is a preorder.
		
		\item If
		$
		a\rightsquigarrow_c b
		$
		and
		$
		b\not\rightsquigarrow_c a,
		$
		then
		$
		{
			c(\{a\})=I.
		}
		$
		
		\item For every nonempty \(S\subseteq I\),
		\[
		{
			c(S)
			=
			\begin{cases}
				I,
				&
				\text{if \(S\) contains a pair \(a\perp_c b\),}
				\\[1ex]
				\displaystyle
				\bigcup_{a\in S}c(\{a\}),
				&
				\text{otherwise.}
			\end{cases}
		}
		\]
		
		\item The closure \(c\) is generated by
		\[
		{
			\begin{aligned}
				\mathcal B_3(c)
				={}&
				\{a\to b:
				a\ne b,\ 
				a\rightsquigarrow_c b\}
				\\
				&\cup
				\{\{a,b\}\to z:
				a\perp_c b,\ 
				z\in I\setminus\{a,b\}\}.
			\end{aligned}
		}
		\]		
	\end{enumerate}
\end{lemma}

\begin{proof}
	Pass to the joint spectral quotient $\Omega$, where $2\le|\Omega|\le3$.
	If $|\Omega|=2$, every nonconstant coordinate induces the discrete partition, so every singleton determines every label.
	If $|\Omega|=3$, a nonconstant coordinate induces either the discrete partition or one of the three $2+1$ partitions.
	A discrete-type label determines everything.
	Equal $2+1$ types determine one another.
	Distinct $2+1$ types are singleton-incomparable and their meet is discrete.
	These three observations give (i)--(iii) immediately, and the unary implications together with the binary rules for incomparable pairs generate exactly the cases described in (iii), proving (iv).
\end{proof}

\begin{lemma}[Closure of a three-type native tuple]
	\label{LEM:PackagingThreeTypeClosureFormula}
	Put
	\[
	P_0=\operatorname{diag}(0,0,1),
	\qquad
	P_1=\operatorname{diag}(0,1,0),
	\qquad
	P_2=\operatorname{diag}(1,0,0).
	\]
	Let \(I\) be finite and let
	\[
	\theta:I\to\{0,1,2\}.
	\]
	For
	\[
	D_a=P_{\theta(a)}
	\qquad(a\in I),
	\]
	let \(c_\theta=c_{\mathcal D}\).
	
	Then
	\[
	c_\theta(\varnothing)=\varnothing,
	\]
	and for \(S\ne\varnothing\),
	\[
	{
		c_\theta(S)
		=
		\begin{cases}
			\theta^{-1}(\theta(S)),
			&
			|\theta(S)|=1,
			\\[1ex]
			I,
			&
			|\theta(S)|\ge2.
		\end{cases}
	}
	\]
	Consequently \(c_\theta\) depends only on the equality relation induced by \(\theta\).
\end{lemma}

\begin{proof}
	The matrices $P_0,P_1,P_2$ realize the three $2+1$ partitions of a three-point spectral set.
	One type leaves its corresponding $2+1$ partition, whereas two distinct types have discrete common refinement.  Proposition~\ref{PROP:PackagingNativeCommutantDependence} therefore gives the displayed formula and shows that only the equality pattern of $\theta$ matters.
\end{proof}

\begin{proposition}[Exact globalization can force a strict realization-dimension increase]
	\label{PROP:PackagingNativeStrictGlobalDimensionJump}
	Let $V$ have four labels, and let $\mathcal B_4=\{V\setminus\{v\}:v\in V\}$ be its simplex-boundary cover.
	Assign to every three-label facet $I$ the uniform-circuit closure
	\[
	 u_I(S)=
	 \begin{cases}
	 S,&|S|\le1,\\
	 I,&|S|\ge2.
	 \end{cases}
	\]
	The resulting atlas is compatible and each local closure has native realization dimension exactly three.
	It has an exact global closure, but \emph{every} exact global closure $C'$ satisfies
	\[
	\boxed{d_{\mathrm{comm}}(C')\ge4.}
	\]
	This bound is attained by an exact global closure $C$ with $d_{\mathrm{comm}}(C)=4$.
	Consequently,
	\[
	\mathsf{UG}_{\mathrm{comm}}(3,\mathcal B_4)
	\quad\text{holds, but}\quad
	\neg\mathsf{UG}^{=3}_{\mathrm{comm}}(3,\mathcal B_4).
	\]
\end{proposition}

\begin{proof}
	The restriction of $u_I$ to any proper subset of $I$ is the identity closure.
	Two distinct facets intersect in two labels, so all local closures agree exactly on overlaps.
	Lemma~\ref{LEM:PackagingNativeUniformCircuit} gives $d_{\mathrm{comm}}(u_I)=3$.

	Define a closure on $V$ by
	\[
	C(S)=
	\begin{cases}
	S,&|S|\le1,\\
	V,&|S|\ge2.
	\end{cases}
	\]
	The displayed operator is extensive, monotone, and idempotent.
	Its restriction to each three-element facet is exactly $u_I$, so it is an exact global extension.

	Suppose $C'$ is any exact global closure realized by a commuting self-adjoint tuple in dimension at most three.
	Retain one row from each joint spectral sector.
	There are at most three rows, and their coordinate equality partitions represent $C'$ by Proposition~\ref{PROP:PackagingNativeCommutantDependence}.
	Exactness on the facets implies that no attribute is constant and no singleton attribute determines any distinct attribute: any two distinct labels lie together in a facet, where both singleton closures contain only their premises.
	With two joint spectral rows, every nonconstant attribute separates the rows and therefore determines every other attribute, a contradiction.
	With three rows, an attribute inducing the discrete three-block partition would likewise determine all other attributes.
	Thus every one of the four labels must induce a nontrivial two-block partition of a three-row set.
	There are exactly three such $2+1$ partitions.
	By the pigeonhole principle two different labels induce the same partition and therefore determine each other, again contradicting the prescribed facet closure.
	A one-row realization is ruled out by empty preservation.
	Hence $d_{\mathrm{comm}}(C')\ge4$ for every exact global closure $C'$.

	For attainability, take four rows $\Omega=\{1,2,3,4\}$.
	Choose four distinct two-element subsets, for instance
	\[
	\{1,2\},\qquad\{1,3\},\qquad\{1,4\},\qquad\{2,3\},
	\]
	and for each label use the partition whose only nonsingleton block is the corresponding pair.
	Realize each partition as a real-valued function on $\Omega$ with distinct values on distinct blocks, and take the corresponding four diagonal self-adjoint matrices.
	Each attribute is nonconstant.
	No one such partition refines another, so no singleton determines a different label.
	The common refinement of any two distinct partitions is discrete, because their nonsingleton pairs differ; therefore any two labels determine every label.
	The resulting closure is exactly $C$ and is realized in dimension four.
	The lower bound just proved gives $d_{\mathrm{comm}}(C)=4$.

	Finally, universal exact globalization for this particular cover and local bound $d=3$ follows independently from Theorem~\ref{THM:PackagingNativeDimensionThreeClassification}, since $\mathcal B_4$ has complete primal graph and guards every three-vertex clique.
	The preceding atlas shows that the stronger same-dimension statement fails.
\end{proof}

\begin{theorem}[Three-dimensional base case]
	\label{THM:PackagingNativeDimensionThreeClassification}
	Let \(\mathcal H\) be a finite cover of a finite label set \(V\).
	The following are equivalent:	
	\begin{enumerate}		
		\item
		\[
		{
			\mathsf{UG}_{\mathrm{comm}}(3,\mathcal H).
		}
		\]
		
		\item The primal graph of \(\mathcal H\) is chordal and every three-vertex clique is guarded.		
	\end{enumerate}	
	When these conditions hold, the union
	\[
	T=
	\bigcup_{I\in\mathcal H}\mathcal B_3(c_I)
	\]
	of the canonical local bases generates the least exact global closure.	
\end{theorem}

\begin{proof}
	For sufficiency, let $T$ be the union of the canonical bases $\mathcal B_3(c_I)$.
	If the target chart $I_0$ is empty, exactness follows from empty preservation because no rule has an empty premise.
	Otherwise choose a perfect elimination ordering with $I_0$ last, as permitted by Lemma~\ref{LEM:PackagingCliqueLastElimination}, and eliminate the labels outside $I_0$.
	
	At each stage retain a Horn theory whose models are exactly the projections of the models of $T$.
	Every nontrivial rule is a canonical consequence of	an original chart, has a unary or singleton-incomparable binary premise, and has guarded support of cardinality at most three.
	The initial theory has these properties by Lemma~\ref{LEM:PackagingNativeDimensionThreeStructure}.
	
	At the next simplicial vertex $v$, Lemma~\ref{LEM:PackagingHornProjection} produces nontrivial raw resolvents from parents
	\[
	P\to v,\qquad Q\cup\{v\}\to b,
	\qquad |P|\in\{1,2\},\quad |Q|\in\{0,1\}.
	\]
	All surviving labels of the parent supports are neighbors of $v$, hence form a clique.
	Each triangle used below is therefore guarded, and compatibility transfers every relevant consequence or singleton incomparability between original charts containing its support.
	
	If $P=\{a,c\}$, then $a\perp c$ by the canonical-premise invariant.
	The guarding chart for $\{a,c,b\}$ gives $\{a,c\}\to b$, which implies the raw resolvent, even when its premise also contains $Q$.
	If $P=\{a\}$ and $Q=\varnothing$, the chart guarding $\{a,v,b\}$ transfers both parents.
	Transitivity of its singleton dependence gives $a\to b$.
	
	In the remaining case $P=\{a\}$ and $Q=\{c\}$, the second parent has $v\perp c$.
	In the chart guarding $\{a,v,c\}$, either $v\to a$, so $a$ and $v$ are singleton-equivalent, or $a\to v$ is strict.
	In the first case, either singleton implication between $a$ and $c$ would contradict $v\perp c$.
	Thus $a\perp c$, and the chart guarding	$\{a,c,b\}$ gives $\{a,c\}\to b$.
	In the second case, the strict implication transfers to the chart guarding $\{a,v,b\}$, where
	Lemma~\ref{LEM:PackagingNativeDimensionThreeStructure} makes $a$ singleton-universal and gives $a\to b$.
	
	Each replacement is a canonical unary or incomparable-binary consequence of an original guarding chart, supported on at most three surviving labels.
	It is therefore satisfied by every projected model of the original theory $T$.
	Conversely, the replacement implies the corresponding raw Horn resolvent.
	Together with Lemma~\ref{LEM:PackagingHornProjection}, these two directions preserve the projected model family exactly at each elimination stage.
	The canonical premise form, original-chart provenance, and guarded-support invariant are therefore maintained.
	
	After eliminating all labels outside $I_0$, every surviving rule is valid in $c_{I_0}$ by compatibility.
	Conversely, the complete canonical basis $\mathcal B_3(c_{I_0})$ survives because it uses only labels in $I_0$.
	The generated global closure is therefore exact on every original chart	and is least among exact global extensions.
	
	For necessity, choose an induced cycle $v_0,v_1,\ldots,v_{q-1},v_0$ with $q\ge4$, or an unguarded triangle with $q=3$, and put $S=\{v_0,\ldots,v_{q-1}\}$.
	Every chart meets $S$ in at most two vertices, and any such pair is a cycle edge.
	For each chart $I$, define a separate map $\theta_I:I\to\{0,1,2\}$.
	All labels in $I\setminus S$ have type $0$.
	An isolated cycle label has type $1$.
	The two endpoints of a path edge $\{v_i,v_{i+1}\}$, $0\le i<q-1$, both have type $1$.
	On the closing edge, $v_{q-1}$ has type $1$ and $v_0$ has type $2$.
	Use the three spectral types of Lemma~\ref{LEM:PackagingThreeTypeClosureFormula} to define $c_I$.
	
	These are local type assignments, not a single type assignment on $V$.
	Their equality patterns agree on every overlap:
	off-cycle labels are all of one type, every cycle label differs from that type, and a pair of cycle labels occurs together only on its specified edge, whose equality or inequality status is fixed.
	The cited lemma therefore gives a compatible atlas of empty-preserving closures in $\mathbf{CommCl}_3$.
	
	Each path edge forces $v_i\to v_{i+1}$.
	Any exact global closure would compose these implications to obtain $v_0\to v_{q-1}$, while the closing-edge chart forbids that consequence because its endpoints have different spectral types.
	This contradiction proves necessity.
\end{proof}

\subsection{Boundary obstructions and weak-filling lemmas}

The all-dimensional argument uses the exact simplex-boundary obstruction and conservative completion of missing clique faces.
A weak filler may have Carath\'eodory number at most \(d-1\) without belonging to \(\mathbf{CommCl}_d\).
Theorem~\ref{THM:PackagingNativeRecoveryOfWeakFillers} recovers its consequences in original charts under the original guarding hypothesis.
Lemma~\ref{THM:PackagingSupportPreservingFillerElimination} eliminates auxiliary filler steps while preserving a specified support.
Both results are proved independently of the native clique-resolution theorem; see Remark~\ref{REM:PackagingNativeFillingNoCircularity}.

\subsubsection{Conservative filling and native recovery}

On a simplex boundary, every proper face lies in an original chart.
The first theorem determines the native obstruction threshold in this geometry.
Subsequent filling results preserve the original global consequence relation and recover auxiliary filler consequences in original charts.

\begin{theorem}[Sharp native simplex-boundary obstruction]
	\label{THM:PackagingNativeSimplexBoundaryGeneral}
	Let
	\[
	d\ge2,
	\qquad
	|V|=n\ge3,
	\]
	and let
	\[
	\partial\Delta_{n-1}
	=
	\{V\setminus\{x\}:x\in V\}
	\]
	be the simplex-boundary cover.	
	Then
	\[
	{
		\mathsf{UG}_{\mathrm{comm}}
		(d,\partial\Delta_{n-1})
		\quad\Longleftrightarrow\quad
		n\ge d+1.
	}
	\]
\end{theorem}

\begin{proof}
	Write
	$
	V_u=V\setminus\{u\}
	~
	(u\in V).
	$	
	\begin{enumerate}
		\item \emph{Sufficiency: $n\ge d+1$.}
		
		Let $(c_u)_{u\in V}$ be a compatible atlas with
		$
		c_u\in\mathbf{CommCl}_d(V_u)
		\quad
		(u\in V).
		$
		Suppose, toward a contradiction, that no exact global extension exists.
		
		On each facet, take the family of all inclusion-minimal nontrivial valid implications.
		These families generate the local closures:
		for every nontrivial consequence, a	minimal determining subset of its premise supplies one of the chosen rules.
		
		By Lemma~\ref{LEM:PackagingTwoRuleObstruction},	there exist distinct $x,y\in V$ and two original-chart implications
		\[
		A\to x,
		\qquad
		B\cup\{x\}\to y,
		\]
		each with inclusion-minimal premise, such that
		\[
		A,B\subseteq D:=V\setminus\{x,y\},
		\qquad
		A\cup B=D,
		\]
		and
		$
		c_x(D)=D.
		$
		
		\begin{itemize}
			\item \emph{Step 1: Dimension drop on $D$.}
			
			Compatibility on the overlap $V_x\cap V_y=D$ gives
			$
			c_D:=c_x|_D=c_y|_D.
			$
			
			Since
			$
			D=V_x\setminus\{y\},
			\quad
			c_x(D)=D\ne V_x,
			$
			the set $D$ is not a key of $c_x$.
			
			Lemma~\ref{LEM:PackagingNativeDimensionDrop} therefore gives
			\[
			d_{\mathrm{comm}}(c_D)
			\le d_{\mathrm{comm}}(c_x)-1
			\le d-1.
			\tag{1}
			\]
			
			\item \emph{Step 2: The restriction to $B$ is the identity closure.}
			
			Let $V_z$ be an original facet containing the second implication
			$
			B\cup\{x\}\to y.
			$
			Because this implication has an	inclusion-minimal premise, Lemma~\ref{LEM:PackagingNativeMinimalPremiseIdentity} gives
			$	
			c_z|_{B\cup\{x\}}
			=
			\operatorname{id}_{B\cup\{x\}}.
			$
			In particular,
			$
			c_z|_B=\operatorname{id}_B.
			$		
			Since
			$
			B\subseteq V_z\cap V_x,
			$
			compatibility yields
			\[
			c_D|_B
			=
			c_x|_B
			=
			c_z|_B
			=
			\operatorname{id}_B.
			\tag{2}
			\]
			
			\item \emph{Step 3: Each label of $A\setminus B$ is indispensable in $D$.}
			
			Put
			$
			P=A\setminus B.
			$
			Because $A\cup B=D$, we have the disjoint decomposition
			$
			D=B\sqcup P.
			$		
			We claim that
			\[
			p\notin c_D(D\setminus\{p\})
			\qquad(p\in P).
			\tag{3}
			\]
			
			Fix $p\in P$ and suppose instead that
			$
			p\in c_D(D\setminus\{p\}).
			$
			Put
			$
			T=D\setminus\{p\}.
			$		
			Since $c_D=c_y|_D$, we have
			$
			p\in c_y(T).
			$		
			The first original implication $A\to x$ is supported in $V_y$, and hence is	valid in $c_y$ by compatibility.		
			Moreover,
			$
			A\setminus\{p\}\subseteq T.
			$
			Thus $A\subseteq c_y(T)$, and the first implication gives
			\[
			x\in c_y(T).
			\tag{4}
			\]
			
			Since
			$
			T\cup\{x\}\subseteq V_y\cap V_p,
			$
			compatibility transfers (4) to
			\[
			x\in c_p(T).
			\tag{5}
			\]
			
			Now $p\notin B$, so
			$
			B\subseteq T.
			$
			The second original implication $B\cup\{x\}\to y$ is supported in $V_p$, because
			$
			p\notin B\cup\{x,y\}.
			$
			Hence it is valid in $c_p$.		
			Combining this implication with (5) gives
			\[
			y\in c_p(T).
			\tag{6}
			\]
			
			Finally,
			$
			T\cup\{y\}\subseteq V_p\cap V_x.
			$
			Compatibility transfers (6) to
			$
			y\in c_x(T)
			\subseteq c_x(D)
			=D.
			$
			This contradicts $y\notin D$.		
			Therefore (3) holds for every $p\in P$.
			
			\item \emph{Step 4: Iterated dimension drop.}
			
			Enumerate
			$
			P=\{p_1,\ldots,p_t\},
			~~
			t=|P|,
			$
			and define
			\[
			S_j=D\setminus\{p_1,\ldots,p_j\}
			\qquad(0\le j\le t).
			\]
			Thus
			\[
			S_0=D,
			\qquad
			S_t=B,
			\qquad
			S_{j-1}=S_j\cup\{p_j\}.
			\]
			
			For every $1\le j\le t$, we have
			$
			S_j\subseteq D\setminus\{p_j\}.
			$
			By monotonicity and (3),
			$
			p_j\notin c_D(S_j).
			$
			
			Consequently,
			$
			p_j\notin
			(c_D|_{S_{j-1}})(S_j),
			$
			so $S_j$ is not a key of the restricted closure $c_D|_{S_{j-1}}$.
			
			Applying Lemma~\ref{LEM:PackagingNativeDimensionDrop} to this restricted closure gives
			\[
			d_{\mathrm{comm}}(c_D|_{S_j})
			\le
			d_{\mathrm{comm}}(c_D|_{S_{j-1}})-1.
			\]
			
			Iterating and using (1), we obtain
			\[
			\begin{aligned}
				d_{\mathrm{comm}}(c_D|_B)
				&\le d_{\mathrm{comm}}(c_D)-t\\
				&\le d-1-t.
			\end{aligned}
			\tag{7}
			\]
			
			But (2) and Theorem~\ref{THM:PackagingNativeIdentityDimension} give
			\[
			d_{\mathrm{comm}}(c_D|_B)
			=
			d_{\mathrm{comm}}(\operatorname{id}_B)
			=
			|B|+1.
			\tag{8}
			\]
			This identity also holds when $B=\varnothing$.
			
			Combining (7) and (8) yields
			$
			|B|+1\le d-1-|P|.
			$
			Hence
			$
			|D|=|B|+|P|\le d-2.
			$			
			Since $|D|=n-2$, it follows that
			$
			n\le d,
			$
			contrary to $n\ge d+1$.			
			Therefore every compatible $\mathbf{CommCl}_d$-atlas on the simplex boundary admits an exact global extension whenever $n\ge d+1$.
		\end{itemize}

		\item \emph{Necessity and sharpness: $3 \le n \le d$.}
		
		Choose pairwise distinct $x,y,z\in V$.
		Define the local atlas by
		\[
		c_y=u_{V_y},
		\qquad
		c_z=u_{V_z},
		\]
		and
		\[
		c_w=\operatorname{id}_{V_w}
		\qquad(w\notin\{y,z\}).
		\]
		
		Every facet has $n-1$ labels.
		By Theorem~\ref{THM:PackagingNativeIdentityDimension}, an identity facet has native dimension
		$
		d_{\mathrm{comm}}(\operatorname{id}_{V_w})
		=n.
		$
		By Lemma~\ref{LEM:PackagingNativeUniformCircuit}, each uniform circuit facet has dimension
		$
		d_{\mathrm{comm}}(u_{V_w})=n-1.
		$		
		Since $n\le d$, all local closures belong to $\mathbf{CommCl}_d$.
		
		The atlas is compatible:
		every intersection of two distinct facets is a proper subset of each, and both the identity closure and the uniform circuit closure restrict to the identity on every proper subset.
		
		We show that no exact global extension exists.		
		Put
		\[
		D=V\setminus\{x,y\},
		\qquad
		E=(D\setminus\{z\})\cup\{x\}.
		\]
		Both sets have cardinality $n-2$.
		
		In the uniform circuit facet $V_y$, the set $D$ is a key.
		Hence
		\[
		x\in c_y(D).
		\tag{9}
		\]
		
		Similarly,
		$
		E=V_z\setminus\{y\}
		$
		is a key of the uniform circuit	closure $c_z$, so
		\[
		y\in c_z(E).
		\tag{10}
		\]
		
		Suppose an exact global closure $C$	existed.		
		By (9) and exactness on $V_y$,
		$
		x\in C(D).
		$
		Since
		$
		E\subseteq D\cup\{x\}
		\subseteq C(D),
		$
		equation (10), exactness on $V_z$, monotonicity, and idempotence imply
		$
		y\in C(E)\subseteq C(D).
		$
		
		On the other hand, $x\notin D$ and the missing-$x$ facet carries the identity closure.
		Therefore
		$
		C(D)\cap V_x
		=
		c_x(D)
		=
		D.
		$
		Since $y\in V_x\setminus D$, this contradicts $y\in C(D)$.	
		Thus the compatible atlas admits no exact global extension whenever	$3\le n\le d$.
	\end{enumerate}	
	Combining sufficiency and necessity proves
	\[
	\mathsf{UG}_{\mathrm{comm}}
	(d,\partial\Delta_{n-1})
	\quad\Longleftrightarrow\quad
	n\ge d+1.
	\]
\end{proof}

\begin{definition}[Weakly fillable clique and least weak filler]
	\label{DEF:PackagingWeaklyFillableClique}
	Let
	$
	(c_I)_{I\in\mathcal H}
	$
	be a compatible closure atlas and let \(M\) be a clique of its primal graph.	
	The induced atlas on \(M\) is
	$
	\left(
	c_I|_{I\cap M}
	\right)_{I\in\mathcal H}.
	$
	\begin{itemize}
		\item We call \(M\) \textbf{weakly fillable} when this induced atlas admits an exact closure extension on \(M\).	
		
		\item If \(M\) is weakly fillable, let
		$
		\lambda_M
		$
		denote the least closure on \(M\) preserving every induced local consequence.	
		By Lemma~\ref{LEM:PackagingClosureRestrictionGeneration}, the existence of any exact extension implies that this least generated closure is itself exact.
		We call
		$
		\lambda_M
		$
		the \textbf{least weak filler} of \(M\).	
		No native realization dimension is imposed on \(\lambda_M\).
	\end{itemize}	
\end{definition}

\begin{lemma}[Premise-bound preservation by an exact face filler]
	\label{LEM:PackagingWeakFillerPremiseBound}	
	Let $M$ be a finite set and let $\mathcal K$ be a cover of $M$ such that every proper subset of $M$ is guarded.	
	Let $(c_J)_{J\in\mathcal K}$ be a closure atlas such that
	$
	\operatorname{car}(c_J) \le k
	~ (J\in\mathcal K),
	$
	where $k\ge1$.	
	Suppose the atlas admits an exact extension, and let $\lambda_M$ be its least exact extension.	
	Then we have
	\[
	\operatorname{car}(\lambda_M)\le k.
	\]
\end{lemma}

\begin{proof}
	For each $J\in\mathcal K$, choose an implicational basis $T_J$ of $c_J$ whose premises have size at most $k$.
	Empty premises are permitted.	
	Such bases exist by the Carath\'eodory bound.	
	Put
	\[
	T=\bigcup_{J\in\mathcal K}T_J.
	\]	
	Since an exact extension exists, Lemma~\ref{LEM:PackagingClosureRestrictionGeneration} shows that $\lambda_M$ is precisely the closure generated by $T$.	
	Let $S \subseteq M$ and $b \in \lambda_M(S)$.	
	If $b \in S$, take $A = \{b\}$.
	Then $|A| = 1 \le k$ and $b \in\lambda_M(A)$.	
	Suppose henceforth that $b\notin S$.	
	\begin{itemize}
		\item \emph{Case 1: $S\cup\{b\}\subsetneq M$.}
		
		By the guarding hypothesis, there is $J\in\mathcal K$ containing $S\cup\{b\}$.		
		Exactness gives
		$
		b\in\lambda_M(S)\cap J=c_J(S).
		$		
		Since $\operatorname{car}(c_J)\le k$, there exists $A\subseteq S$ with
		$
		|A|\le k,
		~
		b\in c_J(A).
		$		
		Exactness again gives
		$
		b\in\lambda_M(A).
		$
		
		\item \emph{Case 2: $S\cup\{b\}=M$.}
		
		Then
		$
		S=M\setminus\{b\}.
		$		
		Since $b\in\lambda_M(S)\setminus S$, the set $S$ is not a model of $T$.		
		Hence some rule $A\to a$ in $T$	has $A\subseteq S$ and $a\notin S$.		
		As $M\setminus S=\{b\}$, necessarily $a=b$.
		Thus
		\[
		A\subseteq S,
		\qquad |A|\le k,
		\qquad A\to b\in T.
		\]		
		Since $\lambda_M$ is generated by $T$, we have
		$
		b\in\lambda_M(A).
		$
	\end{itemize}		
	In both cases every consequence has	a determining subset of size at most $k$.
	Therefore
	\[
	\operatorname{car}(\lambda_M)\le k.
	\]
\end{proof}

\begin{theorem}[Uniform premise bound along minimal weak fillings]
	\label{THM:PackagingWeakFillersPreserveNativeBound}	
	Let $d\ge2$, and let
	$
	\mathcal A_0=(c_I)_{I\in\mathcal H_0}
	$
	be a finite compatible atlas on a finite label set $V$, with
	$
	c_I\in\mathbf{CommCl}_d(I)
	\quad
	(I\in\mathcal H_0).
	$	
	Assume that every primal clique of cardinality at most $d$ is guarded by the original cover.	
	Construct a finite filling sequence
	$
	\mathcal A_0
	\subseteq\mathcal A_1
	\subseteq\cdots\subseteq\mathcal A_r
	$
	as follows.	
	At stage $j$, choose a clique $M_j$ that is inclusion-minimal among all unguarded cliques of the current cover, and whose induced atlas is weakly fillable.
	Obtain $\mathcal A_j$ by adjoining its least weak filler $\lambda_{M_j}$.	
	Then every chart closure on a nonempty chart at every stage is empty-preserving and has	Carath\'eodory number at most $d-1$.	
	In particular,
	\[
	\operatorname{car}(\lambda_{M_j})\le d-1
	\qquad(1\le j\le r).
	\]	
	Moreover, every filled clique satisfies
	\[
	|M_j|\ge d+1.
	\]
\end{theorem}

\begin{proof}
	Put $k=d-1\ge1$.
	We prove the Carath\'eodory bound and compatibility of the successive atlases by induction on the filling stage.	
	\begin{enumerate}
		\item \emph{Initial stage.}
		
		Every original closure belongs to $\mathbf{CommCl}_d$ and is therefore empty-preserving.		
		By Theorem~\ref{THM:PackagingNativeCommutantDimensionBound}, every nonempty original chart closure satisfies
		\[
		\operatorname{car}(c_I)\le d-1=k.
		\]
		The condition on an empty chart is vacuous.		
		The original atlas is compatible by hypothesis.
		
		\item \emph{Inductive step.}
		
		Suppose that $\mathcal A_{j-1}$ is compatible and that every nonempty chart closure in it is empty-preserving and has Carath\'eodory number at most $k$.		
		Let $M=M_j$ be the clique selected at stage $j$.		
		We verify the hypotheses of	Lemma~\ref{LEM:PackagingWeakFillerPremiseBound}	for the induced atlas on $M$.
		
		First, since $M$ is inclusion-minimal among all unguarded cliques, every proper subset
		$
		F\subsetneq M
		$
		is guarded by the current cover.		
		Indeed, every subset of a primal clique is again a primal clique, and a proper unguarded subset would contradict the inclusion-minimality of $M$.		
		Thus, for every such $F$, there exists a current chart $I$ with
		$
		F\subseteq I.
		$
		Since $F\subseteq M$, this gives
		$
		F\subseteq I\cap M.
		$
		Hence every proper subset of $M$ is also guarded by the induced cover on $M$.
		
		Second, the induced atlas consists of the restrictions
		$
		c_I|_{I\cap M}.
		$
		These closures are compatible because restrictions compose and the current atlas is compatible.		
		They also inherit the Carath\'eodory bound $k$.		
		Indeed, if
		\[
		S\subseteq I\cap M,
		\qquad
		b\in(c_I|_{I\cap M})(S),
		\]
		then $b\in c_I(S)$.
		The inductive bound supplies $A\subseteq S$ with
		$
		|A|\le k,
		\quad b\in c_I(A).
		$
		Since $b\in I\cap M$, it follows that
		$
		b\in(c_I|_{I\cap M})(A).
		$		
		Thus every induced local closure has Carath\'eodory number at most $k$.
		Empty preservation is likewise inherited under restriction.
		
		Third, the induced atlas is weakly fillable by the choice of $M$.
		Its least exact extension is precisely $\lambda_M$.		
		Lemma~\ref{LEM:PackagingWeakFillerPremiseBound} therefore applies and gives
		$
		\operatorname{car}(\lambda_M)\le k=d-1.
		$		
		Moreover, $\lambda_M$ is empty-preserving.
		Indeed, if $b\in\lambda_M(\varnothing)$, choose an induced chart $J$ containing $b$.
		Exactness would imply
		\[
		b\in\lambda_M(\varnothing)\cap J
		=c_J(\varnothing)
		=\varnothing,
		\]
		a contradiction.		
		By exactness of the filler, for every current chart $I$,
		\[
		\lambda_M|_{I\cap M}
		=
		c_I|_{I\cap M}.
		\]
		Thus adding the chart \((M,\lambda_M)\) preserves compatibility.		
		The previously existing chart closures are unchanged, and the new closure satisfies the same determination bound.
		This completes the induction.
		
		\item \emph{Size of the filled cliques.}
		
		Every $M_j$ is a clique of the current primal graph.
		Adjoining it as a chart therefore creates no new primal edges.		
		Hence the primal graph remains unchanged throughout the filling sequence.
		
		The original cover guards every primal clique of cardinality at most $d$.
		Since filling only adds charts, this guarding property persists at every stage.		
		But $M_j$ is unguarded when selected.
		Therefore
		$
		|M_j|>d,
		$
		and hence
		$
		|M_j|\ge d+1.
		$
		This proves all assertions.
	\end{enumerate}
\end{proof}

The next lemma records the main conservation property needed in the recovery argument:
adjoining a least weak filler changes neither the primal graph nor the least generated global closure.

\begin{lemma}[Conservative weak filling]
	\label{LEM:PackagingConservativeWeakFilling}	
	Let
	$
	(c_I)_{I\in\mathcal H}
	$
	be a finite compatible closure atlas on \(V\), and let
	$
	C_*
	$
	be its least generated global closure.	
	Let \(M\subseteq V\) be a weakly fillable clique, and let
	$
	\lambda_M
	$
	be its least weak filler.	
	Adjoin
	$
	(M,\lambda_M)
	$
	as one additional chart, obtaining an enlarged atlas, and let
	$
	C_*^{+}
	$
	be its least generated global closure.	
	Then we have:	
	\begin{enumerate}[label=(\roman*),leftmargin=2.5em]		
		\item the enlarged atlas is compatible;
		
		\item adjoining \(M\) creates no new primal edge, so the primal graph is unchanged;
		
		\item
		$
		{
			C_*^{+}=C_*;
		}
		$
		
		\item more precisely, for every \(S\subseteq M\),
		$
		{
			\lambda_M(S)
			\subseteq
			C_*(S)\cap M;
		}
		$
		
		\item if every closure in the original atlas is empty-preserving, then
		$
		{
			\lambda_M(\varnothing)=\varnothing.
		}
		$		
	\end{enumerate}
\end{lemma}

\begin{proof}
	We establish the five assertions, proving the consequence inclusion before conservation.
	
	\noindent\emph{(i) Compatibility.}		
	Since $\lambda_M$ is an exact extension of the induced atlas on $M$, for every original chart $I\in\mathcal H$ we have
	$
	\lambda_M|_{I\cap M}
	=
	c_I|_{I\cap M}.
	$
	Thus the new chart is compatible with every	original chart.
	Compatibility among the	original charts is unchanged.
	
	\par\medskip
	\noindent\emph{(ii) Preservation of the primal graph.}		
	Since $M$ is already a clique of the original primal graph, every edge introduced by adjoining $M$ is already present.
	Hence the primal graph is unchanged.
	
	\par\medskip
	\noindent\emph{(iv) Inclusion of filler consequences.}		
	Consider the restriction $C_*|_M$ of the least generated global closure to $M$.		
	For every original chart $I$ and every $A\subseteq I\cap M$, we have
	\[
	\begin{aligned}
		(c_I|_{I\cap M})(A)
		&=c_I(A)\cap(I\cap M)\\
		&\subseteq C_*(A)\cap M\\
		&=(C_*|_M)(A).
	\end{aligned}
	\]
	Thus $C_*|_M$ preserves every consequence of the induced local atlas.		
	By the least-generated-extension property in Lemma~\ref{LEM:PackagingClosureRestrictionGeneration}, $\lambda_M$ is the least closure on $M$ preserving these consequences.
	Therefore
	$
	\lambda_M\le C_*|_M,
	$
	or equivalently,
	\[
	\lambda_M(S)\subseteq C_*(S)\cap M
	\qquad(S\subseteq M).
	\]
	
	\par\medskip
	\noindent\emph{(iii) Conservation of the global closure.}		
	Since the enlarged atlas contains all original charts, its least generated global closure satisfies
	$
	C_*\le C_*^+.
	$		
	Conversely, assertion \textup{(iv)} shows that every consequence of the newly adjoined closure $\lambda_M$ already holds in $C_*$.		
	Hence $C_*$ preserves every local consequence of the enlarged atlas.
	By the least-generated-extension property,
	$
	C_*^+\le C_*.
	$
	Consequently,
	$
	C_*^+=C_*.
	$
	
	\par\medskip
	\noindent\emph{(v) Empty preservation.}		
	Suppose every original closure is empty-preserving.		
	The induced charts $I\cap M$ cover $M$.
	By exactness of $\lambda_M$, for every original chart $I$,
	\[
	\begin{aligned}
		\lambda_M(\varnothing)\cap(I\cap M)
		&=(c_I|_{I\cap M})(\varnothing)\\
		&=c_I(\varnothing)\cap(I\cap M)\\
		&=\varnothing.
	\end{aligned}
	\]
	Since these intersections cover $M$, it follows that
	$
	\lambda_M(\varnothing)=\varnothing.
	$
	
\end{proof}

\begin{corollary}[Finite conservative filling sequence]
	\label{COR:PackagingConservativeWeakFillingSequence}	
	Start from a finite compatible atlas
	$
	\mathcal A_0
	$
	with least generated closure \(C_*\).	
	Suppose that successively one adjoins finitely many charts
	\[
	(M_1,\lambda_{M_1}),
	\ldots,
	(M_r,\lambda_{M_r}),
	\]
	where at stage \(j\), \(\lambda_{M_j}\) is the least weak filler of the induced atlas on \(M_j\) for the atlas obtained after the preceding fillings.	
	Then, after every stage:	
	\begin{enumerate}[label=(\roman*),leftmargin=2.5em]		
		\item the enlarged atlas remains compatible;
		
		\item the primal graph is unchanged from that of \(\mathcal A_0\);
		
		\item the least generated global closure remains
		$
		{C_*};
		$
		
		\item if \(\mathcal A_0\) is empty-preserving, every adjoined weak filler and every chart of the enlarged atlas is empty-preserving.		
	\end{enumerate}	
	In particular, every closure-theoretic consequence present after the fillings was already a consequence of the original atlas.
\end{corollary}

\begin{proof}
	We argue by induction on the number of fillings.
	
	At stage zero, all assertions hold for the original atlas $\mathcal A_0$.
	
	Suppose they hold after stage $j-1$.
	By hypothesis, $M_j$ is a weakly fillable clique of the current primal graph, and $\lambda_{M_j}$ is its least weak filler.
	
	Lemma~\ref{LEM:PackagingConservativeWeakFilling} therefore applies to the current atlas.
	It shows that adjoining	$(M_j,\lambda_{M_j})$ preserves	compatibility and the primal graph, and leaves the least generated global	closure unchanged.
	
	If $\mathcal A_0$ is empty-preserving, the inductive hypothesis ensures that the current atlas is empty-preserving.
	The same lemma then shows that $\lambda_{M_j}$ is empty-preserving.
	Since the existing chart closures are unchanged, the enlarged atlas is empty-preserving as well.
	
	This proves \textup{(i)}--\textup{(iv)} at every stage.	
	Finally, the least generated global	closure is $C_*$ at every stage.
	Hence every implication valid in the final generated closure was already valid in $C_*$, and therefore was already a consequence of the original atlas.
\end{proof}

\begin{theorem}[Native recovery of a weak filler consequence]
	\label{THM:PackagingNativeRecoveryOfWeakFillers}
	Let \(d\ge2\), and let
	$
	\mathcal A_0
	=
	(c_I)_{I\in\mathcal H}
	$
	be a finite compatible atlas with
	$
	c_I\in\mathbf{CommCl}_d(I).
	$	
	Assume that every clique of the original primal graph of cardinality at most \(d\) is guarded by an original chart.
	Let
	$
	\lambda_M
	$
	be any least weak filler produced by a minimal filling sequence as in Theorem~\ref{THM:PackagingWeakFillersPreserveNativeBound}.	
	Then for every
	$
	S\subseteq M
	$
	and every
	$
	b\in\lambda_M(S)\setminus S,
	$
	there exist	
	\[
	D\subseteq S,
	\qquad
	1\le|D|\le d-1,
	\]
	and an original chart
	$
	I\in\mathcal H
	$
	such that
	$
	{
		D\cup\{b\}\subseteq I
	}
	$
	and
	$
	{
		b\in c_I(D).
	}
	$	
	Thus every consequence of the possibly nonnative filler	\(\lambda_M\) has a bounded-premise witness in the original \(d\)-dimensional native atlas.
\end{theorem}

\begin{proof}
	By Theorem~\ref{THM:PackagingWeakFillersPreserveNativeBound} and conservative filling, $\lambda_M$ is empty-preserving and has Carath\'eodory number at most $d-1$.
	Thus $b\in\lambda_M(S)\setminus S$ has a nonempty witness $D\subseteq S$, $|D|\le d-1$, with $b\in\lambda_M(D)$.
	
	Earlier weak fillings leave the primal graph unchanged.
	Hence $D\cup\{b\}\subseteq M$, being a subset of the current primal clique $M$, is a clique of the original primal graph.
	Since $b\notin D$ and $|D|\le d-1$, its cardinality is at most $d$, so some original chart $I$ guards $D\cup\{b\}$.
	
	By definition, the least weak filler is exact on the induced current atlas on $M$.
	In particular the original chart $I$, which remains a chart of the current atlas, satisfies
	\[
	\lambda_M|_{I\cap M}=c_I|_{I\cap M}.
	\]
	Because $D\cup\{b\}\subseteq I\cap M$, the membership $b\in\lambda_M(D)$ is equivalent to $b\in c_I(D)$, proving the claimed
	original-chart recovery.
\end{proof}

\begin{lemma}[Native recovery in a minimal filling enlargement]
	\label{LEM:PackagingNativeEnlargedChartRecovery}
	Let
	$
	\mathcal A_0\subseteq\mathcal A_1\subseteq\cdots\subseteq\mathcal A_r
	$
	be a minimal weak-filling sequence as in Theorem~\ref{THM:PackagingWeakFillersPreserveNativeBound}, starting from an original compatible \(\mathbf{CommCl}_d\)-atlas whose original primal graph guards every clique of cardinality at most \(d\).	
	Let \(J\) be any chart of the enlarged atlas \(\mathcal A_r\), with closure \(c_J\).	
	If
	\[
	b \in c_J(S)\setminus S,
	\quad S\subseteq J,
	\]
	then there exist
	\[
	D\subseteq S,
	\qquad
	1\le|D|\le d-1,
	\]
	and an original chart \(I\) such that
	\[
	D\cup\{b\}\subseteq I,
	\qquad
	b\in c_I(D).
	\]
\end{lemma}

\begin{proof}
	Let $J$ be a chart of $\mathcal A_r$, and suppose
	\[
	b\in c_J(S)\setminus S,
	\qquad S\subseteq J.
	\]	
	We distinguish two cases.
	
	\emph{Case 1: $J$ is an original chart.}
	
	Then $c_J\in\mathbf{CommCl}_d(J)$.
	By Theorem~\ref{THM:PackagingNativeCommutantDimensionBound}, there exists $D\subseteq S$ such that
	\[
	|D|\le d-1,
	\qquad b\in c_J(D).
	\]	
	Since $c_J$ is empty-preserving, $D$ cannot be empty.
	Hence
	\[
	1\le |D|\le d-1.
	\]	
	Taking $I=J$ gives the required original-chart consequence, with
	\[
	D\cup\{b\}\subseteq I.
	\]
	
	\emph{Case 2: $J$ is an adjoined filler chart.}
	
	Then $c_J=\lambda_J$ is the least weak filler adjoined at some stage of the filling sequence.	
	The hypotheses of Theorem~\ref{THM:PackagingNativeRecoveryOfWeakFillers} are satisfied by the original atlas and the given filling sequence.	
	Applying that theorem to
	\[
	b\in\lambda_J(S)\setminus S
	\]
	gives a nonempty $D\subseteq S$ and an original chart $I$ satisfying
	\[
	|D|\le d-1,
	\qquad
	D\cup\{b\}\subseteq I,
	\qquad
	b\in c_I(D).
	\]	
	Thus the conclusion holds in both cases.
\end{proof}

\begin{lemma}[Native recovery is stable under restriction]
	\label{LEM:PackagingNativeRestrictedEnlargedRecovery}
	In the setting of Lemma~\ref{LEM:PackagingNativeEnlargedChartRecovery}, let \(J\) be any chart of the enlarged atlas, let \(L\subseteq J\), and equip \(L\) with the restricted closure
	$
	\lambda_L=c_J|_L.
	$
	If
	\[
	b\in\lambda_L(S)\setminus S,
	\qquad S\subseteq L,
	\]
	then there exist
	\[
	D\subseteq S,
	\qquad 1\le |D|\le d-1,
	\]
	and an original chart \(I\) such that
	\[
	D\cup\{b\}\subseteq I,
	\qquad
	b\in c_I(D).
	\]
	This shows that we can recover consequences from restricted filler charts in the original native charts.
\end{lemma}

\begin{proof}
	Since
	\[
	\lambda_L=c_J|_L,
	\]
	for every $S\subseteq L$ we have
	\[
	\lambda_L(S)=c_J(S)\cap L.
	\]	
	Thus
	\[
	b\in\lambda_L(S)\setminus S
	\]
	implies
	\[
	b\in c_J(S)\setminus S.
	\]
	
	Lemma~\ref{LEM:PackagingNativeEnlargedChartRecovery} therefore supplies an original chart $I$ and a nonempty set $D\subseteq S$ such that
	\[
	|D|\le d-1,
	\qquad
	D\cup\{b\}\subseteq I,
	\qquad
	b\in c_I(D).
	\]
	
	Since $D\subseteq S\subseteq L$, the recovered premise lies in the restricted source set, as required.
\end{proof}

\begin{lemma}[Reflection across a missing vertex]
	\label{LEM:PackagingMissingVertexReflection}
	Let \(Q\) be a set.
	\begin{itemize}
		\item Let
		$
		v,b\in Q,
		~
		v\ne b,
		$
		and put
		$
		F=Q\setminus\{v\};
		$
		
		\item Let
		$
		A,B\subseteq Q\setminus\{v,b\};
		$	
		
		\item Let \(C\) be a closure on \(Q\) satisfying
		$
		v\in C(A)
		$
		and
		$
		b\in C(B\cup\{v\}).
		$
	\end{itemize}	
	Then
	\begin{itemize}
		\item $b\in C(A\cup B)$.
		
		\item If \(\lambda\) is a closure on \(F\) such that $C|_F=\lambda$, then $b\in\lambda(A\cup B)$.
	\end{itemize}
\end{lemma}

\begin{proof}
	Since
	\[
	A\subseteq A\cup B,
	\]
	monotonicity gives
	\[
	v\in C(A)\subseteq C(A\cup B).
	\]
	
	By extensivity,
	\[
	B\cup\{v\}\subseteq C(A\cup B).
	\]
	
	Therefore monotonicity and idempotence give
	\[
	\begin{aligned}
		b
		&\in C(B\cup\{v\})\\
		&\subseteq C(C(A\cup B))\\
		&=C(A\cup B).
	\end{aligned}
	\]
	
	This proves the first assertion.
	
	For the second, observe that
	\[
	A\cup B\subseteq F,
	\qquad b\in F.
	\]
	If $C|_F=\lambda$, then
	\[
	\lambda(A\cup B)
	=
	C(A\cup B)\cap F.
	\]
	Hence
	\[
	b\in\lambda(A\cup B),
	\]
	as claimed.
\end{proof}

\begin{lemma}[Boundary reflection criterion for a raw Horn bridge]
	\label{THM:PackagingBoundaryReflectionCriterion}
	Let \(Q\) be finite.
	\begin{itemize}
		\item Let
		$
		v,b\in Q,
		~
		v\ne b,
		$
		and put
		$
		F_v=Q\setminus\{v\}.
		$
		
		\item Let
		$
		\mathcal Q
		=
		\{Q\setminus\{q\}:q\in Q\}
		$
		be the simplex-boundary cover of \(Q\), equipped with a compatible family of local closures
		$
		(c_q)_{q\in Q}.
		$
		
		\item Assume that
		$
		{
			c_v=\lambda_{F_v},
		}
		$
		where \(\lambda_{F_v}\) is the declared least weak filler on the missing-\(v\) face.
		
		\item Let
		$
		A,B\subseteq Q\setminus\{v,b\},
		$
		and assume that the boundary atlas contains the local consequences
		$
		A\to v,
		~
		B\cup\{v\}\to b.
		$
		That is, each implication is valid in some boundary chart containing its full support.
	\end{itemize}	
	If the boundary atlas admits an exact global closure
	$
	C_Q
	$
	on \(Q\), then
	$
	b\in
	\lambda_{F_v}(A\cup B).
	$	
	Equivalently,
	$
	b\notin
	\lambda_{F_v}(A\cup B)
	$
	is itself a certificate that this compatible boundary atlas admits no exact global closure.
\end{lemma}

\begin{proof}
	Suppose the boundary atlas admits an exact global closure $C_Q$.
	
	By hypothesis, there are original boundary charts in which the implications
	\[
	A\to v,
	\qquad
	B\cup\{v\}\to b
	\]
	are valid.
	
	Since $C_Q$ is exact on every boundary chart, it preserves both consequences.
	Thus
	\[
	v\in C_Q(A),
	\qquad
	b\in C_Q(B\cup\{v\}).
	\]
	
	Lemma~\ref{LEM:PackagingMissingVertexReflection} therefore gives
	\[
	b\in C_Q(A\cup B).
	\]	
	Now
	\[
	A\cup B\subseteq F_v,
	\qquad
	b\in F_v.
	\]
	
	Exactness on the missing-$v$ face yields
	\[
	\begin{aligned}
		b
		&\in C_Q(A\cup B)\cap F_v\\
		&=c_v(A\cup B)\\
		&=\lambda_{F_v}(A\cup B).
	\end{aligned}
	\]	
	Consequently,
	\[
	b\notin\lambda_{F_v}(A\cup B)
	\]
	rules out every exact global closure of the boundary atlas.
	
	This proves both assertions.
\end{proof}

\begin{definition}[Support-constrained derivability]
	\label{DEF:PackagingSupportConstrainedDerivability}	
	Let $\mathcal H$ be a finite cover of a	finite label set $V$, and let
	$
	\mathcal A=(c_I)_{I\in\mathcal H}
	$
	be a closure atlas.
	\begin{itemize}
		\item Its \textbf{full local implication theory} is
		\[
		\operatorname{Imp}(\mathcal A)
		:=
		\left\{
		A\to b:
		\begin{array}{l}
			\exists I\in\mathcal H,\quad A\subseteq I,\\
			b\in c_I(A)\setminus A
		\end{array}
		\right\}.
		\]
		
		\item For $Q\subseteq V$, define the \textbf{$Q$-supported implication theory}
		by
		\[
		\operatorname{Imp}^{[Q]}(\mathcal A)
		:=
		\left\{
		A\to b\in\operatorname{Imp}(\mathcal A):
		A\cup\{b\}\subseteq Q
		\right\}.
		\]
		
		\item For $S\subseteq Q$, define the \textbf{support-constrained derivability closure} by
		\[
		C_{\mathcal A}^{[Q]}(S)
		:=
		C_{\operatorname{Imp}^{[Q]}(\mathcal A)}(S),
		\]
		where the right-hand side is the closure on $Q$ generated by the $Q$-supported implications.
	\end{itemize}	
	Equivalently,
	$
	b\in C_{\mathcal A}^{[Q]}(S)
	$
	if and only if $b$ admits a finite forward derivation from $S$ using only implications whose premises and conclusions lie in $Q$.
	
	A derivation of length zero is allowed, and implication premises may be empty.
	The source charts of the implications need not themselves be contained in $Q$.	
	More precisely, let
	$
	\mathcal A|_Q
	:=
	\bigl(c_I|_{I\cap Q}\bigr)_{I\in\mathcal H}
	$
	be the induced atlas on $Q$.
	Then
	$
	\operatorname{Imp}^{[Q]}(\mathcal A)
	=
	\operatorname{Imp}(\mathcal A|_Q).
	$
	Consequently, $C_{\mathcal A}^{[Q]}$ is the least closure on $Q$ preserving	every consequence of the induced atlas.
	
	If $C_*$ denotes the global closure	generated by $\operatorname{Imp}(\mathcal A)$, then
	\[
	C_{\mathcal A}^{[Q]}(S)
	\subseteq
	C_*(S)\cap Q
	\qquad(S\subseteq Q).
	\]	
	This inclusion may be strict:
	an unrestricted derivation may leave $Q$ and subsequently return, whereas a $Q$-supported derivation cannot.
\end{definition}

The support constraint is essential in the filling argument.
Equality of the unrestricted generated closures does not by itself imply equality of $C_{\mathcal A}^{[Q]}$: a derivation that leaves $Q$ is unavailable to a subsequent elimination confined to $Q$.
The next two lemmas prove the stronger statement by replacing filler steps, starting with the last filling stage.
We never treat an auxiliary filler as if it were one of the original \(d\)-native charts.

\begin{lemma}[Stagewise expansion of a weak filler]
	\label{LEM:PackagingStagewiseFillerExpansion}	
	Let $\mathcal A_0$ be a finite closure atlas on a finite label set $V$, and fix $Q\subseteq V$.
	
	For $j=1,\ldots,r$, suppose that $\mathcal A_j$ is obtained from $\mathcal A_{j-1}$ by adjoining a clique
	\[
	M_j\subseteq Q
	\]
	with its least weak filler $\lambda_j$.
	
	Then, for every $j=1,\ldots,r$ and every $A\subseteq M_j$,
	\[
	\lambda_j(A)
	\subseteq
	C_{\mathcal A_0}^{[Q]}(A)\cap M_j.
	\]
	
	Equivalently, every	$b\in\lambda_j(A)$ admits a finite derivation from $A$ using only original local implications supported in $Q$.
\end{lemma}

\begin{proof}
	Put
	\[
	K=C_{\mathcal A_0}^{[Q]}.
	\]	
	We prove by induction on $j$ that
	\[
	\lambda_j\le K|_{M_j}.
	\]	
	Fix $j$ and suppose the assertion holds for all earlier filling stages.
	
	By definition, $\lambda_j$ is the least	closure on $M_j$ preserving all consequences of the atlas induced by $\mathcal A_{j-1}$ on $M_j$.	
	We show that $K|_{M_j}$ also preserves all these consequences.
	
	Let $J$ be a chart of $\mathcal A_{j-1}$, with closure $c_J$, and put
	\[
	L=J\cap M_j.
	\]	
	Take $P\subseteq L$ and
	\[
	b\in(c_J|_L)(P)=c_J(P)\cap L.
	\]	
	If $b\in P$, extensivity gives $b\in K(P)$.	
	Suppose $b \notin P$.
	
	\emph{Case 1: $J$ is an original chart.}
	
	Since
	\[
	P\cup\{b\}\subseteq L\subseteq Q,
	\]
	the implication $P\to b$ belongs to	$\operatorname{Imp}^{[Q]}(\mathcal A_0)$.	
	Hence
	$
	b \in K(P).
	$
	
	\emph{Case 2: $J$ was adjoined at an earlier stage $i<j$.}
	
	Then
	\[
	c_J=\lambda_i,
	\qquad
	P\subseteq M_i.
	\]	
	By the induction hypothesis,
	\[
	b\in\lambda_i(P)\subseteq K(P).
	\]	
	Thus in both cases
	\[
	(c_J|_L)(P)
	\subseteq K(P)\cap M_j
	=
	(K|_{M_j})(P).
	\]	
	Therefore $K|_{M_j}$ preserves every consequence of the induced atlas on $M_j$.
	
	By the least-generated-extension property of Lemma~\ref{LEM:PackagingClosureRestrictionGeneration},
	\[
	\lambda_j\le K|_{M_j}.
	\]	
	This completes the induction.	
	Consequently, for every $A\subseteq M_j$,
	\[
	\lambda_j(A)
	\subseteq K(A)\cap M_j.
	\]	
	Finally, $K$ is generated by the finite implication theory $\operatorname{Imp}^{[Q]}(\mathcal A_0)$.
	Hence every membership $b\in K(A)$ has a finite forward derivation from $A$ using only original implications supported in $Q$.
\end{proof}

\begin{lemma}[Support-preserving elimination of conservative fillers]
	\label{THM:PackagingSupportPreservingFillerElimination}	
	Let $\mathcal A_0$ be a finite compatible closure atlas on a finite label set $V$, and fix $Q\subseteq V$.	
	Suppose
	\[
	\mathcal A_0
	\subseteq\mathcal A_1
	\subseteq\cdots\subseteq\mathcal A_r
	\]
	is a finite filling sequence in which $\mathcal A_j$ is obtained from $\mathcal A_{j-1}$ by adjoining a clique $M_j\subseteq Q$ with its least weak filler $\lambda_j$.	
	Then, for every $j=0,\ldots,r$ and every $S\subseteq Q$,
	\[
	C_{\mathcal A_j}^{[Q]}(S)
	=
	C_{\mathcal A_0}^{[Q]}(S).
	\]	
	Equivalently, every finite $Q$-supported derivation in the filled atlas can be replaced by a derivation with the same initial set and conclusion,
	using only original local implications and no labels outside $Q$.
\end{lemma}

\begin{proof}
	Put
	\[
	K=C_{\mathcal A_0}^{[Q]}.
	\]	
	Fix $j\in\{0,\ldots,r\}$.
	
	Since every original chart remains present in $\mathcal A_j$,
	\[
	\operatorname{Imp}^{[Q]}(\mathcal A_0)
	\subseteq
	\operatorname{Imp}^{[Q]}(\mathcal A_j).
	\]	
	Hence
	\[
	K\le C_{\mathcal A_j}^{[Q]}.
	\]	
	For the reverse inequality, consider any implication
	\[
	P\to b
	\in
	\operatorname{Imp}^{[Q]}(\mathcal A_j).
	\]
	
	If it comes from an original chart, then it belongs to $\operatorname{Imp}^{[Q]}(\mathcal A_0)$, so
	\[
	b\in K(P).
	\]	
	Otherwise, it comes from a filler $\lambda_i$ added at some stage $i\le j$.	
	Then
	\[
	P\subseteq M_i,
	\qquad
	b\in\lambda_i(P).
	\]	
	Lemma~\ref{LEM:PackagingStagewiseFillerExpansion} gives
	\[
	b\in K(P).
	\]	
	Thus every implication in $\operatorname{Imp}^{[Q]}(\mathcal A_j)$ is valid in $K$.
	
	By the minimality of the closure generated by these implications,
	\[
	C_{\mathcal A_j}^{[Q]}\le K.
	\]	
	Combining both inequalities yields
	\[
	C_{\mathcal A_j}^{[Q]}
	=
	C_{\mathcal A_0}^{[Q]}.
	\]	
	Finally, the first lemma supplies a finite original-atlas derivation, supported in $Q$, for every filler implication used in a derivation.	
	Replacing each filler-rule application by such a derivation preserves the initial set and the conclusion.
	Monotonicity of forward chaining ensures that subsequent rule applications remain valid.	
	The resulting derivation uses only original local implications and contains no label outside $Q$.
\end{proof}

\subsection{All-dimensional native resolution}
\label{SUBSEC:PackagingGeneralNativeResolution}

Theorem~\ref{THM:PackagingNativeGeneralResolution} recovers the composite of two locally valid implications in an original chart, with premise size at most \(d-1\).
Its proof combines spectral compression, elimination of critical two-rule configurations, and a minimal-counterexample argument, uniformly for \(d\ge3\).

\subsubsection{Minimal premises and relative spectral compression}

A minimal native implication has an identity restriction on its premise, which consumes part of the available spectral dimension.
The compression lemmas bound the remaining capacity for distinctions and shorten determining sets while preserving the required closure state.

\begin{lemma}[Minimal native premises are independent]
	\label{LEM:PackagingNativeMinimalPremiseIdentity}	
	Let
	$
	c\in\mathbf{CommCl}_d(I),
	$
	and let
	$
	A\to b
	$
	be an inclusion-minimal nontrivial implication.	
	Then
	\[
	{A\ne\varnothing,}
	\]
	and
	\[
	{
		c|_A=\operatorname{id}_A.
	}
	\]
	Consequently,
	\[
	{
		1\le |A|\le d-1.
	}
	\]
\end{lemma}

\begin{proof}
	Since \(c\) is empty-preserving and the implication is nontrivial,
	\[
	b\in c(A)\setminus A,
	\]
	one cannot have \(A=\varnothing\).
	Hence
	\[
	A\ne\varnothing.
	\]	
	Suppose that there exist
	\[
	a\in A,
	\qquad
	S\subseteq A\setminus\{a\},
	\]
	such that
	\[
	a\in c(S).
	\]
	Then
	\[
	A\subseteq c(A\setminus\{a\}),
	\]
	because every element of \(A\setminus\{a\}\) is present by extensivity and	\(a\in c(S)\subseteq c(A\setminus\{a\})\).
	Therefore
	\[
	b\in c(A)
	\subseteq
	c(c(A\setminus\{a\}))
	=
	c(A\setminus\{a\}),
	\]
	contradicting inclusion-minimality of the premise \(A\).	
	Thus, for every \(T\subseteq A\),
	\[
	c(T)\cap A=T.
	\]
	Equivalently,
	\[
	c|_A=\operatorname{id}_A.
	\]	
	By Lemma~\ref{LEM:PackagingNativeRestrictionCompression}, the restricted closure \(c|_A\) belongs to \(\mathbf{CommCl}_d(A)\).
	On the other hand, Theorem~\ref{THM:PackagingNativeIdentityDimension} gives
	\[
	d_{\mathrm{comm}}(\operatorname{id}_A)
	=
	|A|+1.
	\]
	Hence
	\[
	|A|+1\le d,
	\]
	and therefore
	\[
	1\le |A|\le d-1.
	\]
\end{proof}

\begin{lemma}[Finite coordinate refinement budget]
	\label{LEM:PackagingPartitionRefinementBudget}
	Let \(\Omega\) be finite, let \(\Pi_0\) be a partition of \(\Omega\), and let
	\[
	(\Pi_u)_{u\in U}
	\]
	be a finite family of partitions.
	For \(T\subseteq U\), put
	\[
	\Pi_T
	=
	\Pi_0\wedge\bigwedge_{u\in T}\Pi_u.
	\]
	If
	\[
	\Pi_T\ne\Pi_U,
	\]
	then there exists
	\[
	u\in U\setminus T
	\]
	such that
	\[
	\Pi_{T\cup\{u\}}
	\text{ strictly refines }\Pi_T.
	\]
	Consequently, if $m\ge1$ and $\beta\ge0$ are integers, $\Pi_U$ has at most $m$ blocks, and $\Pi_0$ has at least $m-\beta$ blocks, then there exists
	\[
	C\subseteq U,
	\qquad
	|C|\le\beta,
	\]
	such that
	\[
	{
		\Pi_C=\Pi_U.
	}
	\]
\end{lemma}

\begin{proof}
	If no remaining partition \(\Pi_u\), \(u\in U\setminus T\), strictly refined \(\Pi_T\), then
	\[
	\Pi_T\wedge\Pi_u=\Pi_T
	\qquad
	(u\in U\setminus T).
	\]
	Taking the meet over all remaining \(u\) would give \(\Pi_U=\Pi_T\), a contradiction.
	Thus some remaining coordinate produces a strict refinement.
	Every strict refinement increases the number of blocks by at least one.
	Starting with at least \(m-\beta\) blocks and ending with at most \(m\) blocks therefore permits at most \(\beta\) strict refinements.
\end{proof}

\begin{lemma}[Relative spectral compression over an identity base]
	\label{LEM:PackagingNativeRelativeSpectralCompression}
	Let
	\[
	d\ge3,
	\qquad
	0\le\beta\le d-2.
	\]
	Let $B,U$ be disjoint finite label sets and let $y$ be a label outside $B\cup U$.
	Assume
	\[
	|B|=d-2-\beta,
	\]
	and let
	\[
	c\in\mathbf{CommCl}_d(B\cup U\cup\{y\}).
	\]
	Assume
	\[
	c|_B=\operatorname{id}_B
	\]
	and
	\[
	y\notin c(B\cup U).
	\]
	Then there exists
	\[
	C\subseteq U,
	\qquad
	|C|\le\beta,
	\]
	such that
	\[
	{
		U\subseteq c(B\cup C).
	}
	\]	
	More precisely, with
	\[
	S=B\cup U,
	\]
	the restricted closure \(c|_S\) admits a native realization on at most \(d-1\) joint spectral sectors for which
	\[
	{
		\Pi_{B\cup C}=\Pi_{B\cup U}.
	}
	\]
	Consequently, the full closure states agree:
	\[
	c(B\cup C)=c(B\cup U)=S.
	\]
\end{lemma}

\begin{proof}
	Put
	$
	S=B\cup U.
	$
	Since \(y\notin c(S)\), the set \(S\) is not a key.
	By Lemma~\ref{LEM:PackagingNativeDimensionDrop}, the restricted closure \(c|_S\) admits a native realization on at most \(d-1\) joint spectral sectors.
	Fix such a realization, let \(\Omega\) be its joint spectral quotient, and write
	\[
	m=|\Omega|\le d-1.
	\]
	For \(T\subseteq S\), let \(\Pi_T\) denote the joint partition of \(\Omega\) induced by the coordinates in \(T\).
	
	Because
	$
	c|_B=\operatorname{id}_B,
	$
	the subtuple indexed by \(B\) realizes the identity closure.
	In the chosen realization of \(c|_S\), the joint spectral sectors of this subtuple are exactly the blocks of \(\Pi_B\).
	By Lemma~\ref{LEM:PackagingNativeRestrictionCompression}, \(\operatorname{id}_B\) therefore admits a native realization on \(|\Pi_B|\) dimensions.
	Theorem~\ref{THM:PackagingNativeIdentityDimension} gives
	\[
	|B|+1
	=
	d_{\mathrm{comm}}(\operatorname{id}_B)
	\le
	|\Pi_B|.
	\]
	Hence
	\[
	|\Pi_B|
	\ge
	|B|+1
	=
	d-1-\beta.
	\]	
	Since \(m\le d-1\), the preceding estimate also gives
	$
	|\Pi_B|\ge d-1-\beta\ge m-\beta.
	$
	Apply Lemma~\ref{LEM:PackagingPartitionRefinementBudget} to the base partition \(\Pi_B\) and the coordinate partitions indexed by \(U\), with the numerical parameter in that lemma equal to the actual sector count \(m=|\Omega|\).
	The final partition \(\Pi_{B\cup U}\) has at most \(m\) blocks, so there exists
	\[
	C\subseteq U,
	\qquad
	|C|\le\beta,
	\]
	such that
	$
	\Pi_{B\cup C}=\Pi_{B\cup U}.
	$	
	By Proposition~\ref{PROP:PackagingNativeCommutantDependence}, every coordinate in \(U\) is constant on each block of \(\Pi_{B\cup U}\).
	The displayed partition equality therefore implies that every coordinate in \(U\) is also constant on each block of \(\Pi_{B\cup C}\).
	Hence
	\[
	U\subseteq c(B\cup C).
	\]
	Finally, $B\cup C\subseteq S$ and the preceding inclusion, together with extensivity, give
	\[
	B\cup C\subseteq S\subseteq c(B\cup C).
	\]
	By monotonicity and idempotence,
	\[
	c(B\cup C)=c(S).
	\]
	Since $c$ is a closure on $S\cup\{y\}$ and $y\notin c(S)$, we have $c(S)=S$, proving the asserted equality.
\end{proof}

\subsubsection{Closure-state compression and critical two-rule bridges}

Closure-state compression preserves the full closure induced on a chosen face, not merely one implication.
The resulting lemmas reduce a hypothetical failure of native recovery to a critical pair of rules, which the bridge theorem excludes under the stated hypotheses.

\begin{lemma}[Global seed with closure-state preservation]
	\label{LEM:PackagingNativeGlobalSeed}
	Let $d\ge3$, and let
	\[
	\mathcal A_0=(c_I)_{I\in\mathcal H},
	\qquad c_I\in\mathbf{CommCl}_d(I),
	\]
	be a compatible atlas on a finite cover $\mathcal H$ of a finite label set $V$ whose original primal graph guards every clique of cardinality at most $d$.
	Let $Q$ be a finite primal clique, let $x,y\in Q$ be distinct, and let $A,B\subseteq Q\setminus\{x,y\}$ satisfy
	\[
	A\to x,\qquad B\cup\{x\}\to y,
	\qquad A\cup B=Q\setminus\{x,y\},
	\]
	where both implications have inclusion-minimal premises in their respective original source closures.	
	Suppose a compatible enlargement of $\mathcal A_0$ guards both $Q\setminus\{x\}$ and $Q\setminus\{y\}$.
	Denote the closures	induced there by $\lambda_x$ and $\lambda_y$, respectively.
	Assume
	\[
	y\notin\lambda_x(A\cup B).
	\]
	Put
	\[
	X=A\setminus B,\qquad \beta=d-2-|B|.
	\]
	Then $\beta\ge0$, and there exists
	\[
	C\subseteq X,\qquad |C|\le\beta,
	\]
	such that
	\[
	\lambda_y(B\cup C)=\lambda_y(A\cup B).
	\]
	In particular,
	\[
	A\subseteq\lambda_y(B\cup C),
	\qquad x\in\lambda_y(B\cup C).
	\]
\end{lemma}

\begin{proof}
	Choose an original source chart for $B\cup\{x\}\to y$.
	By minimality and Lemma~\ref{LEM:PackagingNativeMinimalPremiseIdentity}, its restriction to $B\cup\{x\}$ is the identity.
	Hence
	\[
	|B|+1\le d-1,
	\qquad \beta=d-2-|B|\ge0.
	\]
	The implication $A\to x$ is supported in the missing-$y$ face, so compatibility gives
	\[
	x\in\lambda_y(A).
	\tag{*}
	\]
	Also $A\cup B=B\cup X$.
	
	We process the elements of $X$, maintaining sets $C\subseteq P\subseteq X$ with
	\[
	|C|\le\beta,
	\qquad
	\lambda_y(B\cup C)=\lambda_y(B\cup P).
	\tag{**}
	\]
	Initially $C=P=\varnothing$.
	Suppose the invariant holds and choose $a\in X\setminus P$.
	Put $P'=P\cup\{a\}$.
	
	If $|C|<\beta$, put $C'=C\cup\{a\}$.
	By the invariant,
	\[
	P\subseteq\lambda_y(B\cup C)
	\subseteq\lambda_y(B\cup C'),
	\]
	and $a\in C'$.
	Thus $P'\subseteq\lambda_y(B\cup C')$.
	Since $C'\subseteq P'$, monotonicity and idempotence give
	\[
	\lambda_y(B\cup C')=\lambda_y(B\cup P').
	\]
	
	Suppose instead that $|C|=\beta$.
	Put $U=C\cup\{a\}$ and
	\[
	R=B\cup U\cup\{y\}.
	\]
	Then $|U|=\beta+1$ and $|R|=d$.
	Since $R\subseteq Q$ is a primal clique, an original $d$-native chart $J_R$ guards $R$.
	Write $c_R=c_{J_R}|_R$.
	By compatibility with the original source chart of the second parent,
	\[
	c_R|_B=\operatorname{id}_B.
	\]
	Moreover $B\cup U\subseteq A\cup B$, so
	\[
	y\notin\lambda_x(B\cup U).
	\]
	Because $R\subseteq Q\setminus\{x\}$, compatibility of the enlargement gives
	\[
	y\notin c_R(B\cup U).
	\]
	Lemma~\ref{LEM:PackagingNativeRelativeSpectralCompression} now provides $C'\subseteq U$ with $|C'|\le\beta$ and
	\[
	U\subseteq c_R(B\cup C').
	\]
	Every implication $B\cup C'\to u$, $u\in U$, has support in $Q\setminus\{y\}$.
	Compatibility with the missing-$y$ face therefore gives
	\[
	U\subseteq\lambda_y(B\cup C').
	\tag{***}
	\]
	Since $C\subseteq U$, the displayed inclusion implies
	\[
	\lambda_y(B\cup C)
	\subseteq\lambda_y(B\cup C').
	\]
	By the old invariant, $P\subseteq\lambda_y(B\cup C)$.
	By (***), $a\in\lambda_y(B\cup C')$.
	Hence $P'\subseteq\lambda_y(B\cup C')$.
	Together with $C'\subseteq U\subseteq P'$, this gives
	\[
	\lambda_y(B\cup C')=\lambda_y(B\cup P').
	\]
	Thus the invariant is preserved in both cases.
	
	At termination $P=X$, so (**), together with $A\cup B=B\cup X$, yields
	\[
	\lambda_y(B\cup C)=\lambda_y(A\cup B),
	\qquad |C|\le\beta.
	\]
	Extensivity gives $A\subseteq\lambda_y(A\cup B)=\lambda_y(B\cup C)$.
	By (*), monotonicity and idempotence then give
	\[
	x\in\lambda_y(A)
	\subseteq\lambda_y(B\cup C).
	\]
\end{proof}

\begin{lemma}[Native normalization of a critical two-rule obstruction]
	\label{LEM:PackagingNativeTwoRuleNormalization}
	Let \(d\ge3\), and let
	\[
	(c_I)_{I\in\mathcal H},
	\qquad
	c_I\in\mathbf{CommCl}_d(I),
	\]
	be an original compatible atlas on a finite cover whose original primal graph guards every clique of cardinality at most \(d\).	
	Let \(Q\) be a finite primal clique of the original cover.
	Suppose a compatible enlargement has been obtained from the original atlas by a finite minimal weak-filling sequence.
	No restriction on the supports of earlier filler charts relative to $Q$ is imposed.
	Suppose each coatom face of $Q$ is equipped with the restriction of a chart of this enlargement.	
	Assume the resulting simplex-boundary atlas on \(Q\) has no exact extension, and let
	\[
	x,y\in Q,
	\qquad
	D=Q\setminus\{x,y\},
	\]
	be a two-rule obstruction supplied by Lemma~\ref{LEM:PackagingTwoRuleObstruction}.	
	Then there exist
	$
	A,B\subseteq D
	$
	such that
	\[
	{
		A\to x,
		\qquad
		B\cup\{x\}\to y
	}
	\]
	are inclusion-minimal implications valid in original \(d\)-native charts,
	$
		A \cup B=D,
	$
	and the missing-\(x\) face closure satisfies
	$
		y\notin\lambda_x(D).
	$
\end{lemma}

\begin{proof}
	Let the two boundary consequences supplied by
	Lemma~\ref{LEM:PackagingTwoRuleObstruction} be
	\[
	A_1\to x,
	\qquad
	P_1\to y,
	\]
	with
	\[
	A_1\subseteq D,
	\qquad
	P_1\subseteq D\cup\{x\},
	\]
	and with
	\[
	y\notin\lambda_x(D).
	\]
	
	Because each coatom closure is the restriction of a chart of the compatible minimal weak-filling enlargement, and because the original atlas guards every original primal clique of cardinality at most \(d\), Lemma~\ref{LEM:PackagingNativeRestrictedEnlargedRecovery}	applies to both boundary consequences.
	
	Hence there are original \(d\)-native charts carrying implications
	\[
	A'\to x,
	\qquad
	P'\to y,
	\]
	with
	\[
	A'\subseteq A_1\subseteq D,
	\qquad
	P'\subseteq P_1\subseteq D\cup\{x\},
	\qquad
	1\le |A'|,|P'|\le d-1.
	\]
	Shrink \(A'\) and \(P'\) inside their respective original source closures to inclusion-minimal premises.
	Relabel the resulting premises again by	\(A\) and \(P\).
	
	We claim that
	$
	x\in P.
	$
	Indeed, if \(x\notin P\), then
	$
	P\subseteq D.
	$
	The implication \(P\to y\) is valid in an original chart, and its entire support \(P\cup\{y\}\) lies in the missing-\(x\) face.
	By compatibility with the chart whose restriction is \(\lambda_x\),	Lemma~\ref{LEM:PackagingCompatibleConsequenceTransfer} would give
	\[
	y\in\lambda_x(P)
	\subseteq
	\lambda_x(D),
	\]
	contradicting the defining obstruction.
	Thus
	$
	P=B\cup\{x\}
	$
	for some \(B\subseteq D\).
	
	It remains to prove
	$
	A\cup B=D.
	$
	Suppose instead that there exists
	\[
	z\in D\setminus(A\cup B).
	\]
	Then the missing-\(z\) face
	$
	Q\setminus\{z\}
	$
	contains the complete supports of both original implications
	\[
	A\to x,
	\qquad
	B\cup\{x\}\to y.
	\]
	By compatibility, both implications are valid in its face closure \(\lambda_z\).
	Hence
	\[
	x\in\lambda_z(A)
	\]
	and consequently, by monotonicity and idempotence,
	\[
	y\in\lambda_z(A\cup B).
	\]
	
	Since
	\[
	z\notin A\cup B\cup\{y\},
	\]
	the support
	\[
	A\cup B\cup\{y\}
	\]
	lies in the overlap of the missing-\(z\) and missing-\(x\) faces.
	Compatibility therefore gives
	\[
	y\in\lambda_x(A\cup B)
	\subseteq
	\lambda_x(D),
	\]
	again contradicting
	$
	y\notin\lambda_x(D).
	$
	Thus
	$
	A\cup B=D.
	$
\end{proof}

\begin{theorem}[Nonexistence of native critical two-rule circuits]
	\label{THM:PackagingNoNativeCriticalCircuit}
	Let \(d\ge3\), and let \(Q\) be a finite primal clique of the original cover.
	Start with an original compatible atlas whose local closures lie in \(\mathbf{CommCl}_d\), and form a compatible enlargement by a conservative
	minimal weak-filling sequence supported on proper cliques of \(Q\).
	Suppose	each coatom face \(Q\setminus\{z\}\) is equipped with the restriction \(\lambda_z\) of a chart of this enlargement.
	
	Assume that for distinct \(x,y\in Q\) there are inclusion-minimal implications valid in respective original source closures,
	\[
	A\to x,
	\qquad
	B\cup\{x\}\to y,
	\]
	where
	\[
	A,B\subseteq Q\setminus\{x,y\},
	\qquad
	A\cup B=Q\setminus\{x,y\},
	\]
	and assume
	\[
	y\notin\lambda_x(A\cup B).
	\]
	Assume every clique of cardinality at most \(d\) is guarded by an original native chart.
	Then no such critical two-rule configuration exists.
\end{theorem}

\begin{proof}
	Minimality of the second parent and Lemma~\ref{LEM:PackagingNativeMinimalPremiseIdentity} give
	\[
	|B|+1\le d-1.
	\]
	Hence, with
	\[
	\beta=d-2-|B|,
	\]
	one has \(\beta\ge0\).
	Lemma~\ref{LEM:PackagingNativeGlobalSeed} gives
	\[
	C\subseteq A\setminus B,
	\qquad
	|C|\le\beta,
	\]
	with
	\[
	\lambda_y(B\cup C)=\lambda_y(A\cup B),
	\]
	and in particular
	\[
	x\in\lambda_y(B\cup C).
	\]
	Lemma~\ref{LEM:PackagingNativeRestrictedEnlargedRecovery} supplies an original implication
	\[
	E\to x,
	\qquad
	E\subseteq B\cup C,
	\qquad
	|E|\le d-1.
	\]
	Now
	\[
	K=B\cup C\cup\{x,y\}
	\]
	is a subset of the primal clique \(Q\), hence is itself a clique.
	Its size is	at most
	\[
	(d-2-\beta)+\beta+2=d,
	\]
	so it lies in an original native chart.
	By Lemma~\ref{LEM:PackagingCompatibleConsequenceTransfer}, the implications \(E\to x\) and \(B\cup\{x\}\to y\) are both valid in that guarded chart.
	Starting from \(B\cup C\), the first implication derives \(x\), and the second then derives \(y\).
	Hence
	\[
	y\in c_K(B\cup C),
	\]
	where \(c_K\) denotes the guarded-chart closure restricted to \(K\).
	Since
	\[
	B\cup C\cup\{y\}\subseteq Q\setminus\{x\},
	\]
	Lemma~\ref{LEM:PackagingCompatibleConsequenceTransfer} transfers this consequence to the missing-\(x\) face closure:
	\[
	y\in\lambda_x(B\cup C)\subseteq\lambda_x(A\cup B),
	\]
	contradiction.
\end{proof}

\subsubsection{Fillability of proper faces}

The proper-face filling lemma establishes that every unguarded proper face of a minimal resolution counterexample is weakly fillable.
Its proof uses the preceding boundary and filling results, without invoking native clique resolution.

\begin{lemma}[Proper-face fillability in a minimal native resolution counterexample]
	\label{LEM:PackagingNativeMinimalCounterexampleFillability}
	Let \(d\ge3\), and let
	\[
	(c_I)_{I\in\mathcal H},
	\qquad
	c_I\in\mathbf{CommCl}_d(I),
	\]
	be a compatible atlas on a finite cover whose original primal graph guards every clique of cardinality at most \(d\).	
	Suppose original-chart implications
	\[
	A\to v,
	\qquad
	B\cup\{v\}\to b
	\]
	have primal-clique support
	\[
	Q=A\cup B\cup\{v,b\}
	\]
	and admit no original-chart replacement
	\[
	D\to b,
	\qquad
	D\subseteq A\cup B,
	\qquad
	1\le |D|\le d-1.
	\]
	Assume both parent premises have been shrunk to inclusion-minimal ones, \(v\) remains in the second premise, and \(Q\) has minimum cardinality among all such counterexamples.	
	Let
	\[
	\mathcal A_0
	\subseteq
	\mathcal A_1
	\subseteq\cdots\subseteq
	\mathcal A_t
	\]
	be any finite prefix of a minimal weak-filling sequence supported on proper subsets of \(Q\):
	at every preceding stage, an inclusion-minimal unguarded proper clique is filled by adjoining its least weak filler.	
	Then every inclusion-minimal unguarded proper clique
	\[
	M\subsetneq Q
	\]
	of the current cover \(\mathcal A_t\) is weakly fillable.
\end{lemma}

\begin{proof}
	Assume, towards a contradiction, that the current cover \(\mathcal A_t\) has an inclusion-minimal unguarded proper clique
	\[
	M\subsetneq Q
	\]
	whose induced atlas is not weakly fillable.
	
	By Corollary~\ref{COR:PackagingConservativeWeakFillingSequence}, the current enlargement is compatible and has the same primal graph as the original atlas.
	In particular \(M\), being a current primal clique, is also a primal clique of the original cover.
	Since every original primal clique of size at most \(d\) is guarded, the fact that \(M\) is currently unguarded implies
	$
	|M|\ge d+1.
	$
	
	Every singleton is guarded because the current charts cover the label set.
	Moreover, if \(|M|=2\), then \(M\) is an edge of the current primal graph, hence its two vertices already occur together in a current chart.
	Thus \(M\) would be guarded.
	Therefore
	\[
	{|M|\ge3.}
	\]	
	Because \(M\) is inclusion-minimal unguarded in the current cover, every proper subset of \(M\) is guarded.
	For each
	\[
	z\in M
	\]
	choose a current chart containing the coatom
	\[
	M\setminus\{z\},
	\]
	and let
	\[
	\lambda_z
	\]
	be the restriction of its closure to that coatom.
	Compatibility of the current enlargement makes this restriction independent of the chosen guarding chart.
	The same compatibility shows that the coatom closures \((\lambda_z)_{z\in M}\) agree on every pairwise intersection, so they form a compatible simplex-boundary atlas on \(M\).
	
	We claim that the resulting simplex-boundary atlas
	\[
	(\lambda_z)_{z\in M}
	\]
	has no exact extension on \(M\).
	Indeed, suppose that a closure \(C_M\) on \(M\) extended all coatom closures exactly.
	Let \(I\) be any current chart and put
	$
	T=I\cap M.
	$
	Since \(M\) is unguarded,
	$
	T\subsetneq M.
	$
	Choose \(z\in M\setminus T\).
	Then
	$
	T\subseteq M\setminus\{z\}.
	$
	The current chart \(I\) and the current chart chosen to guard \(M\setminus\{z\}\) belong to the same compatible enlargement, so their restrictions agree on \(T\).
	Exactness of \(C_M\) on the coatom therefore gives
	\[
	C_M|_T
	=
	\lambda_z|_T
	=
	c_I|_T.
	\]
	Thus \(C_M\) would exactly extend the entire induced current atlas on \(M\), contrary to the assumed nonfillability of \(M\).
	
	Hence the simplex-boundary atlas on \(M\) is nonextendible.
	Apply Lemma~\ref{LEM:PackagingTwoRuleObstruction}, followed by Lemma~\ref{LEM:PackagingNativeTwoRuleNormalization}.
	The latter lemma applies to this entire finite prefix:
	its recovery argument does not require earlier filler charts to be supported inside $M$.
	The original guardedness hypothesis and minimal weak-filling conditions are unchanged.
	
	We obtain distinct
	\[
	x,y\in M
	\]
	and inclusion-minimal implications valid in original \(d\)-native charts,
	\[
	A'\to x,
	\qquad
	B'\cup\{x\}\to y,
	\]
	such that
	\[
	A'\cup B'
	=
	M\setminus\{x,y\},
	\]
	while the missing-\(x\) face closure satisfies
	\[
	y\notin\lambda_x(A'\cup B').
	\]
	
	We claim that these two implications admit no original-chart replacement
	\[
	D\to y,
	\qquad
	D\subseteq A'\cup B',
	\qquad
	1\le|D|\le d-1.
	\]
	For if such a replacement were valid in an original chart, then its support would lie in the missing-\(x\) face.
	Since that original chart and the current chart defining \(\lambda_x\) belong to the same compatible enlargement,	Lemma~\ref{LEM:PackagingCompatibleConsequenceTransfer} would give
	\[
	y\in\lambda_x(D)
	\subseteq
	\lambda_x(A'\cup B'),
	\]
	a contradiction.	
	Finally, conservative filling does not change the primal graph.
	Thus \(M\) is also a primal clique of the original cover.
	The pair
	\[
	A'\to x,
	\qquad
	B'\cup\{x\}\to y
	\]
	is therefore a native resolution counterexample supported on the strictly smaller primal clique
	$
	M\subsetneq Q.
	$
	This contradicts the minimal choice of \(Q\).	
	Hence every inclusion-minimal unguarded proper clique of the current cover is weakly fillable.
\end{proof}

\begin{remark}[Independent inputs for proper-face filling]
	\label{REM:PackagingNativeFillingNoCircularity}
	Lemma~\ref{LEM:PackagingNativeMinimalCounterexampleFillability} relies on Theorem~\ref{THM:PackagingNativeRecoveryOfWeakFillers} and Lemma~\ref{LEM:PackagingNativeRestrictedEnlargedRecovery} to recover consequences from previously filled charts.
	Those results use the Carath\'eodory bound for exact weak fillers and the original clique-guarding condition, but not Theorem~\ref{THM:PackagingNativeGeneralResolution}.
\end{remark}

\begin{remark}[Hypotheses used in native clique resolution]
	\label{REM:PackagingNativeResolutionHypothesisLedger}
	The proof of Theorem~\ref{THM:PackagingNativeGeneralResolution} uses the following distinct inputs.
	The finite spectral partition model establishes restriction stability and the dimension-drop bound (Lemmas~\ref{LEM:PackagingNativeRestrictionCompression} and~\ref{LEM:PackagingNativeDimensionDrop}).
	The relative spectral-compression lemma applies only after a missing-face consequence is shown to be absent.
	Its identity-base hypothesis is supplied by minimality of the second parent premise.
	Conservative weak filling changes neither the primal graph nor the globally generated closure.
	We establish the Carath\'eodory bound for each minimal filler before using original-chart recovery.
	The original-chart recovery theorem applies to restrictions of \emph{any} filler produced by such a sequence, not only to fillers supported inside the final clique.
	For a counterexample with minimal support \(Q\), this recovery shows that every proper face can be filled.
	We can then keep all later fillings inside proper subsets of \(Q\).
	Finally, the critical two-rule contradiction uses a guarding chart for $B\cup C\cup\{x,y\}$ of cardinality at most $d$ and transfers the resulting $B\cup C\to y$ consequence to the missing-$x$ face.
	In particular, the argument never assumes that a filler has native dimension at most $d$.
\end{remark}

\subsubsection{All-dimensional native clique resolution}

The preceding lemmas now give the central native result.
Theorem~\ref{THM:PackagingNativeGeneralResolution} takes two locally valid implications with clique support and recovers their combined consequence in an original chart.
Its premise uses at most \(d-1\) labels, for every \(d\ge3\).
We then express the conclusion directly in the language of functional dependencies.

\begin{theorem}[Native clique-resolution theorem]
	\label{THM:PackagingNativeGeneralResolution}
	Let $d\ge3$, and let $(c_I)_{I\in\mathcal H}$ be a compatible atlas on a finite cover $\mathcal H$ of a finite label set $V$, with $c_I\in\mathbf{CommCl}_d(I)$ for every $I\in\mathcal H$.
	Assume that every primal clique of cardinality at most $d$ is guarded by an original chart.
	No chordality assumption is made.	
	Let $v,b\in V$ be distinct, and let $A,B\subseteq V\setminus\{v,b\}$.
	Suppose there are original charts $I_1,I_2\in\mathcal H$ such that
	\[
	\begin{aligned}
		A\cup\{v\}&\subseteq I_1,
		&v&\in c_{I_1}(A),\\
		B\cup\{v,b\}&\subseteq I_2,
		&b&\in c_{I_2}(B\cup\{v\}).
	\end{aligned}
	\]
	If $Q:=A\cup B\cup\{v,b\}$ is a primal clique, then there exist
	an original chart $J\in\mathcal H$ and a set $D$ such that
	\[
	\begin{gathered}
		D\subseteq A\cup B,\qquad 1\le |D|\le d-1,\\
		D\cup\{b\}\subseteq J,\qquad b\in c_J(D).
	\end{gathered}
	\]
	Thus the composite consequence has a bounded-premise witness
	in an original $d$-native chart.
	Empty preservation forces $A\ne\varnothing$; no nonemptiness assumption on $B$ is needed.
\end{theorem}

\begin{proof}
	We split the proof into four steps:
	\begin{enumerate}
		\item \emph{Minimal counterexample.}
		Suppose the conclusion fails for a compatible atlas.
		In the original source charts, shrink $A$ and $B\cup\{v\}$
		to inclusion-minimal premises $A_0\subseteq A$ and
		$P\subseteq B\cup\{v\}$ for $v$ and $b$, respectively.
		If $v\notin P$, then $P\subseteq B$, and
		Lemma~\ref{LEM:PackagingNativeMinimalPremiseIdentity} gives
		$1\le |P|\le d-1$.
		Thus $P\to b$ is already a replacement in the original chart $I_2$, a contradiction.
		Therefore $P=B_0\cup\{v\}$ for some $B_0\subseteq B$.
		Any replacement for the reduced parents would also replace
		the original parents, so the reduced pair remains a counterexample.
		
		Among all counterexamples in this fixed original atlas whose two parent premises are inclusion-minimal and whose	second premise contains the intermediate label, choose one with minimum combined-support cardinality.
		Here the intermediate and conclusion labels, source charts,	and premises are all allowed to vary.		
		Denote the chosen data again by	$v,b,I_1,I_2,A,B$, and put
		\[
		Q=A\cup B\cup\{v,b\}.
		\]
		The support $Q$ is a primal clique by the definition of the counterexample.		
		Since $A,B\subseteq V\setminus\{v,b\}$,
		\[
		A\cup B=Q\setminus\{v,b\}.
		\]
		
		\item \emph{The support $Q$ is unguarded.}
		If an original chart $J$ guarded $Q$, compatibility would transfer
		both parent implications to $c_J$:
		\[
		v\in c_J(A),\qquad b\in c_J(B\cup\{v\}).
		\]
		By monotonicity and idempotence, $b\in c_J(A\cup B)$.
		Theorem~\ref{THM:PackagingNativeCommutantDimensionBound}
		then yields $D\subseteq A\cup B$ with
		$|D|\le d-1$ and $b\in c_J(D)$.
		Empty preservation implies $D\ne\varnothing$.
		This is a forbidden original-chart replacement.
		Hence $Q$ is unguarded, so $|Q|\ge d+1$.
		
		\item \emph{Conservatively fill the proper faces.}
		Starting from the original atlas, while some proper subset of $Q$
		is unguarded, choose an inclusion-minimal unguarded proper clique
		$M\subsetneq Q$ and adjoin its least weak filler.
		Since $M$ is a primal clique, every subset of $M$ is also a primal clique.
		By the minimal choice of $M$ among unguarded proper cliques of $Q$, every proper subset of $M$ is already guarded.		
		Consequently, $M$ is inclusion-minimal among all unguarded cliques of the current cover.
		Lemma~\ref{LEM:PackagingNativeMinimalCounterexampleFillability}
		therefore applies and shows that $M$ is weakly fillable.
		By Lemma~\ref{LEM:PackagingConservativeWeakFilling} and
		Corollary~\ref{COR:PackagingConservativeWeakFillingSequence},
		each filling preserves compatibility, the primal graph, empty preservation, and the least generated global closure.
		Each step guards a previously unguarded proper subset of the finite set $Q$, and guardedness is never lost.
		Thus the process terminates	with every proper subset of $Q$ guarded.
		
		For each $z\in Q$, let $\lambda_z$ be the restriction to
		$F_z:=Q\setminus\{z\}$ of a current chart guarding $F_z$.
		Compatibility makes $\lambda_z$ independent of that chart's choice.
		
		\item \emph{Exclude a native critical circuit.}
		We claim that
		$
		b\notin\lambda_v(A\cup B).
		$
		Otherwise, since $b\notin A\cup B$,
		Lemma~\ref{LEM:PackagingNativeRestrictedEnlargedRecovery} would provide an original chart $J$ and a nonempty
		$D\subseteq A\cup B$ with
		\[
		1\le |D|\le d-1,\qquad
		D\cup\{b\}\subseteq J,\qquad b\in c_J(D),
		\]
		contradicting the choice of the counterexample.
		
		Now apply Theorem~\ref{THM:PackagingNoNativeCriticalCircuit}
		with $x=v$ and $y=b$.
		Its hypotheses hold:
		$Q$ is an original	primal clique;
		the enlargement is a conservative minimal weak-filling sequence on proper cliques of $Q$;
		its coatom closures are the $\lambda_z$;
		both parent implications are inclusion-minimal and valid in original source charts;
		$A\cup B=Q\setminus\{v,b\}$;
		the missing-$v$ face does not derive $b$ from $A\cup B$;
		and all original cliques of size at most $d$ are guarded.
		The cited theorem excludes exactly this configuration, a contradiction.
	\end{enumerate}
\end{proof}

\begin{corollary}[Original-chart resolution for bounded-Armstrong FD theories]
	\label{COR:PackagingArmstrongOriginalChartResolution}	
	Let \(d\ge3\), let \(\mathcal H\) be a finite cover, and suppose every primal clique of cardinality at most \(d\) is guarded.	
	Let
	\[
	(\Sigma_I)_{I\in\mathcal H}
	\]
	be a compatible family of functional-dependency theories whose induced	closures are empty-preserving and satisfy
	\[
	\operatorname{arm}(\Sigma_I)\le d
	\qquad
	(I\in\mathcal H).
	\]	
	Suppose that, in original charts,
	\[
	\Sigma_{I_1}\models A\to v,
	\qquad
	\Sigma_{I_2}\models B\cup\{v\}\to b,
	\qquad
	v\ne b,
	\]
	where
	\[
	A,B\subseteq V\setminus\{v,b\},
	\]
	and suppose
	\[
	Q=A\cup B\cup\{v,b\}
	\]
	is a primal clique.	
	Then there exist an original chart
	\[
	J\in\mathcal H
	\]
	and a set
	\[
	D\subseteq A\cup B,
	\qquad
	1\le |D|\le d-1,
	\]
	such that
	\[
	D\cup\{b\}\subseteq J
	\]
	and
	\[
	{
		\Sigma_J\models D\to b.
	}
	\]
\end{corollary}

\begin{proof}
	Translate each local FD theory to its closure $c_I$.
	On every nonempty chart, the equality $d_{\mathrm{comm}}(c_I)=\operatorname{arm}(c_I)$ from Theorem~\ref{THM:PackagingCommutantDimensionEqualsArmstrongSize}
	places $c_I$ in $\mathbf{CommCl}_d(I)$.
	On an empty chart the unique closure belongs to $\mathbf{CommCl}_d(\varnothing)$ by definition.
	Compatibility of the complete FD theories is compatibility of their closure traces.
	The conclusion now follows from
	Theorem~\ref{THM:PackagingNativeGeneralResolution}.
\end{proof}

\subsection{Fixed-dimensional native globalization}
\label{SUBSEC:PackagingSharpNativeGlobalization}

The sharp native recovery rank gives globalization sufficiency by chordal elimination.
Simplex-boundary and cycle obstructions establish necessity.
Together they yield the exact \(d\)-conformality criterion for \(d\ge3\) in Theorem~\ref{THM:PackagingNativeGeneralClassification}.

\subsubsection{Exact native recovery rank and globalization sufficiency}

The native propagation--recovery rank equals \(d\), giving the exact clique-guarding scale for recovering implications in original charts.
Combined with chordal elimination, this yields an exact global closure; a separate corollary gives clique completion under bounded conformality without requiring chordality.

\begin{theorem}[Exact native propagation--recovery rank]
	\label{THM:PackagingNativePropagationRecoveryRank}
	Let \(d\ge2\), and put
	\[
	\mathscr N_d(I)=\mathbf{CommCl}_d(I).
	\]
	For \(c\in\mathscr N_d(I)\), define
	\[
	\mathscr R_d^{\mathrm{comm}}(c)
	=
	\left\{
	A\to b:
	1\le|A|\le d-1,
	\ b\in c(A)\setminus A
	\right\}.
	\]
	Then \(\mathscr R_d^{\mathrm{comm}}\) is a chart-local rule scheme and
	\[
	\boxed{
		\operatorname{prk}
		\bigl(\mathscr N_d,\mathscr R_d^{\mathrm{comm}}\bigr)
		=d.
	}
	\]
	Thus the native recovery scale is exactly \(d\), whereas forgetting the native structure and retaining only	\(\operatorname{car}\le d-1\) gives the generic scale \(2d-1\).
\end{theorem}

\begin{proof}
	Restriction stability follows from Lemma~\ref{LEM:PackagingNativeRestrictionCompression}.
	For every nonempty native closure in \(\mathbf{CommCl}_d\), Theorem~\ref{THM:PackagingNativeCommutantDimensionBound} gives
	Carath\'eodory number at most \(d-1\).
	Hence \(\mathscr R_d^{\mathrm{comm}}(c)\) generates \(c\).
	Trace locality follows from
	\[
	(c|_J)(A)=c(A)\cap J
	\]
	whenever \(A\cup\{b\}\subseteq J\).
	The empty label set is trivial.
	
	Suppose first that \(d=2\).
	By Proposition~\ref{THM:PackagingNativeDimensionTwoClassification}, every nonempty closure in \(\mathbf{CommCl}_2(I)\) satisfies
	\[
	c(S)=I\qquad(S\ne\varnothing).
	\]
	The rule scheme therefore consists of its valid unary implications.
	In a local propagation step the parent rules necessarily have the form
	\[
	a\to v,
	\qquad
	v\to b,
	\]
	with \(a,v,b\) distinct.
	If their combined support is a primal clique, then \(\{a,b\}\) is an edge and hence is guarded by an original chart.
	On that nonempty chart one has \(a\to b\), giving the required original-chart recovery with \(D=\{a\}\).
	Thus the class is locally \(2\)-propagation-recoverable.
	Since the definition of propagation--recovery rank takes the minimum over \(m\ge2\), its rank is exactly \(2\).
	
	Assume now that \(d\ge3\).
	For the upper bound, take a cover on which every clique of size at most \(d\) is guarded, with no chordality assumption, a compatible \(\mathbf{CommCl}_d\)-atlas, and two chart-local parent rules whose combined support is a primal clique.
	Theorem~\ref{THM:PackagingNativeGeneralResolution} supplies an original-chart consequence
	\[
	D\to b,
	\qquad
	D\subseteq A\cup B,
	\qquad
	1\le|D|\le d-1.
	\]
	Thus
	\[
	D\to b\in\mathscr R_d^{\mathrm{comm}}(c_K)
	\]
	for an original chart \(K\), proving
	\[
	\operatorname{prk}
	\bigl(\mathscr N_d,\mathscr R_d^{\mathrm{comm}}\bigr)
	\le d.
	\]
	
	For sharpness, let \(V\) have \(d\) labels and take the simplex-boundary cover
	\[
	\mathcal H=\{V\setminus\{q\}:q\in V\}.
	\]
	Its primal graph is complete and every clique of cardinality at most \(d-1\) is guarded.
	Choose distinct
	\[
	x,y,z\in V,
	\qquad
	D=V\setminus\{x,y\}.
	\]
	On \(V\setminus\{x\}\) take the identity closure.
	On \(V\setminus\{y\}\) and \(V\setminus\{z\}\) take the corresponding uniform circuit closures.
	On all remaining facets take the identity closure.
	As in the necessity construction of Theorem~\ref{THM:PackagingNativeSimplexBoundaryGeneral}, this is a compatible \(\mathbf{CommCl}_d\)-atlas.
	
	The two implications
	\[
	D\to x,
	\qquad
	(D\setminus\{z\})\cup\{x\}\to y
	\]
	are valid in the two circuit facets.
	Their premises have cardinality	\(d-2\), so both belong to \(\mathscr R_d^{\mathrm{comm}}\), and their combined support is the clique \(V\).
	Suppose there were an original-chart recovery
	\[
	E\to y,
	\qquad
	E\subseteq D,
	\]
	belonging to \(\mathscr R_d^{\mathrm{comm}}(c_J)\).
	The identity facet \(V\setminus\{x\}\) contains \(E\cup\{y\}\).
	Compatibility on the common support would transfer \(E\to y\) to that identity closure, which is impossible.
	Hence the native class is not locally \((d-1)\)-propagation-recoverable.
	Therefore the rank is exactly \(d\).
\end{proof}

\begin{corollary}[Canonical clique completion for bounded native dimension]
	\label{COR:PackagingNativeCanonicalCliqueCompletion}
	Let $d\ge3$, and let $\mathcal H$ be any finite $d$-conformal cover, with no chordality assumption.
	Every compatible atlas of local closures in $\mathbf{CommCl}_d$ admits a canonical least exact filler on each primal clique.
	These fillers agree on all inclusions and pairwise intersections, can be added simultaneously on the maximal cliques without changing the least generated global closure, and are exact on every induced original-chart trace.
	No bound by $d$ on the realizations of the added fillers is asserted.
\end{corollary}

\begin{proof}
	Theorem~\ref{THM:PackagingNativePropagationRecoveryRank} gives local $d$-propagation-recoverability for $\mathbf{CommCl}_d$.
	Apply Theorem~\ref{THM:PackagingCanonicalCliqueCompletion} with $m=d$.
\end{proof}

\begin{theorem}[Universal native gluing with bounded local dimension]
	\label{THM:PackagingNativeGeneralGluingSufficiency}
	Let \(d\ge3\), and let \(\mathcal H\) be a finite cover whose primal graph is chordal and such that every clique of cardinality at most \(d\) is guarded.
	Then
	\[
	{\mathsf{UG}_{\mathrm{comm}}(d,\mathcal H).}
	\]
	Equivalently, every compatible atlas
	\[
	(c_I)_{I\in\mathcal H},
	\qquad
	c_I\in\mathbf{CommCl}_d(I),
	\]
	admits an exact global closure.
	By Proposition~\ref{PROP:PackagingNativeCommutantDependence}, that finite global closure in turn admits some finite-dimensional native commutant realization.
\end{theorem}

\begin{proof}
	By Theorem~\ref{THM:PackagingNativePropagationRecoveryRank}, the native rule scheme
	\[
	\bigl(\mathscr N_d,\mathscr R_d^{\mathrm{comm}}\bigr)
	\]
	is locally \(d\)-propagation-recoverable.
	The hypotheses of the present theorem are therefore exactly the geometric hypotheses required by
	Theorem~\ref{THM:PackagingControlledPropagationGlobalization}.
	It follows that the closure \(C_*\) generated by the union of the local native rule families is the least exact global closure.
	
	Finally, \(V\) is finite and \(C_*\) is empty-preserving.
	Hence Proposition~\ref{PROP:PackagingNativeCommutantDependence} realizes \(C_*\) by	a finite commuting self-adjoint tuple.
	The controlled-propagation theorem only proves weak exact globalization:
	the dimension of this global native realization is not asserted to remain equal to \(d\).
\end{proof}

\subsubsection{Neutral lifts and geometric necessity}

The positive theorem proves that a geometric condition is sufficient.
To show it is necessary, we build compatible atlases that cannot extend when the condition fails.
Neutral lifts add labels without increasing local native dimension.
They let us place small boundary obstructions inside larger covers.

\begin{lemma}[Neutral one-point lifts of native coatoms]
	\label{LEM:PackagingNativeNeutralCoatomLift}
	Let
	\[
	r\ge3,
	\qquad
	|F|=r-1,
	\]
	and let \(q\notin F\).
	For each of the two closures
	\[
	\gamma=\operatorname{id}_F
	\qquad\text{or}\qquad
	\gamma=u_F,
	\]
	there exists an empty-preserving closure
	\[
	\widehat\gamma\in\mathbf{CommCl}_r(F\cup\{q\})
	\]
	such that
	\[
	{\widehat\gamma|_F=\gamma}
	\]
	and, for every proper subset \(J\subsetneq F\),
	$
		\widehat\gamma|_{J\cup\{q\}}
		=
		\operatorname{id}_{J\cup\{q\}}.
	$
\end{lemma}

\begin{proof}
	First suppose \(\gamma=\operatorname{id}_F\).
	Take
	\[
	\widehat\gamma=u_{F\cup\{q\}}.
	\]
	Since \(|F\cup\{q\}|=r\),
	Lemma~\ref{LEM:PackagingNativeUniformCircuit} gives \(\widehat\gamma\in\mathbf{CommCl}_r\), and every proper restriction of a uniform circuit closure is the identity closure.
	Hence both required restriction identities hold.
	
	Now suppose \(\gamma=u_F\).
	Choose \(f_0\in F\), and use the \(r\)-element joint spectral set
	\[
	\Omega=F\sqcup\{*\}.
	\]
	For \(a\in F\), define
	\[
	v_a(f)=\delta_{af}\quad(f\in F),
	\qquad
	v_a(*)=\delta_{a f_0}.
	\]
	Thus the row \(*\), restricted to the coordinates in \(F\), duplicates the row \(f_0\).
	Define the auxiliary coordinate by
	\[
	v_q(f)=0\quad(f\in F),
	\qquad
	v_q(*)=1.
	\]
	Because \(r\ge3\), every coordinate \(v_a\) with \(a\in F\) is nonconstant, and \(v_q\) is nonconstant as well.
	Hence the native closure of these \(r\) rows is empty-preserving.
	Denote it by \(\widehat\gamma\).
	Forgetting \(q\) leaves the standard-basis realization of \(u_F\), together with one duplicate row, so
	\[
	\widehat\gamma|_F=u_F.
	\]
	
	Let \(J\subsetneq F\).
	The coordinate \(q\) is not determined by \(J\), because the rows \(f_0\) and \(*\) agree on every coordinate of \(F\), hence on \(J\), but have different \(q\)-values.
	Next fix \(a\in J\) and choose \(h\in F\setminus J\).
	The rows \(a\) and \(h\) both have \(q\)-value zero, agree on every coordinate in \(J\setminus\{a\}\), and differ at \(a\).
	Thus
	\[
	a\notin
	\widehat\gamma\bigl((J\setminus\{a\})\cup\{q\}\bigr).
	\]
	Hence no label of \(J\cup\{q\}\) is determined by all the others.
	By monotonicity no nontrivial dependency occurs there, so
	\[
	\widehat\gamma|_{J\cup\{q\}}
	=
	\operatorname{id}_{J\cup\{q\}}.
	\]
\end{proof}

\begin{lemma}[Cloning an auxiliary native coordinate]
	\label{LEM:PackagingNativeCloneLift}
	Let
	\[
	c\in\mathbf{CommCl}_d(T\cup\{q\}),
	\]
	and let \(O\) be a finite nonempty set disjoint from \(T\cup\{q\}\).
	There is a closure
	\[
	\operatorname{Clone}_{q\to O}(c)
	\in
	\mathbf{CommCl}_d(T\cup O)
	\]
	obtained by replacing every \(o\in O\) by a copy of the coordinate \(q\).
	Explicitly, for \(S\subseteq T\cup O\), put
	\[
	\phi(S)
	=
	(S\cap T)
	\cup
	\begin{cases}
		\{q\},&S\cap O\ne\varnothing,\\
		\varnothing,&S\cap O=\varnothing.
	\end{cases}
	\]
	Then
	\[
	\operatorname{Clone}_{q\to O}(c)(S)
	=
	\bigl(c(\phi(S))\cap T\bigr)
	\cup
	\begin{cases}
		O,&q\in c(\phi(S)),\\
		\varnothing,&q\notin c(\phi(S)).
	\end{cases}
	\]
	Moreover, cloning commutes with restriction:
	if \(T'\subseteq T\) and \(\varnothing\ne O'\subseteq O\), then
	\[
	{
		\operatorname{Clone}_{q\to O}(c)|_{T'\cup O'}
		=
		\operatorname{Clone}_{q\to O'}
		\bigl(c|_{T'\cup\{q\}}\bigr).
	}
	\]
	If \(O'=\varnothing\), the restriction is simply \(c|_{T'}\).
\end{lemma}

\begin{proof}
	Choose a commuting self-adjoint realization
	\[
	(D_t)_{t\in T},\quad D_q
	\]
	of \(c\) in \(M_d(\mathbb C)\), and put \(D_o=D_q\) for every \(o\in O\).
	The displayed formula and the restriction identities follow immediately from
	Proposition~\ref{PROP:PackagingNativeCommutantDependence}.
\end{proof}

\begin{theorem}[Native lifting from an unguarded clique]
	\label{THM:PackagingNativeUnguardedCliqueLifting}
	Let \(d\ge3\), and let \(\mathcal H\) be a finite cover of a finite label set.
	Suppose the primal graph of \(\mathcal H\) contains an unguarded clique \(K\) with
	\[
	3\le|K|\le d.
	\]
	Then there exists a compatible atlas
	\[
	(c_I)_{I\in\mathcal H},
	\qquad
	c_I\in\mathbf{CommCl}_d(I),
	\]
	having no exact global extension.  Consequently,
	\[
	{\neg\mathsf{UG}_{\mathrm{comm}}(d,\mathcal H).}
	\]
\end{theorem}

\begin{proof}
	Choose an inclusion-minimal unguarded clique \(M\subseteq K\), and put
	$
	r=|M|.
	$
	Every two-element subset of \(M\) is guarded by the definition of the primal graph, so \(r\ge3\); also \(r\le d\).
	Minimal unguardedness implies that every proper subset of \(M\) is guarded.
	
	Choose distinct \(x,y,z\in M\), and for \(m\in M\) put
	\[
	F_m=M\setminus\{m\}.
	\]
	Define a simplex-boundary atlas by
	\[
	\gamma_x=\operatorname{id}_{F_x},
	\qquad
	\gamma_y=u_{F_y},
	\qquad
	\gamma_z=u_{F_z},
	\]
	and
	\[
	\gamma_m=\operatorname{id}_{F_m}
	\qquad
	(m\notin\{x,y,z\}).
	\]
	Every proper restriction of a uniform circuit closure is identity, so the family \((\gamma_m)_{m\in M}\) is compatible.
	
	It is not exactly extendible.
	Put
	\[
	D=M\setminus\{x,y\}.
	\]
	Since \(|D|=r-2\), the circuit closure on \(F_y\) gives
	\[
	x\in\gamma_y(D).
	\]
	Also
	\[
	(D\setminus\{z\})\cup\{x\}
	=
	M\setminus\{y,z\}
	\]
	has cardinality \(r-2\), so the circuit closure on \(F_z\) gives
	\[
	y\in
	\gamma_z\bigl((D\setminus\{z\})\cup\{x\}\bigr).
	\]
	Any exact global closure on \(M\) would therefore derive \(D\to x\to y\).
	But \(\gamma_x=\operatorname{id}_{F_x}\) and \(D\subseteq F_x\), so
	\[
	y\notin\gamma_x(D),
	\]
	a contradiction.
	
	We lift this boundary obstruction to the original cover.
	Introduce an auxiliary symbol \(q\notin M\).
	For each coatom \(F_m\), choose a neutral lift
	\[
	\widehat\gamma_m
	\in
	\mathbf{CommCl}_r(F_m\cup\{q\})
	\]
	from Lemma~\ref{LEM:PackagingNativeNeutralCoatomLift}.
	Thus
	\[
	\widehat\gamma_m|_{F_m}=\gamma_m
	\]
	and
	\[
	\widehat\gamma_m|_{J\cup\{q\}}
	=
	\operatorname{id}_{J\cup\{q\}}
	\qquad
	(J\subsetneq F_m).
	\]	
	For every proper subset \(T\subsetneq M\), define a closure \(\eta_T\) on \(T\cup\{q\}\) by
	\[
	\eta_T
	=
	\begin{cases}
		\operatorname{id}_{T\cup\{q\}},&|T|\le r-2,\\
		\widehat\gamma_m,&T=F_m.
	\end{cases}
	\]
	The neutral-lift lemma gives the restriction rule
	\[
	{
		\eta_T|_{T'\cup\{q\}}
		=
		\eta_{T'}
		\qquad
		(T'\subseteq T\subsetneq M).
	}
	\tag{\(*\)}
	\]
	
	For each original chart \(I\in\mathcal H\), put
	\[
	T_I=I\cap M,
	\qquad
	O_I=I\setminus M.
	\]
	Since \(M\) is unguarded, \(T_I\subsetneq M\).
	If \(O_I=\varnothing\), set
	\[
	c_I=\eta_{T_I}|_{T_I}.
	\]
	If \(O_I\ne\varnothing\), set
	\[
	c_I
	=
	\operatorname{Clone}_{q\to O_I}(\eta_{T_I}).
	\]
	
	These closures have native dimension at most \(r\).
	Indeed, if \(|T_I|\le r-2\), then \(\eta_{T_I}\) is the identity closure on \(|T_I|+1\) labels and is realizable in dimension \(|T_I|+2\le r\) by Theorem~\ref{THM:PackagingNativeIdentityDimension}.
	If \(T_I\) is a	coatom, its neutral lift belongs to \(\mathbf{CommCl}_r\).
	Cloning does not increase dimension.
	Dimension padding therefore gives
	\[
	c_I\in\mathbf{CommCl}_d(I).
	\]
	
	We next prove compatibility.
	Let \(I,J\in\mathcal H\), put \(R=I\cap J\), and write
	\[
	T=R\cap M,
	\qquad
	O=R\setminus M.
	\]
	If \(O=\varnothing\), the restriction of either local closure to \(R=T\) is \(\eta_T|_T\), by \((*)\).
	If \(O\ne\varnothing\),	Lemma~\ref{LEM:PackagingNativeCloneLift} and \((*)\) show that the restriction of either local closure to \(R\) is
	\[
	\operatorname{Clone}_{q\to O}(\eta_T).
	\]
	Hence \(c_I|_R=c_J|_R\).
	
	Finally, every coatom \(F_m\) is guarded, so choose \(I_m\in\mathcal H\) with \(F_m\subseteq I_m\).
	Because no chart contains all of \(M\), one has	\(I_m\cap M=F_m\), and therefore
	\[
	c_{I_m}|_{F_m}=\gamma_m.
	\]
	If an exact global closure \(C\) existed on the full label set, then its restriction
	\[
	(C|_M)(S)=C(S)\cap M\qquad(S\subseteq M)
	\]
	would be a closure on \(M\).
	For each \(m\), exactness of \(C\) on the original chart \(I_m\), followed by restriction to \(F_m\), gives
	\[
	(C|_M)|_{F_m}=c_{I_m}|_{F_m}=\gamma_m.
	\]
	Thus \(C|_M\) would exactly extend every \(\gamma_m\), contradicting the nonextendibility of the boundary atlas above.
	Hence the lifted atlas has no exact global extension.
\end{proof}

\begin{theorem}[Necessary geometry for native gluing with bounded local dimension]
	\label{THM:PackagingNativeGeneralNecessity}
	Let \(d\ge3\), and let \(\mathcal H\) be a finite cover of a finite label set.
	If
	\[
	\mathsf{UG}_{\mathrm{comm}}(d,\mathcal H),
	\]
	then the primal graph of \(\mathcal H\) is chordal and every clique of cardinality at most \(d\) is guarded.
\end{theorem}

\begin{proof}
	By dimension padding,
	\[
	\mathbf{CommCl}_3(I)
	\subseteq
	\mathbf{CommCl}_d(I).
	\]
	Thus \(\mathsf{UG}_{\mathrm{comm}}(d,\mathcal H)\) implies \(\mathsf{UG}_{\mathrm{comm}}(3,\mathcal H)\), and
	Theorem~\ref{THM:PackagingNativeDimensionThreeClassification} forces the primal graph to be chordal.
	
	Suppose now that there were an unguarded clique \(K\) with \(|K|\le d\).
	Singletons are guarded because \(\mathcal H\) covers the label set, and two-element cliques are guarded by the definition of the primal graph.
	Hence \(K\) contains an inclusion-minimal unguarded subclique \(M\) with
	$
	3\le|M|\le d.
	$
	Theorem~\ref{THM:PackagingNativeUnguardedCliqueLifting} then constructs a compatible \(\mathbf{CommCl}_d\)-atlas with no exact extension, a contradiction.
	Thus every clique of size at most \(d\) is guarded.
\end{proof}

\subsubsection{Classification and Armstrong-size consequences}

We now have both directions of the native globalization classification.
The main theorem states exactly which covers work: their primal graphs must be chordal and their small cliques guarded through size \(d\).
We translate the result into minimum Armstrong size and give explicit geometric obstruction families.
Only the local realization dimension is bounded; a global realization may need more dimensions.

\begin{theorem}[Sharp bounded-conformality classification for native globalization]
	\label{THM:PackagingNativeGeneralClassification}
	Let \(d\ge3\), and let \(\mathcal H\) be a finite cover of a finite label set.
	The following are equivalent:
	\begin{enumerate}[label=(\roman*),leftmargin=2.5em]
		\item
		$
		\mathsf{UG}_{\mathrm{comm}}(d,\mathcal H);
		$
		
		\item the primal graph \(G_{\mathcal H}\) is chordal and \(\mathcal H\) is \(d\)-conformal;
		
		\item
		$
		d<\nu(\mathcal H).
		$
	\end{enumerate}
	Equivalently,
	\[
		\mathfrak G_d
		=
		\left\{
		\mathcal H:
		G_{\mathcal H}\text{ is chordal and }
		\mathcal H\text{ is \(d\)-conformal}
		\right\},
		\qquad d\ge3.
	\]
	Thus, for this native class, full conformality is not required: the exact universal threshold is bounded conformality at the native scale \(d\).
\end{theorem}

\begin{proof}
	The implication \(\textup{(ii)}\Longrightarrow\textup{(i)}\) is Theorem~\ref{THM:PackagingNativeGeneralGluingSufficiency}, and
	\(\textup{(i)}\Longrightarrow\textup{(ii)}\) is	Theorem~\ref{THM:PackagingNativeGeneralNecessity}.
	The equivalence of \textup{(ii)} and \textup{(iii)} is Lemma~\ref{LEM:PackagingObstructionRankDictionary}.
\end{proof}

\begin{remark}[Relation to classical acyclicity and local consistency]
	\label{REM:PackagingNativeRecoveryLiteratureScope}	
	The geometric ingredients in Theorem~\ref{THM:PackagingNativeGeneralClassification}
	are classical.
	Chordality together with full conformality characterizes the usual acyclicity condition for finite database hypergraphs
	\cite{BeeriFaginMaierYannakakis1983Acyclic}.
	Bounded conformality is also standard hypergraph terminology \cite{Otto2012HighlyAcyclic}.
	
	Functional-dependency acyclicity and join consistency have been studied for database schemes equipped with fixed dependency theories
	\cite{SagivShmueli1989FDAcyclicity}.
	Related local-to-global methods arise in relational and constraint-network consistency \cite{DechterVanBeek1997RelationalConsistency}.
	Locally consistent $K$-relations provide another dependency semantics in which ordinary FD transitivity need not hold
	\cite{BarlagEtAl2026LocallyConsistentKRelations}.
	
	The quantifiers and extension requirement of the present theorem are different.
	The cover is fixed, while the complete local dependency closures vary over $\mathbf{CommCl}_d$.
	Compatibility means equality of their restrictions on overlaps, and exact globalization requires a single global closure with precisely the prescribed local restrictions.
	
	Thus the theorem is not a reformulation of classical hypergraph acyclicity.
	For $d\ge3$, full conformality is replaced by the sharp $d$-conformality condition, while full chordality remains necessary.
	The accompanying recovery theorem identifies $d$ as the exact original-chart propagation--recovery scale.
\end{remark}

\begin{corollary}[Unified obstruction-rank form]
	\label{COR:PackagingNativeUnifiedObstructionRank}
	For every integer \(d\ge2\) and every finite cover \(\mathcal H\),
	\[
	\boxed{
		\mathsf{UG}_{\mathrm{comm}}(d,\mathcal H)
		\quad\Longleftrightarrow\quad
		d<\nu(\mathcal H).
	}
	\]
	For \(d=2\) both sides hold for every finite cover.
	For \(d\ge3\) the criterion is exactly chordality together with \(d\)-conformality.
\end{corollary}

\begin{proof}
	If \(d=2\), Proposition~\ref{THM:PackagingNativeDimensionTwoClassification} gives universal exact globalization on every finite cover, while
	\(\nu(\mathcal H)\ge3\) gives \(2<\nu(\mathcal H)\).
	For \(d\ge3\), apply Theorem~\ref{THM:PackagingNativeGeneralClassification}.
\end{proof}

\begin{corollary}[Armstrong-size form of the native classification]
	\label{COR:PackagingNativeArmstrongSizeClassification}
	Let \(d\ge3\), and let \(\mathcal H\) be a finite cover of a finite label set.
	The following are equivalent:
	\begin{enumerate}[label=(\roman*),leftmargin=2.5em]
		\item every compatible atlas of empty-preserving functional-dependency closure operators \((c_I)_{I\in\mathcal H}\) satisfying
		\[
		\operatorname{arm}(c_I)\le d
		\qquad(I\in\mathcal H)
		\]
		admits an exact global closure;
		
		\item the primal graph of \(\mathcal H\) is chordal and every primal clique of cardinality at most \(d\) is guarded.
	\end{enumerate}
	Equivalently, the native theorem classifies covers using a uniform bound on the minimum number of rows in an Armstrong relation for each local dependency theory.
\end{corollary}

\begin{proof}
	If \(I=\varnothing\), then \(c_I\) is the unique closure on the empty set and belongs to \(\mathbf{CommCl}_d(\varnothing)\) by definition.
	For \(I\ne\varnothing\), Theorem~\ref{THM:PackagingCommutantDimensionEqualsArmstrongSize}, together with dimension padding, gives
	\[
	\operatorname{arm}(c_I)\le d\quad\Longleftrightarrow\quad c_I\in\mathbf{CommCl}_d(I).
	\]
	Therefore the asserted universal property is exactly \(\mathsf{UG}_{\mathrm{comm}}(d,\mathcal H)\).
	Apply Theorem~\ref{THM:PackagingNativeGeneralClassification}.
\end{proof}

\begin{lemma}[Constant-coordinate reduction for functional dependencies]
	\label{LEM:PackagingFDConstantReduction}	
	Let \(I\) be finite and let
	\[
	c:2^I\longrightarrow 2^I
	\]
	be a closure operator.
	Put
	\[
	K_c:=c(\varnothing),
	\qquad
	I^\circ:=I\setminus K_c,
	\]
	and define
	\[
	c^\circ(X):=c(X)\setminus K_c,
	\qquad
	X\subseteq I^\circ.
	\]	
	Then the following hold.	
	\begin{enumerate}[label=(\roman*),leftmargin=2.5em]		
		\item \(c^\circ\) is an empty-preserving closure operator on \(I^\circ\), and for every \(X\subseteq I\),
		\[
		{
			c(X)
			=
			K_c
			\cup
			c^\circ(X\setminus K_c).
		}
		\]
		
		\item With the zero-row Armstrong convention, one has
		\[
		\operatorname{arm}(c^\circ)=\operatorname{arm}(c),
		\]
		including the case $K_c=I$.
		
		\item Let
		\[
		(c_I)_{I\in\mathcal H}
		\]
		be a compatible atlas of finite closure operators on a finite cover \(\mathcal H\) of \(V\).
		Then constant status is globally well defined: there is a unique set
		\[
		K\subseteq V
		\]
		such that
		\[
		{
			K\cap I=c_I(\varnothing)
			\qquad
			(I\in\mathcal H).
		}
		\]		
		Writing
		\[
		V^\circ=V\setminus K,
		\qquad
		I^\circ=I\setminus K,
		\qquad
		\mathcal H^\circ
		=
		\{I^\circ:I\in\mathcal H\},
		\]
		the reduced closures
		\[
		c_I^\circ(X)
		=
		c_I(X)\setminus K
		\qquad
		(X\subseteq I^\circ)
		\]
		form a compatible empty-preserving atlas on \(\mathcal H^\circ\).
		
		If distinct original charts have the same reduced chart \(I^\circ\), their reduced closures coincide by compatibility, so duplicate reduced charts may be identified without ambiguity.
		
		\item The primal graph of \(\mathcal H^\circ\) is the induced subgraph
		\[
		G_{\mathcal H}[V^\circ].
		\]
		Consequently, chordality and guardedness of all cliques of cardinality at most \(d\) pass from \(\mathcal H\) to \(\mathcal H^\circ\).
		
		\item If \(C^\circ\) is an exact global extension of the reduced atlas, then
		\[
		{
			C(X)
			=
			K\cup C^\circ(X\setminus K)
		}
		\qquad
		(X\subseteq V)
		\]
		is an exact global extension of the original atlas.		
	\end{enumerate}
\end{lemma}

\begin{proof}
	Write $K=c(\varnothing)$.
	
	\emph{(i) Reduction of constants.}
	For every $X\subseteq I$, monotonicity gives $K\subseteq c(X)$.
	Hence
	\[
	c(X\cup K)=c(X).
	\]
	It follows directly that $c^\circ$ is extensive, monotone, idempotent and empty-preserving on $I^\circ$.
	Moreover,
	\[
	c(X)=K\cup c^\circ(X\setminus K).
	\]
	
	\emph{(ii) Armstrong-size preservation.}
	If $K=I$, both $c$ and the unique closure on $I^\circ=\varnothing$ have minimum Armstrong size zero by the zero-row convention.	
	Suppose $K\ne I$.
	In any Armstrong relation for $c$, the columns indexed by $K$ are constant, because $\varnothing\to a$ holds for every $a\in K$.
	Deleting these columns preserves the number of distinct rows and represents exactly $c^\circ$.
	Conversely, adjoining fixed constant columns to an Armstrong relation for $c^\circ$ preserves its row count and represents $c$ by part~\textup{(i)}.
	Taking minima proves
	\[
	\operatorname{arm}(c^\circ)=\operatorname{arm}(c).
	\]
	
	\emph{(iii) Compatibility after reduction.}
	If $a\in I\cap J$, compatibility gives
	\[
	a\in c_I(\varnothing)
	\quad\Longleftrightarrow\quad
	a\in c_J(\varnothing).
	\]
	Thus the constant status of each attribute is independent of its original chart, and the set $K$ is uniquely determined by $K\cap I=c_I(\varnothing)$.
	For $X\subseteq I^\circ\cap J^\circ$, compatibility of $c_I$ and $c_J$, followed by deletion of $K$, yields
	\[
	(c_I^\circ|_{I^\circ\cap J^\circ})(X)
	=(c_J^\circ|_{I^\circ\cap J^\circ})(X).
	\]
	Each reduced closure is empty-preserving by	part~\textup{(i)}.
	If two charts have the same	reduction, this identity also makes their reduced closures identical.
	
	\emph{(iv) Geometry.}
	Two labels in $V^\circ$ occur together in a reduced chart exactly when they occur together in the corresponding original chart.
	Therefore
	\[
	G_{\mathcal H^\circ}=G_{\mathcal H}[V^\circ].
	\]
	Chordality passes to induced subgraphs.
	Any clique guarded by an original chart remains guarded after intersecting that chart with $V^\circ$.
	
	\emph{(v) Reconstruction.}
	Define
	\[
	C(X)=K\cup C^\circ(X\setminus K).
	\]
	This is a closure on $V$.
	For $A\subseteq I$,	exactness of $C^\circ$ on the reduced chart gives
	\[
	\begin{aligned}
		C(A)\cap I
		&=(K\cap I)\cup
		\bigl(C^\circ(A\setminus K)\cap I^\circ\bigr)\\
		&=c_I(\varnothing)\cup c_I^\circ(A\setminus K)\\
		&=c_I(A),
	\end{aligned}
	\]
	where the last equality is part~\textup{(i)} applied to $c_I$.
	Thus $C$ is an exact extension of the original atlas.
\end{proof}

\begin{theorem}[Exact globalization theorem for bounded-Armstrong FD theories]
	\label{THM:PackagingClassicalBoundedArmstrongFDGlobalization}	
	Let
	$
	\mathcal H=\{I_1,\ldots,I_m\}
	$
	be a finite cover of a finite attribute set \(V\), and let \(d\ge3\).	
	For every \(I\in\mathcal H\), let \(\Sigma_I\) be a functional-dependency theory on \(I\), and write
	\[
	\operatorname{Th}_J(\Sigma_I)
	=
	\{X\to a:
	X\cup\{a\}\subseteq J,\;
	\Sigma_I\models X\to a\}.
	\]
	Functional dependencies with empty left-hand side are allowed throughout.
	Thus local theories may contain or imply constant attributes.
	All right-hand sides are written as single attributes, with no loss of generality by the usual decomposition of \(X\to Y\) into
	\(X\to y\), \(y\in Y\).	
	Say that the family
	\[
	(\Sigma_I)_{I\in\mathcal H}
	\]
	is compatible if
	\[
	\operatorname{Th}_{I\cap J}(\Sigma_I)
	=
	\operatorname{Th}_{I\cap J}(\Sigma_J)
	\qquad
	(I,J\in\mathcal H).
	\]	
	Then the following are equivalent.	
	\begin{enumerate}[label=(\roman*),leftmargin=2.5em]		
		\item For every compatible family
		\[
		(\Sigma_I)_{I\in\mathcal H}
		\]
		such that
		\[
		\operatorname{arm}(\Sigma_I)\le d
		\qquad
		(I\in\mathcal H),
		\]
		there exists a functional-dependency theory \(\Sigma\) on \(V\) satisfying
		\[
		{
			\operatorname{Th}_I(\Sigma)
			=
			\operatorname{Th}_I(\Sigma_I)
			\qquad
			(I\in\mathcal H).
		}
		\]
		No bound is imposed on the minimum Armstrong size of the global theory \(\Sigma\).
		The parameter \(d\) constrains the local theories only.
		The theorem concerns weak exact globalization rather than resource-preserving globalization at Armstrong budget \(d\).
		
		\item The primal graph \(G_{\mathcal H}\) is chordal and every primal clique	of cardinality at most \(d\) is contained in some member of \(\mathcal H\).
		
		\item $d < \nu(\mathcal H)$.
	\end{enumerate}	
	When these conditions fail, failure is witnessed by a finite compatible family of FD theories on the original cover whose nonextendibility is supported by a primitive geometric obstruction:
	either an induced primal cycle or an inclusion-minimal unguarded clique of cardinality at most \(d\).
	Equivalently, the counterexample is obtained by lifting an explicit primitive obstruction to the ambient cover.
\end{theorem}

\begin{proof}
	Assume first that \textup{(i)} holds.
	Then it holds in particular for compatible families whose induced FD closures are empty-preserving.
	Corollary~\ref{COR:PackagingNativeArmstrongSizeClassification} therefore implies that \(G_{\mathcal H}\) is chordal and that every primal clique of cardinality at most \(d\) is guarded.
	Thus \textup{(i)} implies \textup{(ii)}.
	
	Conversely, assume \textup{(ii)}, and let
	\[
	(\Sigma_I)_{I\in\mathcal H}
	\]
	be an arbitrary compatible family satisfying
	\[
	\operatorname{arm}(\Sigma_I)\le d.
	\]
	Let \(c_I=c_{\Sigma_I}\) be the corresponding FD closures.
	
	Apply Lemma~\ref{LEM:PackagingFDConstantReduction}.
	There is a globally determined constant set
	\[
	K \subseteq V
	\]
	and a compatible empty-preserving reduced atlas
	\[
	(c_I^\circ)_{I^\circ\in\mathcal H^\circ}
	\]
	on
	\[
	V^\circ=V\setminus K.
	\]
	Every nonempty reduced chart satisfies
	\[
	\operatorname{arm}(c_I^\circ)
	\le
	\operatorname{arm}(c_I)
	=
	\operatorname{arm}(\Sigma_I)
	\le d.
	\]
	If \(I^\circ=\varnothing\), then \(c_I^\circ\) is the unique closure on the empty set and belongs to \(\mathbf{CommCl}_d(\varnothing)\) by definition.
	
	Moreover,
	\[
	G_{\mathcal H^\circ}
	=
	G_{\mathcal H}[V^\circ],
	\]
	so \(G_{\mathcal H^\circ}\) is chordal and every primal clique of cardinality at most \(d\) is guarded.
	Corollary~\ref{COR:PackagingNativeArmstrongSizeClassification}
	therefore gives an exact global closure
	\[
	C^\circ
	\]
	for the reduced atlas.
	
	Lemma~\ref{LEM:PackagingFDConstantReduction} lifts it to
	\[
	C(X)
	=
	K\cup C^\circ(X\setminus K),
	\]
	which is an exact global closure for the original family.
	Since every finite closure operator is represented by a finite functional-dependency theory, choose a theory \(\Sigma\) whose closure is \(C\).
	Then
	\[
	\operatorname{Th}_I(\Sigma)
	=
	\operatorname{Th}_I(\Sigma_I)
	\qquad
	(I\in\mathcal H),
	\]
	proving \textup{(i)}.
	
	Finally, Lemma~\ref{LEM:PackagingObstructionRankDictionary} identifies \textup{(ii)} with
	\[
	d<\nu(\mathcal H),
	\]
	because \(d\ge3\).
	Thus \textup{(ii)} and \textup{(iii)} are equivalent.
	
	Finally, suppose the geometric condition fails.
	If the primal graph is nonchordal, the three-dimensional cycle obstruction constructed in
	Theorem~\ref{THM:PackagingNativeDimensionThreeClassification}, followed by dimension padding, gives a compatible nonextendible
	\(\mathbf{CommCl}_d\)-atlas.
	If instead there is an unguarded primal clique of cardinality at most \(d\), Theorem~\ref{THM:PackagingNativeUnguardedCliqueLifting}
	gives such an atlas directly.
	For each local closure \(c_I\), choose a finite	FD theory \(\Sigma_I\) with complete closure \(c_I\).
	Then \(\operatorname{arm}(\Sigma_I)=\operatorname{arm}(c_I)\le d\), and compatibility of the closures is exactly compatibility of the complete local FD theories.
	Any exact global FD theory would induce an exact global closure, contradicting the chosen native obstruction.
	Thus the asserted finite FD obstruction witnesses exist.
\end{proof}

The obstruction-rank notation used above was fixed at the beginning of Section~\ref{SEC:NativeRecoveryA}.
The later spectrum section develops its presentation-independent role further.

\begin{theorem}[Primitive native obstruction basis]
	\label{THM:PackagingNativePrimitiveObstructionBasis}
	Let \(d\ge3\) and let \(\mathcal H\) be a finite cover.
	Then
	\[
	\neg\mathsf{UG}_{\mathrm{comm}}(d,\mathcal H)
	\]
	if and only if at least one of the following geometric witnesses exists:
	\begin{enumerate}[label=(\roman*),leftmargin=2.5em]
		\item the primal graph contains an induced cycle of length at least four;
		
		\item there is an unguarded primal clique \(K\) with
		\[
		3\le|K|\le d.
		\]
	\end{enumerate}
	Moreover these witnesses carry native obstructions at their intrinsic geometric size:
	\begin{enumerate}[label=(\alph*),leftmargin=2.5em]
		\item an induced cycle supports a three-dimensional obstruction, which pads to all dimensions \(d\ge3\);
		
		\item if \(M\subseteq K\) is inclusion-minimal unguarded and \(|M|=r\), then the neutral-lift/cloning construction of
		Theorem~\ref{THM:PackagingNativeUnguardedCliqueLifting} gives an obstruction already in \(\mathbf{CommCl}_r\), hence in every \(\mathbf{CommCl}_d\) with \(d\ge r\).
	\end{enumerate}
	Thus every failure of universal fixed-dimensional native globalization has a primitive certificate of one of these two forms.
\end{theorem}

\begin{proof}
	If neither witness exists, the cover is chordal and $d$-conformal, so 	Theorem~\ref{THM:PackagingNativeGeneralClassification} gives universal exact globalization.
	Conversely, an induced cycle carries the explicit three-dimensional obstruction from Theorem~\ref{THM:PackagingNativeDimensionThreeClassification}, which pads to all $d\ge3$.
	An unguarded clique of size at most $d$ contains a minimal unguarded subclique $M$.
	Theorem~\ref{THM:PackagingNativeUnguardedCliqueLifting}	provides the obstruction already in dimension $|M|$, hence in dimension $d$.
\end{proof}

\begin{corollary}[Strict native globalization hierarchy]
	\label{COR:PackagingNativeStrictHierarchy}
	The weak fixed-dimensional native gluing classes satisfy
	\[
	{
		\mathfrak G_2
		\supsetneq
		\mathfrak G_3
		\supsetneq
		\mathfrak G_4
		\supsetneq
		\mathfrak G_5
		\supsetneq\cdots.
	}
	\]
	More precisely, for every \(d\ge3\),
	\[
	{\mathfrak G_{d+1}\subsetneq\mathfrak G_d.}
	\]
\end{corollary}

\begin{proof}
	Dimension padding gives the inclusions.
	For strictness, the simplex boundary on $d+1$ vertices is chordal and $d$-conformal but not $(d+1)$-conformal, so
	Theorem~\ref{THM:PackagingNativeGeneralClassification} places it in $\mathfrak G_d\setminus\mathfrak G_{d+1}$.
	The triangle boundary similarly	separates $\mathfrak G_2$ from $\mathfrak G_3$.
\end{proof}

\section{Cycle-free local-to-global amalgamation}
\label{SEC:PackagingCycleGeometry}

The preceding results concern implications between labels.
Here we test the same cover geometry on a different kind of data: transitive binary relations.
The main theorem, Theorem~\ref{THM:PackagingTransitiveUniversalAmalgamation}, characterizes when every compatible atlas of such relations extends exactly.
The same answer holds separately for preorders, partial orders, strict partial orders, and equivalence relations: cycle-freedom, which for finite covers is equivalent to chordality and \(3\)-conformality.
The proof uses directed paths and explicit cycle obstructions, rather than closure-recovery arguments.
We also give obstruction certificates and a quotient interpretation for equivalence relations.

\subsection{Kernel traces and minimal overlap obstructions}
\label{SUBSEC:PackagingExactProfileGlobalization}

For equivalence relations, a global extension must include every equality already required locally.
Taking the transitive closure of these equalities gives the smallest candidate global relation.
We test whether it conflicts with any pair that a chart keeps distinct.
The resulting theorems describe exact extension and reduce failures to small overlap obstructions.

\begin{theorem}[Kernel-level overlap extension and cycle obstruction]
	\label{THM:PackagingKernelOverlapExtension}
	Let \((I_K)_{K\in\mathcal K}\) be a set-indexed cover of a set \(E\), and let \(E_K\) be an equivalence relation on \(I_K\).
	Empty or repeated charts are permitted.
	Pairwise agreement is not assumed.
	Put
	\[
	E^{\mathrm{gen}}=\left(\bigcup_K E_K\right)^{\mathrm{eq}}
	\subseteq E^2.
	\]
	The following are equivalent:
	\begin{enumerate}[label=(\roman*),leftmargin=2em]
		\item A global equivalence relation restricts to exactly \(E_K\) on every \(I_K\).
		
		\item \(E^{\mathrm{gen}}\cap I_K^2=E_K\) for every \(K\).
		
		\item There is no finite chain \(x=x_0,\ldots,x_m=y\), each consecutive	pair belonging to some local \(E_{K_j}\), whose endpoints lie in one \(I_K\) but satisfy \((x,y)\notin E_K\).
	\end{enumerate}
	When these conditions hold, \(E^{\mathrm{gen}}\) is the unique least exact global extension under inclusion.
	Other exact extensions may exist.
	The conditions imply pairwise agreement on overlaps.
\end{theorem}

\begin{proof}
	Every extension contains the local relations and therefore their equivalence closure.
	Exact restriction then gives
	\[
	E_K\subseteq E^{\mathrm{gen}}\cap I_K^2
	\subseteq E^{\mathrm{glob}}\cap I_K^2=E_K.
	\]
	This proves (i) implies (ii).
	Conversely, (ii) makes \(E^{\mathrm{gen}}\) an exact extension and proves leastness.
	
	The relation of being joined by a finite local-equivalence chain, including a length-zero chain, is reflexive, symmetric by reversing a chain, and transitive by concatenation.
	It contains every \(E_K\), and every equivalence containing the \(E_K\) contains every such chain.
	It is thus exactly \(E^{\mathrm{gen}}\), proving (ii) equivalent to (iii).
	If an exact extension exists, restricting it to	\((I_K\cap I_L)^2\) gives agreement of \(E_K\) and \(E_L\).
	All arguments include \(E=\varnothing\).
\end{proof}

\begin{definition}[Primal graph, guarded sets, and chordality]
	\label{DEF:PackagingPrimalOverlapGraph}	
	Let
	$
	\mathcal I=(I_K)_{K\in\mathcal K}
	$
	be a cover of $E$, and regard it as the hypergraph
	$
	\mathcal H_{\mathcal I}
	=
	(E,\{I_K:K\in\mathcal K\}),
	$
	ignoring repeated hyperedges in its geometry.	
	\begin{enumerate}[label=(\roman*),leftmargin=2.5em]		
		\item \textbf{Primal graph.}
		The primal graph (or Gaifman graph)	$G_{\mathcal I}$ has vertex set $E$, with distinct vertices $x,y\in E$ adjacent if and only if
		$
		\{x,y\}\subseteq I_K
		$
		for some $K\in\mathcal K$.		
		Its vertices are units, not package labels.
		It is distinct from the package-overlap graph of Definition~\ref{DEF:PackagingIncidenceHypergraph}.
		
		\item \textbf{Guarded sets and $3$-conformality.}
		A set $S\subseteq E$ is guarded if
		$
		S\subseteq I_K
		$
		for some $K\in\mathcal K$.		
		The cover is $3$-conformal if every	nonempty clique of $G_{\mathcal I}$ of cardinality at most three is guarded.		
		Since every vertex and primal edge is guarded, this is equivalent to requiring every primal triangle to be guarded.
		
		\item \textbf{Chordality.}
		The primal graph $G_{\mathcal I}$ is chordal if it has no induced cycle of finite length at least four.
		
		This definition also applies when $E$ is infinite.		
	\end{enumerate}
\end{definition}

\begin{definition}[Pairwise-compatible local equivalence atlas]
	\label{DEF:PackagingPairwiseCompatibleKernelAtlas}
	A \textbf{local equivalence atlas} is \(\mathcal A=(E,(I_K,E_K)_{K\in\mathcal K})\), with \((I_K)_K\) a cover and each \(E_K\) an equivalence on \(I_K\).
	It is \textbf{pairwise compatible} when
	\[
	E_K\cap(I_K\cap I_L)^2=E_L\cap(I_K\cap I_L)^2
	\qquad(K,L\in\mathcal K).
	\]
\end{definition}

\begin{definition}[Signed primal graph and support-minimal obstructions]
	\label{DEF:PackagingSignedPrimalObstruction}	
	Let
	$
	\mathcal A
	=
	(E,(I_K,E_K)_{K\in\mathcal K})
	$
	be a pairwise-compatible equivalence atlas, and let $G_{\mathcal I}$ be the primal graph of its underlying cover.	
	\begin{enumerate}[label=(\roman*),leftmargin=2.5em]		
		\item \textbf{Signed primal graph.}		
		A primal edge $\{x,y\}$ is \textbf{positive} if
		$
		(x,y)\in E_K
		$
		for some chart containing both vertices, and \textbf{negative} otherwise.		
		Pairwise compatibility makes this sign independent of the chosen chart.		
		Denote the positive and negative edge sets by
		$
		\mathcal E_+(\mathcal A),
		~
		\mathcal E_-(\mathcal A),
		$
		respectively.
		Thus
		\[
		E(G_{\mathcal I})
		=
		\mathcal E_+(\mathcal A)
		\,\dot\cup\,
		\mathcal E_-(\mathcal A).
		\]		
		The \textbf{signed primal graph} is	$G_{\mathcal I}$ equipped with this edge-sign partition.		
		Its \textbf{positive subgraph} is
		$
		G_+(\mathcal A)
		=
		(E,\mathcal E_+(\mathcal A)).
		$
		
		\item \textbf{Support-minimal obstructions.}		
		For $S\subseteq E$, define the restricted atlas
		\[
		\mathcal A|_S
		=
		\bigl(
		S,
		(I_K\cap S,\,
		E_K\cap S^2)_{K\in\mathcal K}
		\bigr).
		\]		
		A finite set $S\subseteq E$ is a \textbf{support-minimal obstruction} if $\mathcal A|_S$ admits no exact	extension, but $\mathcal A|_T$ admits an exact extension for every proper subset $T\subsetneq S$.		
		The primal and positive graphs of $\mathcal A|_S$ are the corresponding induced subgraphs on $S$.		
		Define the \textbf{obstruction support number} by
		\[
		\begin{aligned}
			g_{\mathrm{obs}}(\mathcal A)
			:=
			1+\min\bigl\{
			&d_{G_+(\mathcal A)}(x,y):
			\\
			&\{x,y\}\in\mathcal E_-(\mathcal A),
			\\
			&d_{G_+(\mathcal A)}(x,y)<\infty
			\bigr\},
		\end{aligned}
		\]
		with value $\infty$ when the set is empty.		
		Here distances are lengths of shortest finite positive paths.		
		Thus $g_{\mathrm{obs}}(\mathcal A)$ is the minimum number of vertices on a positive path whose endpoints form a negative primal edge.		
	\end{enumerate}
\end{definition}

\begin{theorem}[Canonical minimal-obstruction theorem]
	\label{THM:PackagingCanonicalMinimalObstruction}	
	Let $\mathcal A$ be a pairwise-compatible equivalence atlas.	
	\begin{enumerate}[label=(\roman*),leftmargin=2em]
		\item \emph{Exact extension criterion.}
		
		An exact extension exists if and only if
		$
		g_{\mathrm{obs}}(\mathcal A)=\infty.
		$
		Equivalently, no negative primal edge has both endpoints in the same positive connected component.
		
		\item \emph{Canonical minimal obstruction.}
		
		Suppose
		$
		g_{\mathrm{obs}}(\mathcal A)=q<\infty.
		$
		Choose a minimizing negative edge $\{x,y\}$ and a positive geodesic
		$
		x=v_0,v_1,\ldots,v_{q-1}=y.
		$		
		Then $q\ge3$, and closing this path with the negative edge $\{x,y\}$ gives an induced primal cycle with exactly one negative edge.		
		If $q=3$, this triangle is necessarily unguarded.		
		Its vertex set
		$
		S=\{v_0,\ldots,v_{q-1}\}
		$
		is a support-minimal obstruction.
		
		\item \emph{Classification of all minimal obstructions.}
		
		A finite set $S\subseteq E$ is a support-minimal obstruction if and only if its induced primal graph is a cycle of length at least three with exactly one negative edge.
		All remaining cycle edges are positive.		
	\end{enumerate}	
	Consequently,
	\[
		g_{\mathrm{obs}}(\mathcal A)
		=
		\min\bigl\{
		|S|:
		S\text{ is a support-minimal obstruction}
		\bigr\},
	\]
	with $\min\varnothing=\infty$.	
	Thus every support-minimal failure has the signed form
	\[
	\underbrace{+\;\cdots\;+}_{q-1}\;-.
	\]	
	Here ``canonical'' describes the normal form, not a unique obstruction support or a preferred cyclic labeling.
\end{theorem}

\begin{proof}
	\leavevmode
	\begin{enumerate}[label=(\roman*),leftmargin=2em]	
		\item \emph{Exact extension criterion.}
		
		Let $R_+$ be the equivalence relation on $E$ whose classes are the connected components of $G_+(\mathcal A)$.		
		Every positive edge belongs to $R_+$, so $R_+$ contains every local equivalence relation $E_K$.		
		Moreover, any equivalence relation extending all $E_K$ must contain $R_+$, since it must identify the endpoints of every finite positive path.
		
		Suppose first that an exact extension $R$ exists.		
		If a negative edge $\{x,y\}$ had its endpoints in one positive component, then $(x,y)\in R$.		
		But some chart contains both $x$ and $y$, and their negative sign means that they are not equivalent in that chart.
		This contradicts exactness of $R$.		
		Conversely, suppose no negative edge has its endpoints in one positive connected component.
		
		We claim that $R_+$ is an exact extension.		
		For every chart $I_K$, we already have
		$
		E_K\subseteq R_+|_{I_K}.
		$
		Conversely, let $x,y\in I_K$ satisfy $(x,y)\in R_+$.		
		If $x=y$, then $(x,y)\in E_K$ by reflexivity.		
		If $x\ne y$, the pair $\{x,y\}$ is a primal edge.
		By assumption it cannot be negative, since its endpoints lie in the same positive component.
		Hence it is positive, giving $(x,y)\in E_K$.		
		Therefore
		$
		R_+|_{I_K}=E_K
		\quad(K\in\mathcal K).
		$		
		Thus $R_+$ is an exact extension.
		
		Finally, two vertices lie in the same positive component precisely when their positive graph distance is finite.
		By the definition of $g_{\mathrm{obs}}$, this proves
		\[
		R_+\text{ is exact}
		\quad\Longleftrightarrow\quad
		g_{\mathrm{obs}}(\mathcal A)=\infty.
		\]
		
		\item \emph{Canonical minimal obstruction.}
		
		Suppose
		$
		g_{\mathrm{obs}}(\mathcal A)=q<\infty.
		$		
		Choose a negative edge $\{x,y\}$ realizing this minimum and a shortest positive path
		$
		x=v_0,v_1,\ldots,v_{q-1}=y.
		$
		
		Since $\{x,y\}$ is negative, it is not a positive edge.
		Hence
		$
		d_{G_+(\mathcal A)}(x,y)=q-1\ge2,
		$
		and therefore $q\ge3$.
		
		A shortest path has distinct vertices.
		Closing it with $\{x,y\}$ gives a primal cycle with $q-1$ positive edges and one negative edge.		
		We show that this cycle is induced.
		
		For $q=3$, this is automatic.
		Moreover, if the triangle were guarded by a chart, its two positive edges would belong to the same local equivalence relation.		
		Transitivity would then force $(x,y)$ to belong to that relation, contradicting the negative sign of $\{x,y\}$.		
		Hence the triangle is unguarded.
		
		Now let $q\ge4$ and suppose the cycle has a chord
		\[
		\{v_i,v_j\},
		\qquad
		0\le i<j\le q-1,
		\]
		where $j-i\ge2$ and $(i,j)\ne(0,q-1)$.
				
		If this chord is positive, replacing the positive subpath from $v_i$ to $v_j$ by the chord gives a shorter positive path from $x$ to $y$.		
		This contradicts the choice of the geodesic.
		
		If the chord is negative, the positive subpath from $v_i$ to $v_j$ gives
		\[
		d_{G_+(\mathcal A)}(v_i,v_j)
		\le j-i
		\le q-2.
		\]
		
		Thus this negative edge would yield
		\[
		g_{\mathrm{obs}}(\mathcal A)
		\le 1+d_{G_+(\mathcal A)}(v_i,v_j)
		\le q-1,
		\]
		contradicting the minimality of $q$.		
		Therefore no chord exists, and the cycle is induced.
		
		Put
		$
		S=\{v_0,\ldots,v_{q-1}\}.
		$		
		The restricted atlas $\mathcal A|_S$ has no exact extension: the positive path joins the endpoints of the negative edge $\{x,y\}$.
		
		Now let $T\subsetneq S$.		
		Since the induced primal graph on $S$ is a cycle with exactly one negative edge, deleting at least one vertex either removes that negative edge or breaks its unique positive connecting path.		
		Hence the restricted atlas on $T$ contains no negative edge whose endpoints are positively connected.		
		By part~\textup{(i)}, $\mathcal A|_T$ admits an exact extension.		
		Thus $S$ is a support-minimal obstruction, proving \textup{(ii)}.
		
		\item \emph{Classification of all minimal obstructions.}
		
		Let $S\subseteq E$ be a finite support-minimal obstruction.		
		By part~\textup{(i)}, the restricted atlas $\mathcal A|_S$ has a negative edge whose endpoints are connected by a positive path.
		
		Apply part~\textup{(ii)} to $\mathcal A|_S$.		
		It produces a subset $T \subseteq S$ whose induced primal graph is a cycle with exactly one negative edge, and for which $\mathcal A|_T$ has no exact extension.		
		Since $S$ is support-minimal, necessarily
		$
		T=S.
		$		
		Therefore the induced primal graph on $S$ has the asserted form.
				
		Conversely, suppose the induced primal graph on a finite set $S$ is a cycle of length at least three with exactly one negative edge.		
		The remaining edges form a positive path joining the endpoints of that negative edge.		
		Part~\textup{(i)} therefore shows that $\mathcal A|_S$ admits no exact extension.		
		Every proper induced subgraph of this cycle breaks the positive connecting path or removes the negative edge.		
		Hence every proper restriction admits an exact extension, again by part~\textup{(i)}.		
		Thus $S$ is support-minimal, proving \textup{(iii)}.
		
		\medskip
		Finally, part~\textup{(ii)} produces a support-minimal obstruction of cardinality $g_{\mathrm{obs}}(\mathcal A)$
		whenever this number is finite.
		
		Conversely, by part~\textup{(iii)}, every support-minimal obstruction $S$ contains a positive path of	length $|S|-1$ joining the endpoints of a negative edge.		
		Therefore
		$
		g_{\mathrm{obs}}(\mathcal A)\le |S|.
		$		
		Taking the minimum over all such $S$ proves the displayed equality, including the case where no obstruction exists.
	\end{enumerate}
\end{proof}

\begin{proposition}[Cycle-freedom is chordality plus $3$-conformality]
	\label{PROP:PackagingBraultBaronCycleFree}	
	Let $\mathcal I=(I_K)_{K\in\mathcal K}$ be a cover of $E$.	
	Call the cover \textbf{cycle-free} if there is no finite set
	$$
	S=\{v_0,\ldots,v_{q-1}\},
	\quad q\ge3,
	$$
	whose inclusion-maximal nonempty traces among
	$$
	\{I_K\cap S:K\in\mathcal K\}
	$$
	are exactly the cycle edges
	$
	\{v_0,v_1\},\ldots,
	\{v_{q-2},v_{q-1}\},
	\{v_{q-1},v_0\}.
	$	
	For finite hypergraphs, this is Brault--Baron's cycle-free notion
	\cite{BraultBaron2016HypergraphAcyclicity}.
	The same finite-obstruction definition is used here for arbitrary set-indexed covers.	
	The following are equivalent:
	\begin{enumerate}[label=(\roman*),leftmargin=2em]		
		\item The cover hypergraph $\mathcal H_{\mathcal I}$ is cycle-free.
		
		\item The primal graph $G_{\mathcal I}$ is chordal and the cover is $3$-conformal.		
	\end{enumerate}	
	Equivalently, cycle-freedom means that every primal triangle is guarded and the primal graph has no induced finite cycle of length at least four.	
	Empty and repeated charts do not affect these conditions.
\end{proposition}

\begin{proof}
	We prove the two implications.
	\begin{enumerate}
		\item \emph{(ii) $\Longrightarrow$ (i).}
		
		Suppose, towards a contradiction, that the cover is not cycle-free.		
		Then there exists a finite set
		$
		S=\{v_0,\ldots,v_{q-1}\},
		~~
		q\ge3,
		$
		whose inclusion-maximal nonempty traces are exactly the edges of the displayed cycle.
		
		Since $S$ is finite, every trace $I_K\cap S$ is contained in an inclusion-maximal trace.		
		Consequently, two distinct vertices	of $S$ occur together in a chart if and only if they are consecutive on the cycle.		
		Thus
		$
		G_{\mathcal I}[S]\cong C_q.
		$
		
		If $q=3$, then $S$ is a primal triangle, but no chart contains all of $S$, contradicting $3$-conformality.
		
		If $q\ge4$, then $G_{\mathcal I}[S]$ is an induced cycle of length at least four, contradicting chordality.
		
		Hence the cover is cycle-free.
		
		\item \emph{(i) $\Longrightarrow$ (ii).}
		
		We argue by contraposition.		
		Suppose first that the cover is not $3$-conformal.		
		Since every vertex and primal edge is guarded, there exists an unguarded primal triangle
		$
		S=\{v_0,v_1,v_2\}.
		$		
		Every pair of vertices of $S$ occurs together in some chart, but no chart contains all three.		
		Therefore the inclusion-maximal nonempty traces on $S$ are exactly
		\[
		\{v_0,v_1\},
		\quad
		\{v_1,v_2\},
		\quad
		\{v_2,v_0\}.
		\]
		This is a forbidden cycle restriction.
		
		Now suppose instead that the primal graph is not chordal.		
		Then it contains an induced cycle on a finite vertex set
		\[
		S=\{v_0,\ldots,v_{q-1}\},
		\qquad q\ge4.
		\]		
		Every chart trace $I_K\cap S$ is a clique in the induced graph
		$G_{\mathcal I}[S]$.		
		Since this induced graph is a chordless cycle of length at least four, every such clique has at most two vertices, and any two-vertex trace must be a cycle edge.		
		Conversely, every cycle edge is	contained in some chart, by the	definition of the primal graph.		
		Hence the inclusion-maximal nonempty traces on $S$ are exactly the cycle edges.		
		This again contradicts cycle-freedom.		
		Thus cycle-freedom implies both chordality and $3$-conformality.
	\end{enumerate}
	Finally, empty charts contribute no	nonempty traces, and repeated charts do not change the set of traces.
	Therefore neither affects the stated conditions.
\end{proof}

\begin{definition}[Amalgamation cycle rank]
	\label{DEF:PackagingAmalgamationCycleRank}
	For a set-indexed cover \(\mathcal H=(I_K)_{K\in\mathcal K}\) of \(E\), define the \textbf{amalgamation cycle rank}
	\[
	{
		\lambda_{\mathrm{am}}(\mathcal H)
		=
		\min\left\{
		q\ge3:
		\begin{array}{l}
			\text{there is }S=\{v_0,\ldots,v_{q-1}\}\subseteq E\text{ such that}\\
			\text{the inclusion-maximal nonempty traces }I_K\cap S\\
			\text{are exactly the cycle edges }\{v_i,v_{i+1}\}\ (i\bmod q)
		\end{array}
		\right\},
	}
	\]
	with value \(\infty\) when no such finite support exists.
	Thus \(\lambda_{\mathrm{am}}(\mathcal H)=\infty\) exactly when the cover is cycle-free.
\end{definition}

\subsection{Universal cycle-free amalgamation}
\label{SUBSEC:PackagingUniversalCycleFreeAmalgamation}

We now ask when \emph{every} compatible atlas of transitive relations on a fixed cover extends exactly.
Theorem~\ref{THM:PackagingTransitiveUniversalAmalgamation} answers this question for transitive relations and four related classes.
Path elimination proves the positive direction, while triangles and longer induced cycles provide counterexamples when the cover condition fails.

\begin{definition}[Pairwise-compatible transitive atlas]
	\label{DEF:PackagingTransitiveAtlas}	
	Let
	$
	\mathcal I=(I_K)_{K\in\mathcal K}
	$
	be a cover of $E$.	
	\begin{enumerate}[label=(\roman*),leftmargin=2.5em]		
		\item \textbf{Transitive atlas.}
		A local transitive atlas on $\mathcal I$ is	a family
		\[
		\mathcal A=(E,(I_K,R_K)_{K\in\mathcal K}),
		\]
		where each
		$
		R_K\subseteq I_K^2
		$
		is a transitive relation.
		
		\item \textbf{Pairwise compatibility.}
		The atlas $\mathcal A$ is \textbf{pairwise compatible} if
		\[
		R_K\cap(I_K\cap I_L)^2
		=
		R_L\cap(I_K\cap I_L)^2
		\]
		for every $K,L\in\mathcal K$.		
	\end{enumerate}	
	A transitive relation $R\subseteq E^2$ is an \textbf{exact global extension} of $\mathcal A$ if
	\[
	R\cap I_K^2=R_K
	\qquad(K\in\mathcal K).
	\]	
	Define the \textbf{generated transitive relation} by
	\[
	R_{\mathcal A}^{\mathrm{tr}}
	:=
	\left(\bigcup_{K\in\mathcal K}R_K\right)^+,
	\]
	where $T^+$ denotes reachability by finite directed paths of positive length in $T$.	
	Thus $R_{\mathcal A}^{\mathrm{tr}}$ is the least transitive relation on $E$ containing all local relations $R_K$.
	It need not be an exact extension.	
	The use of positive-length paths means that reflexivity is not imposed:
	diagonal pairs arise only when supported by positive-length directed closed walks.
\end{definition}

\begin{lemma}[Generated transitive extension criterion]
	\label{LEM:PackagingTransitiveExactExtension}
	For any local transitive atlas on a cover, an exact transitive extension exists if and only if
	\[
	R^{\mathrm{tr}}\cap I_K^2=R_K\qquad(K\in\mathcal K).
	\]
	When it exists, \(R^{\mathrm{tr}}\) is its least such extension.
	Failure has a positive-length local-edge path whose endpoints form a forbidden ordered pair in one chart.
	The endpoints may coincide.
	Neither pairwise compatibility nor a geometric hypothesis is needed for this criterion.
\end{lemma}

\begin{proof}
	Positive-length reachability contains the generating relation and is transitive by concatenation.
	Every transitive relation containing the generators contains every finite positive-length path, by induction on its length.
	It is therefore the least transitive extension.
	If an exact	extension \(H\) exists, then \(R_K\subseteq R^{\mathrm{tr}}\cap I_K^2\subseteq H\cap I_K^2=R_K\).
	Conversely, the displayed equalities make \(R^{\mathrm{tr}}\) exact.
	A pair in their difference has the stated path witness.
	On the diagonal that witness must be a nonempty closed walk, not a length-zero path.
\end{proof}

\begin{lemma}[One-pivot elimination for positive-length transitive closure]
	\label{LEM:PackagingTransitiveOnePivotElimination}
	Let $W$ be finite, let $T\subseteq W^2$, let $v\in W$, and put $J=W\setminus\{v\}$.
	Define
	\[
	\operatorname{Comp}_v(T)
	=\{(x,y)\in J^2:(x,v)\in T,\ (v,y)\in T\}.
	\]
	Writing $T^+$ for positive-length transitive closure, one has
	\[
	T^+\cap J^2
	=\bigl((T\cap J^2)\cup\operatorname{Comp}_v(T)\bigr)^+.
	\]
	Diagonal pairs are retained in both the composition set and the closure.
\end{lemma}

\begin{proof}
	Put
	$
	U=(T\cap J^2)\cup\operatorname{Comp}_v(T).
	$
	We prove both inclusions.
	\begin{itemize}
		\item \emph{First inclusion: $U^+\subseteq T^+\cap J^2$.}
		
		Every edge in $T\cap J^2$ is already an edge of $T$.		
		If
		$
		(x,y)\in\operatorname{Comp}_v(T),
		$
		then
		\[
		(x,v)\in T,
		\qquad
		(v,y)\in T.
		\]
		Thus
		$
		x\longrightarrow v\longrightarrow y
		$
		is a $T$-walk of length two.
		
		Consequently, every edge of $U$	can be replaced by a positive-length $T$-walk with the same endpoints.		
		Expanding the edges of any positive	$U$-walk therefore produces a positive $T$-walk.		
		Since all vertices of $U$ lie in $J$, this proves
		$
		U^+\subseteq T^+\cap J^2.
		$
		
		\item \emph{Reverse inclusion: $T^+\cap J^2\subseteq U^+$.}
		
		Let
		$
		(x,y)\in T^+\cap J^2.
		$
		By definition of positive-length transitive closure, there exists a finite $T$-walk
		\[
		x=w_0,w_1,\ldots,w_n=y,
		\qquad n\ge1,
		\]
		such that
		$
		(w_{i-1},w_i)\in T
		\quad(1\le i\le n).
		$		
		List, in increasing order, all indices at which this walk visits $J$:
		\[
		0=i_0<i_1<\cdots<i_m=n.
		\]
		Such indices exist because $x,y\in J$, and $m\ge1$ since $n\ge1$.		
		We claim that
		\[
		(w_{i_{t-1}},w_{i_t})\in U
		\qquad(1\le t\le m).
		\]		
		Fix $t$.		
		If
		$
		i_t=i_{t-1}+1,
		$
		the corresponding vertices are consecutive in the original walk.
		Both lie in $J$, so
		\[
		(w_{i_{t-1}},w_{i_t})
		\in T\cap J^2\subseteq U.
		\]		
		Otherwise,
		$
		i_t\ge i_{t-1}+2.
		$
		By the choice of consecutive indices of visits to $J$, every intermediate vertex lies outside $J$.		
		Since $W \setminus J = \{v\}$, we have
		\[
		w_{i_{t-1}+1}
		=\cdots=
		w_{i_t-1}
		=v.
		\]		
		In particular, the first and last edges of this segment give
		\[
		(w_{i_{t-1}},v)\in T,
		\qquad
		(v,w_{i_t})\in T.
		\]		
		Therefore
		\[
		(w_{i_{t-1}},w_{i_t})
		\in\operatorname{Comp}_v(T)
		\subseteq U.
		\]		
		The claim follows in both cases.		
		Hence
		$
		x=w_{i_0},w_{i_1},\ldots,w_{i_m}=y
		$
		is a $U$-walk of length $m\ge1$.
		Thus
		$
		(x,y)\in U^+.
		$
		This proves
		\[
		T^+\cap J^2\subseteq U^+.
		\]
	\end{itemize}
	Combining the two inclusions gives
	$
	T^+\cap J^2
	=
	\bigl((T\cap J^2)
	\cup\operatorname{Comp}_v(T)\bigr)^+.
	$	
	The argument also covers $x=y$:
	the compressed walk still has positive length, and a segment leaving and returning to the same vertex of $J$ produces a diagonal edge in $\operatorname{Comp}_v(T)$.	
	Thus diagonal pairs are preserved without inserting reflexive pairs by convention.
\end{proof}

For two charts, strong amalgamation into the union for these relation classes is covered by Lipparini \cite[Proposition~3.2 and Corollary~3.3]{Lipparini2021SAPU}.
Otto's amalgamation-pattern framework gives finite global realizations of locally specified finite templates under a different, highly general realization setup \cite{Otto2024AmalgamationPatterns}.
Here the cover itself is fixed, may have arbitrarily many charts, and exact restriction to every original chart is required.
We prove both directions of the resulting cycle-free universal criterion directly.

\begin{remark}[Classical ingredients and scope of the universal theorem]
	\label{REM:PackagingTransitiveAmalgamationNoveltyScope}
	The cover geometry is classical.
Cycle-freedom is standard in hypergraph acyclicity, while \(\alpha\)-acyclicity additionally requires full conformality \cite{BeeriFaginMaierYannakakis1983Acyclic,BraultBaron2016HypergraphAcyclicity}.
	Strong amalgamation of two structures in these relation classes is also known \cite{Lipparini2021SAPU}, and finite preorder sufficiency on chordal constraint graphs is closely related to chordal path-consistency results
	\cite{BliekSamHaroud1999Path,LecoutreCardonVion2011Consistency}.
	
	Accordingly, the role of the theorem here is not to claim these classical ingredients as new.
	Our theorem fixes a cover, allows any finite number of charts, and requires exact agreement on every original chart.
The same criterion works for five classes of transitive relations.
We also give small, explicit obstructions when the criterion fails.
	Its importance for packaging is comparative:
	the same obstruction geometry that governs the closure specializations reappears without using Horn resolution, Carath\'eodory bounds, Armstrong relations, or spectral compression.
\end{remark}

\begin{theorem}[Universal cycle-free amalgamation for transitive binary structures]
	\label{THM:PackagingTransitiveUniversalAmalgamation}
	Let $\mathcal I=(I_K)_{K\in\mathcal K}$ be a set-indexed cover of a set $E$.
	The following are equivalent:
	\begin{enumerate}[label=(\roman*),leftmargin=2em]
		\item The cover hypergraph $\mathcal H_{\mathcal I}$ is cycle-free.
		
		\item Every pairwise-compatible transitive atlas on $\mathcal I$ admits an exact global transitive relation.
		
		\item For every pairwise-compatible transitive atlas $\mathcal A=(E,(I_K,R_K)_{K\in\mathcal K})$,
		\[
		R_{\mathcal A}^{\mathrm{tr}}
		=\left(\bigcup_{K\in\mathcal K}R_K\right)^+
		\]
		satisfies
		\[
		R_{\mathcal A}^{\mathrm{tr}}\cap I_K^2=R_K
		\qquad(K\in\mathcal K),
		\]
		and is the unique least exact global transitive extension.
		
		\item Every pairwise-compatible local preorder atlas admits	an exact global preorder.
		
		\item Every pairwise-compatible local partial-order atlas admits an exact global partial order.
		
		\item Every pairwise-compatible local strict-partial-order atlas admits an exact global strict partial order.
		
		\item Every pairwise-compatible local equivalence atlas admits an exact global equivalence relation.
	\end{enumerate}
	Whenever these conditions hold, the generated relation $R_{\mathcal A}^{\mathrm{tr}}$ is also the unique least exact extension within each of the four indicated subclasses.
	No finiteness, connectedness, nonempty-chart, or distinct-chart	hypothesis is required.
\end{theorem}

\begin{proof}
	Put $U=\bigcup_{K\in\mathcal K}R_K$ and $R^+=U^+$ for a	pairwise-compatible transitive atlas $\mathcal A$.	
	\begin{itemize}
		\item \emph{(i) $\Longrightarrow$ (iii): exactness of the generated transitive relation.}
		Assume the cover is cycle-free.
		By Proposition~\ref{PROP:PackagingBraultBaronCycleFree}, its primal graph is chordal and every primal triangle is guarded.
		Certainly $R_K\subseteq R^+\cap I_K^2$ for every $K$.
		
		Suppose first that an off-diagonal violation occurs:
		for some chart $I_K$ and distinct $a,b\in I_K$,
		\[
		(a,b)\in R^+\setminus R_K.
		\]
		Among all such violations, choose one for which the length $n$ of a shortest $U$-walk from $a$ to $b$ is minimum.
		Thus there are
		\[
		a=v_0\longrightarrow v_1\longrightarrow\cdots
		\longrightarrow v_n=b,
		\qquad (v_{i-1},v_i)\in U.
		\]
		This walk is a simple path: deleting the segment between repeated vertices would produce a shorter $U$-walk witnessing the same violation.
		The case $n=1$ is impossible, since pairwise compatibility transfers the edge from its source chart to $R_K$.
		If $n=2$, the three vertices form a primal triangle, since its path edges and the endpoint pair $\{a,b\}$ are primal edges.
		A chart guarding the triangle contains both path edges by compatibility.
		Local transitivity then gives $(a,b)$, which transfers to $R_K$, a contradiction.
		
		If $n\ge3$, the path together with the primal edge $\{a,b\}$ forms a cycle of length at least four.
		Chordality supplies a chord	$\{v_i,v_j\}$, where
		\[
		0\le i<j\le n,\qquad j-i\ge2,\qquad (i,j)\ne(0,n).
		\]
		Choose a chart $I_L$ containing $v_i,v_j$.
		The directed subpath from $v_i$ to $v_j$ has length $j-i<n$.
		If $(v_i,v_j)\in R_L$, replace this subpath by that one local edge, obtaining a shorter $U$-walk for the original violation $(a,b)$.
		Otherwise the subpath gives an off-diagonal violation $(v_i,v_j)\in R^+\setminus R_L$ with a shorter witness.
		Both alternatives contradict the choice of $n$.
		Hence no off-diagonal violation exists.
		
		Next suppose a diagonal violation occurs: for some $a\in I_K$,
		\[
		(a,a)\in R^+\setminus R_K.
		\]
		Choose a positive-length $U$-closed walk $a=v_0\to v_1\to\cdots\to v_n=a$.
		If its first edge is a loop at $a$, compatibility with its source chart already gives $(a,a)\in R_K$.
		We may therefore take $v_1\ne a$.
		Let $I_L$ be a source chart for $a\to v_1$.
		The remaining positive walk gives $(v_1,a)\in R^+$.
		Since off-diagonal violations have been excluded and $a,v_1\in I_L$, we obtain $(v_1,a)\in R_L$.
		Together with $(a,v_1)\in R_L$, local transitivity gives $(a,a)\in R_L$, and compatibility transfers the loop to $R_K$.
		This is a contradiction.
		Thus $R^+$ is exact on every chart.
		
		Since $R^+$ is the least transitive relation containing $U$, every exact global transitive extension contains it.
		Consequently $R^+$ is the unique least exact transitive extension, proving~\textup{(iii)}.
		
		\item \emph{(ii) $\Longleftrightarrow$ (iii).}
		Assertion~\textup{(iii)} plainly implies~\textup{(ii)}.
		Conversely, if an exact global transitive relation $R$ exists, minimality of positive-length transitive closure gives $R^+\subseteq R$.
		Therefore, for every $K$,
		\[
		R_K\subseteq R^+\cap I_K^2
		\subseteq R\cap I_K^2=R_K.
		\]
		Hence $R^+$ is exact, and the same minimality makes it the unique least exact transitive extension.
		Thus these two assertions are equivalent for any cover.
		
		\item \emph{(i) $\Longrightarrow$ (iv)--(vii): preservation of the four subclasses.}
		Continue to assume~\textup{(i)}, so that $R^+$ is exact.		
		For~\textup{(iv)}, if each $R_K$ is reflexive, every $x\in E$ belongs to some chart, so $(x,x)\in U\subseteq R^+$.
		Thus $R^+$ is a preorder.
		
		For~\textup{(v)}, suppose every $R_K$ is a partial order.
		The preceding argument gives reflexivity.
		If distinct $x,y$ satisfied both $(x,y),(y,x)\in R^+$, choose a $U$-walk from $x$ to $y$ and discard any initial loops.
		Its first edge $x\to z$ has $z\ne x$.
		Let $I_L$ contain this edge.
		The remainder of that walk, followed by a $U$-walk from $y$ to $x$, gives $(z,x)\in R^+$.
		Exactness on $I_L$ gives $(z,x)\in R_L$, while $(x,z)\in R_L$.
		This contradicts the antisymmetry of $R_L$.
		Hence $R^+$ is a partial order.
		
		For~\textup{(vi)}, if every $R_K$ is irreflexive, any loop $(x,x)\in R^+$ would belong, by exactness, to $R_K$ for a chart containing $x$.
		This is impossible.
		Thus $R^+$ is irreflexive and transitive, hence a strict partial order.
		
		For~\textup{(vii)}, if every $R_K$ is an equivalence relation, $U$ is reflexive and symmetric.
		Reversing a positive $U$-walk shows that $R^+$ is symmetric.
		It is therefore an equivalence relation.
		
		In each case, $R^+$ is already the least transitive relation containing the local data, and consequently is also the unique least exact extension within the specified subclass.
		
		\item \emph{(ii), (iv), (v), (vi), or (vii) $\Longrightarrow$ (i): common geometric counterexample.}
		We prove all five implications by contraposition.
		Suppose the cover is not cycle-free.
		Then there is a finite set
		\[
		S=\{v_0,\ldots,v_{q-1}\},\qquad q\ge3,
		\]
		whose inclusion-maximal nonempty chart traces are exactly the consecutive edges of the cycle $C_q$.
		Consequently every $I_K\cap S$ has at most two vertices, and any two-vertex trace is one of those edges.
		
		For each $K$, put
		\[
		P_K=\{(v_i,v_{i+1}):0\le i<q-1,
		\ \{v_i,v_{i+1}\}\subseteq I_K\},
		\qquad \Delta_K=\{(x,x):x\in I_K\}.
		\]
		Each $P_K$ contains at most one ordered pair.
		On the original	cover define five families of local relations:
		\[
		\begin{aligned}
			R_K^{\mathrm{tr}}=R_K^{\mathrm{spo}}&=P_K,\\
			R_K^{\mathrm{pre}}=R_K^{\mathrm{po}}&=\Delta_K\cup P_K,\\
			R_K^{\mathrm{eq}}&=\Delta_K\cup P_K\cup P_K^{-1}.
		\end{aligned}
		\]
		Because each $P_K$ is empty or consists of a single nonloop edge, these are, respectively, transitive relations, strict partial orders, preorders, partial orders, and equivalence relations.
		The definitions depend only on the complete support of each off-diagonal pair, so all five atlases are pairwise compatible.
		
		Every path edge $v_i\to v_{i+1}$, $0\le i<q-1$, is present in its local charts.
		Any global extension in any of the five transitive classes must therefore contain $(v_0,v_{q-1})$.
		But some chart contains the closing primal edge $\{v_{q-1},v_0\}$, and its prescribed local relation omits $(v_0,v_{q-1})$ in each of the five atlases.
		Exactness is therefore impossible in every class.
		This disproves~\textup{(ii)} and each of~\textup{(iv)}--\textup{(vii)} whenever~\textup{(i)} fails.
	\end{itemize}	
	Combining the implications proves the equivalence of \textup{(i)}--\textup{(vii)}.
\end{proof}

\begin{corollary}[Rank-three obstruction-spectrum form]
	\label{COR:PackagingTransitiveObstructionRankForm}
	Let $\mathcal H$ be a finite cover.
	The following are equivalent:
	\begin{enumerate}[label=(\roman*),leftmargin=2.5em]
		\item $3<\nu(\mathcal H)$;
		
		\item the primal graph is chordal and the cover is $3$-conformal;
		
		\item the cover is cycle-free;
		
		\item every pairwise-compatible transitive atlas admits an exact global transitive relation;
		
		\item universal exact amalgamation holds separately for each of the classes of preorders, partial orders, strict partial
		orders, and equivalence relations.
	\end{enumerate}
	Thus transitive-relation amalgamation realizes the rank-three member of the finite-cover obstruction spectrum shared with the closure and bounded-Armstrong specializations.
\end{corollary}

\begin{proof}
	\emph{(i) $\Longleftrightarrow$ (ii).}
	This is Lemma~\ref{LEM:PackagingObstructionRankDictionary} with $m=3$.
	
	\emph{(ii) $\Longleftrightarrow$ (iii).}
	This is Proposition~\ref{PROP:PackagingBraultBaronCycleFree}.
	
	\emph{(iii) $\Longleftrightarrow$ (iv).}
	This is the equivalence of~\textup{(i)} and~\textup{(ii)} in Theorem~\ref{THM:PackagingTransitiveUniversalAmalgamation}.
	
	\emph{(iii) $\Longleftrightarrow$ (v).}
	This follows from the equivalence of~\textup{(i)} with each of \textup{(iv)}--\textup{(vii)} in that theorem.
	Hence all five conditions are equivalent.
\end{proof}

\begin{theorem}[Sharp amalgamation obstruction rank and primitive basis]
	\label{THM:PackagingAmalgamationPrimitiveBasis}
	Let $\mathcal H=(I_K)_{K\in\mathcal K}$ be an arbitrary set-indexed cover of $E$, with amalgamation cycle rank $\lambda_{\mathrm{am}}(\mathcal H)$ as in Definition~\ref{DEF:PackagingAmalgamationCycleRank}.
	Then:
	\begin{enumerate}[label=(\roman*),leftmargin=2.5em]
		\item The following are equivalent:
		\[
		\lambda_{\mathrm{am}}(\mathcal H)=\infty
		\quad\Longleftrightarrow\quad
		\text{universal exact amalgamation holds}
		\]
		in each of the classes of transitive relations, preorders, partial orders, strict partial orders, and equivalence relations.
		
		\item If there is an unguarded primal triangle, then $\lambda_{\mathrm{am}}(\mathcal H)=3$.
		Otherwise, if the primal graph is nonchordal, $\lambda_{\mathrm{am}}(\mathcal H)$ is the length of a shortest induced primal cycle.
		If the primal graph is chordal and every primal triangle is guarded, then $\lambda_{\mathrm{am}}(\mathcal H)=\infty$.
		
		\item Suppose $q=\lambda_{\mathrm{am}}(\mathcal H)<\infty$.
		There is a set $S\subseteq E$ with $|S|=q$ such that, for each of the five indicated relation classes, the restricted cover $\mathcal H|_S=(I_K\cap S)_{K\in\mathcal K}$ admits a pairwise-compatible atlas without an exact global extension, but whose restriction to every proper $T\subsetneq S$ is exactly extendible in the same class.
		In the transitive and strict-order examples, each restricted chart has at most one nontrivial directed edge.
For preorders and partial orders we add only the diagonal.
In the equivalence example, each chart has at most one nontrivial two-point equivalence block.
		No nonextendible compatible atlas in any of these classes can occur on a finite vertex restriction having fewer than $q$ vertices.
	\end{enumerate}
	Consequently, up to relabeling and deletion of empty, repeated, or nonmaximal traces, the \emph{vertex-restriction-minimal forbidden covers} for universal exact amalgamation are exactly the cycle hypergraphs
	\[
	C_q=\left(\{0,\ldots,q-1\},
	\bigl\{\{i,(i+1)\bmod q\}:0\le i<q\bigr\}\right),
	\qquad q\ge3.
	\]
\end{theorem}

\begin{proof}
	\emph{(i) Universal amalgamation criterion.}
	By Definition~\ref{DEF:PackagingAmalgamationCycleRank}, $\lambda_{\mathrm{am}}(\mathcal H)=\infty$ holds exactly when $\mathcal H$ is cycle-free.
	Theorem~\ref{THM:PackagingTransitiveUniversalAmalgamation} identifies cycle-freedom with universal exact amalgamation in each of the five relation classes.
	This proves~\textup{(i)}.
	
	\medskip
	\noindent\emph{(ii) Exact value of the obstruction rank.}
	A three-vertex cycle-trace obstruction is precisely an unguarded primal triangle.
	Thus if one exists, the minimum in the definition of $\lambda_{\mathrm{am}}$ equals $3$.
	Otherwise every triangle is guarded.
	By Proposition~\ref{PROP:PackagingBraultBaronCycleFree}, a finite cycle-trace obstruction must then be an induced primal cycle of length at least four.
	Conversely, every such induced cycle has exactly its cycle edges as maximal nonempty traces.
	Minimizing the vertex count therefore gives the shortest induced cycle length if the graph is nonchordal, and gives $\infty$ if it is chordal.
	This proves~\textup{(ii)}.
	
	\medskip
	\noindent\emph{(iii) Sharp witnesses on a minimum support.}
	Let $q=\lambda_{\mathrm{am}}(\mathcal H)<\infty$.
	Choose a witnessing set $S=\{v_0,\ldots,v_{q-1}\}$ whose maximal nonempty traces are exactly the edges of the cyclic order.
	Put $J_K=I_K\cap S$ and
	\[
	P_K=\{(v_i,v_{i+1}):0\le i<q-1,
	\ \{v_i,v_{i+1}\}\subseteq J_K\},
	\qquad
	\Delta_{J_K}=\{(x,x):x\in J_K\}.
	\]
	Every $J_K$ is contained in one cycle edge, so $P_K$ contains at most one nonloop pair.
	Define the five local atlases on the induced cover by
	\[
	\begin{aligned}
		R_K^{\mathrm{tr}}=R_K^{\mathrm{spo}}&=P_K,\\
		R_K^{\mathrm{pre}}=R_K^{\mathrm{po}}&=\Delta_{J_K}\cup P_K,\\
		R_K^{\mathrm{eq}}&=\Delta_{J_K}\cup P_K\cup P_K^{-1}.
	\end{aligned}
	\]
	Each local relation belongs to its stated class: a restricted chart has at most one nonloop path edge.
	Moreover, membership of	any off-diagonal pair is determined solely by its two endpoints, and all local diagonal requirements agree.
	Thus the five atlases are pairwise compatible.
	Every proposed exact extension would propagate along $v_0\to v_1\to\cdots\to v_{q-1}$, contradicting the local relation on the closing edge $\{v_{q-1},v_0\}$.
	Hence all five atlases are nonextendible.
	
	Let $T\subsetneq S$.
	If the restricted cover $\mathcal H|_T$	contained any cycle-trace obstruction on $P\subseteq T$, then $P$ would also be such an obstruction in $\mathcal H$, with $|P|<q$, contrary to the definition of $q$.
	Thus $\mathcal H|_T$ is cycle-free.
	The universal theorem supplies an exact extension, in the corresponding relation class, for the restriction of each constructed atlas to $T$.
	The same argument applies to every finite support of cardinality less than $q$, proving the last assertion in~\textup{(iii)}.
	
	\medskip
	\noindent\emph{Primitive forbidden-cover basis.}
	Let $\mathcal H$ be a vertex-restriction-minimal cover for which universal exact amalgamation fails.
	By part~\textup{(i)}, the cover has a finite cycle-trace obstruction on some $S\subseteq E$.
	If $S\subsetneq E$, the restricted cover $\mathcal H|_S$ would itself fail universal exact amalgamation, contradicting vertex
	minimality.
	Therefore $S=E$, so $E$ is finite and the maximal nonempty traces of $\mathcal H$ are exactly the edges of a cycle $C_q$ with $q\ge3$.
	
	Conversely, each $C_q$ contains a cycle-trace obstruction on its full vertex set, whereas deletion of any vertex destroys every cycle-trace obstruction.
	Thus $C_q$ fails universal exact amalgamation but every proper vertex restriction satisfies it.
	Removing empty, repeated, or nonmaximal chart traces does not change this characterization. This proves the primitive basis assertion.
\end{proof}

\begin{remark}[Support size versus globalization threshold rank]
	\label{REM:PackagingAmalgamationRankVersusNu}
	For a finite cover,
	\[
	\lambda_{\mathrm{am}}(\mathcal H)=\infty
	\quad\Longleftrightarrow\quad
	3<\nu(\mathcal H).
	\]
	The finite values of these two invariants need not agree.
	The first records the smallest vertex support of a cycle-trace amalgamation obstruction.
	The second is a threshold rank:
	nonchordality contributes the value $3$, regardless of the length of its shortest induced cycle.
	Hence, if every triangle is guarded and the shortest induced primal cycle has length $q\ge4$, then
	\[
	\lambda_{\mathrm{am}}(\mathcal H)=q,
	\qquad \nu(\mathcal H)=3.
	\]
	This distinction separates minimal obstruction support size from the presentation-independent globalization threshold.
\end{remark}

\subsection{Equivalence specialization, sharpness, and obstruction certificates}
\label{SUBSEC:PackagingEquivalenceAmalgamationConsequences}

The cycle-free classification also applies to equivalence relations.
We compare this result with path consistency and with the stronger full-conformality condition used in database acyclicity.
We also describe explicit certificates when exact extension fails.
These comparisons clarify what is shared with classical criteria and what differs in our exact-trace problem.

\noindent\textbf{Comparison with chordal path consistency.}
Finite preorder sufficiency is also a consequence of the equivalence of conservative and partial path consistency on chordal constraint graphs
\cite{BliekSamHaroud1999Path,LecoutreCardonVion2011Consistency}.
Here is the reduction, with the constraint graph retained explicitly.
Give each unit \(x\) two Boolean coordinates \((p_x,a_x)\).
On a primal edge \(xy\), require \(a_x=a_y\), together with \(p_x\le p_y\) if \(xR_Ky\) and \(p_y\le p_x\) if \(yR_Kx\), in any chart containing
the pair.
Compatibility makes these conditions well defined.
The auxiliary equality makes every primal edge a nonuniversal constraint, even when the two units are locally incomparable.

Every permitted assignment on an edge of a guarded triangle extends to the third vertex.
Indeed, give the first coordinates value one on the upward closure, in the guarding preorder, of the assigned one-vertices, and value zero elsewhere.
This retains each assigned zero, since an implication from an assigned one to that zero would already violate the endpoint constraint.
Give the auxiliary coordinates their common value.
The resulting network is conservatively path consistent, and the cited chordal theorem makes it partially path consistent, namely consistent
along all paths with adjacent endpoints.
If a positive directed path forced a locally forbidden \(xRy\), the permitted endpoint assignment \(p_x=1,p_y=0\), with equal auxiliary coordinates, could not extend along that path.
This contradiction gives exactness.
Reflexivity already settles the diagonal for a preorder.
The direct proof above avoids dependence on this translation and includes nonreflexive relations and diagonal violations.
Its finite-witness argument does not require an infinite perfect elimination ordering.

\begin{remark}[Classical multicut comparison]
	Mark a positive primal edge zero and a negative edge one.
	Exact equivalence extension holds precisely when no cycle has exactly one one-edge.
	The zero-edge connected components give the least equivalence.
	This is the classical multicut criterion, including its reduction to chordless cycles \cite{ChopraRao1993Partition,HornakovaLangeAndres2017Multicuts}.
	Theorem~\ref{THM:PackagingCanonicalMinimalObstruction} additionally states its vertex-support-minimal form for the given compatible atlas.
	The	directed transitive theorem has different local data and also treats loops.
\end{remark}

\begin{proposition}[Separation from the classical database threshold]
	\label{PROP:PackagingDatabaseAcyclicitySeparation}
	For a finite cover hypergraph:
	\begin{enumerate}[label=(\roman*),leftmargin=2em]
		\item Every family of database relations with equal projections on pairwise overlaps has a relation with those exact local projections if and only if the cover is \(\alpha\)-acyclic.
		Here the universal property	ranges over ordinary set-valued relations and value domains on the specified schemas, without extra dependencies \cite{BeeriFaginMaierYannakakis1983Acyclic}.
		
		\item \(\alpha\)-acyclicity is conformality together with cycle-freeness, where conformality requires every nonempty primal clique to lie in a chart \cite{BraultBaron2016HypergraphAcyclicity}.
		
		\item Every compatible local transitive atlas, or every compatible local equivalence atlas, extends exactly if and only if the cover is cycle-free, equivalently if and only if its primal graph is chordal and the cover is $3$-conformal.
	\end{enumerate}
	Consequently the class of finite \(\alpha\)-acyclic covers is a proper subclass of the finite cycle-free covers.
\end{proposition}

\begin{proof}
	The first two clauses are the cited classical characterization and	hypergraph criterion.
	The third is the theorem just proved.
	Full conformality implies $3$-conformality.
	Since $\alpha$-acyclicity is conformality plus cycle-freedom, this gives the inclusion.
	The following explicit example proves strictness.
	The conditions include empty/repeated charts by the usual redundant-chart conventions.
	On an empty vertex set the extension assertion is immediate.
\end{proof}

\begin{example}[Tetrahedral separation of cycle-freeness from \(\alpha\)-acyclicity]
	\label{EX:PackagingTetrahedralAcyclicitySeparation}
	Take \(E=\{1,2,3,4\}\) and \(I_i=E\setminus\{i\}\).
	The primal graph is \(K_4\) and every triangle is a chart.
	Thus the primal graph is chordal, the cover is $3$-conformal, and hence the cover is cycle-free.
	The four-vertex clique lies in no chart, so the	cover is not conformal or \(\alpha\)-acyclic.
	This geometric example is also recorded in \cite{BraultBaron2016HypergraphAcyclicity}.
	Every compatible atlas in each of the five relation classes of the main theorem nevertheless extends exactly on this cover.
	
	For the contrasting database problem, give each attribute the domain \(\mathbb F_2\), and set
	\[
	Q_i=\{t\in\mathbb F_2^{I_i}:\sum_{j\in I_i}t_j=0\}.
	\]
	Each projection to a two-attribute overlap is full:
	arbitrary values on the two coordinates determine the third.
	Hence the relations agree on all pairwise projections.
	In their natural join, subtracting any two of the equations \(\sum_{j\ne i}t_j=0\) shows that all four coordinates are equal.
	Each equation then forces their common value to be zero.
	The join is therefore the singleton zero tuple.
	Its projection has one tuple, whereas each \(Q_i\) has four.
	Any proposed exact global relation would be a subset of the natural join and could not have those projections.
	Thus global consistency fails, despite the nonempty join.
	Projection of relations on attribute assignments and restriction of binary relations on units are distinct extension problems.
\end{example}

\begin{remark}[Dependence of gluing on premise size]
	\label{REM:PackagingTransitiveNotArbitraryHorn}
	Empty-preserving unary dependency closures are preorder closures.
	Their simplex-boundary threshold is the unary case of Theorem~\ref{THM:PackagingSharpDependencyGluing}, whose general basis bound
	is \(n\ge2k+2\).
	Corollary~\ref{COR:PackagingDependencyHierarchy} separates successive bounds.
	The present transitive argument must not be used to	assert extension of arbitrary Horn dependencies on every cycle-free cover.
\end{remark}

\begin{example}[Pairwise-compatible data with a cycle obstruction]
	\label{EX:PackagingKernelCycleObstruction}
	Let \(E=\{a,b,c\}\), with charts \(\{a,b\},\{b,c\},\{c,a\}\).
	Use the discrete equivalence on the first chart and the indiscrete equivalences on the other two.
	All overlaps are singletons, hence the atlas is compatible.
	But \(b\sim c\sim a\) forces the pair distinguished	on \(\{a,b\}\).
	There is no exact extension.
	This is the smallest unguarded-triangle obstruction.
	On the same non-cycle-free cover, choosing every local equivalence discrete instead gives the discrete exact extension.
	The geometric threshold is universal over atlases.
\end{example}

\subsection{Presentation invariance and quotient realization}
\label{SUBSEC:PackagingGlobalizationPresentationInvariance}

The exact-extension conditions depend on which units belong to which charts and on the local relations, not on how the units and charts are named.
We make this invariance precise, including repeated and empty charts.
When an equivalence atlas extends exactly, its global relation also gives an ordinary quotient presentation.
This does not turn quotient formation into a primitive packaging axiom.

\begin{proposition}[Presentation invariance and degenerate atlases]
	\label{PROP:PackagingKernelOverlapPresentationInvariance}
	The kernel-extension criteria are invariant under a unit bijection and chart reindexing transporting all local equivalences.
	Repeating or deleting an identical chart--relation pair does not change the criterion or the least extension.  Feasibility can be tested componentwise on disconnected primal graphs.
	Empty covered sets,	empty charts, and singleton charts are allowed.
	Equal chart interiors can admit exact extension only if their local relations agree.
	The analogous assertions hold for the transitive generated-extension test, with the local diagonal retained as supplied.
\end{proposition}

\begin{proof}
	For a transporting bijection \(h:E\to E'\), direct image under \(h\times h\) commutes with unions and with either closure:
	a finite path transports in both directions, with its length unchanged.
	Thus it transports generated relations and restriction equalities.
	Identical duplicates leave the generating union unchanged.
	
	Every local edge lies in a primal component, and so does every generating path.
	Restriction exactness holds globally exactly when it holds in each component, and the least extension is the disjoint union of the component
	least relations.
	This does not require all other exact global extensions to split by components: absent cross-component local conditions can permit additional identifications or ordered pairs.
	
	An empty chart has its unique empty relation.
	On a singleton the local equivalence is the diagonal, whereas a transitive relation may be empty or diagonal.
	The test preserves that distinction.
	On \(E=\varnothing\) the unique global relation is empty.
	Exact restriction to equal interiors forces equality of their relations.
	After equality holds, the duplicate	argument applies.
\end{proof}

\begin{corollary}[Global quotient realization of locally specified refinements]
	\label{COR:PackagingKernelOverlapQuotientRealization}
	Under the kernel-extension criterion, the quotient \(q:E\twoheadrightarrow E/E^{\mathrm{gen}}\) restricts to a descriptor with kernel \(E_K\) on every chart.
	More precisely,
	\[
	b_K:I_K/E_K\xrightarrow{\cong}q(I_K),
	\qquad [x]_{E_K}\longmapsto[x]_{E^{\mathrm{gen}}},
	\]
	is a canonical bijection, and \((q|_{I_K})^{\mathrm{im}}=b_K\circ q_K\), where \(q_K:I_K\twoheadrightarrow I_K/E_K\) is the local quotient.
	Conversely, the kernels of restrictions of any one global descriptor satisfy the kernel-extension criterion.
	This concerns unlabeled refinement classes.
	It neither glues prescribed output labels literally nor asserts	admissibility of a global descriptor.
\end{corollary}

\begin{proof}
	The restriction has kernel \(E^{\mathrm{gen}}\cap I_K^2=E_K\).
	This equality makes \(b_K\) well defined and injective.
	Its attained-image codomain makes it surjective and proves the displayed factorization, including empty charts.
	Conversely, the kernel of a global descriptor is an exact extension of all its restricted kernels, so the kernel theorem
	applies.
	For example, two constant descriptors on one shared singleton have identical kernels even when their prescribed values are different.
	Kernel compatibility alone does not reconcile those labels.
\end{proof}

\begin{example}[Cycle-free geometry does not select compatible candidates]
	\label{EX:PackagingCycleFreeSelectionLimit}
	Take unit set \(E=\{0,1,2\}\), distinct carrier labels \(A,B\), and interiors \(I_A=\{0,1\}\) and \(I_B=E\).
	On each interior take all refinements as local candidates, but select only the indiscrete refinement on $I_A$ and only the discrete refinement on $I_B$.
	The larger chart guards every primal triangle, so the cover is cycle-free, and both selected candidate families are nonempty.
	Nevertheless their only possible selections disagree on $I_A\cap I_B=\{0,1\}$.
	No global equivalence relation can have these two exact restrictions.
	Hence the intrinsic completion space $\mathcal G(\mathfrak P)$ is empty.
	Admit the two carriers in the structural domain and prescribe a singleton formation datum whose only output is this target.
	This realizes a	packaging system without repairing the completion failure.
\end{example}

For selected equivalence candidates \(r_K\) that do agree on overlaps, cycle-freeness gives \(\Gamma(r)=r\) and the least extension \(G(r)\).
The fixed-point and extension assertions are instances of the intrinsic completion theory in Theorem~\ref{THM:PackagingIntrinsicRefinementCompletionA}.
Restrictions on selected refinements beyond overlap compatibility require separate hypotheses.
Cycle-free geometry alone neither selects local candidates nor imposes additional admissibility conditions.
The transitive theorem concerns binary relations, not arbitrary algebraic operations on package carriers.

For finite covers, the transitive classification is equivalent to $3<\nu(\mathcal H)$.
The comparison with closure and Armstrong globalization is developed in Section~\ref{SEC:SpectrumA}.
This extension theorem does not identify the global relation or its quotient presentation with the native formation \(P:\Xi\rightsquigarrow\mathfrak P\).

\section{Intrinsic recovery and the globalization obstruction spectrum}
\label{SEC:SpectrumA}

The previous sections give three sharp globalization results with different local proof methods.
Here we compare them using a common geometric obstruction rank.
The main theorem, Theorem~\ref{THM:PackagingGlobalizationSpectrum}, combines the closure, bounded-Armstrong, and transitive classifications into one spectrum.
We also define intrinsic recovery ranks and show which parts of the argument hold for general closure classes.
Theorem~\ref{THM:PackagingIntrinsicRecoverySufficiency} proves a sufficient recovery condition without transport assumptions.
Theorem~\ref{THM:PackagingAbstractRecoveryObstructionDuality} gives the converse when obstructions can be transported.
The counterexample in Theorem~\ref{THM:PackagingTransportNecessityCounterexample} explains why this extra assumption is necessary.

\subsection{Packaging-native recovery rank}
\label{SUBSEC:PackagingNativeRecoveryRank}

Closures and transitive relations use different kinds of elementary consequences.
To compare them, we define a common recovery structure that records assertions, their supports, how they restrict, and how they compose.
The associated recovery rank measures how much clique guarding is needed to recover a composed assertion inside an original chart.
Theorem~\ref{THM:PackagingTransitiveRecoveryRank} computes this rank for transitive relations.
Compatible atlases of these local states can be represented as realized packaging systems by Theorem~\ref{THM:PackagingAtlasRealizationBridge}; their native recovery signatures and elimination laws are additional specialization data, not properties supplied by the primitive packaging axioms.

\begin{definition}[Native recovery signature]
	\label{DEF:PackagingNativeRecoverySignature}	
	Let $\mathscr P=(\mathscr P(I))_I$ be a class of local packaged states indexed by finite label sets, equipped with restriction maps
	\[
	X\longmapsto X|_J\in\mathscr P(J)
	\qquad(J\subseteq I)
	\]
	satisfying
	\[
	X|_I=X,
	\qquad
	(X|_J)|_K=X|_K
	\quad(K\subseteq J\subseteq I).
	\]	
	A \textbf{native recovery signature} on	$\mathscr P$ consists of the following data and axioms:	
	\begin{enumerate}[label=(N\arabic*),leftmargin=3em]		
		\item \textbf{Supported native assertions.}
		
		For every finite label set $I$, there is a finite set $\mathsf A(I)$ of native assertions, equipped with a support map
		\[
		\operatorname{supp}:
		\mathsf A(I)\longrightarrow 2^I.
		\]		
		For every $J\subseteq I$, there is a canonical identification
		\[
		\mathsf A(J)
		\cong
		\{\alpha\in\mathsf A(I):
		\operatorname{supp}(\alpha)\subseteq J\}.
		\]
		These identifications preserve supports and are coherent under iterated restriction and enlargement.
		
		\item \textbf{Trace-local native core.}
		
		Every $X\in\mathscr P(I)$ has a canonically determined native core
		$
		\mathsf N(X)\subseteq\mathsf A(I).
		$		
		For every $J\subseteq I$ and every assertion $\alpha$ supported in $J$,
		\[
		\alpha\in\mathsf N(X)
		\quad\Longleftrightarrow\quad
		\alpha\in\mathsf N(X|_J),
		\]
		under the canonical identification.
		
		\item \textbf{Support-local one-pivot composition.}
		
		For every $v\in I$, there is a partially defined operation
		\[
		\star_v:
		\mathsf A(I)\times\mathsf A(I)
		\dashrightarrow
		\mathsf A(I\setminus\{v\}).
		\]		
		If $\alpha\star_v\beta$ is defined, then
		\[
		v\in
		\operatorname{supp}(\alpha)
		\cap\operatorname{supp}(\beta)
		\]
		and
		\[
		\operatorname{supp}(\alpha\star_v\beta)
		\subseteq
		\bigl(
		\operatorname{supp}(\alpha)
		\cup\operatorname{supp}(\beta)
		\bigr)\setminus\{v\}.
		\]		
		The operation is support-local:
		if $J\subseteq I$ contains $\operatorname{supp}(\alpha)
		\cup\operatorname{supp}(\beta)$, then composition is defined in ambient set $I$ if and only if it is defined in ambient set $J$, and the resulting assertions agree	under canonical identification.		
		We call such a pair \textbf{natively composable across $v$}.
		
		For $T\subseteq\mathsf A(I)$, put
		\[
		\operatorname{Comp}_v(T)
		=
		\{\alpha\star_v\beta:
		\alpha,\beta\in T,\;
		\alpha\star_v\beta
		\text{ is defined}\}.
		\]
		
		\item \textbf{Recovery dominance.}
		
		For every finite $I$, there is a binary relation $\trianglerighteq$ on $\mathsf A(I)$, compatible with the canonical identifications.		
		Thus the comparison of two assertions supported in $J\subseteq I$ is independent of whether it is made in $\mathsf A(J)$ or $\mathsf A(I)$.		
		The relation specifies which native assertions are admissible replacements for other assertions.
		
		\item \textbf{Local recovery soundness.}
		
		If $X\in\mathscr P(I)$, $\alpha,\beta\in\mathsf N(X)$, and $\alpha\star_v\beta$ is defined, then there exists
		$
		\gamma\in
		\mathsf N(X|_{I\setminus\{v\}})
		$
		such that
		\[
		\operatorname{supp}(\gamma)
		\subseteq
		\bigl(
		\operatorname{supp}(\alpha)
		\cup\operatorname{supp}(\beta)
		\bigr)\setminus\{v\}
		\]
		and
		\[
		\gamma\trianglerighteq
		\alpha\star_v\beta.
		\]
		
		\item \textbf{Separation.}
		
		For $X,Y\in\mathscr P(I)$,
		\[
		\mathsf N(X)=\mathsf N(Y)
		\quad\Longrightarrow\quad
		X=Y.
		\]		
	\end{enumerate}	
	A class equipped with a specified native recovery signature is called a \textbf{native-recovery packaging class}.
\end{definition}

\begin{definition}[Admissible native recovery]
	\label{DEF:PackagingAdmissibleNativeRecovery}	
	Let $\mathscr P$ be a native-recovery packaging class, let $v\in I$, and let $\alpha,\beta\in\mathsf A(I)$ be natively composable across $v$.	
	Put
	\[
	Q(\alpha,\beta)
	=
	\operatorname{supp}(\alpha)
	\cup\operatorname{supp}(\beta).
	\]	
	An assertion
	$
	\gamma\in\mathsf A(I\setminus\{v\})
	$
	is an \textbf{admissible $v$-recovery} of $(\alpha,\beta)$ if
	\[
	\operatorname{supp}(\gamma)
	\subseteq Q(\alpha,\beta)\setminus\{v\}
	\]
	and
	\[
	\gamma\trianglerighteq
	\alpha\star_v\beta.
	\]	
	Denote the set of admissible recoveries by
	\[
	\operatorname{Rec}_v(\alpha,\beta)
	=
	\left\{
	\gamma\in\mathsf A(I\setminus\{v\}):
	\begin{array}{l}
		\operatorname{supp}(\gamma)
		\subseteq Q(\alpha,\beta)\setminus\{v\},\\
		\gamma\trianglerighteq
		\alpha\star_v\beta
	\end{array}
	\right\}.
	\]	
	By local recovery soundness, whenever $\alpha,\beta\in\mathsf N(X)$ for some $X\in\mathscr P(I)$,
	\[
	\operatorname{Rec}_v(\alpha,\beta)
	\cap
	\mathsf N(X|_{I\setminus\{v\}})
	\ne\varnothing.
	\]
\end{definition}

\begin{definition}[Packaging-native recovery rank]
	\label{DEF:PackagingNativeRecoveryRank}
	Let $\mathscr P$ be a native-recovery packaging class with a specified native recovery signature.
	An atlas on a finite cover $\mathcal H$ of a finite label set $V$ is a family
	\[
	\mathcal A=(X_I)_{I\in\mathcal H},
	\qquad X_I\in\mathscr P(I).
	\]
	It is \textbf{compatible} if
	\[
	X_I|_{I\cap J}=X_J|_{I\cap J}
	\qquad (I,J\in\mathcal H).
	\]	
	Fix an integer $m\ge2$.
	We call $\mathscr P$ \textbf{packaging-$m$-recoverable} if the following condition holds.
	For every finite cover $\mathcal H$ on which every nonempty primal clique of cardinality at most $m$ is guarded, every compatible
	$\mathscr P$-atlas, all original charts $I,J\in\mathcal H$, and every pair
	\[
	\alpha\in\mathsf N(X_I),
	\qquad \beta\in\mathsf N(X_J),
	\]
	identify $\alpha$ and $\beta$ canonically in $\mathsf A(V)$.
	Whenever the two assertions are natively composable across a pivot $v$ and their combined support
	\[
	Q(\alpha,\beta)
	=\operatorname{supp}(\alpha)\cup
	\operatorname{supp}(\beta)
	\]
	is a primal clique, there must exist an original chart $K\in\mathcal H$ and an admissible recovery
	\[
	\gamma\in\operatorname{Rec}_v(\alpha,\beta)
	\]
	such that
	\[
	\operatorname{supp}(\gamma)\subseteq K,
	\qquad \gamma\in\mathsf N(X_K),
	\]
	under canonical identification.
	No chordality assumption is made.	
	The \textbf{packaging-native recovery rank} is
	\[
	\rho_{\mathrm{pkg}}(\mathscr P)
	=\min\bigl\{m\ge2:\mathscr P
	\text{ is packaging-$m$-recoverable}\bigr\},
	\]
	with value $\infty$ if there is no such $m$.
	The rank is relative to the specified native recovery signature.
\end{definition}

\begin{theorem}[Packaging-native recovery rank of transitive structures]
	\label{THM:PackagingTransitiveRecoveryRank}	
	For each finite label set $I$, let $\mathscr T(I)$ be the class of transitive relations $R\subseteq I^2$.
	Equip it with the native recovery signature
	\[
	\mathsf A(I)=\{[x\to y]:x,y\in I\},
	\qquad
	\operatorname{supp}([x\to y])=\{x,y\},
	\]
	\[
	\mathsf N(R)=\{[x\to y]:(x,y)\in R\},
	\qquad
	[x\to v]\star_v[v\to y]=[x\to y]
	\qquad
	(x,y\in I\setminus\{v\}),
	\]
	with no other one-pivot composites and with recovery dominance equal to equality.
	Then
	\begin{equation}\label{EQ:PackagingNativeRecoveryRank}
		\rho_{\mathrm{pkg}}(\mathscr T)=3.
	\end{equation}
	The same equality holds for the subclasses of preorders, partial orders, strict partial orders, and equivalence relations with their inherited native recovery signatures.
\end{theorem}

\begin{proof}
	We first verify the native recovery signature, then prove the upper and	lower rank bounds.	
	\begin{enumerate}
		\item[(a)] \emph{Verification of the native recovery signature.}
		
		Restrictions are ordinary relation restrictions:
		$
		R|_J=R\cap J^2.
		$
		They preserve transitivity and satisfy the restriction identities.		
		The assertion
		$
		[x\to y]
		$
		depends only on the ordered pair $(x,y)$ and has support $\{x,y\}$.
		Thus assertions are canonically identified under restriction whenever their supports remain present.		
		Trace locality follows from
		\[
		(x,y)\in R
		\quad\Longleftrightarrow\quad
		(x,y)\in R|_J
		\qquad(x,y\in J).
		\]
		
		The specified one-pivot composition is support-local and satisfies
		\[
		\operatorname{supp}([x\to y])
		\subseteq
		\{x,v,y\}\setminus\{v\}.
		\]
		
		Equality as recovery dominance is independent of the ambient label set, verifying axiom~\textup{(N4)}.
		
		To verify local recovery soundness, suppose
		\[
		[x\to v],[v\to y]\in\mathsf N(R).
		\]
		Then
		$
		(x,v) \in R,
		\quad (v,y) \in R.
		$
		Transitivity gives
		$
		(x,y)\in R,
		$
		and hence
		$
		[x\to y]
		\in\mathsf N(R|_{I\setminus\{v\}}).
		$
		
		Since recovery dominance is equality, $[x\to y]$ is an admissible recovery.		
		This also covers $x=y$: transitivity then gives the diagonal pair $(x,x)$.
		No reflexive pair has been inserted by convention.		
		Finally,
		\[
		\mathsf N(R)
		=
		\{[x\to y]:(x,y)\in R\}
		\]
		determines $R$ uniquely, proving separation.		
		Thus all native recovery signature axioms hold.		
		The same verification applies to preorders, partial orders, strict partial orders, and equivalence relations, which are restriction-stable subclasses of transitive relations.
		
		\item[(b)] \emph{Upper bound: $\rho_{\mathrm{pkg}}\le3$.}
		
		Let $\mathcal H$ be a finite cover in which every primal clique of cardinality at most three is guarded, and let $(R_I)_{I\in\mathcal H}$
		be a compatible transitive atlas.		
		Suppose original-chart assertions
		\[
		[x\to v]\in\mathsf N(R_I),
		\qquad
		[v\to y]\in\mathsf N(R_J)
		\]
		are composable across $v$, and their combined support
		$
		Q=\{x,v,y\}
		$
		is a primal clique.		
		Since $|Q| \le 3$, there exists an original chart $K$ containing $Q$.		
		Compatibility on overlaps transfers	both original assertions to $R_K$:
		\[
		(x,v)\in R_K,
		\qquad
		(v,y)\in R_K.
		\]		
		By transitivity of $R_K$, we have
		$
		(x,y)\in R_K.
		$
		Therefore
		$
		\gamma=[x\to y]\in\mathsf N(R_K).
		$		
		Moreover,
		$
		\operatorname{supp}(\gamma)
		=\{x,y\}\subseteq Q\setminus\{v\},
		$
		and
		$
		\gamma=[x\to v]\star_v[v\to y].
		$		
		Thus $\gamma$ is an admissible original-chart recovery.		
		Hence
		\[
		\rho_{\mathrm{pkg}}(\mathscr T)\le3.
		\]		
		Since all four indicated subclasses	are transitive and inherit the same signature, the identical argument proves the upper bound three for each subclass.
		
		\item[(c)] \emph{Sharpness: packaging-$2$-recoverability fails.}
		
		Let
		$
		V=\{x,v,y\}
		$
		and take the triangle-boundary cover
		$
		\mathcal H
		=
		\bigl\{
		\{x,v\},\{v,y\},\{x,y\}
		\bigr\}.
		$		
		Every primal clique of cardinality at most two is guarded, while the three-vertex clique is unguarded.
		
		First define local strict partial orders by
		\[
		R_{\{x,v\}}=\{(x,v)\},
		\qquad
		R_{\{v,y\}}=\{(v,y)\},
		\qquad
		R_{\{x,y\}}=\varnothing.
		\]		
		These relations are transitive and irreflexive.
		Their restrictions to singleton intersections are empty, so the atlas is pairwise compatible.		
		The two parent assertions
		\[
		[x\to v],\qquad[v\to y]
		\]
		are valid in original charts, and their combined support is the	primal clique $V$.		
		Their unique admissible recovery, since dominance is equality, is
		$
		[x\to y].
		$		
		But the only original chart containing $\{x,y\}$ is $\{x,y\}$, whose local relation is empty.		
		Thus no original chart contains the required recovery.
		Packaging-$2$-recoverability fails for transitive relations and strict partial orders.
		
		For preorders and partial orders, adjoin the diagonal
		$
		\Delta_F=\{(a,a):a\in F\}
		$
		to each local relation on its chart	$F$.		
		The resulting local relations are partial orders, hence also preorders.
		They remain compatible and still do not contain $(x,y)$ on the chart $\{x,y\}$.		
		Thus packaging-$2$-recoverability fails for both classes.
		
		Finally, for equivalence relations, take
		\[
		\begin{aligned}
			E_{\{x,v\}}
			&=
			\Delta_{\{x,v\}}
			\cup\{(x,v),(v,x)\},\\
			E_{\{v,y\}}
			&=
			\Delta_{\{v,y\}}
			\cup\{(v,y),(y,v)\},\\
			E_{\{x,y\}}
			&=
			\Delta_{\{x,y\}}.
		\end{aligned}
		\]		
		These local equivalence relations are pairwise compatible, since their restrictions to singleton intersections are identities.		
		Again,
		\[
		[x\to v]\star_v[v\to y]
		=
		[x\to y],
		\]
		but
		$
		(x,y)\notin E_{\{x,y\}}.
		$		
		Hence packaging-$2$-recoverability also fails for equivalence relations.		
		Therefore, in every one of the five relation classes,
		\[
		\rho_{\mathrm{pkg}}\ge3.
		\]		
	\end{enumerate}
	Combining the lower bound of (c) with part~\textup{(b)} proves Eq. \eqref{EQ:PackagingNativeRecoveryRank},
	$
	\rho_{\mathrm{pkg}} = 3,
	$
	for transitive relations, preorders, partial orders, strict partial orders, and equivalence relations.
\end{proof}

\subsection{Packaging-native recovery and elimination}
\label{SUBSEC:PackagingNativeRecoveryElimination}

Recovering one composed consequence does not automatically prove that an entire elimination process preserves local traces.
We state the additional conditions needed to repeat the recovery step safely.
The main theorem of this subsection then proves exact global extension under those conditions.
A separate obstruction theorem explains what is needed to obtain a converse.

\begin{definition}[Recovery-equipped packaging class]
	\label{DEF:PackagingRecoveryEquippedClass}
	A \textbf{recovery-equipped packaging class}
	\[
	\mathscr P=(\mathscr P(I))_I
	\]
	is a native-recovery packaging class together with an \textbf{ambient class}
	\[
	\widehat{\mathscr P}
	=(\widehat{\mathscr P}(I))_I,
	\qquad \mathscr P(I)\subseteq\widehat{\mathscr P}(I),
	\]
	and the following additional data and axioms.
	\begin{enumerate}[label=(E\arabic*),leftmargin=3em]
		\item \textbf{Ambient restrictions.}
		The ambient class admits restriction maps
		\[
		X\longmapsto X|_J\in\widehat{\mathscr P}(J)
		\qquad(J\subseteq I),
		\]
		extending the local restrictions, with
		\[
		X|_I=X,
		\qquad (X|_J)|_K=X|_K
		\quad(K\subseteq J\subseteq I).
		\]
		
		\item \textbf{Generated ambient states.}
		For every $T\subseteq\mathsf A(I)$ there is a generated	ambient state
		\[
		\langle T\rangle_I\in\widehat{\mathscr P}(I).
		\]
		The native core generates each local state:
		\[
		\langle\mathsf N(X)\rangle_I=X
		\qquad(X\in\mathscr P(I)).
		\]
		
		\item \textbf{Redundancy of dominance.}
		If $\gamma\in T\subseteq\mathsf A(I)$ and $\gamma\trianglerighteq\delta$, then
		\[
		\langle T\cup\{\delta\}\rangle_I
		=\langle T\rangle_I.
		\]
		
		\item \textbf{Exact one-pivot elimination.}
		For $v\in I$, put $J=I\setminus\{v\}$ and
		\[
		T[J]=\{\alpha\in T:
		\operatorname{supp}(\alpha)\subseteq J\}.
		\]
		Then
		\[
		\langle T\rangle_I|_J
		=\left\langle
		T[J]\cup\operatorname{Comp}_v(T)
		\right\rangle_J.
		\tag{E}
		\]
	\end{enumerate}
	The local class restricts the permitted original states.
	The ambient class contains generated states and need not satisfy the same local complexity bounds.
\end{definition}

\begin{remark}[Native meaning of the elimination axioms]
	\label{REM:PackagingNativeEliminationSpecializations}
	The formalism separates the ambient assertions, the native core belonging to each original chart, and the state generated by an arbitrary assertion family.
	For closures, assertions are implications, one-pivot composition is Horn resolution, and dominance is premise inclusion.
	In the bounded Armstrong/native class, recovery returns consequences with premise size at most $d-1$ to an original chart.
	A raw resolvent need not itself be a rooted circuit or satisfy the native premise bound.
	It belongs first to the ambient generated calculus.
	
	For transitive structures, assertions are ordered pairs, dominance is equality, and composition is $[x\to v]\star_v[v\to y]=[x\to y]$.
	Lemma~\ref{LEM:PackagingTransitiveOnePivotElimination} verifies exact one-pivot elimination.
	In particular, $[x\to v]\star_v[v\to x]=[x\to x]$ records a positive closed walk.
	It must not be discarded as though it were the tautological Horn rule $x\to x$.
	
	Thus the shared interface expresses exact elimination and native recovery without identifying the concrete assertion semantics or requiring the	generated global state to retain the original local complexity bound.
\end{remark}

\begin{definition}[Primitive-obstruction portability]
	\label{DEF:PackagingPrimitiveObstructionPortability}
	Let $m\ge3$ and let $\mathscr P$ be a recovery-equipped packaging class.
	We say that $\mathscr P$ is	\textbf{primitive-obstruction portable through rank $m$} if, for every finite cover $\mathcal H$ of a finite label set $V$, both of the following conditions hold.
	\begin{enumerate}[label=(P\arabic*),leftmargin=3em]
		\item \textbf{Induced-cycle portability.}
		If the primal graph of $\mathcal H$ is nonchordal, then $\mathcal H$ supports a compatible $\mathscr P$-atlas with no exact global extension in $\widehat{\mathscr P}(V)$.
		
		\item \textbf{Unguarded-clique portability.}
		If $\mathcal H$ has an inclusion-minimal unguarded primal clique $M$ with
		\[
		3\le|M|\le m,
		\]
		then $\mathcal H$ supports a compatible $\mathscr P$-atlas with no exact global extension in $\widehat{\mathscr P}(V)$.
	\end{enumerate}
\end{definition}

\begin{remark}[Role of obstruction portability]
	\label{REM:PackagingPortabilityInterface}
	Conditions \textup{(P1)}--\textup{(P2)} specify the two counterexample constructions needed for the converse globalization implication.
	Neither condition follows solely from the numerical recovery ranks \(\rho_{\mathrm{pkg}}\) or \(\operatorname{ipr}\).
	For the generic closure and Armstrong/native classes, portability is established through relabeling stability and neutral-cone transport of explicit obstructions.
\end{remark}

For a recovery-equipped packaging class \(\mathscr P\) and a finite cover \(\mathcal H\) of \(V\), write
$
\operatorname{UG}_{\mathscr P}(\mathcal H)
$
when every compatible \(\mathscr P\)-atlas \((X_I)_{I\in\mathcal H}\) has an exact global extension
\[
X\in\widehat{\mathscr P}(V),
\qquad
X|_I=X_I
\quad(I\in\mathcal H).
\]

\begin{theorem}[Packaging-native recovery sufficiency]
	\label{THM:PackagingNativeRecoverySufficiency}
	Let $\mathscr P$ be a recovery-equipped packaging class with
	\[
	p:=\rho_{\mathrm{pkg}}(\mathscr P)<\infty,
	\qquad p\ge3.
	\]
	For every finite cover $\mathcal H$ of a finite label set $V$,
	\[
	p<\nu(\mathcal H)
	\quad\Longrightarrow\quad
	\operatorname{UG}_{\mathscr P}(\mathcal H).
	\]
	Thus finite packaging-native recovery rank yields a sufficient condition for universal exact globalization,	without requiring primitive-obstruction portability.
\end{theorem}

\begin{proof}
	Let $\mathcal H$ be a finite cover of $V$ such that $p<\nu(\mathcal H)$, and fix a compatible $\mathscr P$-atlas $(X_I)_{I\in\mathcal H}$.
	
	\medskip
	\noindent\emph{Geometry and the generated state.}
	By Lemma~\ref{LEM:PackagingObstructionRankDictionary}, the primal graph $G_{\mathcal H}$ is chordal and every primal clique of size at most $p$ is guarded.
	Use the canonical identifications of assertions to put
	\[
	T=\bigcup_{I\in\mathcal H}\mathsf N(X_I)
	\subseteq\mathsf A(V),
	\qquad Y=\langle T\rangle_V
	\in\widehat{\mathscr P}(V).
	\]
	For $W\subseteq V$, write
	\[
	T[W]=\{\alpha\in T:
	\operatorname{supp}(\alpha)\subseteq W\}.
	\]
	
	\medskip
	\noindent\emph{Elimination invariant.}
	Fix an original chart $I_0\in\mathcal H$.
	Because $I_0$ is a primal clique and $G_{\mathcal H}$ is chordal, choose a perfect elimination ordering	of $V\setminus I_0$ that leaves $I_0$ last (Lemma~\ref{LEM:PackagingCliqueLastElimination}).
	We prove, for each successive remaining vertex set $W$, that
	\[
	Y|_W=\langle T[W]\rangle_W.
	\tag{1}
	\]
	This holds initially at $W=V$.
	
	Suppose (1) holds at $W$, and let $v\in W\setminus I_0$ be the next vertex to eliminate.
	It is simplicial in $G_{\mathcal H}[W]$.
	Put $J=W\setminus\{v\}$.
	Axiom~\textup{(E4)} gives
	\[
	Y|_J
	=\left\langle
	T[J]\cup\operatorname{Comp}_v(T[W])
	\right\rangle_J.
	\tag{2}
	\]
	
	Let $\delta=\alpha\star_v\beta$ be any composite occurring in (2), with $\alpha,\beta\in T[W]$.
	Each parent belongs to the native core of some original chart.
	Its support is therefore a clique of $G_{\mathcal H}$.
	Both parent supports contain $v$.
	Every other vertex in their union is thus a neighbor of $v$ in $G_{\mathcal H}[W]$.
	Since $v$ is simplicial, their combined support
	\[
	Q=\operatorname{supp}(\alpha)
	\cup\operatorname{supp}(\beta)
	\]
	is a primal clique.
	
	Packaging-$p$-recoverability supplies an original chart $K$ and an admissible recovery $\gamma\in\mathsf N(X_K)$ such that
	\[
	\operatorname{supp}(\gamma)
	\subseteq Q\setminus\{v\}\subseteq J,
	\qquad \gamma\trianglerighteq\delta.
	\]
	Therefore $\gamma\in T[J]$.
	Axiom~\textup{(E3)} makes $\delta$ redundant over $T[J]$.
	Since $\mathsf A(W)$ is finite,	every composite in (2) can be deleted in this way.
	Consequently,
	\[
	Y|_J=\langle T[J]\rangle_J.
	\]
	This proves (1) by induction.
	
	\medskip
	\noindent\emph{Exactness on every original chart.}
	At the end of the elimination, $W=I_0$, so
	\[
	Y|_{I_0}=\langle T[I_0]\rangle_{I_0}.
	\]
	We claim that
	\[
	T[I_0]=\mathsf N(X_{I_0}).
	\]
	One inclusion is immediate.
	For the other, suppose $\alpha\in\mathsf N(X_I)$ for an original chart $I$ and $\operatorname{supp}(\alpha)\subseteq I_0$.
	The support lies in $I\cap I_0$.
	Compatibility and trace locality, axiom~\textup{(N2)}, therefore transfer $\alpha$ to $\mathsf N(X_{I_0})$.
	The claim follows.
	
	By the native-core generation axiom~\textup{(E2)},
	\[
	Y|_{I_0}
	=\langle\mathsf N(X_{I_0})\rangle_{I_0}
	=X_{I_0}.
	\]
	Since $I_0$ was arbitrary, $Y$ is an exact global extension of the atlas in $\widehat{\mathscr P}(V)$.
\end{proof}

\begin{theorem}[Packaging-native recovery--obstruction theorem]
	\label{THM:PackagingNativeRecoveryObstruction}
	Let $\mathscr P$ be a recovery-equipped packaging class with
	\[
	p:=\rho_{\mathrm{pkg}}(\mathscr P)<\infty,
	\qquad p\ge3.
	\]
	Assume $\mathscr P$ is primitive-obstruction portable through rank $p$.
	Then, for every finite cover $\mathcal H$ of a finite label set $V$,
	\[
	\boxed{
		\operatorname{UG}_{\mathscr P}(\mathcal H)
		\quad\Longleftrightarrow\quad
		p<\nu(\mathcal H).
	}
	\]
	Thus universal exact globalization is governed by the packaging-native recovery rank of the local class and the obstruction rank of the cover.
\end{theorem}

\begin{proof}
	Fix a finite cover $\mathcal H$ of $V$.
	We prove both implications.
	
	\medskip
	\noindent\emph{Sufficiency: $p<\nu(\mathcal H)$.}
	This is exactly Theorem~\ref{THM:PackagingNativeRecoverySufficiency}.
	No portability assumption is needed in this direction.
	
	\medskip
	\noindent\emph{Necessity (by contraposition): $\nu(\mathcal H)\le p$.}
	If the primal graph is nonchordal, portability axiom~\textup{(P1)} gives a compatible $\mathscr P$-atlas with no exact global extension
	in $\widehat{\mathscr P}(V)$.
	
	Otherwise the primal graph is chordal.
	By the definition of $\nu(\mathcal H)$, there is an unguarded primal clique $K$ with $|K|\le p$.
	Choose an inclusion-minimal unguarded subclique	$M\subseteq K$.
	Singletons are guarded because the charts cover $V$, and primal edges are guarded by their definition.
	Hence
	\[
	3\le |M|\le p.
	\]
	Portability axiom~\textup{(P2)} again gives a compatible nonextendible $\mathscr P$-atlas.
	
	Therefore $\nu(\mathcal H)\le p$ implies $\neg\operatorname{UG}_{\mathscr P}(\mathcal H)$.
	This is the contrapositive of necessity and	completes the equivalence.
\end{proof}

\subsection{Intrinsic recovery rank}
\label{SUBSEC:PackagingAbstractNecessity}

Earlier, recovery ranks depended on a chosen set of local rules.
We now use \emph{rooted circuits}---the implications with minimal premises---to define a rank that does not depend on that choice.
The main characterization shows how these intrinsic consequences control original-chart recovery.
This rank becomes the parameter in the abstract globalization results below.

\begin{definition}[Rooted circuits and intrinsic propagation-recovery rank]
	\label{DEF:PackagingIntrinsicRecoveryRank}
	\begin{enumerate}[label=(\roman*),leftmargin=2.5em]		
		\item \textbf{Rooted circuits.}		
		Let $c$ be an empty-preserving closure operator on a finite label set $I$.		
		A pair \( (A,b) \) with \(\varnothing\ne A\subseteq I\setminus\{b\}\) is a \textbf{rooted circuit} of $c$ if
		$
		b\in c(A)
		$
		but
		$
		b \notin c(A')
		\quad
		\text{for every }A'\subsetneq A.
		$		
		The set $A$ is called the \textbf{stem}, and $b$ the \textbf{root}.
		
		Write
		$
		A\rightsquigarrow_c b
		$
		for this relation, and define the \textbf{rooted-circuit rule family} by
		\[
		\mathscr R^{\mathrm{rc}}(c)
		=
		\{A\to b:A\rightsquigarrow_c b\}.
		\]
		
		\item \textbf{Intrinsic propagation-recovery rank.}
		
		Let
		$
		\mathscr C=(\mathscr C(I))_I
		$
		be a restriction-stable family of finite empty-preserving closures.		
		Fix an integer $m\ge2$.		
		We say that $\mathscr C$ is	\textbf{intrinsically $m$-recoverable} if the following property holds.
		
		For every finite cover $\mathcal H$ of a finite label set $V$ in which every nonempty primal clique of cardinality at most $m$ is guarded, and for every compatible atlas
		$
		(c_I)_{I\in\mathcal H},
		\quad
		c_I\in\mathscr C(I),
		$
		consider any two rooted circuits in original charts $I,J\in\mathcal H$:
		\[
		A\rightsquigarrow_{c_I}v,
		\qquad
		B\cup\{v\}\rightsquigarrow_{c_J}b,
		\]
		where
		$
		v\ne b,
		\quad
		A,B\subseteq V\setminus\{v,b\}.
		$
		
		Suppose their combined support
		$
		Q=A\cup B\cup\{v,b\}
		$
		is a primal clique.		
		Then there exist an original chart $K\in\mathcal H$ and a nonempty set
		$
		D\subseteq A\cup B
		$
		such that
		$
		D\cup\{b\}\subseteq K
		$
		and
		$
		D\rightsquigarrow_{c_K}b.
		$		
		No chordality assumption is imposed on the primal graph.
		
		The \textbf{intrinsic propagation-recovery rank} of	$\mathscr C$ is
		\[
		\operatorname{ipr}(\mathscr C)
		=
		\min\bigl\{
		m\ge2:
		\mathscr C
		\text{ is intrinsically }
		m\text{-recoverable}
		\bigr\},
		\]
		with value $\infty$ if no such finite $m$ exists.		
	\end{enumerate}
\end{definition}

For comparison with admissible implicational presentations, write \(\operatorname{Adm}(\mathscr C)\) for the family of chart-local rule schemes
for \(\mathscr C\) satisfying local generation \textup{(R1)} and trace locality \textup{(R2)} of Definition~\ref{DEF:PackagingPropagationRecoveryScheme}.

\begin{proposition}[Rooted-circuit closure classes satisfy the recovery-equipped axioms]
	\label{PROP:PackagingClosureRecoveryEquipped}	
	Let $\mathscr C=(\mathscr C(I))_I$ be a	restriction-stable family of finite	empty-preserving closures.	
	For each finite label set $I$, let $\mathsf A(I)$ be the set of all	syntactically nontrivial implications
	\[
	A\to b,
	\qquad
	\varnothing\ne A\subseteq I\setminus\{b\},
	\]
	with support
	\[
	\operatorname{supp}(A\to b)=A\cup\{b\}.
	\]	
	For $c\in\mathscr C(I)$, define its	native core by
	$
	\mathsf N(c)=\mathscr R^{\mathrm{rc}}(c).
	$
		
	Define native one-pivot composition by
	\[
	(A\to v)\star_v
	(B\cup\{v\}\to b)
	=
	A\cup B\to b
	\]
	whenever
	$
	v,b\notin A\cup B
	$
	and the displayed resolvent is nontrivial.
	No other composites are defined.

	Define recovery dominance by
	\[
	D\to b\trianglerighteq C\to b
	\quad\Longleftrightarrow\quad
	D\subseteq C.
	\]
	
	Let $\widehat{\mathscr C}(I)$ be the class of all empty-preserving closures on $I$, and for	$T\subseteq\mathsf A(I)$ let
	$
	\langle T\rangle_I
	$
	be the closure generated by $T$.	
	Then the following hold:	
	\begin{enumerate}[label=(\roman*),leftmargin=2.5em]
		\item \textbf{Native recovery signature.}
		
		The assertion family, support map, rooted-circuit core, one-pivot composition, and recovery dominance satisfy axioms \textup{(N1)}--\textup{(N6)} of Definition~\ref{DEF:PackagingNativeRecoverySignature}.		
		In particular, every composable pair of rooted circuits in one local closure admits an admissible rooted-circuit recovery in that closure.
		
		\item \textbf{Ambient generation and recovery redundancy.}
		
		The ambient family $\widehat{\mathscr C}$ and the generated-state operation	$T\mapsto\langle T\rangle_I$ satisfy axioms \textup{(E1)}--\textup{(E3)} of Definition~\ref{DEF:PackagingRecoveryEquippedClass}.		
		In particular,
		\[
		\langle\mathsf N(c)\rangle_I=c
		\qquad(c\in\mathscr C(I)),
		\]
		and whenever
		\[
		\gamma\in T,
		\qquad
		\gamma\trianglerighteq\delta,
		\]
		one has
		\[
		\langle T\cup\{\delta\}\rangle_I
		=
		\langle T\rangle_I.
		\]
		
		\item \textbf{Exact one-pivot elimination.}
		
		For every finite $I$, every $T\subseteq\mathsf A(I)$, and every	$v\in I$, putting
		\[
		J=I\setminus\{v\},
		\qquad
		T[J]
		=
		\{\alpha\in T:
		\operatorname{supp}(\alpha)\subseteq J\},
		\]
		gives
		\[
		\langle T\rangle_I|_J
		=
		\left\langle
		T[J]\cup\operatorname{Comp}_v(T)
		\right\rangle_J.
		\]
		Thus axiom \textup{(E4)} holds.		
	\end{enumerate}	
	Consequently, $\mathscr C$, equipped with the specified native signature, ambient class, and generation operation, is a \textbf{recovery-equipped packaging class} in the sense of Definition~\ref{DEF:PackagingRecoveryEquippedClass}.
\end{proposition}

\begin{proof}
	We verify the three assertions.
	\begin{enumerate}[label=(\roman*),leftmargin=2.5em]
		\item \emph{Native recovery signature.}
		
		We check axioms \textup{(N1)}--\textup{(N6)}.
		\begin{enumerate}
			\item[(N1):] \emph{Supported assertions.}
			Every assertion has a uniquely determined finite support $A\cup\{b\}$.
			For $J\subseteq I$, the assertions supported in $J$ are exactly	$\mathsf A(J)$, under the canonical	identification.
			These identifications are coherent under iterated restriction.
			
			\item[(N2):] \emph{Trace locality.}
			Let $J\subseteq I$ and suppose $A\cup\{b\}\subseteq J$.
			For every $A'\subseteq A$,
			\[
			b\in(c|_J)(A')
			\quad\Longleftrightarrow\quad
			b\in c(A').
			\]		
			Consequently,
			\[
			A\rightsquigarrow_c b
			\quad\Longleftrightarrow\quad
			A\rightsquigarrow_{c|_J}b.
			\]
			Thus the rooted-circuit core is trace-local.
			
			\item[(N3):] \emph{One-pivot composition.}
			Whenever
			\[
			(A\to v)\star_v
			(B\cup\{v\}\to b)
			\]
			is defined, the pivot $v$ belongs to both parent supports.		
			The resolvent has support
			\[
			(A\cup B)\cup\{b\}
			\subseteq
			\bigl(
			\operatorname{supp}(A\to v)
			\cup
			\operatorname{supp}(B\cup\{v\}\to b)
			\bigr)\setminus\{v\}.
			\]
			
			The definition depends only on the parent implications, not on their ambient label set.
			Hence composition is support-local.
			
			\item[(N4):] \emph{Recovery dominance.}
			The condition
			\[
			D\to b\trianglerighteq C\to b
			\quad\Longleftrightarrow\quad
			D\subseteq C
			\]
			is determined entirely by the two implications.
			It is therefore preserved under canonical restriction and enlargement.
			
			\item[(N5):] \emph{Local recovery soundness.}
			Suppose
			\[
			A\rightsquigarrow_c v,
			\qquad
			B\cup\{v\}\rightsquigarrow_c b
			\]
			are composable rooted circuits.		
			Since
			$
			v\in c(A)
			\subseteq c(A\cup B),
			$
			extensivity gives
			\[
			B\cup\{v\}\subseteq c(A\cup B).
			\]		
			By monotonicity and idempotence,
			\[
			\begin{aligned}
				b
				&\in c(B\cup\{v\})\\
				&\subseteq c(c(A\cup B))\\
				&=c(A\cup B).
			\end{aligned}
			\]		
			Since $I$ is finite, choose an inclusion-minimal subset
			$
			D\subseteq A\cup B
			$
			such that
			$
			b\in c(D).
			$		
			Empty preservation gives $D\ne\varnothing$.
			By minimality,
			$
			D \rightsquigarrow_c b.
			$		
			Moreover,
			\[
			D\to b
			\trianglerighteq
			A\cup B\to b
			\]
			because $D\subseteq A\cup B$.		
			The recovered assertion is supported in $I\setminus\{v\}$.
			By trace locality, it belongs to the native core of	$c|_{I\setminus\{v\}}$.		
			This proves local recovery soundness.
			
			\item[(N6):] \emph{Separation.}
			Every nontrivial consequence
			$
			b\in c(S)\setminus S
			$
			has, by finiteness, an inclusion-minimal determining subset $D\subseteq S$.		
			Since $c$ is empty-preserving, $D$ is nonempty, and
			$
			D\rightsquigarrow_c b.
			$		
			Hence the rooted circuits determine all consequences of $c$.
			Therefore
			\[
			\mathsf N(c)=\mathsf N(c')
			\quad\Longrightarrow\quad
			c=c'.
			\]		
		\end{enumerate}
		All six native recovery axioms hold.
		
		\item \emph{Ambient generation and recovery redundancy.}
		
		We verify axioms \textup{(E1)}--\textup{(E3)}.
		\begin{enumerate}
			\item[(E1):] \emph{Ambient restrictions.}
			The restriction of an empty-preserving closure is again an empty-preserving closure, and closure restrictions satisfy the required identities.		
			Thus $\widehat{\mathscr C}$ is restriction-stable and extends the restrictions on $\mathscr C$.
			
			\item[(E2):] \emph{Generation.}
			Every implication in $\mathsf A(I)$ has a nonempty premise.
			Therefore the closure generated by any $T\subseteq\mathsf A(I)$ is empty-preserving and belongs to $\widehat{\mathscr C}(I)$.		
			By the determining-subset argument in \textup{(N6)}, every consequence of $c$ follows from a rooted circuit belonging to $\mathsf N(c)$.		
			Conversely, every rooted circuit is valid in $c$.
			Hence
			$
			\langle\mathsf N(c)\rangle_I=c.
			$
			
			\item[(E3):] \emph{Recovery redundancy.}
			Suppose
			\[
			\gamma=D\to b\in T,
			\qquad
			\delta=C\to b,
			\qquad
			D\subseteq C.
			\]		
			Since $\gamma$ belongs to $T$, we have
			$
			b\in\langle T\rangle_I(D).
			$		
			By monotonicity, we have
			$
			b\in\langle T\rangle_I(C).
			$		
			Thus $\delta$ is already valid in the closure generated by $T$.
			Adding it changes no consequence:
			$
			\langle T\cup\{\delta\}\rangle_I
			=
			\langle T\rangle_I.
			$		
		\end{enumerate}		
		This proves \textup{(E1)}--\textup{(E3)}.
		
		\item \emph{Exact one-pivot elimination.}		
		Fix $v\in I$ and put
		$
		J=I\setminus\{v\}.
		$		
		The operation $\star_v$ is exactly definite-Horn one-pivot resolution:
		its composites are the nontrivial resolvents
		$
		A\cup B\to b
		$
		formed from
		\[
		A\to v,
		\qquad
		B\cup\{v\}\to b.
		\]		
		Consequently, $\operatorname{Comp}_v(T)$ coincides with the resolvent family $\operatorname{Res}_v(T)$ of Lemma~\ref{LEM:PackagingNativeOnePivotHorn}.		
		That lemma gives
		$
		\langle T\rangle_I|_J
		=
		\left\langle
		T[J]\cup\operatorname{Comp}_v(T)
		\right\rangle_J.
		$	
		This is precisely axiom \textup{(E4)}.
	\end{enumerate}	
	Combining \textup{(i)}--\textup{(iii)},	all axioms of Definition~\ref{DEF:PackagingRecoveryEquippedClass} are satisfied.
\end{proof}

\begin{proposition}[Closure specialization of packaging-native recovery rank]
	\label{PROP:PackagingNativeRankClosureSpecialization}
	Let \(\mathscr C\) be a restriction-stable family of finite empty-preserving closures, equipped with the rooted-circuit recovery
	signature of Proposition~\ref{PROP:PackagingClosureRecoveryEquipped}.
	Then
	\[
	\rho_{\mathrm{pkg}}(\mathscr C)=\operatorname{ipr}(\mathscr C).
	\]
\end{proposition}

\begin{proof}
	Fix \(m\ge2\).
	In this signature, native assertions range over	syntactically nontrivial implications, while the native core of each closure consists precisely of its rooted circuits.
	Two core assertions are composable exactly when they have the form
	\[
	A\rightsquigarrow_{c_I}v,
	\qquad B\cup\{v\}\rightsquigarrow_{c_J}b
	\]
	with \(v,b\notin A\cup B\).
	Their composite is \(A\cup B\to b\).
	By the definition of recovery dominance, an admissible original-chart recovery is exactly a rooted circuit \(D\rightsquigarrow_{c_K}b\) with
	\[
	D\subseteq A\cup B,
	\qquad D\cup\{b\}\subseteq K.
	\]
	Hence packaging-\(m\)-recoverability and intrinsic \(m\)-recoverability are equivalent for every \(m\ge2\).
	Taking their minima, including the value \(\infty\), proves the claim.
\end{proof}

\begin{theorem}[Canonical rooted-circuit characterization of intrinsic recovery]
	\label{THM:PackagingRootedCircuitCharacterization}
	Let $\mathscr C$ be a restriction-stable family of finite empty-preserving closures.
	Then:
	\begin{enumerate}[label=(\roman*),leftmargin=2.5em]
		\item the rooted-circuit scheme $\mathscr R^{\mathrm{rc}}$ is an admissible chart-local generating scheme;
		
		\item it is unavoidable among admissible trace-local generating schemes:
		for every $\mathscr R\in\operatorname{Adm}(\mathscr C)$ and every $c \in \mathscr C(I)$,
		\[
		{
			\mathscr R^{\mathrm{rc}}(c)\subseteq\mathscr R(c);
		}
		\]
		
		\item the presentation-minimized rank is exactly the directly defined rooted-circuit rank:
		\[
		{
			\operatorname{ipr}(\mathscr C)
			=
			\operatorname{prk}(\mathscr C,\mathscr R^{\mathrm{rc}})
			=
			\min_{\mathscr R\in\operatorname{Adm}(\mathscr C)}
			\operatorname{prk}(\mathscr C,\mathscr R).
		}
		\]
	\end{enumerate}
	Thus the minimization over implication presentations is not part of the	intrinsic invariant.
	It is a derived variational characterization of the rooted-circuit recovery rank.
\end{theorem}

\begin{proof}
	\begin{enumerate}[label=(\roman*),leftmargin=2.5em]
		\item \emph{Admissibility of the rooted-circuit scheme.}
		Every rule \(A\to b\in\mathscr R^{\mathrm{rc}}(c)\) is nontrivial and valid in \(c\).
		Conversely, if \(b\in c(S)\setminus S\), finiteness and empty preservation yield an inclusion-minimal nonempty \(A\subseteq S\) with
		\(b\in c(A)\).
		Thus \(A\rightsquigarrow_c b\), and the	rooted-circuit rules generate \(c\), proving \textup{(R1)}.
		
		For \(J\subseteq I\) with \(A\cup\{b\}\subseteq J\),
		\[
		b\in(c|_J)(A')\quad\Longleftrightarrow\quad
		b\in c(A')\qquad(A'\subseteq A).
		\]
		Hence \(A\rightsquigarrow_c b\) if and only if \(A\rightsquigarrow_{c|_J}b\), proving \textup{(R2)}.
		Therefore \(\mathscr R^{\mathrm{rc}} \in \operatorname{Adm}(\mathscr C)\).
		
		\item \emph{Unavoidability under trace locality.}
		Let \(\mathscr R\in\operatorname{Adm}(\mathscr C)\), let \(c\in\mathscr C(I)\), and fix a rooted circuit \(A\rightsquigarrow_c b\).
		Put \(J=A\cup\{b\}\) and \(c_J=c|_J\in\mathscr C(J)\).
		
		By \textup{(R1)}, the scheme \(\mathscr R(c_J)\) generates \(c_J\).
		Since \(b\in c_J(A)\), there is a finite forward derivation of \(b\) from \(A\) using these rules.
		The initial set \(A\) contains every label of \(J\) except \(b\).
		Hence the first rule that adds a new label has the form	\(E\to b\in\mathscr R(c_J)\) with \(E\subseteq A\).
		This rule is valid in \(c_J\), so \(b\in c(E)\).
		Rooted minimality forces \(E=A\).
		Thus \(A\to b\in\mathscr R(c_J)\).
		By trace locality \textup{(R2)}, \(A\to b\in\mathscr R(c)\).
		Since the circuit was arbitrary,
		\[
		\mathscr R^{\mathrm{rc}}(c)\subseteq\mathscr R(c).
		\]
		
		\item \emph{Equality of the presentation-minimized ranks.}
		Fix \(m\ge2\).
		By definition, intrinsic \(m\)-recoverability is precisely local \(m\)-propagation-recoverability for the rooted-circuit scheme.
		Consequently,
		\[
		\operatorname{ipr}(\mathscr C)
		=\operatorname{prk}(\mathscr C,\mathscr R^{\mathrm{rc}}).
		\]
		
		Suppose now that some \(\mathscr R\in\operatorname{Adm}(\mathscr C)\) is locally \(m\)-propagation-recoverable.
		Take any two original-chart	rooted circuits
		\[
		A\rightsquigarrow_{c_I}v,
		\qquad B\cup\{v\}\rightsquigarrow_{c_J}b
		\]
		whose combined support is a primal clique in a compatible atlas on an \(m\)-guarded cover.
		By part~\textup{(ii)}, both parent rules belong to \(\mathscr R(c_I)\) and \(\mathscr R(c_J)\), respectively.
		Recovery for \(\mathscr R\) supplies an original chart \(K\) and a rule \(E\to b\in\mathscr R(c_K)\) with
		\[
		E\subseteq A\cup B,
		\qquad E\cup\{b\}\subseteq K.
		\]
		Because the recovered rule is valid and \(c_K\) is finite and empty-preserving, there exists an inclusion-minimal nonempty \(D\subseteq E\) with \(b\in c_K(D)\).
		Then \(D\rightsquigarrow_{c_K}b\) is the required intrinsic recovery.
		Therefore every guarding scale admissible for \(\mathscr R\) is admissible for the rooted-circuit scheme, and
		\[
		\operatorname{ipr}(\mathscr C)
		\le\operatorname{prk}(\mathscr C,\mathscr R).
		\]
		On the other hand, part~\textup{(i)} places the	rooted-circuit scheme in \(\operatorname{Adm}(\mathscr C)\).
		Taking the minimum over that nonempty family proves
		\[
		\operatorname{ipr}(\mathscr C)
		=\operatorname{prk}(\mathscr C,\mathscr R^{\mathrm{rc}})
		=\min_{\mathscr R\in\operatorname{Adm}(\mathscr C)}
		\operatorname{prk}(\mathscr C,\mathscr R),
		\]
		including when the common value is \(\infty\).
	\end{enumerate}
\end{proof}

\begin{remark}[Relation to classical rooted-circuit theory]
	\label{REM:PackagingRootedCircuitNoveltyScope}
	Rooted circuits and their relation to implicational representations of finite closure systems are classical
	\cite{NakamuraKashiwabara2013RootedCircuits,CaspardMonjardet2003ClosureSystems}.
	Likewise, least generated closures and obstruction tests for a fixed closure atlas belong to established closure-system and recent closure-atlas theory \cite{Kim2026ClosureAtlases}.
	
	The invariant introduced here is different:
	it asks for the least guarding scale at which every elementary composition of local rooted circuits with provenance can be recovered inside an original chart, uniformly over all compatible atlases in the class.
	Theorem~\ref{THM:PackagingRootedCircuitCharacterization} shows that this scale is invariant under admissible chart-local, trace-local generating rule schemes.
	It does not assert that arbitrary implicational bases contain every rooted circuit.
	The subsequent duality theorem then compares this intrinsic local recovery scale with the independently defined cover obstruction rank.
\end{remark}

\subsection{Provenance and intrinsic recovery}
\label{SUBSEC:PackagingProvenanceRecovery}

A global closure records which implications hold, but not necessarily which chart supplied each implication.
That missing information matters when we require recovery in an original chart.
The first theorem shows that generated closures and guard geometry can lose this provenance.
The next theorem gives a finite critical-support description of intrinsic recovery rank that retains it.

\begin{definition}[Guard complex and local circuit provenance]
	\label{DEF:PackagingCircuitProvenanceShadow}
	Let $\mathcal A=(c_I)_{I\in\mathcal H}$ be a compatible closure atlas on a finite label set $V$.
	Its \textbf{guard complex} is the downward-closed family
	\[
	\mathsf G_{\mathcal A}
	=\{S\subseteq V:\exists I\in\mathcal H,\ S\subseteq I\}.
	\]
	Its \textbf{local rooted-circuit provenance} is
	\[
	\mathsf L_{\mathcal A}
	=\{(A,b):\exists I\in\mathcal H,\
	A\rightsquigarrow_{c_I}b\}.
	\]
	Let $C_{\mathcal A}$ be the closure generated on $V$ by the implications $A\to b$ for $(A,b)\in\mathsf L_{\mathcal A}$.
	The \textbf{provenance-marked circuit shadow} is
	\[
	\mathsf P(\mathcal A)
	=(C_{\mathcal A},\mathsf G_{\mathcal A},
	\mathsf L_{\mathcal A}).
	\]
	Thus $C_{\mathcal A}$ is determined by $\mathsf L_{\mathcal A}$.
	It is displayed separately to distinguish global derivability from original-chart provenance.
	
	A set $Q\subseteq V$ is a \textbf{primal clique} of the shadow if every two-element subset of $Q$ lies in $\mathsf G_{\mathcal A}$.
\end{definition}

\begin{theorem}[Provenance is necessary for local recovery]
	\label{THM:PackagingBareClosureRecoveryNoGo}
	The existence of an original-chart rooted-circuit recovery for a specified pair of local parent circuits is not determined by the generated global closure together with the complete guard complex.
	More precisely, there exist compatible atlases $\mathcal A_0$ and $\mathcal A_1$ on the same three-label edge cover such that
	\[
	C_{\mathcal A_0}=C_{\mathcal A_1},\qquad
	\mathsf G_{\mathcal A_0}=\mathsf G_{\mathcal A_1},
	\]
	and the same two parent rooted circuits occur locally in both, but their rooted recovery occurs in an original chart only in $\mathcal A_1$.
	Both atlases consist of finite empty-preserving	closures generated by unary implications.
\end{theorem}

\begin{proof}
	Put $V=\{a,v,b\}$ and take the triangle-boundary cover
	\[
	\mathcal H=\bigl\{\{a,v\},\{v,b\},\{a,b\}\bigr\}.
	\]
	On $\{a,v\}$ take the closure generated by $a\to v$, and on $\{v,b\}$ take the closure generated by $v\to b$, in both atlases.
	On $\{a,b\}$ take the identity closure in $\mathcal A_0$ and the closure generated by $a\to b$ in $\mathcal A_1$.
	
	Every local closure is empty-preserving, so its restriction to any singleton is the identity.
	As distinct edge charts meet in	singletons, both atlases are compatible.
	They also have the same guard complex, namely the subsets of $V$ of size at most two.
	
	In either atlas, the global closure generated by local rooted circuits is the closure generated by
	\[
	a\to v,\qquad v\to b.
	\]
	Indeed, $a\to b$ already follows from these two implications.
	Adjoining it in $\mathcal A_1$ therefore changes no global consequence.
	Hence $C_{\mathcal A_0}=C_{\mathcal A_1}$.
	
	Both atlases have the local rooted circuits $\{a\}\rightsquigarrow v$ and $\{v\}\rightsquigarrow b$.
	The only possible nonempty recovered stem contained in $\{a\}$ is $\{a\}$ itself.
	The circuit	$\{a\}\rightsquigarrow b$ occurs in the $\{a,b\}$ chart of $\mathcal A_1$, but not in that chart of $\mathcal A_0$.
	No other original chart contains its support.
	Thus original-chart recovery differs, despite equality of the generated closures and guard complexes.
\end{proof}

\begin{definition}[$\mathscr C$-provenance circuit systems]
	\label{DEF:PackagingProvenanceCircuitSystem}
	Let $\mathscr C=(\mathscr C(I))_I$ be a restriction-stable family of finite empty-preserving closures.
	A \textbf{$\mathscr C$-provenance circuit system} is a triple $\mathfrak S=(V,\mathsf G,\mathsf L)$ satisfying:
	\begin{enumerate}[label=(P\arabic*),leftmargin=3em]
		\item $V$ is finite, $\mathsf G\subseteq2^V$ is downward closed, and $\bigcup\mathsf G=V$.
		
		\item $\mathsf L$ consists of pairs $(A,b)$ with $\varnothing\ne A\subseteq V\setminus\{b\}$ and $A\cup\{b\}\in\mathsf G$.
		
		\item For each inclusion-maximal face $M$ of $\mathsf G$, put
		\[
		\mathsf L_M=
		\{(A,b)\in\mathsf L:A\cup\{b\}\subseteq M\},
		\]
		and let $c_M$ be the closure on $M$ generated by $\{A\to b:(A,b)\in\mathsf L_M\}$.
		Then
		\[
		c_M\in\mathscr C(M),\qquad
		\{A\to b:(A,b)\in\mathsf L_M\}
		=\mathscr R^{\mathrm{rc}}(c_M).
		\]
	\end{enumerate}
	Denote the closure generated on $V$ by the marked implications by $C_{\mathfrak S}$.
	Condition \textup{(P3)} says that the marked implications supported in each maximal guarded face are exactly its rooted-circuit rule family.
\end{definition}

\begin{lemma}[Facet compatibility is encoded by global circuit provenance]
	\label{LEM:PackagingProvenanceCircuitCompatibility}
	Let $\mathscr C$ be restriction-stable.
	\begin{enumerate}[label=(\roman*),leftmargin=2.5em]
		\item For every $\mathscr C$-provenance circuit system $\mathfrak S=(V,\mathsf G,\mathsf L)$, its facet closures $(c_M)_{M\in\max\mathsf G}$ are compatible on intersections.
		
		\item Conversely, every finite compatible $\mathscr C$-atlas, after discarding repeated and nonmaximal charts, determines a $\mathscr C$-provenance circuit system having the same guard complex and the same locally realized rooted circuits.
	\end{enumerate}
\end{lemma}

\begin{proof}
	\emph{(i) Compatibility.}
	Let $M,N\in\max\mathsf G$, put $J=M\cap N$, and let $S\subseteq J$.
	We show $c_M(S)\cap J=c_N(S)\cap J$.
	Take $b\in c_M(S)\cap J$.
	If $b\in S$, then $b\in c_N(S)$	by extensivity.
	Otherwise, finiteness and empty preservation give an inclusion-minimal nonempty $A\subseteq S$ with
	$A\rightsquigarrow_{c_M}b$.
	By \textup{(P3)}, $(A,b)\in\mathsf L_M\subseteq\mathsf L$.
	Since $A\cup\{b\}\subseteq J\subseteq N$, it also lies in $\mathsf L_N$.
	Again by \textup{(P3)}, $A\rightsquigarrow_{c_N}b$, and hence $b\in c_N(S)$.
	The reverse inclusion is symmetric.
	
	\medskip
	\noindent\emph{(ii) Reconstruction from an atlas.}
	Start with a finite compatible atlas.
	Every nonmaximal chart is contained in a maximal chart, and compatibility makes its closure the restriction of the latter.
	Repeated chart supports carry identical closures, so discarding either kind changes neither the guard complex nor the locally realized rooted circuits.
	
	Let $\mathsf G$ be the downward closure of the remaining charts, and let $\mathsf L$ be the union of their rooted-circuit pairs.
	Conditions \textup{(P1)} and \textup{(P2)} are immediate.
	Fix a maximal chart $M$.
	Every rooted circuit of its closure is included in $\mathsf L_M$.
	Conversely, if $(A,b)\in\mathsf L_M$ originates in another chart $N$, its full support lies in $M\cap N$.
	Compatibility identifies the closure traces there.
	In particular, rootedness transfers to $M$.
	Hence $\mathsf L_M$ consists exactly of the rooted circuits of the original closure on $M$.
	As rooted circuits generate that closure, \textup{(P3)} holds.
\end{proof}

\begin{lemma}[Restriction stability of provenance circuit systems]
	\label{LEM:PackagingProvenanceSystemRestriction}
	Let $\mathfrak S=(V,\mathsf G,\mathsf L)$ be a $\mathscr C$-provenance circuit system, with $\mathscr C$ restriction-stable, and let $Q\subseteq V$.
	Put
	\[
	\mathsf G_Q=\mathsf G\cap2^Q,\qquad
	\mathsf L_Q=\{(A,b)\in\mathsf L:A\cup\{b\}\subseteq Q\}.
	\]
	Then $\mathfrak S|_Q=(Q,\mathsf G_Q,\mathsf L_Q)$ is again a $\mathscr C$-provenance circuit system.
\end{lemma}

\begin{proof}
	Conditions \textup{(P1)}--\textup{(P2)} pass to the restriction.
	Let $N$ be a maximal face of $\mathsf G_Q$.
	Because $\mathsf G$ is finite and downward closed, there is a maximal face $M$ of $\mathsf G$ containing $N$.
	Then $N\subseteq M\cap Q\in\mathsf G_Q$, so maximality yields $N=M\cap Q$.
	
	By \textup{(P3)}, the marked circuits supported in $N$ are precisely those of $c_M$ supported there.
	Rooted-circuit membership is trace-local:
	for
	$A\cup\{b\}\subseteq N$,
	\[
	A\rightsquigarrow_{c_M}b
	\quad\Longleftrightarrow\quad
	A\rightsquigarrow_{c_M|_N}b.
	\]
	Thus the marked implications on $N$ are exactly $\mathscr R^{\mathrm{rc}}(c_M|_N)$ and generate $c_M|_N$.
	Restriction stability gives	$c_M|_N\in\mathscr C(N)$.
	This verifies \textup{(P3)} for every maximal face $N$.
\end{proof}

\begin{definition}[Provenance-critical circuit rank]
	\label{DEF:PackagingProvenanceCriticalRank}
	Let $\mathscr C$ be restriction-stable and let $\mathfrak S=(V,\mathsf G,\mathsf L)$ be a $\mathscr C$-provenance circuit system.
	A pair of marked rooted circuits
	\[
	(A,v),\qquad (B\cup\{v\},b)\qquad\text{in }\mathsf L
	\]
	is \textbf{provenance-critical} if $v,b\notin A\cup B$, its combined support
	\[
	Q=A\cup B\cup\{v,b\}
	\]
	is a primal clique of the one-skeleton of $\mathsf G$, and there is no marked recovery
	\[
	(D,b)\in\mathsf L,
	\qquad\varnothing\ne D\subseteq A\cup B.
	\]
	
	For a primal clique $Q$ in a finite guard complex $\mathsf G$, define its \textbf{first unguarded support size} by
	\[
	\gamma_{\mathsf G}(Q)
	=\min\{|S|:\varnothing\ne S\subseteq Q,\ S\notin\mathsf G\},
	\]
	with $\min\varnothing=\infty$.
	Since $Q$ is a primal clique, all its subsets are primal cliques.
	
	Define the \textbf{provenance-critical circuit rank} of $\mathscr C$ by
	\[
	\boxed{
		\operatorname{pcr}(\mathscr C)
		=\sup\!\left(
		\{2\}\cup
		\left\{
		\gamma_{\mathsf G}(V):
		\begin{array}{l}
			(V,\mathsf G,\mathsf L)
			\text{ is a }\mathscr C
			\text{-provenance circuit system}\\
			\text{containing a provenance-critical pair}\\
			\text{whose combined support equals }V
		\end{array}
		\right\}
		\right).
	}
	\]
	The supremum is allowed to be infinite.
	The definition ranges over single-support provenance systems.
	It does not require a separate choice of a cover and a compatible atlas.
\end{definition}

\begin{lemma}[Critical supports are necessarily unguarded]
	\label{LEM:PackagingCriticalSupportUnguarded}
	Let $\mathfrak S=(V,\mathsf G,\mathsf L)$ be a $\mathscr C$-provenance circuit system.
	If a provenance-critical pair has combined support $Q=A\cup B\cup\{v,b\}$, then $Q\notin\mathsf G$.
	In particular,
	\[
	3\le\gamma_{\mathsf G}(Q)\le|Q|<\infty.
	\]
\end{lemma}

\begin{proof}
	Suppose $Q\in\mathsf G$ and choose a maximal face $M\supseteq Q$.
	By \textup{(P3)}, both marked parents are rooted circuits of $c_M$.
	Consequently,
	\[
	v\in c_M(A),\qquad b\in c_M(B\cup\{v\}),
	\qquad\text{so }b\in c_M(A\cup B).
	\]
	By finiteness and empty preservation, there is a nonempty inclusion-minimal $D\subseteq A\cup B$ with $D\rightsquigarrow_{c_M}b$.
	Condition \textup{(P3)} places $(D,b)$ in $\mathsf L$, contradicting criticality.
	Thus $Q\notin\mathsf G$.
	Every singleton of $Q$ is guarded by \textup{(P1)}, and every pair is guarded because $Q$ is a primal clique.
	Therefore the smallest unguarded subset of $Q$ has size at least three and at most $|Q|$.
\end{proof}

\begin{theorem}[Exact provenance-critical characterization of intrinsic recovery]
	\label{THM:PackagingProvenanceCriticalRank}
	For every restriction-stable family $\mathscr C$ of finite empty-preserving closures,
	\[
	\boxed{\operatorname{pcr}(\mathscr C)
		=\operatorname{ipr}(\mathscr C).}
	\]
	Thus the intrinsic rank can be expressed as the supremum of the	first unguarded support sizes of critical provenance systems normalized to their combined supports.
	The provenance marking cannot in general be replaced by the generated closure and guard complex alone, by
	Theorem~\ref{THM:PackagingBareClosureRecoveryNoGo}.
\end{theorem}

\begin{proof}
	We establish both inequalities, with the convention that either rank may be infinite.
	
	\medskip
	\noindent\emph{First inequality: $\operatorname{ipr}(\mathscr C)\ge	\operatorname{pcr}(\mathscr C)$.}
	Let $\mathfrak S=(Q,\mathsf G,\mathsf L)$ be any system	occurring in the defining supremum, and put	$r=\gamma_{\mathsf G}(Q)$.
	By Lemma~\ref{LEM:PackagingCriticalSupportUnguarded}, $3\le r\le|Q|$.
	
	Take the compatible $\mathscr C$-atlas on the maximal faces of $\mathsf G$ supplied by Lemma~\ref{LEM:PackagingProvenanceCircuitCompatibility}.
	Its guard complex is $\mathsf G$, so it guards every primal clique of size at most $r-1$.
	The marked parent circuits occur in original facet charts by \textup{(P3)}, and their combined support is the primal clique $Q$.
	Any original-chart rooted recovery would contribute a pair $(D,b)\in\mathsf L$ with $\varnothing\ne D\subseteq A\cup B$, contrary to criticality.
	
	The atlas therefore witnesses failure of intrinsic $(r-1)$-recoverability.
	As recoverability is monotone in the guarding scale, $\operatorname{ipr}(\mathscr C)\ge r$.
	Taking the supremum over all such systems, and noting $\operatorname{ipr}(\mathscr C)\ge2$, yields $\operatorname{ipr}(\mathscr C)\ge \operatorname{pcr}(\mathscr C)$.
	
	\medskip
	\noindent\emph{Second inequality: $\operatorname{pcr}(\mathscr C)\ge \operatorname{ipr}(\mathscr C)$.}
	Let $m\ge2$ be an integer with $m<\operatorname{ipr}(\mathscr C)$.
	By the definition of intrinsic rank, there is a finite cover guarding every primal clique of size at most $m$, and a compatible $\mathscr C$-atlas with parent rooted circuits
	\[
	A\rightsquigarrow v,\qquad
	B\cup\{v\}\rightsquigarrow b
	\]
	of primal-clique combined support $Q=A\cup B\cup\{v,b\}$ for which no original-chart rooted recovery with stem contained in $A\cup B$ exists.
	
	Restrict the atlas to $Q$.
	The restricted chart closures remain in $\mathscr C$ by restriction stability, and compatibility and the parent circuits survive by trace locality.
	The induced cover of $Q$ still guards every	clique of size at most $m$.
	By Lemma~\ref{LEM:PackagingProvenanceCircuitCompatibility}, its maximal charts determine a provenance system $\mathfrak S_Q=(Q,\mathsf G_Q,\mathsf L_Q)$.
	The two parent circuits lie in $\mathsf L_Q$.
	If $(D,b)\in\mathsf L_Q$ with $\varnothing\ne D\subseteq A\cup B$, then its source chart would provide a recovered rooted circuit in the original atlas by compatibility and restriction.
	Thus the pair remains provenance-critical.

	Since all primal cliques of size at most $m$ are guarded in $\mathsf G_Q$,
	\[
	\gamma_{\mathsf G_Q}(Q)>m.
	\]
	This system occurs in the defining supremum, so	$\operatorname{pcr}(\mathscr C)>m$.

	If $\operatorname{ipr}(\mathscr C)=2$, the desired inequality is immediate from the defining lower bound $\operatorname{pcr}(\mathscr C)\ge2$.
	If $\operatorname{ipr}(\mathscr C)$ is finite and at least three, take $m=\operatorname{ipr}(\mathscr C)-1$.
	If $\operatorname{ipr}(\mathscr C)=\infty$, let $m$	range over all integers $m\ge2$.
	Both cases give $\operatorname{pcr}(\mathscr C)\ge \operatorname{ipr}(\mathscr C)$.
	Combining the two inequalities proves the theorem.
\end{proof}

\subsection{Boundary obstructions and one-point extension calculus}
\label{SUBSEC:PackagingBoundaryObstructions}

A recovery bound gives a sufficient condition for globalization.
To prove a sharp converse, we also need compatible atlases that cannot extend.
Simplex boundaries provide basic examples.
We then show how to add labels while controlling the local complexity bound, so these examples can be used inside larger covers.

For \(r\ge3\), put \(V_r=\{1,\ldots,r\}\) and let
$
\mathcal B_r
=
\{V_r\setminus\{i\}:i\in V_r\}
$
be the simplex-boundary cover.
Its primal graph is \(K_r\).
Every proper subset of \(V_r\) is guarded, while \(V_r\) itself is not.

\begin{definition}[Boundary-obstruction rank]
	\label{DEF:PackagingBoundaryObstructionRank}
	Let $\mathscr C=(\mathscr C(I))_I$ be a family of classes of finite empty-preserving closures.
	An integer $r\ge3$ is \textbf{boundary-obstructed for $\mathscr C$} if the simplex-boundary cover $\mathcal B_r$ supports a compatible atlas
	\[
	(c_F)_{F\in\mathcal B_r},\qquad c_F\in\mathscr C(F),
	\]
	with no exact global closure on $V_r$.
	
	The \textbf{boundary-obstruction rank} is
	\[
	\operatorname{bor}(\mathscr C)
	:=\sup\left\{m\in\mathbb Z_{\ge2}:
	\begin{array}{l}
		\mathcal B_r\text{ is boundary-obstructed for }\mathscr C\\
		\text{for every integer }3\le r\le m
	\end{array}\right\},
	\]
	where the supremum may be $\infty$.
	
	Thus $\operatorname{bor}(\mathscr C)=2$ if $\mathcal B_3$ is not boundary-obstructed, and $\operatorname{bor}(\mathscr C)=\infty$ if every $\mathcal B_r$, $r\ge3$, is boundary-obstructed.
	This is an \emph{initial-segment} obstruction rank: it need not equal the largest individual obstructed boundary size.
	Since the local closures are empty-preserving and the facets cover $V_r$, any exact global closure would automatically be empty-preserving.
\end{definition}

\begin{definition}[Universal inert one-point extension]
	\label{DEF:PackagingUniversalInertExtension}
	Let $c$ be a finite empty-preserving closure on a finite set $T$, and let $q\notin T$.
	Define
	\[
	{
		\operatorname{Inert}_q(c)(S)
		=
		c(S\cap T)\cup(S\cap\{q\})
		\qquad(S\subseteq T\cup\{q\}).
	}
	\]
	Thus the new coordinate is present exactly when it is observed and neither creates nor receives any additional consequence.
\end{definition}

\begin{theorem}[Universal one-point extension calculus]
	\label{THM:PackagingUniversalOnePointExtensionCalculus}
	Let $T$ be finite, let $c$ be an empty-preserving closure on $T$, and let $q\notin T$.
	Then:
	\begin{enumerate}[label=(\roman*),leftmargin=2.5em]
		\item \textbf{Closure and empty preservation.}
		$\operatorname{Inert}_q(c)$ is an empty-preserving closure on $T\cup\{q\}$.
		
		\item \textbf{Restriction naturality.}
		\[
		\operatorname{Inert}_q(c)|_T=c,
		\qquad
		\operatorname{Inert}_q(c)|_{J\cup\{q\}}
		=\operatorname{Inert}_q(c|_J)
		\quad(J\subseteq T).
		\]
		
		\item \textbf{Carath\'eodory complexity.}
		If $T\ne\varnothing$, then
		\[
		\operatorname{car}(\operatorname{Inert}_q(c))
		=\operatorname{car}(c).
		\]
		
		\item \textbf{Armstrong complexity.}
		If $T\ne\varnothing$, then, viewing $c$ as a finite	functional-dependency closure,
		\[
		\operatorname{arm}(\operatorname{Inert}_q(c))
		=\operatorname{arm}(c)+1.
		\]
	\end{enumerate}
\end{theorem}

\begin{proof}
	\emph{(i) Closure and empty preservation.}
	For $S\subseteq T\cup\{q\}$, the defining formula
	\[
	\operatorname{Inert}_q(c)(S)
	=c(S\cap T)\cup(S\cap\{q\})
	\]
	gives extensivity and monotonicity. Applying it twice gives
	\[
	\begin{aligned}
		\operatorname{Inert}_q(c)^2(S)
		&=c(c(S\cap T))\cup(S\cap\{q\})\\
		&=\operatorname{Inert}_q(c)(S),
	\end{aligned}
	\]
	and $\operatorname{Inert}_q(c)(\varnothing)=\varnothing$.
	
	\medskip
	\noindent\emph{(ii) Restriction naturality.}
	For $S\subseteq T$ the formula gives $\operatorname{Inert}_q(c)(S)\cap T=c(S)$.
	For $J\subseteq T$ and $S\subseteq J\cup\{q\}$,
	\[
	\begin{aligned}
		\operatorname{Inert}_q(c)(S)\cap(J\cup\{q\})
		&=(c(S\cap J)\cap J)\cup(S\cap\{q\})\\
		&=\operatorname{Inert}_q(c|_J)(S),
	\end{aligned}
	\]
	where $c(S\cap J)\cap J=(c|_J)(S\cap J)$.
	
	\medskip
	\noindent\emph{(iii) Carath\'eodory complexity.}
	Restriction to $T$ gives $\operatorname{car}(c)\le \operatorname{car}(\operatorname{Inert}_q(c))$.
	Every consequence with an old conclusion has a determining subset inside $T$ of size at most $\operatorname{car}(c)$.
	The new conclusion $q$ occurs only when $q$ is already present, and is witnessed by $\{q\}$.
	Thus the reverse bound is $\max\{\operatorname{car}(c),1\}$.
	Since $T\ne\varnothing$, empty preservation gives $\operatorname{car}(c)\ge1$, proving equality.
	
	\medskip
	\noindent\emph{(iv) Armstrong complexity.}
	Put $r=\operatorname{arm}(c)$ and choose an $r$-row	Armstrong relation $R$ for $c$.
	Adjoin a column $q$ with value $0$ on all old rows.
	Duplicate one old row and assign value $1$ to $q$ only on this duplicate.
	The resulting relation has $r+1$ rows.
	Its projection to $T$ realizes exactly $c$.
	The two rows with identical old coordinates and different $q$-values show that no old premise determines $q$.
	Conversely, any failure of an old dependency remains witnessed by two original rows, both with $q=0$, even after $q$ is	added to the premise.
	Thus $q$ creates no additional old consequences, and the relation realizes precisely $\operatorname{Inert}_q(c)$.
	This proves
	\[
	\operatorname{arm}(\operatorname{Inert}_q(c))\le r+1.
	\]
	
	For the reverse inequality, let an $n$-row Armstrong relation realize $\operatorname{Inert}_q(c)$.
	Because $q\notin\operatorname{Inert}_q(c)(T)$, two rows agree on all columns in $T$ but differ at $q$.
	Upon projecting to $T$ and deleting duplicate projected rows, at most $n-1$ rows remain, and the projection still realizes $c$.
	Hence $r\le n-1$, giving the reverse inequality.
\end{proof}

\begin{definition}[Trace-neutral one-point suspension and suspension cost]
	\label{DEF:PackagingTightOnePointSuspension}
	Let $c$ be an empty-preserving closure on a nonempty finite set $F$, and let $q\notin F$.
	\begin{enumerate}[label=(\roman*),leftmargin=2.5em]
		\item \textbf{Trace-neutral suspension.}
		A closure $\Sigma$ on $F\cup\{q\}$ is a \textbf{trace-neutral one-point suspension} of $c$ if
		\[
		\Sigma|_F=c,
		\qquad
		\Sigma|_{J\cup\{q\}}
		=\operatorname{Inert}_q(c|_J)
		\quad(J\subsetneq F).
		\]
		Such a suspension is automatically empty-preserving:
		its restrictions to $F$ and to $\{q\}$ are empty-preserving.
		The new coordinate may interact with the full old face, but	on each proper old face it is inert.
		
		\item \textbf{Cost, defect, and tightness.}
		Let $\mu$ be a finite-valued numerical complexity of finite	closures, invariant under relabeling and monotone under restriction.
		Define the \textbf{$\mu$-suspension cost} and \textbf{$\mu$-suspension defect} by
		\[
		\operatorname{susp}_{\mu}(c)
		:=\min_{\Sigma\text{ trace-neutral}}\mu(\Sigma),
		\qquad
		\delta_{\mu}(c)
		:=\operatorname{susp}_{\mu}(c)-\mu(c).
		\]
		The minimum exists because there are only finitely many closures on $F\cup\{q\}$ and $\operatorname{Inert}_q(c)$ is a trace-neutral suspension.
		Monotonicity gives $\delta_\mu(c)\ge0$.
		
		A particular suspension $\Sigma$ is \textbf{$\mu$-tight at budget $p$} when
		\[
		\mu(c)=\mu(\Sigma)=p.
		\]
		The equality $\delta_\mu(c)=0$ asserts that	some $\mu$-tight suspension exists.
		It does not say that every suspension is tight.
	\end{enumerate}
\end{definition}

\begin{theorem}[Exact one-point suspension-cost dichotomy]
	\label{THM:PackagingOnePointSuspensionCostDichotomy}
	Let $F\ne\varnothing$ be finite and let $c$ be a finite empty-preserving closure on $F$.
	\begin{enumerate}[label=(\roman*),leftmargin=2.5em]
		\item \textbf{Carath\'eodory cost.}
		\[
		\operatorname{susp}_{\operatorname{car}}(c)
		=\operatorname{car}(c),\qquad
		\delta_{\operatorname{car}}(c)=0.
		\]
		
		\item \textbf{Armstrong cost.}
		Viewing $c$ as a finite dependency closure and putting $p=\operatorname{arm}(c)$, one has
		\[
		\operatorname{susp}_{\operatorname{arm}}(c)
		=\begin{cases}
			p,&c=\operatorname{id}_F,\\
			p+1,&c\ne\operatorname{id}_F.
		\end{cases}
		\]
		Equivalently,
		\[
		\delta_{\operatorname{arm}}(c)=0
		\iff c=\operatorname{id}_F,
		\qquad
		\delta_{\operatorname{arm}}(c)=1
		\quad(c\ne\operatorname{id}_F).
		\]
	\end{enumerate}
\end{theorem}

\begin{proof}
	\emph{(i) Carath\'eodory cost.}
	Every trace-neutral suspension $\Sigma$ restricts to $c$, so $\operatorname{car}(\Sigma)\ge\operatorname{car}(c)$.
	By Theorem~\ref{THM:PackagingUniversalOnePointExtensionCalculus}, the inert suspension attains this lower bound.
	Thus $\operatorname{susp}_{\operatorname{car}}(c) = \operatorname{car}(c)$ and the defect is zero.
	
	\medskip
	\noindent\emph{(ii) Armstrong cost.}
	Put $p=\operatorname{arm}(c)$.
	Restriction and the inert-extension theorem give
	\[
	p\le\operatorname{susp}_{\operatorname{arm}}(c)
	\le p+1.
	\]
	Suppose some trace-neutral suspension $\Sigma$ admits a $p$-row Armstrong relation $R$.
	The projection of $R$ to $F$ realizes $c$.
	It must have $p$ distinct projected rows, since otherwise $c$ would have Armstrong size less than $p$.
	
	Fix $a\in F$ and put $J=F\setminus\{a\}$.
	Trace neutrality gives
	\[
	\Sigma|_{J\cup\{q\}}
	=\operatorname{Inert}_q(c|_J),
	\qquad q\notin\Sigma(J).
	\]
	Thus two rows of $R$ agree on $J$ and differ at $q$.
	Their projections to $F$ are distinct, so these rows must differ at $a$.
	Consequently
	\[
	a\notin c(F\setminus\{a\})
	\qquad(a\in F).
	\]
	Any nontrivial dependency $A\to a$ would have $A\subseteq F\setminus\{a\}$ and would contradict this conclusion by monotonicity.
	Therefore $c=\operatorname{id}_F$.
	
	Conversely, if $c=\operatorname{id}_F$, set	$X=F\cup\{q\}$ and $\Sigma=u_X$.
	Every proper restriction of the uniform circuit is an identity closure.
	Hence
	\[
	\Sigma|_F=\operatorname{id}_F,
	\qquad
	\Sigma|_{J\cup\{q\}}
	=\operatorname{Inert}_q(\operatorname{id}_J)
	\quad(J\subsetneq F),
	\]
	so $\Sigma$ is trace-neutral.
	By Theorem~\ref{THM:PackagingNativeIdentityDimension}, Lemma~\ref{LEM:PackagingNativeUniformCircuit}, and $d_{\mathrm{comm}}=\operatorname{arm}$,
	\[
	\operatorname{arm}(\Sigma)=|F|+1
	=\operatorname{arm}(c)=p.
	\]
	Thus cost $p$ is attained exactly for identity closures.
	For every other closure the integer-valued cost is $p+1$.
\end{proof}

\begin{theorem}[Rigidity and uniqueness of tight Armstrong suspension]
	\label{THM:PackagingArmstrongSuspensionRigidity}
	Let $F$ be finite and nonempty, let $c$ be an empty-preserving dependency closure on $F$, and put $p=\operatorname{arm}(c)$.
	The following are equivalent:
	\begin{enumerate}[label=(\roman*),leftmargin=2.5em]
		\item $c$ admits an Armstrong-tight trace-neutral suspension at budget $p$;
		
		\item $c=\operatorname{id}_F$;
		
		\item $p=|F|+1$, and the uniform circuit $u_{F\cup\{q\}}$ is an Armstrong-tight trace-neutral suspension of $c$.
	\end{enumerate}
	Moreover, if $\Sigma$ is any trace-neutral suspension with $\operatorname{arm}(\Sigma)=p$, then
	\[
	\Sigma=u_{F\cup\{q\}}.
	\]
	Thus the uniform circuit is the unique closure-theoretic rank-preserving Armstrong suspension.
\end{theorem}

\begin{proof}
	\emph{(i)$\Longrightarrow$(ii).}
	A tight suspension has $\operatorname{arm}(\Sigma)=p$, so the Armstrong-cost dichotomy
	(Theorem~\ref{THM:PackagingOnePointSuspensionCostDichotomy}) forces $c=\operatorname{id}_F$.
	
	\medskip
	\noindent\emph{(ii)$\Longrightarrow$(iii).}
	For the identity closure, the native identity-dimension formula gives $p=|F|+1$.
	The uniform circuit $u_{F\cup\{q\}}$ is a trace-neutral suspension with Armstrong size $|F|+1$ by the construction in the suspension-cost theorem.
	It is therefore tight.
	
	\medskip
	\noindent\emph{(iii)$\Longrightarrow$(i).}
	The suspension explicitly supplied in (iii) is tight, so (i) follows.
	
	\medskip
	\noindent\emph{Uniqueness.}
	By the established equivalence, $c=\operatorname{id}_F$.
	Put $X=F\cup\{q\}$; then $p=|X|$.
	Every coatom $X\setminus\{z\}$ carries the identity trace:
	this is the old face $F$ when $z=q$, and follows from trace neutrality when $z\in F$.
	Its minimum Armstrong size is $p$.
	
	Let $R$ be any $p$-row Armstrong relation for $\Sigma$.
	The projection of $R$ to each coatom realizes the identity closure there, hence has $p$ distinct rows.
	Therefore that coatom is a key of $\Sigma$:
	\[
	z\in\Sigma(X\setminus\{z\})
	\qquad(z\in X).
	\]
	On the other hand, if $S\subseteq X$, $|S|\le p-2$, and $b\in X\setminus S$, choose $z\in X\setminus(S\cup\{b\})$.
	Since $\Sigma|_{X\setminus\{z\}}$ is the identity,
	\[
	b\notin\Sigma(S).
	\]
	Thus every coatom is a key and no smaller set has a nontrivial consequence.
	These are exactly the dependencies of $u_X$, so $\Sigma=u_X$.
\end{proof}

\begin{theorem}[Inert-or-trace-neutral cone assembly]
	\label{THM:PackagingInertTightConeAssembly}
	Let $\mathscr C$ be a restriction-stable family of finite empty-preserving closures.
	Let $(M,\mathcal F)$ be a primitive core in the sense of Definition~\ref{DEF:PackagingPrimitiveCoreNeutralCone}:
	$\mathcal F$ is an antichain of proper subsets covering $M$.
	Put
	\[
	\Delta(\mathcal F)
	=\{T\subseteq M:T\subseteq F\text{ for some }F\in\mathcal F\}.
	\]
	Let $\Gamma=(\gamma_F)_{F\in\mathcal F}$ be a compatible $\mathscr C$-atlas, and write $\gamma_T=\gamma_F|_T$ for $T\subseteq F\in\mathcal F$.
	Compatibility makes this independent of $F$.	
	Suppose a fresh symbol $q\notin M$ can be chosen such that:
	\begin{enumerate}[label=(C\arabic*),leftmargin=3em]
		\item \textbf{Inert lifts on nonmaximal faces.}
		For every $T\in\Delta(\mathcal F)\setminus\mathcal F$,
		\[
		\operatorname{Inert}_q(\gamma_T)
		\in\mathscr C(T\cup\{q\}).
		\]
		
		\item \textbf{Inert or trace-neutral lifts on facets.}
		For each $F\in\mathcal F$, either
		\[
		\operatorname{Inert}_q(\gamma_F)
		\in\mathscr C(F\cup\{q\}),
		\]
		or there is a trace-neutral suspension $\Sigma_F\in\mathscr C(F\cup\{q\})$ of $\gamma_F$, namely
		\[
		\Sigma_F|_F=\gamma_F,
		\qquad
		\Sigma_F|_{J\cup\{q\}}
		=\operatorname{Inert}_q(\gamma_J)
		\quad(J\subsetneq F).
		\]
	\end{enumerate}
	Then $\Gamma$ admits a neutral cone with cone coordinate $q$.
\end{theorem}

\begin{proof}
	We construct the lifted closures on every face $T\in\Delta(\mathcal F)$ and verify both neutral-cone identities.
	
	\emph{Nonmaximal faces (C1).}
	For $T\in\Delta(\mathcal F)\setminus\mathcal F$, set
	\[
	\widehat\gamma_T=\operatorname{Inert}_q(\gamma_T).
	\]
	These closures belong to $\mathscr C(T\cup\{q\})$ by (C1), and restrict to $\gamma_T$ on $T$.
	If $T'\subseteq T$ are both nonmaximal,	restriction naturality gives
	\[
	\widehat\gamma_T|_{T'\cup\{q\}}
	=\operatorname{Inert}_q(\gamma_{T'})
	=\widehat\gamma_{T'}.
	\]
	
	\medskip
	\noindent\emph{Facets (C2).}
	For each $F\in\mathcal F$, set $\widehat\gamma_F=\operatorname{Inert}_q(\gamma_F)$ if that lift belongs to $\mathscr C(F\cup\{q\})$.
	Otherwise take a trace-neutral $\Sigma_F$ as in (C2).
	Both choices restrict to $\gamma_F$ on $F$.
	For every proper $J\subsetneq F$, either restriction naturality or trace neutrality gives
	\[
	\widehat\gamma_F|_{J\cup\{q\}}
	=\operatorname{Inert}_q(\gamma_J)
	=\widehat\gamma_J.
	\]	
	Since $\mathcal F$ is an antichain, these are all possible strict face inclusions.
	Thus $\widehat\gamma_T|_{T'\cup\{q\}}=\widehat\gamma_{T'}$ whenever $T'\subseteq T$ in $\Delta(\mathcal F)$.
	In particular, two facet lifts agree over their intersection $J=F\cap F'$, since both restrict there to $\operatorname{Inert}_q(\gamma_J)$.
	Therefore $(\widehat\gamma_T)_{T\in\Delta(\mathcal F)}$ is a neutral cone over $\Gamma$.
\end{proof}

\subsection{Obstruction transport and portability}
\label{SUBSEC:PackagingObstructionTransport}

A counterexample on a small simplex boundary is not automatically a counterexample on a larger cover.
We therefore identify operations that carry an obstruction into a larger cover without exceeding the permitted local complexity.
These transport results provide the extra hypotheses needed for the converse of the abstract recovery theorem.

\begin{definition}[Relabeling stability of a closure class]
	\label{DEF:PackagingClosureRelabelingStability}
	Let $\mathscr C=(\mathscr C(I))_I$ be a family of closure operators on finite label sets.
	For a bijection $\sigma:I\to J$ and a closure $c$ on $I$, define the \textbf{transported closure} $\sigma_*c$ on $J$ by
	\[
	(\sigma_*c)(A)
	:=\sigma\bigl(c(\sigma^{-1}(A))\bigr)
	\qquad(A\subseteq J).
	\]
	The class $\mathscr C$ is \textbf{relabeling-stable} if
	\[
	c\in\mathscr C(I)\quad\Longrightarrow\quad
	\sigma_*c\in\mathscr C(J)
	\]
	for every such bijection and every $c\in\mathscr C(I)$.
	This is a property of the specialized closure class, not an axiom of general packaging or of the primitive realization judgement.
\end{definition}

\begin{definition}[Canonical coordinate cloning and clone-stability]
	\label{DEF:PackagingCanonicalCoordinateCloning}
	Let $T$ be finite, $q\notin T$, and let $O$ be a nonempty finite set disjoint from $T\cup\{q\}$.
	\begin{enumerate}[label=(\roman*),leftmargin=2.5em]
		\item \textbf{Canonical coordinate cloning.}
		For $S\subseteq T\cup O$ put
		\[
		\kappa_{q\to O}(S)
		=(S\cap T)\cup
		\begin{cases}
			\{q\},&S\cap O\ne\varnothing,\\
			\varnothing,&S\cap O=\varnothing.
		\end{cases}
		\]
		For a closure $c$ on $T\cup\{q\}$ define
		\[
		\begin{aligned}
			\operatorname{Clone}_{q\to O}(c)(S)
			:={}&\bigl(c(\kappa_{q\to O}(S))\cap T\bigr)\\
			&\cup
			\begin{cases}
				O,&q\in c(\kappa_{q\to O}(S)),\\
				\varnothing,&q\notin c(\kappa_{q\to O}(S)).
			\end{cases}
		\end{aligned}
		\]
		Thus all coordinates in $O$ are indistinguishable copies of $q$.
		
		\item \textbf{Clone-stability.}
		A restriction-stable closure class $\mathscr C$ is \textbf{clone-stable} if, for every such $T,q,O$,
		\[
		c\in\mathscr C(T\cup\{q\})
		\quad\Longrightarrow\quad
		\operatorname{Clone}_{q\to O}(c)\in\mathscr C(T\cup O).
		\]
	\end{enumerate}
\end{definition}

\begin{lemma}[Universal clone calculus]
	\label{LEM:PackagingUniversalCloneCalculus}
	For every closure operator $c$ on $T\cup\{q\}$, the following hold.
	\begin{enumerate}[label=(\roman*),leftmargin=2.5em]
		\item $\operatorname{Clone}_{q\to O}(c)$ is a closure operator on $T\cup O$.
		
		\item If $c(\varnothing)=\varnothing$, then $\operatorname{Clone}_{q\to O}(c)(\varnothing)=\varnothing$.
		
		\item For $T'\subseteq T$ and $O'\subseteq O$,
		\[
		\operatorname{Clone}_{q\to O}(c)|_{T'\cup O'}
		=\begin{cases}
			\operatorname{Clone}_{q\to O'}
			\bigl(c|_{T'\cup\{q\}}\bigr),
			&O'\ne\varnothing,\\[1ex]
			c|_{T'},&O'=\varnothing.
		\end{cases}
		\]
	\end{enumerate}
\end{lemma}

\begin{proof}
	Let $\pi:T\cup O\to T\cup\{q\}$ fix $T$ and send every element of $O$ to $q$.
	Because $O\ne\varnothing$, the map $\pi$ is surjective.
	By the cloning formula,
	\[
	\widetilde c(S)
	:=\operatorname{Clone}_{q\to O}(c)(S)
	=\pi^{-1}\bigl(c(\pi(S))\bigr).
	\]
	
	\emph{(i) Closure axioms.}
	Extensivity and monotonicity follow from those of $c$ and the inclusion $S\subseteq\pi^{-1}(\pi(S))$.
	Surjectivity gives $\pi(\pi^{-1}(A))=A$ for all $A\subseteq T\cup\{q\}$.
	Hence
	\[
	\widetilde c\bigl(\widetilde c(S)\bigr)
	=\pi^{-1}\bigl(c(c(\pi(S)))\bigr)
	=\widetilde c(S).
	\]
	Thus $\widetilde c$ is a closure.
	
	\emph{(ii) Empty preservation.}
	If $c(\varnothing)=\varnothing$, then
	\[
	\widetilde c(\varnothing)
	=\pi^{-1}(c(\varnothing))=\varnothing.
	\]
	
	\emph{(iii) Restriction naturality.}
	Fix $S\subseteq T'\cup O'$.
	If $O'\ne\varnothing$, the restriction $\pi':T'\cup O'\to T'\cup\{q\}$ is surjective and $\pi'(S)=\pi(S)$.
	Consequently,
	\[
	\begin{aligned}
		\widetilde c(S)\cap(T'\cup O')
		&= (\pi')^{-1}
		\bigl(c(\pi'(S))\cap(T'\cup\{q\})\bigr)\\
		&=\operatorname{Clone}_{q\to O'}
		(c|_{T'\cup\{q\}})(S).
	\end{aligned}
	\]
	If $O'=\varnothing$, then $S\subseteq T'$ and
	\[
	\widetilde c(S)\cap T'=c(S)\cap T'=(c|_{T'})(S).
	\]
	These are precisely the two asserted restriction formulas.
\end{proof}

\begin{definition}[Primitive cores, neutral cones, and embeddings]
	\label{DEF:PackagingPrimitiveCoreNeutralCone}
	Let $\mathscr C$ be a restriction-stable family of finite empty-preserving closures.
	\begin{enumerate}[label=(\roman*),leftmargin=2.5em]
		\item \textbf{Primitive core and its traces.}
		A \textbf{primitive core} $(M,\mathcal F)$ consists of a finite set $M$ and an antichain $\mathcal F$ of proper subsets of $M$ whose union is $M$.
		Put
		\[
		\Delta(\mathcal F)
		=\{T\subseteq M:T\subseteq F
		\text{ for some }F\in\mathcal F\}.
		\]
		For a compatible $\mathscr C$-atlas $\Gamma=(\gamma_F)_{F\in\mathcal F}$, define
		\[
		\gamma_T:=\gamma_F|_T
		\qquad(T\subseteq F\in\mathcal F).
		\]
		Compatibility makes $\gamma_T$ independent of the selected facet, and restriction stability gives $\gamma_T\in\mathscr C(T)$.
		
		\item \textbf{Neutral cone.}
		A \textbf{neutral cone} over $\Gamma$ consists of a fresh label	$q\notin M$ and a family
		\[
		\widehat\Gamma
		=(\widehat\gamma_T)_{T\in\Delta(\mathcal F)},
		\qquad
		\widehat\gamma_T\in\mathscr C(T\cup\{q\}),
		\]
		satisfying
		\[
		\widehat\gamma_T|_T=\gamma_T
		\]
		for every $T\in\Delta(\mathcal F)$, and
		\[
		\widehat\gamma_T|_{T'\cup\{q\}}
		=\widehat\gamma_{T'}
		\qquad(T'\subseteq T\text{ in }\Delta(\mathcal F)).
		\]
		The cone coordinate need not be inert.
		Neutrality here means preservation of the original traces and coherent restriction throughout the face complex.
		
		\item \textbf{Embedding in an ambient cover.}
		The primitive core $(M,\mathcal F)$ \textbf{embeds} in a cover $\mathcal H$ of $V\supseteq M$ if
		\begin{enumerate}[label=(Emb\arabic*),leftmargin=3em]
			\item $I\cap M\in\Delta(\mathcal F)$ for every $I\in\mathcal H$;
			
			\item every facet $F\in\mathcal F$ is guarded by $\mathcal H$.
		\end{enumerate}
	\end{enumerate}
\end{definition}

\begin{lemma}[Arbitrarily fresh cone coordinates]
	\label{LEM:PackagingNeutralConeFreshness}
	Let $\mathscr C$ be relabeling-stable, and suppose a compatible $\mathscr C$-atlas on a primitive core $(M,\mathcal F)$ admits a neutral cone with coordinate $q$.
	Then, for every finite set $Z$, it admits a	neutral cone with a coordinate
	\[
	q'\notin M\cup Z.
	\]
	Explicitly, transport each $\widehat\gamma_T$ along the bijection
	\[
	\sigma_T:T\cup\{q\}\longrightarrow T\cup\{q'\}
	\]
	fixing $T$ and sending $q$ to $q'$.
\end{lemma}

\begin{proof}
	Choose $q'\notin M\cup Z\cup\{q\}$ and put $\widehat\gamma_T'= (\sigma_T)_*\widehat\gamma_T$.
	Relabeling stability gives $\widehat\gamma_T'\in\mathscr C(T\cup\{q'\})$.
	Since $\sigma_T$ fixes $T$,
	\[
	\widehat\gamma_T'|_T
	=\widehat\gamma_T|_T=\gamma_T.
	\]
	For $T'\subseteq T$, the bijections agree on $T'\cup\{q\}$, so transport commutes with the relevant restriction:
	\[
	\widehat\gamma_T'|_{T'\cup\{q'\}}
	=(\sigma_{T'})_*
	(\widehat\gamma_T|_{T'\cup\{q\}})
	=(\sigma_{T'})_*\widehat\gamma_{T'}
	=\widehat\gamma_{T'}'.
	\]
	Thus the transported family is again a neutral cone.
\end{proof}

\begin{theorem}[Neutral-cone obstruction transport]
	\label{THM:PackagingNeutralConeTransport}
	Let $\mathscr C$ be a restriction-stable, relabeling-stable, and clone-stable family of finite empty-preserving closures.
	Let $(M,\mathcal F)$ be a primitive core, and let $\Gamma=(\gamma_F)_{F\in\mathcal F}$ be a compatible $\mathscr C$-atlas satisfying:
	\begin{enumerate}[label=(\alph*),leftmargin=2.5em]
		\item $\Gamma$ has no exact global closure on $M$;
		
		\item $\Gamma$ admits a neutral cone $\widehat\Gamma$.
	\end{enumerate}
	If $(M,\mathcal F)$ embeds in a finite cover $\mathcal H$ of $V$, then there is a compatible $\mathscr C$-atlas $(c_I)_{I\in\mathcal H}$ satisfying
	\[
	c_I|_{I\cap M}=\gamma_{I\cap M}
	\qquad(I\in\mathcal H),
	\]
	and admitting no exact global closure on $V$.
\end{theorem}

\begin{proof}
	\emph{Construction and admissibility.}
	By Lemma~\ref{LEM:PackagingNeutralConeFreshness}, choose the cone coordinate $q$ outside $V$.
	For each ambient chart $I\in\mathcal H$ put
	\[
	T_I=I\cap M,\qquad O_I=I\setminus M,
	\]
	and define
	\[
	c_I=
	\begin{cases}
		\gamma_{T_I},&O_I=\varnothing,\\[1ex]
		\operatorname{Clone}_{q\to O_I}
		(\widehat\gamma_{T_I}),&O_I\ne\varnothing.
	\end{cases}
	\]
	By embedding condition \textup{(Emb1)}, $T_I\in\Delta(\mathcal F)$.
	If $O_I=\varnothing$, restriction stability gives $c_I=\gamma_{T_I}\in\mathscr C(I)$.
	Otherwise clone-stability gives $c_I\in\mathscr C(I)$.
	The defining cone trace and the cloning restriction identity give
	\[
	c_I|_{T_I}=\gamma_{T_I}
	\qquad(I\in\mathcal H).
	\]
	
	\emph{Compatibility.}
	Fix $I,J\in\mathcal H$, and write
	\[
	T=I\cap J\cap M,
	\qquad O=(I\cap J)\setminus M.
	\]
	If $O=\varnothing$, restriction naturality (including the case of an empty cloned fibre) gives
	\[
	c_I|_{I\cap J}=\gamma_T=c_J|_{I\cap J}.
	\]
	If $O\ne\varnothing$, both charts have a nonempty off-core part.
	Lemma~\ref{LEM:PackagingUniversalCloneCalculus} and cone coherence give
	\[
	c_I|_{I\cap J}
	=\operatorname{Clone}_{q\to O}(\widehat\gamma_T)
	=c_J|_{I\cap J}.
	\]
	Thus the constructed atlas is compatible.
	
	\emph{Persistence of the obstruction.}
	Suppose an exact global closure $C$ on $V$ existed.
	By embedding condition \textup{(Emb2)}, for each $F\in\mathcal F$ choose a chart $I_F$ containing $F$.
	Since $I_F\cap M\in\Delta(\mathcal F)$ and $F$ is a maximal face, $I_F\cap M=F$.
	Therefore
	\[
	(C|_M)|_F
	=C|_F
	=c_{I_F}|_F
	=\gamma_F.
	\]
	Hence $C|_M$ would be an exact global closure for $\Gamma$, contradicting hypothesis~\textup{(a)}.
\end{proof}

In the native commutant specialization, cloning simply repeats an existing matrix coordinate and introduces no additional joint spectral sectors.
Thus the local dimension of the transported atlas is controlled by the chosen cone lifts rather than by the number of off-core labels.

For $q\ge4$, let $\mathcal C_q$ be the edge cover of the chordless cycle on $q$ vertices.
Its edges form an antichain primitive core.
The simplex-boundary cover $\mathcal B_r$ provides the other family of primitive cores used below.

\begin{definition}[Neutral transport completeness]
	\label{DEF:PackagingNeutralTransportCompleteness}
	Let $m\ge3$ and let $\mathscr C$ be a restriction-stable family of finite empty-preserving closures. We call $\mathscr C$ \textbf{neutral-transport complete through rank $m$} if:
	\begin{enumerate}[label=(NT\arabic*),leftmargin=3em]
		\item $\mathscr C$ is clone-stable;
		
		\item for every $q\ge4$, $\mathcal C_q$ supports a compatible nonextendible $\mathscr C$-atlas admitting a neutral cone;
		
		\item for every $3\le r\le m$, $\mathcal B_r$ supports a compatible nonextendible $\mathscr C$-atlas admitting a neutral cone.
	\end{enumerate}
	These are conditions on the canonical primitive cores, not on arbitrary ambient covers.
\end{definition}

\begin{definition}[One-point transport completeness]
	\label{DEF:PackagingOnePointTransportCompleteness}
	Let $m\ge3$ and let $\mathscr C$ be a restriction-stable family of finite empty-preserving closures.
	We call $\mathscr C$ \textbf{one-point transport complete through rank $m$} if:
	\begin{enumerate}[label=(O\arabic*),leftmargin=3em]
		\item $\mathscr C$ is clone-stable;
		
		\item for every $q\ge4$, $\mathcal C_q$ supports a compatible nonextendible $\mathscr C$-atlas satisfying conditions \textup{(C1)}--\textup{(C2)} of
		Theorem~\ref{THM:PackagingInertTightConeAssembly};
		
		\item for every $3\le r\le m$, $\mathcal B_r$ supports a compatible nonextendible $\mathscr C$-atlas satisfying the same conditions \textup{(C1)}--\textup{(C2)}.
	\end{enumerate}
	These hypotheses are weaker than neutral-transport completeness.
They require inert extensions on nonmaximal faces.
For each facet, they require either an inert extension or an allowed trace-neutral one-point suspension.
	They do not presuppose a coherent cone family.
\end{definition}

\begin{proposition}[One-point extension data generate neutral transport]
	\label{PROP:PackagingOnePointImpliesNeutralTransport}
	Let $m\ge3$ and let $\mathscr C$ be a restriction-stable family of finite empty-preserving closures.
	If $\mathscr C$ is one-point transport complete through rank $m$, then it is neutral-transport complete through rank $m$.
\end{proposition}

\begin{proof}
	We verify conditions \textup{(NT1)}--\textup{(NT3)}.
	\emph{(NT1):} Clone-stability is precisely condition \textup{(O1)}.
	
	\emph{(NT2):} For each $q\ge4$, condition \textup{(O2)} and Theorem~\ref{THM:PackagingInertTightConeAssembly} assemble a neutral cone over the stipulated nonextendible cycle atlas.
	
	\emph{(NT3):} The same theorem applied to the boundary witnesses in \textup{(O3)} supplies neutral cones on every $\mathcal B_r$, $3\le r\le m$.
	Thus all three neutral-transport conditions hold.
\end{proof}

\begin{proposition}[Neutral transport generates primitive portability]
	\label{PROP:PackagingNeutralTransportImpliesPortability}
	Let $m\ge3$ and let $\mathscr C$ be a restriction-stable,
	relabeling-stable family of finite empty-preserving closures that is neutral-transport complete through rank $m$.
	Then the following hold.
	\begin{enumerate}[label=(\roman*),leftmargin=2.5em]
		\item Every finite cover whose primal graph contains an induced	cycle of length at least four supports a compatible nonextendible $\mathscr C$-atlas.
		
		\item Every finite cover containing an inclusion-minimal unguarded primal clique $M$ with $3\le|M|\le m$ supports a compatible nonextendible $\mathscr C$-atlas.
	\end{enumerate}
	These facts show that \(\mathscr C\), with its rooted-circuit recovery structure and the ambient class of all empty-preserving closures, is primitive-obstruction portable through rank \(m\) as defined in Definition~\ref{DEF:PackagingPrimitiveObstructionPortability}.
\end{proposition}

\begin{proof}
	\emph{(i) Induced-cycle obstruction.}
	Let $M$ be the vertex set of an induced cycle of length $q\ge4$.
	Every chart intersection $I\cap M$ is a clique in the induced cycle, hence is contained in a cycle edge.
	Every cycle edge is guarded by the definition of the primal graph.
	Thus the cycle core $\mathcal C_q$ embeds in the ambient cover.
	By neutral-transport completeness \textup{(NT2)}, the cycle core has a compatible nonextendible $\mathscr C$-atlas admitting a neutral cone.
	Relabeling stability permits identification with the given cycle, and Theorem~\ref{THM:PackagingNeutralConeTransport} transports its obstruction to the ambient cover.
	
	\emph{(ii) Unguarded-clique obstruction.}
	Let $M$ be inclusion-minimal unguarded, with $3\le|M|\le m$.
	No chart contains $M$, so every $I\cap M$ is a proper subset of $M$ and belongs to the downward closure of the simplex-boundary facets.
	Every coatom $M\setminus\{x\}$ is guarded by inclusion-minimality.
	Hence the simplex-boundary core $\mathcal B_{|M|}$ embeds.
	Condition \textup{(NT3)}, relabeling stability, and Theorem~\ref{THM:PackagingNeutralConeTransport} give the required compatible nonextendible ambient atlas.
	
	Finally, these are precisely the two conditions for primitive-obstruction portability through rank $m$.
	The recovery-equipped closure structure is supplied by Proposition~\ref{PROP:PackagingClosureRecoveryEquipped}.
\end{proof}

\subsection{Abstract recovery--obstruction duality}
\label{SUBSEC:PackagingRecoveryObstructionDuality}

The abstract globalization theorem has two directions, and they require different assumptions.
Finite intrinsic recovery rank is enough for the sufficient direction.
The converse requires transporting primitive obstructions through arbitrary finite covers.
Theorem~\ref{THM:PackagingTransportNecessityCounterexample} gives a class for which the converse fails without that transport property.
We therefore state the full equivalence only under the required hypothesis.

For a finite cover $\mathcal H$ of a finite label set $V$, write
$
\operatorname{UG}_{\mathscr C}(\mathcal H)
$
if every compatible $\mathscr C$-atlas on $\mathcal H$ admits an exact global empty-preserving closure on $V$.

\begin{theorem}[Obstruction transport is indispensable for abstract duality]
	\label{THM:PackagingTransportNecessityCounterexample}
	There exists a restriction-stable and relabeling-stable family $\mathscr C$ of finite empty-preserving closures for which
	\[
	\operatorname{ipr}(\mathscr C)=4,
	\qquad
	\operatorname{bor}(\mathscr C)=2,
	\]
	and yet
	\[
	\operatorname{UG}_{\mathscr C}(\mathcal B_3)
	\quad\text{while}\quad
	\nu(\mathcal B_3)=3.
	\]
	In particular, restriction stability, relabeling stability, and finite intrinsic recovery rank \emph{do not imply} the converse to intrinsic recovery sufficiency.
	Some additional obstruction-existence or transport condition is required for a uniform rank--geometry equivalence.
\end{theorem}

\begin{proof}
	For each finite set $I$, let $\operatorname{id}_I(S)=S$ for all $S\subseteq I$.
	When $|I|=3$, also let $u_I$ be the uniform-circuit closure
	\[
	 u_I(S)=
	\begin{cases}
	S,&|S|\le1,\\
	I,&|S|\ge2.
	\end{cases}
	\]
	Set
	\[
	\mathscr C(I)=
	\begin{cases}
	\{\operatorname{id}_I,u_I\},&|I|=3,\\
	\{\operatorname{id}_I\},&|I|\ne3.
	\end{cases}
	\]
	Every proper restriction of $u_I$ is an identity closure.
	Thus $\mathscr C$ is restriction-stable, and it is manifestly relabeling-stable.

	On the three-vertex simplex boundary $\mathcal B_3$, all original charts have two labels and hence must carry identity closures.
	Their atlas extends by $\operatorname{id}_{V_3}$.
	Therefore $\operatorname{UG}_{\mathscr C}(\mathcal B_3)$ holds, and the initial-segment definition gives $\operatorname{bor}(\mathscr C)=2$.
	The obstruction rank of this boundary is $\nu(\mathcal B_3)=3$.

	We compute the intrinsic recovery rank.
	Let $V=\{a,b,v,y\}$ and take its four simplex-boundary charts.
	Put $u$ on $\{a,b,v\}$ and $\{b,v,y\}$, and identity closures on $\{a,b,y\}$ and $\{a,v,y\}$.
	All overlaps are proper two-element subsets of the three-element charts, and both types of closure restrict to the identity there.
	Hence the atlas is compatible.
	The rooted circuits
	\[
	\{a,b\}\rightsquigarrow v,
	\qquad
	\{b,v\}\rightsquigarrow y
	\]
	are valid in the two uniform-circuit charts.
	Their combined support is the primal clique $V$, and $A\cup B=\{a,b\}$.
	No original chart contains a rooted circuit \(D\rightsquigarrow y\) with \(D\subseteq\{a,b\}\).
Indeed, the chart \(\{a,b,y\}\) has the identity closure.
Compatibility would transfer any such implication to that chart, which is impossible.
	The cover guards every clique of size at most three.
	Therefore $\mathscr C$ is not intrinsically $3$-recoverable, and $\operatorname{ipr}(\mathscr C)\ge4$.
	The same two rules derive $y$ from $\{a,b\}$ in the globally generated closure, contradicting the identity trace on $\{a,b,y\}$; thus $\mathcal B_4$ is in fact boundary-obstructed.

	For the reverse inequality, consider any compatible $\mathscr C$-atlas whose cover guards all primal cliques of size at most four.
	Every rooted circuit in such an atlas comes from a three-element uniform-circuit chart, and has a two-element stem and the third label as its root.
	Suppose two such circuits are composable across a pivot and their union of supports $Q$ is a primal clique.
	If $|Q|\ge4$, choose a four-element subset of $Q$ containing the full three-element support $I$ of the first circuit.
	By four-conformality, an original chart $J$ contains these four labels.
	Since $|J|\ge4$, it carries an identity closure, while its overlap with $I$ contains all of $I$.
	Compatibility then forces the circuit on $I$ to be valid in an identity closure, a contradiction.
	Thus $|Q|=3$.
	But then the first circuit's two-element stem equals $Q\setminus\{v\}$ and contains the other root $b\ne v$, contradicting the composability condition $b\notin A\cup B$.
	There are therefore no composable rooted-circuit pairs with clique support in any such four-conformal atlas.
	Intrinsic $4$-recoverability holds, so $\operatorname{ipr}(\mathscr C)\le4$.
	This completes the rank computation and the counterexample.
\end{proof}

\begin{proposition}[Boundary obstructions lower-bound intrinsic recovery]
	\label{PROP:PackagingBoundaryRankLowerBoundsRecovery}
	Let $\mathscr C$ be a restriction-stable family of finite empty-preserving closures.
	Then
	$
		\operatorname{bor}(\mathscr C)
		\le\operatorname{ipr}(\mathscr C).
	$
	No relabeling or obstruction-transport hypothesis is required.
\end{proposition}

\begin{proof}
	If $\operatorname{ipr}(\mathscr C)=\infty$, the inequality is immediate.
	Otherwise put $p=\operatorname{ipr}(\mathscr C)$, where $p\ge2$.
	
	By Theorem~\ref{THM:PackagingRootedCircuitCharacterization}, the rooted-circuit scheme $\mathscr R^{\mathrm{rc}}$ is admissible and
	\[
	\operatorname{prk}
	(\mathscr C,\mathscr R^{\mathrm{rc}})=p.
	\]
	Fix an integer $r>p$.
	The simplex-boundary cover $\mathcal B_r$ has complete, hence chordal, primal graph.
	Its only unguarded nonempty primal clique is $V_r$ itself.
	Since $|V_r|=r>p$, every primal clique of cardinality at most $p$ is guarded.
	
	Theorem~\ref{THM:PackagingControlledPropagationGlobalization} therefore gives an exact global closure for every compatible $\mathscr C$-atlas on $\mathcal B_r$.
	In particular, $\mathcal B_{p+1}$ is not boundary-obstructed.
	By the initial-segment definition of $\operatorname{bor}(\mathscr C)$, no admissible initial segment extends to $p+1$.
	Hence
	\[
	\operatorname{bor}(\mathscr C)\le p
	=\operatorname{ipr}(\mathscr C).
	\]
\end{proof}

\begin{theorem}[Intrinsic recovery sufficiency without transport]
	\label{THM:PackagingIntrinsicRecoverySufficiency}
	Let $\mathscr C$ be a restriction-stable family of finite empty-preserving closures, and suppose
	\[
	p:=\operatorname{ipr}(\mathscr C)<\infty,
	\qquad p\ge3.
	\]
	Then, for every finite cover $\mathcal H$,
	\[
		p<\nu(\mathcal H)
		\quad\Longrightarrow\quad
		\operatorname{UG}_{\mathscr C}(\mathcal H).
	\]
	Neither relabeling stability, clone-stability, neutral-transport completeness, nor primitive-obstruction portability is assumed.
\end{theorem}

\begin{proof}
	Equip $\mathscr C$ with its canonical rooted-circuit recovery signature and take as ambient states all finite empty-preserving closures. Proposition~\ref{PROP:PackagingClosureRecoveryEquipped}	verifies the recovery-equipped axioms, while Proposition~\ref{PROP:PackagingNativeRankClosureSpecialization} gives
	\[
	\rho_{\mathrm{pkg}}(\mathscr C)
	=\operatorname{ipr}(\mathscr C)=p.
	\]
	Theorem~\ref{THM:PackagingNativeRecoverySufficiency} therefore implies that, whenever $p<\nu(\mathcal H)$, every compatible $\mathscr C$-atlas has an exact global state in the chosen ambient class.
	This is precisely an exact global empty-preserving closure, proving the assertion.
\end{proof}

\begin{theorem}[Transport-complete closure recovery--obstruction duality]
	\label{THM:PackagingAbstractRecoveryObstructionDuality}
	Let $\mathscr C$ be a restriction-stable and relabeling-stable family of finite empty-preserving closures.
	Suppose
	\[
	p:=\operatorname{ipr}(\mathscr C)<\infty,
	\qquad p\ge3,
	\]
	and assume that $\mathscr C$ is neutral-transport complete through rank $p$.
	Then the following hold.
	\begin{enumerate}[label=(\roman*),leftmargin=2.5em]
		\item \textbf{Equality of recovery and obstruction ranks.}
		\[
			\rho_{\mathrm{pkg}}(\mathscr C)
			=\operatorname{ipr}(\mathscr C)
			=\operatorname{bor}(\mathscr C)=p.
		\]
		
		\item \textbf{Exact universal-globalization classification.}
		For every finite cover $\mathcal H$,
		\[
			\operatorname{UG}_{\mathscr C}(\mathcal H)
			\quad\Longleftrightarrow\quad
			p<\nu(\mathcal H).
		\]
	\end{enumerate}
\end{theorem}

\begin{proof}
	\leavevmode
	\begin{enumerate}[label=(\roman*),leftmargin=2.5em]
		\item \emph{Equality of the ranks.}
		Proposition~\ref{PROP:PackagingNativeRankClosureSpecialization} yields
		\[
		\rho_{\mathrm{pkg}}(\mathscr C)
		=\operatorname{ipr}(\mathscr C)=p.
		\]
		Neutral-transport completeness through rank $p$ supplies a compatible nonextendible $\mathscr C$-atlas on every simplex boundary $\mathcal B_r$ with $3\le r\le p$.
		By the definition of the boundary-obstruction rank,
		\[
		\operatorname{bor}(\mathscr C)\ge p.
		\]
		Proposition~\ref{PROP:PackagingBoundaryRankLowerBoundsRecovery}	gives the reverse inequality
		\[
		\operatorname{bor}(\mathscr C)
		\le\operatorname{ipr}(\mathscr C)=p.
		\]
		Thus all three ranks equal $p$.
		
		\item \emph{Universal-globalization classification.}
		Fix a finite cover $\mathcal H$.
		
		\emph{Sufficiency.}
		If $p<\nu(\mathcal H)$, then Theorem~\ref{THM:PackagingIntrinsicRecoverySufficiency} supplies an exact global empty-preserving closure for every compatible $\mathscr C$-atlas on $\mathcal H$.
		Hence $\operatorname{UG}_{\mathscr C}(\mathcal H)$ holds.
		
		\emph{Necessity (by contraposition).}
		Suppose $\nu(\mathcal H)\le p$.
		By Proposition~\ref{PROP:PackagingNeutralTransportImpliesPortability}, relabeling stability and neutral-transport completeness give primitive-obstruction portability through rank $p$.
		
		If $G_{\mathcal H}$ is nonchordal, portability supplies a compatible nonextendible $\mathscr C$-atlas on $\mathcal H$.
		If $G_{\mathcal H}$ is chordal, the definition of $\nu(\mathcal H)$ supplies an unguarded primal clique of minimum cardinality $r=\nu(\mathcal H)$.
		This clique is inclusion-minimal unguarded, and
		$
		3\le r\le p.
		$
		Portability again supplies a compatible nonextendible $\mathscr C$-atlas on $\mathcal H$.
		In either case $\operatorname{UG}_{\mathscr C}(\mathcal H)$ fails.
		Thus universal exact globalization implies $p<\nu(\mathcal H)$, proving the equivalence.
	\end{enumerate}
\end{proof}

\begin{corollary}[One-point sufficient criterion for recovery--obstruction duality]
	\label{COR:PackagingOnePointRecoveryObstructionDuality}
	Let $\mathscr C$ be a restriction-stable and relabeling-stable family of finite empty-preserving closures, with
	\[
	p:=\operatorname{ipr}(\mathscr C)<\infty,
	\qquad p\ge3.
	\]
	If $\mathscr C$ is one-point transport complete through rank $p$, then
	\[
	\rho_{\mathrm{pkg}}(\mathscr C)
	=\operatorname{ipr}(\mathscr C)
	=\operatorname{bor}(\mathscr C)=p,
	\]
	and, for every finite cover $\mathcal H$,
	\[
	\operatorname{UG}_{\mathscr C}(\mathcal H)
	\quad\Longleftrightarrow\quad
	p<\nu(\mathcal H).
	\]
\end{corollary}

\begin{proof}
	Proposition~\ref{PROP:PackagingOnePointImpliesNeutralTransport} converts the one-point transport hypotheses into neutral-transport completeness through rank $p$.
	All hypotheses of Theorem~\ref{THM:PackagingAbstractRecoveryObstructionDuality} are therefore satisfied.
	Its conclusions give both assertions.
\end{proof}

\begin{corollary}[Rank-transfer principle for universal globalization]
	\label{COR:PackagingUniversalRankTransfer}
	Let $\mathscr C$ and $\mathscr D$ each satisfy the hypotheses of Theorem~\ref{THM:PackagingAbstractRecoveryObstructionDuality}, and put
	\[
	p=\operatorname{ipr}(\mathscr C),
	\qquad q=\operatorname{ipr}(\mathscr D),
	\qquad p,q\ge3.
	\]
	Write $\mathfrak G_{\mathscr C}$ and $\mathfrak G_{\mathscr D}$ for their respective classes of finite covers admitting universal exact globalization.
	Then:
	\begin{enumerate}[label=(\roman*),leftmargin=2.5em]
		\item \textbf{Threshold descriptions.}
		\[
		\mathfrak G_{\mathscr C}
		=\{\mathcal H:p<\nu(\mathcal H)\},
		\qquad
		\mathfrak G_{\mathscr D}
		=\{\mathcal H:q<\nu(\mathcal H)\}.
		\]
		
		\item \textbf{Rank comparison.}
		If $p=q$, then $\mathfrak G_{\mathscr C}=\mathfrak G_{\mathscr D}$.
		If $p<q$, then
		$
		\mathfrak G_{\mathscr C}
		\supsetneq\mathfrak G_{\mathscr D}.
		$
		If $q<p$, the strict inclusion is reversed.
	\end{enumerate}
\end{corollary}

\begin{proof}
	\emph{(i) Threshold descriptions.}
	Apply Theorem~\ref{THM:PackagingAbstractRecoveryObstructionDuality} separately to $\mathscr C$ and $\mathscr D$.
	
	\emph{(ii) Rank comparison.}
	The assertion for $p=q$ follows immediately from~\textup{(i)}.
	Suppose $p<q$. Then $q<\nu(\mathcal H)$ implies $p<\nu(\mathcal H)$, so $\mathfrak G_{\mathscr D}\subseteq\mathfrak G_{\mathscr C}$.
	To prove strictness, take the simplex boundary $\mathcal B_{p+1}$.
	Its primal graph is complete and its only unguarded clique is the full $(p+1)$-vertex set.
	Thus
	$
	\nu(\mathcal B_{p+1})=p+1.
	$
	Hence $p<p+1$, but $q\ge p+1$, so $q<\nu(\mathcal B_{p+1})$ fails.
	Consequently
	\[
	\mathcal B_{p+1}
	\in\mathfrak G_{\mathscr C}
	\setminus\mathfrak G_{\mathscr D}.
	\]
	The case $q<p$ follows by interchanging the two classes.
\end{proof}

\subsection{Concrete recovery ranks and the globalization obstruction spectrum}
\label{SUBSEC:PackagingConcreteSpectrum}

For bounded-Carath\'eodory and bounded-Armstrong closures, the recovery and obstruction-transport hypotheses are verified explicitly.
Combining these results with the transitive classification yields the three sharp thresholds in Theorem~\ref{THM:PackagingGlobalizationSpectrum}.
The distinct recovery ranks produce distinct classes of covers with universal exact globalization.

\begin{lemma}[Relabeling stability of the concrete closure classes]
	\label{LEM:PackagingConcreteRelabelingStability}
	For every \(k\ge1\) and \(d\ge2\), the classes
	$
	\mathscr C_k
	$
	of finite empty-preserving closures of Carath\'eodory number at most \(k\), and
	\[
	\mathscr N_d=\mathbf{CommCl}_d
	\]
	of native closures of dimension at most \(d\), are relabeling-stable.
\end{lemma}

\begin{proof}
	\emph{Bounded-Carath\'eodory class.}
	Let $\sigma:I\to J$ be a bijection, let $c\in\mathscr C_k(I)$, and put $c'=\sigma_*c$.
	If $b\in c'(S)$, then $\sigma^{-1}(b)\in c(\sigma^{-1}(S))$.
	There exists $A\subseteq\sigma^{-1}(S)$ of size at most $k$ with $\sigma^{-1}(b)\in c(A)$.
	Consequently, $\sigma(A)\subseteq S$, $|\sigma(A)|\le k$, and $b\in c'(\sigma(A))$.
	Empty preservation is also invariant under relabeling, so $c'\in\mathscr C_k(J)$.
	
	\emph{Native class.}
	If $c\in\mathbf{CommCl}_d(I)$ is realized by commuting self-adjoint matrices $(D_a)_{a\in I}$, assign the matrix $D_a$ to the new label $\sigma(a)$.
	This realizes $\sigma_*c$ in the same dimension.
	Hence $\mathbf{CommCl}_d$ is relabeling-stable.
\end{proof}

\begin{lemma}[Generic one-point extension and clone stability]
	\label{LEM:PackagingGenericNeutralTransport}
	Fix $k\ge1$ and let $\mathscr C_k$ be the class of finite empty-preserving closures of Carath\'eodory number at most $k$.
	Then:
	\begin{enumerate}[label=(\roman*),leftmargin=2.5em]
		\item for every $c\in\mathscr C_k(T)$ and every $q\notin T$,
		\[
		\operatorname{Inert}_q(c)\in\mathscr C_k(T\cup\{q\});
		\]
		
		\item $\mathscr C_k$ is clone-stable;
		
		\item every compatible $\mathscr C_k$-atlas on a primitive core	satisfies the hypotheses of Theorem~\ref{THM:PackagingInertTightConeAssembly} by taking the inert extension on every face.
	\end{enumerate}
\end{lemma}

\begin{proof}
	\emph{(i) Inert extensions.}
	If $T\ne\varnothing$, the exact Carath\'eodory formula in Theorem~\ref{THM:PackagingUniversalOnePointExtensionCalculus}	gives the conclusion.
	If $T=\varnothing$, the inert extension	is $\operatorname{id}_{\{q\}}$, whose Carath\'eodory number	is $1\le k$.
	
	\emph{(ii) Cloning.}
	Let $c\in\mathscr C_k(T\cup\{q\})$ and $c'=\operatorname{Clone}_{q\to O}(c)$, where $O$ is finite, nonempty, and disjoint from $T\cup\{q\}$.
	For $S\subseteq T\cup O$ put $P=\kappa_{q\to O}(S)$.
	If $b\in c'(S)\cap T$, then $b\in c(P)$.
	If $b\in c'(S)\cap O$, then $q\in c(P)$.
	In either case the Carath\'eodory bound on $c$ provides $A_0\subseteq P$, $|A_0|\le k$, determining $b$ or $q$,	respectively.
	If $q\notin A_0$, take $A=A_0\subseteq S$.
	If $q\in A_0$, then $S\cap O\ne\varnothing$.
	Choose	$o\in S\cap O$ and replace $q$ in $A_0$ by $o$.
	In both cases $A\subseteq S$, $|A|\le k$, and the cloning formula gives $b\in c'(A)$.
	The clone is empty-preserving by Lemma~\ref{LEM:PackagingUniversalCloneCalculus}, so $c'\in\mathscr C_k(T\cup O)$.
	
	\emph{(iii) Primitive cones.}
	Fix one fresh label $q$ outside the primitive core.
	On every face $T$ of a compatible primitive atlas, take $\widehat\gamma_T=\operatorname{Inert}_q(\gamma_T)$.
	Part~\textup{(i)} keeps these closures in $\mathscr C_k$, while restriction naturality from Theorem~\ref{THM:PackagingUniversalOnePointExtensionCalculus} gives coherence on all nested faces.
	Thus conditions	\textup{(C1)}--\textup{(C2)} of the cone-assembly theorem hold using only inert extensions.
\end{proof}

\begin{lemma}[Native clone stability and primitive one-point transport]
	\label{LEM:PackagingNativePrimitiveNeutralCones}
	For every $d\ge3$, the native class $\mathscr N_d=\mathbf{CommCl}_d$ has the following properties:
	\begin{enumerate}[label=(\roman*),leftmargin=2.5em]
		\item $\mathscr N_d$ is clone-stable;
		
		\item for every $\ell\ge4$, the cycle core $\mathcal C_\ell$ supports a compatible nonextendible $\mathscr N_d$-atlas satisfying the inert-or-trace-neutral cone-assembly conditions \textup{(C1)}--\textup{(C2)};
		
		\item for every $3\le r\le d$, the simplex-boundary core $\mathcal B_r$ supports such an atlas.
	\end{enumerate}
	In particular, the witnesses in \textup{(ii)}--\textup{(iii)} admit neutral cones by Theorem~\ref{THM:PackagingInertTightConeAssembly}.
\end{lemma}

\begin{proof}
	\emph{(i) Clone stability.}
	Lemma~\ref{LEM:PackagingNativeCloneLift} realizes each cloned coordinate by copying its original commuting self-adjoint matrix.
	The matrix dimension is unchanged, so the cloned closure remains in $\mathbf{CommCl}_d$.
	
	\emph{(ii) Cycle cores.}
	Fix a fresh auxiliary label $q$ outside the primitive core.
	On $\mathcal C_\ell$, $\ell\ge4$, assign the two-label uniform circuit closure $u_2$ to every edge of a path through all $\ell$ vertices, and the identity closure $\operatorname{id}_2$ to the closing edge.
	The atlas is compatible because its	singleton traces are identities.
	It is nonextendible:
	composing the path implications forces its endpoints to	determine each other, contrary to the closing identity edge.
	The path-edge closure $u_2$ has Armstrong size $2$,	whereas $\operatorname{id}_2$ has Armstrong size $3$.
	
	If $d=3$, give every path edge its inert lift, of Armstrong	size $3$, and give the closing edge the uniform circuit	$u_3$, a trace-neutral lift of $\operatorname{id}_2$ of Armstrong size $3$.
	If $d\ge4$, use inert lifts on every edge.
	The largest Armstrong size is $4\le d$.
	Every proper edge trace is an identity on at most one label.
	Its inert extension has Armstrong size at most $3$, and each chosen edge lift restricts to that same inert extension.
	Thus both \textup{(C1)} and \textup{(C2)} hold.
	
	\emph{(iii) Simplex boundaries.}
	Choose a fresh auxiliary label $q$ for the boundary lift.
	Fix $3\le r\le d$ and choose the compatible nonextendible boundary atlas from Theorem~\ref{THM:PackagingNativeSimplexBoundaryGeneral}:
	two facets carry $u_{r-1}$ and all others carry	$\operatorname{id}_{r-1}$.
	Every proper subfacet trace	is an identity closure.
	The uniform circuit facet has Armstrong size $r-1$.
	Its inert lift has size $r\le d$.
	The identity facet has Armstrong size $r$.
	If $r<d$, its inert lift has size $r+1\le d$.
	If $r=d$, instead use its trace-neutral uniform-circuit suspension $u_r$, of Armstrong size $r=d$.
	
	For every nonmaximal face $J$, the canonical lift is $\operatorname{Inert}_q(\operatorname{id}_J) =\operatorname{id}_{J\cup\{q\}}$,	of Armstrong size $|J|+2\le r\le d$.
	This also covers $J=\varnothing$.
	Inert facet lifts restrict naturally, and the uniform-circuit suspension restricts to the identity on every proper subfacet together with $q$.
	Thus the prescribed facet lifts agree with every canonical lower-face lift and satisfy \textup{(C1)}--\textup{(C2)}.
	
	Theorem~\ref{THM:PackagingInertTightConeAssembly} therefore supplies coherent neutral cones for both types of primitive witness.
\end{proof}

\begin{theorem}[Intrinsic sharpness and one-point transport of the generic and Armstrong/native classes]
	\label{THM:PackagingIntrinsicSharpnessConcreteClasses}
	Let $\mathscr C_k$ be the class of finite empty-preserving closures of Carath\'eodory number at most $k$, $k\ge1$.
	Let $\mathscr N_d=\mathbf{CommCl}_d$ be the class of finite	empty-preserving dependency closures of minimum Armstrong size at most $d$, where $d\ge3$.
	Then:
	\begin{enumerate}[label=(\roman*),leftmargin=2.5em]
		\item \emph{Generic closures.}
		\[
		\operatorname{bor}(\mathscr C_k)
		=\operatorname{ipr}(\mathscr C_k)=2k+1,
		\]
		and $\mathscr C_k$ is one-point transport complete through rank $2k+1$.
		
		\item \emph{Armstrong/native closures.}
		\[
		\operatorname{bor}(\mathscr N_d)
		=\operatorname{ipr}(\mathscr N_d)=d,
		\]
		and $\mathscr N_d$ is one-point transport complete through rank $d$.
	\end{enumerate}
\end{theorem}

\begin{proof}
	\leavevmode
	\begin{enumerate}[label=(\roman*),leftmargin=2.5em]
		\item \emph{Generic closures.}
		Theorem~\ref{THM:PackagingBoundedDeterminationPropagationRank} and the rooted-circuit characterization give $\operatorname{ipr}(\mathscr C_k)\le2k+1$.
		On the simplex boundary with $2k+1$ vertices choose the	two-rule witness of Theorem~\ref{THM:PackagingSharpDependencyGluing}.
		Each facet closure is generated by at most one nontrivial implication, whose premise has cardinality at most $k$.
		
		Consequently, every nontrivial consequence in that facet has a determining subset of cardinality at most $k$.
		Thus each local closure belongs to $\mathscr C_k$.
		The two displayed parent implications are rooted circuits, and no original facet can recover their composite:
		the missing-pivot facet is an identity chart containing	the support of every proposed recovery.
		Since every clique of size at most $2k$ is guarded on this boundary, $\operatorname{ipr}(\mathscr C_k)\ge2k+1$.
		
		The same construction gives a compatible nonextendible $\mathscr C_k$-atlas on every boundary $\mathcal B_r$ with $3\le r\le2k+1$:
		the charts still carry at most one rule, so the hereditary Carath\'eodory bound holds.
		Thus $\operatorname{bor}(\mathscr C_k)\ge2k+1$.
		Proposition~\ref{PROP:PackagingBoundaryRankLowerBoundsRecovery} gives the reverse inequality.
		Finally, Lemmas~\ref{LEM:PackagingConcreteRelabelingStability} and~\ref{LEM:PackagingGenericNeutralTransport} supply relabeling stability, clone stability, and the required inert one-point lifts on every primitive core.
		Hence $\mathscr C_k$ is one-point transport complete through rank $2k+1$.
		
		\item \emph{Armstrong/native closures.}
		Theorem~\ref{THM:PackagingNativePropagationRecoveryRank} and the rooted-circuit characterization give $\operatorname{ipr}(\mathscr N_d)\le d$.
		On $\mathcal B_d$, the sharpness atlas of Theorem~\ref{THM:PackagingNativeSimplexBoundaryGeneral} contains two rooted parent implications whose composite cannot be recovered in any original chart:
		the	missing-pivot identity facet forbids every such consequence.
		Every clique of size at most $d-1$ is guarded, so $\operatorname{ipr}(\mathscr N_d)\ge d$.
		The same simplex-boundary theorem supplies nonextendible $\mathscr N_d$-atlases for every $3\le r\le d$, giving	$\operatorname{bor}(\mathscr N_d)\ge d$.
		The reverse inequality follows from the boundary-rank proposition.
		Relabeling stability and Lemma~\ref{LEM:PackagingNativePrimitiveNeutralCones} establish all three one-point transport conditions through rank $d$.
	\end{enumerate}
\end{proof}

\begin{remark}[Relation to classical acyclicity and bounded acyclicity]
	\label{REM:PackagingSpectrumAcyclicityScope}
	The geometric notions of chordality, full conformality, and bounded acyclicity are classical \cite{BeeriFaginMaierYannakakis1983Acyclic,Otto2012HighlyAcyclic}.
	For each integer $m\ge3$, write
	\[
	\mathfrak U_m
	:=\{\mathcal H:\mathcal H\text{ is a finite cover and }
	m<\nu(\mathcal H)\}.
	\]
	By the obstruction-rank dictionary, this means full chordality together with $m$-conformality.
	The contribution is the identification of a class-dependent conformality threshold from independently calculated recovery ranks; the underlying notion of hypergraph acyclicity is unchanged.
\end{remark}

The preceding sections have established the three ingredients needed for the final synthesis:
\[
\operatorname{ipr}(\mathscr C_k)=2k+1,
\qquad
\operatorname{ipr}(\mathscr N_d)=d,
\qquad
\rho_{\mathrm{pkg}}(\mathscr T)=3.
\]

For the closure classes, the recovery--obstruction duality converts the first two intrinsic ranks into exact geometric thresholds.
The transitive theory reaches the same obstruction scale through its independent cycle-composition calculus.
These identities yield the globalization obstruction spectrum.

\begin{theorem}[Main theorem: exact globalization obstruction spectrum]
	\label{THM:PackagingGlobalizationSpectrum}
	Let $\mathcal H$ be a finite cover.
	The following statements hold.
	\begin{enumerate}[label=(\roman*),leftmargin=2.5em]
		\item \emph{Bounded Carath\'eodory closures.}
		For every $k\ge1$, all compatible atlases of finite	empty-preserving closures with $\operatorname{car}\le k$ have exact global extensions if and only if
		\[
		2k+1<\nu(\mathcal H).
		\]
		
		\item \emph{Bounded-Armstrong FD theories.}
		For every $d\ge2$, all compatible families of finite FD theories $(\Sigma_I)_{I\in\mathcal H}$ satisfying $\operatorname{arm}(\Sigma_I)\le d$ admit global FD theories with exactly the prescribed local consequence theories if and only if
		$
		d<\nu(\mathcal H).
		$
		Empty-premise dependencies are permitted.
		No Armstrong-size bound is imposed on the global theory.
		For empty-preserving local theories this is the native statement for $\mathbf{CommCl}_d$.
		
		\item \emph{Transitive binary structures.}
		Universal exact amalgamation holds for each of the five classes of transitive relations, preorders, partial orders,	strict partial orders, and equivalence relations if and only if
		$
		3<\nu(\mathcal H).
		$
		
		\item \emph{Rank-three coincidence.}
		Writing $\operatorname{UG}_{\operatorname{car}\le1}$, $\operatorname{UG}_{\operatorname{arm}\le3}$, and	$\operatorname{UG}_{\mathrm{trans}}$ for the corresponding universal properties, one has
		\[
			\operatorname{UG}_{\operatorname{car}\le1}(\mathcal H)
			\Longleftrightarrow
			\operatorname{UG}_{\operatorname{arm}\le3}(\mathcal H)
			\Longleftrightarrow
			\operatorname{UG}_{\mathrm{trans}}(\mathcal H)
			\Longleftrightarrow
			\mathcal H\text{ is cycle-free}.
		\]
		
		\item \emph{Classical acyclic endpoint.}
		On finite covers,
		\[
			\bigcap_{k\ge1}\mathfrak U_{2k+1}
			=\bigcap_{d\ge3}\mathfrak U_d
			=\{\mathcal H:\mathcal H\text{ is chordal and conformal}\}.
		\]
		This is the classical $\alpha$-acyclic cover class.
	\end{enumerate}
\end{theorem}

\begin{proof}
	\leavevmode
	\begin{enumerate}[label=(\roman*),leftmargin=2.5em]
		\item \emph{Bounded Carath\'eodory closures.}
		The equivalence is Theorem~\ref{THM:PackagingGeneralCaratheodoryGluing} together with Lemma~\ref{LEM:PackagingObstructionRankDictionary}.
		Alternatively, the concrete intrinsic-rank calculation and Theorem~\ref{THM:PackagingAbstractRecoveryObstructionDuality} give the same sharp threshold.
		
		\item \emph{Bounded-Armstrong FD theories.}
		For $d\ge3$, apply Theorem~\ref{THM:PackagingClassicalBoundedArmstrongFDGlobalization} and the obstruction-rank dictionary.
		For $d=2$, Lemma~\ref{LEM:PackagingFDConstantReduction}	identifies the globally constant attributes and reduces the	compatible atlas to empty-preserving closures.
		Every nonempty reduced chart belongs to $\mathbf{CommCl}_2$, by the Armstrong/native dimension identity, and therefore carries the unique two-dimensional closure of Proposition~\ref{THM:PackagingNativeDimensionTwoClassification}.
		These closures admit an exact global extension on any cover.
		If the reduced label set is empty, take its unique closure.
		The constant-coordinate reduction reconstructs an exact global FD closure in either case, represented by a finite FD theory.
		Since $\nu(\mathcal H)\ge3$, the condition $2<\nu(\mathcal H)$ always holds, as required.
		
		\item \emph{Transitive binary structures.}
		By Lemma~\ref{LEM:PackagingObstructionRankDictionary}, $3<\nu(\mathcal H)$ is equivalent to chordality together	with $3$-conformality.
		By Proposition~\ref{PROP:PackagingBraultBaronCycleFree}, that is cycle-freedom.
		Apply Theorem~\ref{THM:PackagingTransitiveUniversalAmalgamation} for all five relation classes.
		
		\item \emph{Rank-three coincidence.}
		In parts~\textup{(i)}--\textup{(iii)}, set $k=1$ and $d=3$. Each resulting universal property is equivalent to $3<\nu(\mathcal H)$, hence to cycle-freedom.
		
		\item \emph{Classical acyclic endpoint.}
		A finite cover belongs to all $\mathfrak U_{2k+1}$,	$k\ge1$, if and only if its obstruction rank is infinite.
		The same is true for all $\mathfrak U_d$, $d\ge3$.
		By Lemma~\ref{LEM:PackagingObstructionRankDictionary}, $\nu(\mathcal H)=\infty$ is equivalent to chordality and full conformality, the classical $\alpha$-acyclicity criterion.
	\end{enumerate}
\end{proof}

\begin{theorem}[Bounded-complexity semantic hierarchy of covers]
	\label{THM:PackagingSemanticAcyclicityHierarchy}
	Let $\mathfrak U_m$, $m\ge3$, be the finite-cover classes defined above.
	Then:
	\begin{enumerate}[label=(\roman*),leftmargin=2.5em]
		\item $\mathfrak U_m$ consists precisely of finite covers with chordal primal graph which are $m$-conformal;
		
		\item
		$
		\mathfrak U_3\supsetneq\mathfrak U_4\supsetneq\cdots,
		$
		and $\bigcap_{m\ge3}\mathfrak U_m$ is the class of finite $\alpha$-acyclic covers;
		
		\item for every $d\ge3$ and $k\ge1$, $\mathfrak U_d$ is the universal exact-globalization class for local FD theories of Armstrong size at most $d$, and $\mathfrak U_{2k+1}$ is the corresponding class for empty-preserving closures with Carath\'eodory number at most $k$;
		
		\item $\mathfrak U_3$ is the universal exact-amalgamation class for transitive relations, preorders, partial orders, strict partial orders, and equivalence relations;
		
		\item let $\mathfrak O_m$ be the class of finite covers which are $m$-acyclic in Otto's sense, namely $m$-conformal and $m$-chordal (no chordless primal cycle of length $4,\ldots,m$) \cite{Otto2012HighlyAcyclic}.
		Then
		$
		\mathfrak U_m\subsetneq\mathfrak O_m.
		$
	\end{enumerate}
\end{theorem}

\begin{proof}
	\emph{(i) Geometric description.}
	Apply Lemma~\ref{LEM:PackagingObstructionRankDictionary}.
	
	\emph{(ii) Strict hierarchy and intersection.}
	The simplex boundary on $m+1$ vertices has $\nu=m+1$, so it belongs to $\mathfrak U_m\setminus\mathfrak U_{m+1}$.
	Membership in every $\mathfrak U_m$ is equivalent to $\nu=\infty$, hence to chordality and full conformality.
	
	\emph{(iii) Closure and FD semantics.}
	These are exactly parts~\textup{(i)}--\textup{(ii)} of Theorem~\ref{THM:PackagingGlobalizationSpectrum}.
	
	\emph{(iv) Transitive semantics.}
	Apply part~\textup{(iii)} of that theorem.
	
	\emph{(v) Comparison with bounded acyclicity.}
	Every chordal primal graph is $m$-chordal, and every cover in $\mathfrak U_m$ is $m$-conformal.
	Thus $\mathfrak U_m\subseteq\mathfrak O_m$.
	The edge cover of a chordless cycle on $m+1$ vertices is $m$-conformal, because its only nonempty cliques have at most two vertices.
	It is also $m$-chordal,	because its only chordless cycle has length $m+1$.
	It is not chordal, so it lies in $\mathfrak O_m\setminus\mathfrak U_m$.
\end{proof}

\begin{corollary}[Strict geometric advantage of Armstrong/native recovery]
	\label{COR:PackagingStrictRecoveryAdvantage}
	For every $d\ge3$ there exists a finite cover $\mathcal H_d$ on which every compatible atlas of empty-preserving FD closures with minimum Armstrong size at most $d$ admits an exact global extension, but some compatible atlas of empty-preserving closures with Carath\'eodory number at most $d-1$ does not.
	In particular, although $\mathscr N_d\subseteq\mathscr C_{d-1}$, the corresponding universal-globalization cover classes are strictly different.
\end{corollary}

\begin{proof}
	By Theorem~\ref{THM:PackagingNativeCommutantDimensionBound}, $\mathscr N_d\subseteq\mathscr C_{d-1}$.
	Let $\mathcal H_d=\mathcal B_{d+1}$ be the boundary of a simplex on $d+1$ vertices.
	Its primal graph is	complete, and its only unguarded clique is the full $(d+1)$-vertex set.
	Therefore
	$
	\nu(\mathcal H_d)=d+1.
	$
	
	Since $d<d+1$, Theorem~\ref{THM:PackagingGlobalizationSpectrum} \textup{(ii)} gives universal exact globalization for the Armstrong/native class on this cover.
	On the other hand, for $d\ge3$,
	\[
	2(d-1)+1=2d-1\ge d+1,
	\]
	so Theorem~\ref{THM:PackagingGlobalizationSpectrum} \textup{(i)} provides a compatible nonextendible atlas in $\mathscr C_{d-1}$ on the same cover.
	This proves strictness.
\end{proof}

\subsection{Explicit fixed-cover separations}
\label{SUBSEC:PackagingFixedCoverDistinctions}

An individual compatible atlas may extend even when another atlas on the same cover does not.
Also, two classes of local structures can behave differently on the same cover.
Proposition~\ref{PROP:PackagingIntroTwoSeparations} proves both facts with explicit examples.
They explain why our main theorems fix a cover and quantify over all compatible atlases in a specified class.

\begin{proposition}[Individual-atlas and semantic-class separation]
	\label{PROP:PackagingIntroTwoSeparations}
	\leavevmode
	\begin{enumerate}[label=(\roman*),leftmargin=2.5em]
		\item There is a finite cover carrying both an exactly extendible and a nonextendible compatible atlas of empty-preserving closures of Carath\'eodory number at most one.
		
		\item On the simplex-boundary cover of a four-element set, every compatible atlas in $\mathbf{CommCl}_3$ admits an exact global closure, whereas some compatible atlas of empty-preserving closures of Carath\'eodory number at most two does not.
	\end{enumerate}
\end{proposition}

\begin{proof}
	\emph{(i)} Let $V=\{a,b,c\}$ and take the edge cover $\mathcal H_3=\{\{a,b\},\{b,c\},\{a,c\}\}$.
	Assigning the identity closure to all three charts produces a compatible atlas extended by the identity closure on $V$.
	For a second atlas, let the closure on $\{a,b\}$ be generated by $a\to b$ and $b\to a$, let the closure on $\{b,c\}$ be generated by $b\to c$ and $c\to b$, and take the identity closure on $\{a,c\}$.
	Each local closure preserves the empty set and has Carath\'eodory number at most one.
	All overlap restrictions agree because every overlap is a singleton and all three closures are empty-preserving.
	In fact, each local closure belongs to $\mathbf{CommCl}_3$:
	the two-way dependency closures admit two spectral rows, while the identity closure on two labels has native dimension three by Theorem~\ref{THM:PackagingNativeIdentityDimension}.
	But a global exact extension would have $b\in C(\{a\})$ and $c\in C(\{b\})$, hence $c\in C(\{a\})$ by monotonicity and idempotence, contradicting exactness on $\{a,c\}$.
	
	\emph{(ii)} Let $V=\{a,b,x,y\}$ and	$\mathcal H_4=\{V\setminus\{t\}:t\in V\}$.
	The primal graph is complete, and the only unguarded clique is $V$ itself.
	Thus $\nu(\mathcal H_4)=4$.
	By Theorem~\ref{THM:PackagingNativeGeneralClassification} with $d=3$, every compatible $\mathbf{CommCl}_3$-atlas on this cover has an exact extension.
	For the contrasting atlas, distribute the two implications
	\[
	\{a,b\}\to x,\qquad \{x\}\to y
	\]
	to exactly those original charts containing their entire supports, and let each chart closure be generated by the rules assigned to that chart.
	No chart contains both supports, so every chart receives at most one nontrivial rule.
	The resulting atlas is empty-preserving and its local closures have Carath\'eodory number at most two.
	It is compatible:
	any nontrivial consequence visible on an overlap is generated by a rule whose full support lies in that overlap, and that rule is assigned to both charts.
	In the union of the chart rules, $\{a,b\}$ determines $x$ and then $y$.
	Yet the missing-$x$ chart $\{a,b,y\}$ has the identity closure, which forbids $\{a,b\}\to y$.
	This identity closure has native dimension four by Theorem~\ref{THM:PackagingNativeIdentityDimension}, so the displayed obstruction atlas is not an atlas in $\mathbf{CommCl}_3$.
	Therefore no exact global extension exists.
\end{proof}

The two examples distinguish individual-atlas extension from a universal fixed-cover property and demonstrate dependence on the class of local structures.
By Theorem~\ref{THM:PackagingAtlasRealizationBridge}, even a nonextendible compatible atlas defines a realized packaging system; exact global completion is a separate property.

\section{Discussion and conclusion}
\label{SEC:DiscussionA}
\label{SEC:ConclusionA}

The packaging theory separates a packaged output from the context and formation datum that realize it, and supplies intrinsic completion invariants and comparison-relative non-collapse tests. Its primitive coordinates also admit a faithful relational presentation (Proposition~\ref{PROP:PackagingRelationalPresentation}).

The atlas realization bridge (Theorem~\ref{THM:PackagingAtlasRealizationBridge}) expresses compatible local data and exact extension as a completion problem for realized packaging systems. The recovery-equipped specializations developed in Section~\ref{SEC:SpectrumA} introduce additional restriction, generation, and elimination laws; under these laws, native recovery gives a general sufficient globalization criterion, while obstruction portability supports a sharp converse. The concrete thresholds of Theorem~\ref{THM:PackagingGlobalizationSpectrum} require class-specific recovery arguments and explicit finite obstructions. The primitive packaging axioms organize these problems but do not by themselves imply the sharp extension criteria.

The same cover obstruction rank controls three classes of local structures, but the classes have different recovery mechanisms.
For \(d\ge3\), the Armstrong/native and generic closure classes have different intrinsic ranks:
\[
\operatorname{ipr}(\mathscr N_d)=d<2d-1=\operatorname{ipr}(\mathscr C_{d-1}).
\]
These different ranks give different classes of covers on which universal globalization holds (Corollary~\ref{COR:PackagingStrictRecoveryAdvantage}).
Transitive relations instead use path composition to reach their rank-three criterion.
For the two closure classes, recovery rank alone gives only sufficiency; sharp necessity also uses transport of explicit obstructions.
Throughout the paper, the complexity bound applies to local data.
A global extension can be more complex.

The strengthened propagation theorem also gives exact localization on each primal clique.
Under bounded conformality, the cliques have compatible canonical fillers even when the full graph is not chordal.
Adding those fillers does not change the closure generated by the original atlas.
A separate example shows that local dimension three can force minimum global dimension four.
Finally, the stable-class counterexample proves that the abstract converse needs an obstruction-existence or transport assumption.
Together, these results connect an axiomatic theory of formation and intrinsic completion with a precise recovery calculus and sharp extension theorems in specified specializations.

\section*{Disclosure of AI assistance}

The author used ChatGPT in developing and refining certain
mathematical ideas and arguments, as well as for language
editing and comparisons with prior literature.
In particular, the AI-assisted mathematical constructions are as below:
\begin{itemize}[leftmargin=2.1em,nosep]
	\item Section~\ref{SEC:FoundationsPackagingSystems}:
	Selection freedom
	(Theorem~\ref{THM:PackagingSelectionFreedom}).

	\item Section~\ref{SEC:PackagingIntrinsicStructureA}:
	Intrinsic refinement completion
	(Theorem~\ref{THM:PackagingIntrinsicRefinementCompletionA}).
	
	\item Section~\ref{SEC:PackagingControlledGlobalization}:
	Controlled propagation, canonical clique completion, and bounded-Carath\'eodory globalization
	(Theorems~\ref{THM:PackagingControlledPropagationGlobalization},
	\ref{THM:PackagingCanonicalCliqueCompletion},
	and~\ref{THM:PackagingGeneralCaratheodoryGluing}).
	
	\item Section~\ref{SEC:ArmstrongNativeA}:
	Commutant realization dimension and minimum Armstrong size
	(Theorem~\ref{THM:PackagingCommutantDimensionEqualsArmstrongSize}).
	
	\item Section~\ref{SEC:NativeRecoveryA}:
	Native clique resolution, bounded-Armstrong globalization,
	and strict realization-dimension separation
	(Theorems~\ref{THM:PackagingNativeGeneralResolution}
	and~\ref{THM:PackagingNativeGeneralClassification};
	Proposition~\ref{PROP:PackagingNativeStrictGlobalDimensionJump}).
	
	\item Section~\ref{SEC:PackagingCycleGeometry}:
	Transitive universal amalgamation
	(Theorem~\ref{THM:PackagingTransitiveUniversalAmalgamation}).
	
	\item Section~\ref{SEC:SpectrumA}:
	Obstruction-transport counterexample, abstract recovery--obstruction duality, and globalization spectrum
	(Theorems~\ref{THM:PackagingTransportNecessityCounterexample},
	\ref{THM:PackagingAbstractRecoveryObstructionDuality},
	and~\ref{THM:PackagingGlobalizationSpectrum}).	
\end{itemize}
The author takes full responsibility for the content of this
paper, including the correctness of its mathematical arguments
and the integrity and accuracy of its citations.

\bibliographystyle{unsrturl}
\bibliography{references_packaging}

\end{document}